\documentclass[11pt]{amsart}
\usepackage{amscd,amsmath,amssymb,amsfonts,verbatim}
\usepackage[cmtip, all]{xy}
\usepackage{MnSymbol}
\usepackage[a4 paper, left=3cm, right=3cm, top = 2.5cm, bottom = 2.5cm ]{geometry}
\usepackage{comment}
\usepackage[colorlinks, 
linkcolor=black,anchorcolor=black, citecolor=black, 
linkcolor=black]{hyperref}

\newtheorem{thm}{Theorem}[section]
\newtheorem{prop}[thm]{Proposition}
\newtheorem{lem}[thm]{Lemma}
\newtheorem{cor}[thm]{Corollary}

\theoremstyle{definition}
\newtheorem{ques}[thm]{Question}

\newtheorem{defn}[thm]{Definition}
\newtheorem{remk}[thm]{Remark}
\newtheorem{remks}[thm]{Remarks}

\newtheorem{exm}[thm]{Example}
\newtheorem{exms}[thm]{Examples}
\newtheorem{notat}[thm]{Notation}
\numberwithin{equation}{section}

{\hfill$\square$\end{defn}}
{\hfill$\square$\end{remk}}
{\hfill$\square$\end{remks}}
{\hfill$\square$\end{exm}}
{\hfill$\square$\end{exms}}
{\hfill$\square$\end{notat}}

\newcommand{\thmref}{Theorem~\ref}
\newcommand{\propref}{Proposition~\ref}
\newcommand{\corref}{Corollary~\ref}

\newcommand{\lemref}{Lemma~\ref}

\newcommand{\sA}{{\mathcal A}}

\newcommand{\sC}{{\mathcal C}}
\newcommand{\sD}{{\mathcal D}}
\newcommand{\sE}{{\mathcal E}}
\newcommand{\sF}{{\mathcal F}}
\newcommand{\sG}{{\mathcal G}}

\newcommand{\sI}{{\mathcal I}}

\newcommand{\sK}{{\mathcal K}}
\newcommand{\sL}{{\mathcal L}}

\newcommand{\sO}{{\mathcal O}}
\newcommand{\sP}{{\mathcal P}}

\newcommand{\sS}{{\mathcal S}}

\newcommand{\sW}{{\mathcal W}}
\newcommand{\sX}{{\mathcal X}}

\newcommand{\sZ}{{\mathcal Z}}
\newcommand{\A}{{\mathbb A}}

\newcommand{\F}{{\mathbb F}}
\newcommand{\G}{{\mathbb G}}
\renewcommand{\H}{{\mathbb H}}

\newcommand{\N}{{\mathbb N}}
\renewcommand{\P}{{\mathbb P}}
\newcommand{\Q}{{\mathbb Q}}

\newcommand{\Z}{{\mathbb Z}}

\newcommand{\fm}{{\mathfrak m}}

\newcommand{\fp}{{\mathfrak p}}
\newcommand{\fq}{{\mathfrak q}}

\newcommand{\ff}{{\mathfrak f}}

\newcommand{\esp}{{\rm sp}}

\newcommand{\Ker}{{\rm Ker}}
\newcommand{\gr}{{\rm gr}}

\newcommand{\Alb}{{\rm Alb}}
\newcommand{\CH}{{\rm CH}}

\newcommand{\surj}{\twoheadrightarrow}
\newcommand{\inj}{\hookrightarrow}
\newcommand{\red}{{\rm red}}

\newcommand{\Pic}{{\rm Pic}}
\newcommand{\Div}{{\rm Div}}
\newcommand{\Hom}{{\rm Hom}}

\newcommand{\Spec}{{\rm Spec \,}}
\newcommand{\sing}{{\rm sing}}
\newcommand{\Char}{{\rm char}}

\newcommand{\ms}{{\rm ms}}

\newcommand{\Tr}{{\rm Tr}}
\newcommand{\ur}{\rm ur}

\newcommand{\Ev}{\rm Ev}
\newcommand{\nr}{\rm nr}
\newcommand{\Gal}{{\rm Gal}}
\newcommand{\divf}{{\rm div}}

\newcommand{\Et}{{\rm {\bf{Et}}}}
\newcommand{\Zars}{{\rm {\bf{Zar}}}}

\newcommand{\NS}{{\operatorname{NS}}}
\newcommand{\id}{{\operatorname{id}}}

\newcommand{\Sch}{{\operatorname{\mathbf{Sch}}}}

\newcommand{\op}{{\text{\rm op}}}
\newcommand{\<}{\langle}
\renewcommand{\>}{\rangle}

\newcommand{\can}{{\operatorname{\rm can}}}

\newcommand{\fin}{{\operatorname{\rm fin}}}

\newcommand{\Wedge}{{\Lambda}}
\newcommand{\eff}{{\operatorname{\rm eff}}}

\newcommand{\et}{{\textnormal{\'et}}}

\newcommand{\ds}{{/\kern-3pt/}}

\newcommand{\res}{{\operatorname{res}}}

\renewcommand{\log}{{\operatorname{log}}}

\newcommand{\Tor}{{\operatorname{Tor}}}

\newcommand{\Br}{{\operatorname{Br}}}

\newcommand{\Picc}{{\mathbf{Pic}}}

\newcommand{\tr}{{\operatorname{tr}}}

\newcommand{\Cor}{{\operatorname{Cor}}}

\newcommand{\un}{\underline}
\newcommand{\ov}{\overline}

\renewcommand{\dim}{\text{\rm dim}}

\newcommand{\tuborg}{\left\{\begin{array}{ll}}
\newcommand{\sluttuborg}{\end{array}\right.}

\newcommand{\pr}{{\rm pr}}

\newcommand{\zar}{{\rm zar}}
\newcommand{\nis}{{\rm nis}}

\newcommand{\homg}{{\rm homog}}

\newcommand{\reg}{{\rm reg}}

\newcommand{\tor}{{\rm tor}}
\newcommand{\ns}{{\rm NS}}
\newcommand{\dlog}{{\rm dlog}}

\newcommand{\ev}{{\rm ev}}
\newcommand{\inv}{{\rm inv}}
\newcommand{\alb}{{\rm alb}}

\newcommand{\spc}{{\rm sp}}

\newcommand{\bn}{{\rm BN}}
\newcommand{\rsw}{{\rm rsw}}
\newcommand{\Rsw}{{\rm Rsw}}

\newcommand{\Sw}{{\rm Sw}}

\newcommand{\wt}{\widetilde}
\newcommand{\wh}{\widehat}
\newcommand{\cont}{{\rm cont}}

\newcommand{\coker}{{\rm Coker}}
\newcommand{\Fil}{{\rm fil}}

\newcommand{\n}{{\un{n}}}

\newcommand{\etl}{{\acute{e}t}}

\newcounter{elno}

\newcounter{elno-abc}   
\newenvironment{listabc}{
                         \begin{list}{\alph{elno-abc})
                                     }{\usecounter{elno-abc}}
                      }{
                         \end{list}}
\newcounter{elno-abc-prime}   

\begin{document}
\title[Brauer group of varieties]{Brauer group of varieties
over local fields of finite characteristic}
\author{Amalendu Krishna, Subhadip Majumder}
\address{Department of Mathematics, South Hall, Room 6607, University of California
  Santa Barbara, CA, 93106-3080, USA.}
\email{amalenduk@math.ucsb.edu}
\address{Institute of Mathematical Sciences, A CI of Homi Bhabha National Institute,
  4th Cross St., CIT Campus, Tharamani, Chennai,
  600113, India.} 
\email{subhamajumdar.sxc@gmail.com}




\maketitle
\vskip .4cm

\begin{quote}\emph{Abstract.}We show that the non-logarithmic version of Kato’s
  ramification filtration on the Brauer group of a separated regular scheme of
  finite type over a henselian discrete valuation field of positive
  characteristic with finite residue field coincides with the evaluation
  filtration. This extends a recent result of Bright--Newton to positive
  characteristic. As applications, we extend several results of Ieronymou,
  Saito--Sato, and Kai to positive characteristic.
\end{quote}

\setcounter{tocdepth}{1}
\tableofcontents

\vskip .4cm

\section{Introduction}\label{sec:Intro}
The Chow group of zero-cycles and the Brauer group of an algebraic variety are
fundamental invariants to study in arithmetic geometry. The objective of this
paper is to study these groups for smooth varieties over a local field
(i.e., a complete discrete valuation field with finite residue field) of positive
characteristic, and to establish several fundamental results that are known over
local fields of characteristic zero but remain unknown in positive
characteristic. In this section, we provide the necessary background and state
the main results.

\subsection{The Bright--Newton theorem}\label{sec:BN-thm}
Let $p > 0$ be a prime number and let $k$ be a henselian discrete
valuation field whose ring of integers $\sO_k$ is
an excellent ring and residue field $\F$ is finite of characteristic $p$.
Let $\fm = (\pi)$ denote the maximal ideal of $\sO_k$.
Let $\sX$ be a  finite type, separated, faithfully flat $\sO_k$-scheme which
is regular. We write $X = \sX \times_S \Spec(k), \ \sX_s = \sX \times_S \Spec(\F)$
and $Y = (\sX_s)_\red$, where $S = \Spec(\sO_k)$. We assume that
$Y$ is a simple normal crossing divisor on $\sX$ with irreducible components
$Y_1, \ldots , Y_r$. We do not assume that $k$ is complete or
$X$ is smooth over $k$.

In order to study the ramification in the
cohomology groups $H^q_\et(X, {\Q}/{\Z}(q-1))$ for $q \ge 1$, 
Kato \cite[\S~7, 8]{Kato-89} introduced a filtration on these
cohomology groups. This filtration for $q =2$ induces a filtration
on the Brauer group of $X$, given by
\begin{equation}\label{eqn:Kato-fil-0}
  \Fil_n \Br(X) = \{\chi \in \Br(X)| \chi_i \in \Fil_n \Br(K_i) \ \forall \ i\}
\end{equation}
for $n \ge 0$. Here, 
$K_i$ is the henselization of the function field of $\sX$ at the generic point
of $Y_i$, $\chi_i$ is the restriction of $\chi$ to $\Br(K_i)$ and
$\Fil_n \Br(K_i)$ is the set of elements $\chi' \in \Br(K_i)$ such that
cup product $\{\chi'_L, 1 + \fm_{K_i}^{n+1}\sO_L\}$ is zero in
$H^3_\et(L, {\Q}/{\Z}(2))$ for every henselian discrete valuation field $L$ over
$K_i$ such that $\sO_{K_i} \subset \sO_L$ and $\fm_L = \fm_{K_i} \sO_L$
(cf \S~\ref{sec:Kato-fil}).
We shall refer to $\Fil_\bullet \Br(X)$ as the Kato filtration of $\Br(X)$.

The Kato filtration of $\Br(X)$ determines the ramification of a Brauer class
on $X$ along $Y$. For instance, $\Fil_0\Br(X)$ can be viewed as
the set of those Brauer classes on $X$ which are `tamely ramified along $Y$'
(cf. \propref{prop:Kato-basic}).

When $\sX$ is smooth over $S$, Bright--Newton introduced another filtration on
$\Br(X)$ in \cite{Bright-Newton}, called the evaluation filtration.
This filtration is defined in terms of the size of the adic neighborhoods of closed
points of $X$ over which a Brauer class on $X$ has constant evaluation  (cf.
\S~\ref{sec:Ev-Br}). Consequently, the evaluation filtration plays a
significant role in determining the relevance of Brauer classes in
the Brauer--Manin obstruction to Hasse principle and weak approximation.

In general, it is very difficult to measure the
size of the open subsets where a Brauer class takes a constant value.
It is therefore natural to ask if the Brauer classes that have constant evaluation on
an open subset can be described in terms of some cohomology groups.
Bright--Newton showed that, if $k$ is a $p$-adic field and $\sX$ is smooth over
$\sO_k$, then the evaluation filtration on $\Br(X)$ coincides with the non-log version
(for $H^1_\et(X, {\Q}/{\Z})$, this version was defined by Matsuda \cite{Matsuda})
of the Kato filtration.

The result of Bright--Newton (cf. \cite[Thm.~A]{Bright-Newton})
has turned out to be of considerable importance in arithmetic geometry,
and has led to the solutions of several important questions about the Brauer group
of smooth projective varieties over local and global fields of characteristic zero.
The objective of this paper is to extend the work of Bright--Newton
to positive characteristic, and give applications
to some important open problems about zero-cycles and
Brauer groups of smooth projective varieties over local fields.

We assume from now on that $\Char(k) = p$.
We introduce a modified version of the Bright--Newton filtration which
works in more general situations (that we work in) than when $\sX/S$ is smooth
(however, see Remark~\ref{remk:BN-mod}). 
Let $X_{(0)}$ be the set of all closed points of $X$ and let
$X^o_\fin$ be the set of those points $P \in X_{(0)}$ whose 
scheme-theoretic closure $\ov{\{P\}}$ in $\sX$ is regular, is finite over $S$ and
is disjoint from the singular locus of $Y$.
For any $n \ge 1$ and $P \in X_{(0)}$, we define a subset $B(P,n) \subset X^o_\fin$
which consists of all points $Q \in X^o_\fin$ such that
$\ov{\{Q\}} \times_{\sX} nY = \ov{\{P\}} \times_S nY$ (cf. \S~\ref{sec:Ev-Br}).
We let $B(P,0) = \{Q\in X^o_{\fin}|[k(Q) : k] = [k(P) : k]\}$.
For $\sA \in \Br(X)$, we let
$\ev_\sA \colon X_{(0)} \to {\Q}/{\Z}$ be the evaluation map
obtained by restricting
$\chi$ to closed points and applying the Hasse invariant maps for residue fields.

For $n \ge -1$, we let
\begin{equation}\label{eqn:Kato-fil-1}
  {\Ev}_n \Br(X) = \left\{\sA \in \Br(X)| \ev_\sA \ \text{is constant on} \
  B(P,n+1) \
\text{ for all } P \in X^o_\fin \right\};
\end{equation}
\[
  {\Ev}_{-2}\Br(X) = \left\{\sA \in \Br(X)|  \ev_\sA \ \text{is zero on} \
  X^o_\fin\right\}.
  \]
We shall refer to ${\Ev}_{\bullet}\Br(X)$ as the evaluation filtration on
$\Br(X)$.
We let $y_i$ denote the generic point of $Y_i$ and let
$\rsw^2_{K_i, n} \colon \frac{\Fil_n \Br(X)}{\Fil_{n-1} \Br(X)} \to
\Omega^2_{k(y_i)} \bigoplus \Omega^1_{k(y_i)}$ be the refined Swan conductor map
for $n \ge 1$ (cf. \S~\ref{sec:RSC}).
Let $g \colon Y \to \Spec(\F)$ be the projection.

\vskip.2cm

The first main result of this paper is the following extension of the theorem
of Bright--Newton in positive characteristic.

\begin{thm}\label{thm:Main-1}
  Let $\sX$ be as above and let $Y^o$ denote the regular locus of $Y$.
  Then there exists a canonical map
  $\partial_X \colon \Fil_0 \Br(X) \to H^1_\et(Y^o, {\Q}/{\Z})$.
  Moreover, we have the following.
  \begin{enumerate}
  \item
    If $n \ge 1$, then
    \[
      {\Ev}_n \Br(X) = \left\{\sA \in \Fil_{n+1} \Br(X)| \rsw^2_{K_i, n+1}(\sA) =
      (\alpha, 0) \in \Omega^2_{k(y_i)} \bigoplus \Omega^1_{k(y_i)}
      \ \forall \ i \right\}.
      \]
    \item
      ${\Ev}_0 \Br(X) = \Fil_0 \Br(X)$.
    \item
      If $\sX_s$ is reduced with geometrically connected irreducible components,
      then
      \[
        {\Ev}_{-1} \Br(X) = \left\{\sA \in \Fil_0 \Br(X)| \partial_X(\sA) \in
        \ {\rm Image}\left(H^1_\et(\F, {\Q}/{\Z}) \xrightarrow{{g}^*}
        H^1_\et(Y^o, {\Q}/{\Z})\right)\right\}.
          \]
\item
  ${\Ev}_{-2} \Br(X) = \Br(\sX)$.
  \end{enumerate}
  \end{thm}

\vskip.2cm

We make some remarks about \thmref{thm:Main-1}.
\begin{remk}\label{remk:BN-mod}
  (1) The item (3) of \thmref{thm:Main-1} remains true even when $\sX_s$ is not
  reduced provided all its components occur with the same multiplicity in its
  presentation as a Weil divisor on $\sX$ (cf. \lemref{lem:K-Ev-8}).

  (2) The above definition of the evaluation filtration is motivated by that of
  Bright--Newton but is slightly different. The modification in the new definition
  is necessary in order to handle the case (as we do in this paper) where the
  special fiber of $\sX$ is not reduced or not irreducible.
  Although it not clear from the
  definition that the evaluation filtration defined above coincides with that of
  \cite{Bright-Newton} when
  $\sX$ is smooth over $\sO_k$ (the case considered in \op. cit.), we show
  that the two filtrations do coincide under the smoothness assumption
  (cf. \thmref{thm:K-Ev-11}).

  (3) While the study of the relation between Kato and evaluation filtrations
  in characteristic zero is done in \cite{Bright-Newton} under the assumption
  that $\sX$ is a smooth $\sO_k$-scheme, no such assumption is made in
  \thmref{thm:Main-1}. It does not also assume that $X$ is a smooth $k$-scheme.
  Note that smoothness of $X$ in characteristic zero is an automatic
  consequence of its regularity but that is not the case over positive
  characteristic local fields. \thmref{thm:Main-1} assumes only that
  $\sX$ is regular and $Y \subset \sX$ is a simple normal crossing divisor.

  (4) In the course of proving \thmref{thm:Main-1}, we show that the inclusions
${\rm RHS} \subset {\rm LHS}$ in parts (1) and (2) of \thmref{thm:Main-1}
hold even when $Y$ is not a simple normal crossing divisor in $\sX$.
In fact, we expect that, with a suitable formulation, the reverse inclusions
in \thmref{thm:Main-1} should also hold without this assumption.
We will return to this question elsewhere.

(5) We show that the evaluation filtration can be completely determined by
evaluating Brauer classes on $X$ at only those closed points of $X$ whose
closures in $\sX$ are transverse to $Y$. In particular, 
this filtration can be completely determined by evaluating Brauer classes on
$X$ at only its unramified closed points, rather than at all points of $X^o_\fin$,
when $\sX_s$ is reduced.

(6) When $X$ is a smooth, projective, geometrically connected curve over $k$,
Yamazaki \cite{Yamazaki} defined a filtration on $\Br(X)$ using the Brauer--Manin
pairing and showed that this filtration coincides with the Kato filtration.
Combining Yamazaki's result with \thmref{thm:Main-1}, it follows that his
filtration also coincides with the evaluation filtration.
 \end{remk}

Now, we explain the applications of \thmref{thm:Main-1} to  Brauer--Manin obstruction,
and to two open problems on
Brauer groups and Chow groups of zero-cycles on projective varieties over local
fields. For more applications of \thmref{thm:Main-1} to Brauer--Manin obstruction
to weak approximation, we refer the reader to \cite{Lazda-Skorobogatov}.

\subsection{Theorems of Ieronymou}\label{sec:Ier**}
One of the applications of \thmref{thm:Main-1} and its mixed
characteristic analogue in  \cite{Bright-Newton} is to determine conditions
on a variety $X$ over a local field $K$ under which all Brauer classes
on $X$ have constant values on the $K$-points of $X$.
This problem is important in the study of the Brauer--Manin obstruction to
Hasse principle and weak approximation. In \cite{Ieronymou},
Ieronymou proved several
results of this kind in characteristic zero.

One of these says that
if $K$ is a $p$-adic field, and $\sX$ is a smooth
projective scheme over $\sO_K$ with generic fiber $X$ such that the special
fiber is separably rationally connected (cf. \cite{Kollar-Szabo}),
then every Brauer class in $\Br(X)$
induces the constant evaluation map on $X(K)$. He proved a similar result if the
generic fiber of $\sX$ is an Enriques surface (cf. \cite{CDL}).

Among several applications of \thmref{thm:Main-1},
we prove the following extensions of
the theorems of Ieronymou to positive characteristic. 

\begin{thm}\label{thm:Main-2}
Let $\sX$ be as in \thmref{thm:Main-1}, and assume additionally that $\sX$ is
smooth and projective over $S$ such that the special fiber $Y$ is separably
rationally connected. Then
\[
  {\ev}_{\sA} \colon X(k) \to {\Q}/{\Z}
  \]
  is a constant map for every $\sA \in \Br(X)$.
\end{thm}

\begin{thm}\label{thm:Main-6}
  Let $\sX$ be as in \thmref{thm:Main-1}, and assume additionally that $\sX$
  is smooth and projective over $S$ such that the generic fiber $X$
  is an Enriques surface.
  Let $\sA \in \Br(X)$ have odd order. Then the evaluation
  map ${\ev}_{\sA} \colon X(k) \to {\Q}/{\Z}$ is constant.
\end{thm}

\begin{remk}\label{remk:Ierony}
  In the characteristic zero version of \thmref{thm:Main-2} proven in
  \cite{Ieronymou}, it is assumed that $k$ is complete,
  $X$ is geometrically integral over $k$ and the
  residue characteristic of $k$ is an odd prime.
  But we do not need any of these hypotheses.
\end{remk}

\vskip.2cm

\subsection{The theorem of Saito--Sato}\label{sec:BM-thm}
If $K$ is any henselian discrete valuation field,
and $X$ is a smooth, connected and complete $K$-variety,
there is a natural pairing (known as the Brauer--Manin pairing)   
\begin{equation}\label{eqn:BMP-00}
  \CH_0(X) \times \Br(X) \to {\Q}/{\Z}
\end{equation}
which is functorial with respect to proper maps of complete $K$-varieties.
Here, $\CH_0(X)$ is the Chow group of 0-cycles on $X$.

When $X$ is geometrically integral of dimension one, Lichtenbaum
\cite{Lichtenbaum} (in characteristic zero) and Saito \cite{Saito-Invent}
(in positive characterstic) showed that ~\eqref{eqn:BMP-00} is a perfect
pairing. In particular, it induces an isomorphism
\begin{equation}\label{eqn:BMP-01}
  A_0(X) \xrightarrow{\cong} \Hom({\Br(X)}/{\Br_0(X)}, {\Q}/{\Z}),
\end{equation}
where $A_0(X) \subset \CH_0(X)$ is the subgroup of degree zero 0-cycles and
$\Br_0(X)$ is the image of the pull-back map $\Br(k) \to \Br(X)$.
It remained an open question whether the results of Lichtenbaum and Saito extend
to higher dimensions. This question was answered negatively by
Parimala--Suresh \cite{Parimala-Suresh}, who showed that ~\eqref{eqn:BMP-00} can
have a nonzero left kernel in general when $\dim(X)=2$.

To say something about the right kernel, let us suppose that $X$ admits a regular
model $\sX$ which is finite type, proper, and flat over $\sO_K$.
In this setting, one observes that ~\eqref{eqn:BMP-00} induces natural
homomorphisms (this uses the fact that $\Br(\sO_K) = 0$)
\begin{equation}\label{eqn:BMP-02}
  \theta_X \colon \Br(X) \to \Hom(\CH_0(X), {\Q}/{\Z});
\end{equation}
\begin{equation}\label{eqn:BMP-03}
  \theta^\star_X \colon A_0(X) \to \Hom({\Br(X)}/{\Br(\sX) + \Br_0(X)}, {\Q}/{\Z}).
\end{equation}

A major problem concerning the pairing ~\eqref{eqn:BMP-00} is to describe the
kernel of $\theta_X$ and the cokernel of $\theta_X^\star$. In
\cite{Saito-Sato-ENS}, Saito--Sato made a major breakthrough in the study of the
Brauer--Manin pairing by resolving this problem. They showed that,
if $K$ is a local field of characteristic zero and $X$ is geometrically integral
over $K$, then $\Ker(\theta_X) = \Br(\sX)$ and $\coker(\theta^\star_X) = 0$.
This was conjectured by Colliot-Th{\'e}l{\`e}ne
\cite[Conj.~1.4]{CT-Bordeaux}, and its  prime-to-$p$ part was
solved by Colliot-Th{\'e}l{\`e}ne--Saito \cite{CTSaito}.
However, the problem has remained open in positive characteristic.

The second main result of this paper is the following extension of the theorem
of Saito--Sato which verifies Colliot-Th{\'e}l{\`e}ne's prediction
in positive characteristic, when
$X$ admits a semistable (but not necessarily strictly semistable) reduction.
This is our second application of \thmref{thm:Main-1}.

\begin{thm}\label{thm:Main-3}
  Let $\sX$ be as in \thmref{thm:Main-1}. Assume additionally the following.
  \begin{enumerate}
  \item
    $k$ is complete.
  \item
    $\sX$ is a projective $\sO_k$-scheme and $X$ is a smooth $k$-scheme.
  \end{enumerate}
  Then $\Ker(\theta_X) = \Br(\sX)$ and $\theta^\star_X$ is
  surjective.
\end{thm}

\vskip.2cm

\begin{remk}\label{remk:Main-2-0}
  (1) We emphasize that, although the first part of \thmref{thm:Main-3} follows
  immediately from \thmref{thm:Main-1}(4), the proof of the second part requires
  additional arguments. One of these involves reducing to the case of curves and
  studying the specialization of the Swan conductor under the
  restriction of a Brauer class to curves
  (cf. \thmref{thm:SC-change-main}).

  (2) We also remark that the methods of \cite{Saito-Sato-ENS} can not be applied
  to prove the second part of \thmref{thm:Main-3}. One reason is that Mattuck's
  theorem \cite{Mattuck}, which is an essential ingredient of op. cit., does not
  hold in positive characteristic. Another point worth noting is that, unlike 
  op. cit., we do not assume $X$ to be geometrically connected in 
  \thmref{thm:Main-3}. 
\end{remk}

As a byproduct of the proof of \thmref{thm:Main-3}, we obtain 
the following result. The analogous result over characteristic zero local
fields was proven by Colliot-Th{\'e}l{\`e}ne-Saito \cite{CTSaito}
(prime-to-$p$ case) and Saito--Sato \cite{Saito-Sato-ENS} (general case).

\begin{thm}\label{thm:Main-5}
  Let $X$ be as in \thmref{thm:Main-3}. There is an inclusion of
  abelian groups
  \begin{equation}\label{eqn:Main-5-0}
  \frac{\Br(X)}{\Br(\sX) + \Br_0(X)} \inj
  \left(\stackrel{\infty}{\underset{i =1}\bigoplus} {\Q_p}/{\Z_p}\right)
  \bigoplus T,
  \end{equation}
  where $T$ is finite. In particular, we have the following.
  \begin{enumerate}
  \item
    $\frac{\Br(X)}{\Br(\sX) + \Br_0(X)}\{p'\}$ is finite.
  \item
    $\frac{\Br(X)}{\Br(\sX) + \Br_0(X)} \otimes_\Z {\Z}/n$ is finite for every
    integer $n$ prime to $p$.
    \end{enumerate}
\end{thm}

\begin{remk}\label{remk:Main-5-1}
  (1) The inclusion of ~\eqref{eqn:Main-5-0} is known to be an isomorphism when
  $\dim(X) =1$ and $g(X) > 0$.
  This shows that item (1) fails for the $p$-primary torsion
  group in general. However, we expect item (2) to hold for every integer $n$.

  (2) The inclusion in ~\eqref{eqn:Main-5-0} may not be an isomorphism (cf.
  Example~\ref{exm:Trivial-Br-p}).
  \end{remk}

\vskip.2cm

\subsection{Theorems of Lichtenbaum and Kai}\label{sec:Kai-thm}
The Albanese map is known to play an important role in the study of the Chow
group of $0$-cycles on smooth projective varieties over fields. If $X$ is a
smooth projective variety over a field $K$, and $\Alb_X$ denotes its
Albanese variety (which exists by \cite[p.~45--46]{Lang}), then the Albanese map
is unfortunately not always surjective (unless $K$ is finite or algebraically
closed), and describing its cokernel is, in general, a difficult problem.

We consider this problem over local fields.
Suppose first that $K$ is a local field of characteristic zero and
$\ov{K}$ is an algebraic closure of $K$ with Galois group $\Gamma$. 
Let $X$ be a geometrically integral smooth projective variety over $K$ and
$\ov{X}$ be the base change of $X$ to $\ov{K}$.
If $\dim(X) =1$, Lichtenbaum proved a duality theorem
in \cite{Lichtenbaum} which says that $\coker(\alb_X)$ is finite, and is
canonically Pontryagin dual to the cokernel of the degree map
$\Pic(\ov{X})^{\Gamma} \to \Z$.
A characteristic free proof of Lichtenbaum's duality was given 
by Saito \cite{Saito-Invent}.

When $k$ is a local field of characteristic zero, Kai \cite{Kai}
proposed a conjectural
formula for the Albanese cokernel in higher dimensions. He proved
this formula under the assumption that the Picard scheme of a smooth model is
smooth. Kai's result, in turn, provided a higher-dimensional extension of
Lichtenbaum's duality theorem.

The higher-dimensional extension of Lichtenbaum's theorem has remained an open
problem in positive characteristic,
and there is no conjectural formula for the Albanese cokernel.
Using \thmref{thm:Main-3}, we extend Kai's duality theorem to positive
characteristic as follows. This shows that his proposed formula for the
Albanese cokernel may be true in positive characteristic too.

\begin{thm}\label{thm:Main-4}
  Let $X$ be as in \thmref{thm:Main-3}. Assume additionally that
  $\sX$ is smooth over $S$ and $X$ is geometrically integral over $k$. Let $k_s$
  denote a separable closure of $k$ and let $\Gamma$ denote the Galois group of
  ${k_s}/k$. Let $X^s$ denote the base change of $X$ by $k_s$.
  Assume that the Picard scheme $\Picc(\sX)$ is smooth over $S$.
  Then we have the following.
  \begin{enumerate}
  \item
    The cokernel of the Albanese map
    \[
    \alb_X \colon A_0(X) \to \Alb_X(k)
    \]
  is a finite $p$-group.
 \item
$\coker(\alb_X)$ is canonically Pontryagin dual to the cokernel of the map
  \[
  \Pic(X^s)^\Gamma \to \ns(X^s)^\Gamma.
  \]
  \end{enumerate}
\end{thm}

 \begin{remk}\label{remk:Main-4-0}
   \thmref{thm:Main-4} implies that the Albanese cokernel of $X$ is finite.
   The finiteness of Albanese cokernel in higher dimensions was not known
   in positive characteristic.
   In characteristic zero, Saito--Sujatha \cite{Saito-Sujatha} proved without any
   assumption on $\Picc(\sX)$ that the Albanese cokernel is finite. But their
   proof does not extend to positive characteristic due to the failure of
   Mattuck's theorem \cite{Mattuck}.
\end{remk}

\vskip .2cm

\subsection{Outline of proofs}\label{sec:Outline}
The proof of \thmref{thm:Main-1} is the most challenging part of this work. Since
the proof is both lengthy and technically involved, we provide an abridged
overview of its key steps. The main ingredients are \thmref{thm:H^1-fil} and
\thmref{thm:RSW-gen}, which were among the main results of \cite{KM-1}, together
with the Kato complex \eqref{eqn:Kato com} and \thmref{thm:SC-change-main}.

The strategy we use for proving
\thmref{thm:Main-1} is fundamentally different from the one employed in
characteristic zero in \cite{Bright-Newton}. The main strategy of op. cit. is an
inductive procedure in which the authors reduce to the case where the
Swan conductor of the given Brauer
class is zero; this case is then handled directly via a residue map. The
reduction of the Swan conductor is achieved by a blow-up argument using Kato's
formula for the Swan conductor along the exceptional divisor of a blow-up.

In contrast, our strategy proceeds by induction on the relative dimension of
$\sX$ over $S$, reducing to the case of relative dimension one, which we handle
directly. This reduction relies on Lefschetz-type hyperplane theorems and on a
functoriality theorem for the generalized refined Swan conductor proved in
\cite{KM-1}. We give below some details of the steps involved.

The first step is to establish the existence of a filtration
$\Fil'_\bullet H^q_\et(X, {\Q}/{\Z}(q-1))$ of $H^q_\et(X, {\Q}/{\Z}(q-1))$ for
$q \ge 1$, which we refer to as the non-logarithmic version of Kato's filtration.
This filtration is an analogue of Matsuda's filtration (cf. \cite{Matsuda}) on
$H^1_{\et}(X,\Q/\Z)$. It induces a corresponding filtration on the Brauer group,
which we denote by $\Fil'_\bullet\Br(X)$. We then compare this filtration with the
evaluation filtration; this comparison is carried out in \S~\ref{sec:MAF}.

To show that ${\Ev}_n \Br(X)$ is contained in $\Fil'_n \Br(X)$
for $n \ge 0$, we use \thmref{thm:SC-change-main} to reduce the proof to the case
when $\sX$ is a relative curve over $S$ (with notation as in \thmref{thm:Main-1})
in which case the Matsuda and Kato filtrations on $\Br(X)$ coincide and all
Brauer classes have type I along the components of the special fiber $Y$.
This coincidence is non-trivial, and is a consequence of  
Theorems~\ref{thm:H^1-fil} and ~\ref{thm:RSW-gen}, which are among the main results
of \cite{KM-1}.

The key argument now is the following.
Suppose that $\chi\in{\Ev}_n\Br(X)$ does not belong to
$\Fil'_n\Br(X)$. Then \thmref{thm:RSW-gen} implies that, after restricting $\chi$
to the local rings at enough number of closed points $x\in Y_i$ for some $i$,
its refined Swan conductor remains nonzero in 
$\Omega^1_{Y_i} \otimes_{\sO_{Y_i}} k(x)$. 
Using the Kato complex \eqref{eqn:Kato com} and \corref{cor:w-3}, we can then
find a rational function $f$ on $\sX$, and two distinct points
$P,Q\in B(P,n+1)\subset X^o_\fin$ such that
$(\partial'_P \circ \partial_P)(\{\chi, f\})
\neq  (\partial'_Q \circ \partial_Q)(\{\chi, f\})$.
It follows that ${\ev}{\chi}(P)\neq{\ev}_{\chi}(Q)$, contradicting the assumption
that $\chi\in{\Ev}_n\Br(X)$.

The proof of the inclusion $\Fil'_n\Br(X)\subset{\Ev}_n\Br(X)$ for $n\geq 0$
proceeds along similar lines. Suppose that $\chi\in\Fil'_n\Br(X)$ and that
$P,Q\in B(P,n+1)$. Since $P$ and $Q$ have a common specialization, we first
localize $\sX$ at this specialization, thereby reducing the problem to a local one.

The next natural step is to reduce further to the case of relative curves.
However, an obstacle arises here: in general, we cannot find a sufficiently nice
relative curve in $\sX$ containing the closures of both $P$ and $Q$. To
circumvent this difficulty, we construct a chain of points
$P_0, P_1, \ldots P_{d-1}, P_d \in B(P,n+1)$,
where $d$ is the relative dimension of $\sX$ over $S$, such that
$P_0=P$, $P_d=Q$, and the closures of each pair of successive points
$P_i,P_{i+1}$ lie on a sufficiently nice relative curve, possibly after further
localization. This reduces the problem to the case where $\sX$ is a relative curve.

We now proceed in the same way as in the previous case. Using corollaries of
\thmref{thm:RSW-gen}, we evaluate the boundary maps $\partial_{P'}$ and
$\partial'_{P'}$, for $P'\in \{P,Q\}$, on the cup product of $\chi$ with a suitable
rational function in the Kato complex. This allows us to conclude that
${\ev}_{\chi}(P) = {\ev}_{\chi}(Q)$. The key ingredient in this final step is a
result of Kato (cf. \cite[Thm.~5.1]{Kato-89}).

We deduce the case $n=-2$ of \thmref{thm:Main-1} using \propref{prop:Kato-Br-0},
together with a theorem from the class field theory of smooth varieties over
finite fields stating that Frobenius elements are dense in the abelianized
{\'e}tale fundamental group. Finally, we treat the case $n=-1$ by reducing to a
suitable intermediate step, after which we can follow the strategy used by
Bright--Newton in \cite{Bright-Newton} to prove the analogous result in
characteristic zero. The argument in our setting, however, is more involved
because of the greater generality covered by \thmref{thm:Main-1}.

The first part of \thmref{thm:Main-3} follows directly from \thmref{thm:Main-1}.
The more difficult second part, however, requires several additional steps. The
key result is \lemref{lem:SC-change-6}, which allows us to proceed by induction
on the dimension of $X$, thereby reducing the problem to the case where $X$ is a
curve. Geometric arguments then further reduce the problem to the case where
$X$ is geometrically integral over $k$. This latter case is
classical (cf. \cite[\S~9]{Saito-Invent}, \cite[Thm.~1.6]{KRS}).

The main work lies in proving \lemref{lem:SC-change-6}. The proof is delicate and
relies on several of the same tools used in the proof of \thmref{thm:Main-1},
namely, the generalized refined Swan conductor and its functoriality, the Kato
complex, and Bertini theorems.

To prove \thmref{thm:Main-5}, we follow the characteristic zero strategy of Kai
\cite{Kai}, although each step of the argument becomes more delicate in positive
characteristic. The proof has three main ingredients: \thmref{thm:Main-3},
\thmref{thm:Milne-4}, and Milne's positive characteristic generalization of Tate
duality for abelian varieties over $p$-adic fields.

Theorem~\ref{thm:Milne-4} allows us to establish a direct relationship between
the Brauer--Manin and Tate--Milne pairings. Our proof of this theorem is
completely different from the characteristic zero argument in \cite{Kai}. Once
this relationship is established, \thmref{thm:Main-3} reduces the proof of
\thmref{thm:Main-5} to an injectivity statement for a certain relative Brauer
group upon passing from the model to its special fiber. This injectivity is then
obtained using the smoothness of the Picard scheme, together with a classical
result of Grothendieck on the Brauer group of the model.

In \S~\ref{sec:KDRW}, we review Kato's filtration for henselian discrete
valuation fields and establish some basic results. In \S~\ref{sec:Kat-Ev*}, we
define the global Kato and evaluation filtrations on the Brauer group and prove
the key result \propref{prop:Kato-Br-0}. In \S~\ref{sec:RSCKC}, we recall the
generalized refined Swan conductor and define the non-logarithmic version of
Kato's filtration. In \S~\ref{sec:Kato-complex}, we give several applications of
the Kato complex, which we then use in \S~\ref{sec:Esp} and \S~\ref{sec:Esp-0} to
prove a specialization theorem for the Swan conductor. 

We complete the proof of \thmref{thm:Main-1} in \S~\ref{sec:KEF} and prove the
comparison between our evaluation filtration and the filtration of
\cite{Bright-Newton} in \S~\ref{sec:BN**}. We begin the proof of
\thmref{thm:Main-3} in \S~\ref{sec:BMP} and complete it in \S~\ref{sec:BMP-2}.
The latter section also contains the proofs of \thmref{thm:Main-2} and
\thmref{thm:Main-6}. Finally, we begin the proof of \thmref{thm:Main-4} in
\S~\ref{sec:Kai} and complete it in \S~\ref{sec:Kai-2}.
Section~\ref{sec:Bertini}
contains the proof of the Bertini theorem, which is a key ingredient in the
proofs of our main results.

\subsection{Notation}\label{sec:Notn}
We shall follow the following common notation throughout this paper.

We fix a prime $p$ for all of this paper and shall assume all schemes to be Noetherian
and separated $\F_p$-schemes. For a field $k$, we shall let $k_s$
(resp. $\ov{k}$)
denote a chosen separable (resp. algebraic) closure of $k$ such that
$k_s \subset \ov{k}$.
If $A \surj B$ is a surjective ring homomorphism and $a \in A$, we shall usually
denote the image of $a$ in $B$ by $\ov{a}$ if it creates no confusion in a
given context. For any prime ideal $\fp$ of $A$, let $A_{\fp}$ denote the localization of $A$ at $\fp$.

For a commutative ring $A$, an $A$-scheme will mean a separated scheme over $\Spec(A)$
and a subscheme will mean a locally closed subscheme.
We shall let $\Sch_A$ denote the category of separated Noetherian $A$-schemes.
The product $X \times_{\Spec(A)} Y$ will be written as $X \times_A Y$.
For a ring homomorphism $A \to A'$ and $X \in \Sch_A$, we shall write
$X \times_{\Spec(A)} \Spec(A')$ as $X_{A'}$.
For a scheme $X$, we let $X^{(q)}$ (resp. $X_{(q)}$) denote the set of points on $X$
having codimension (resp. dimension) $q$. We let $X_\reg$ (resp. $X_\sing$)
denote the regular (resp. singular) locus of $X$. All intersections of subschemes
of a scheme in this paper will be assumed to be scheme-theoretic.
For a point $x \in X$, we let
$k(x)$ denote the residue field of $x$. If $X$ is reduced, we shall let $X^N$ denote
the normalization of $X$.
We let $\Sch_X$ denote the category of separated schemes over $X$.

Given a commutative ring $A$, we let $Q(A)$ denote the total ring of fractions of
$A$. We let $\Sch_{A/\zar}$ (resp. $\Sch_{A/\nis}$, resp. $\Sch_{A/ \etl})$
denote the Zariski (resp. Nisnevich, resp. {\'e}tale) site of $\Sch_{A}$.
We let $\epsilon \colon
\Sch_{A/ \et} \to \Sch_{A/ \nis}$ denote the canonical morphism of sites.
If $\sF$ is a sheaf on $\Sch_{A/\nis}$, 
we shall write $\epsilon^* \sF$ also as $\sF$ as long as the usage of
the {\'e}tale topology is clear in a context.
For a local ring $A$, we let $A^h$ (resp. $A^{sh}$, resp. $\wh{A}$) denote the
henselization (resp. strict henselization, resp. completion) of $A$ with respect to
its maximal ideal.

For an abelian group $A$, we shall write $\Tor^1_{\Z}(A, {\Z}/n)$ as
$A[n]$ and $A/{nA}$ as $A/n$.
The tensor product $A \otimes_{\Z} B$ will be written as $A \otimes B$.
We shall let $A\{p\}$ (resp. $A\{p'\}$) denote the
subgroup of elements of $A$ which are annihilated by a power of (resp.
an integer prime to) $p$. If $\sF$ is a sheaf on some site of a scheme, we shall write
the kernel (resp. cokernel) of the map $\sF \xrightarrow{n} \sF$ as $_n\sF$
(resp. ${\sF}/n$).

We shall let $\N_0$ (resp. $\N$) denote the set of non-negative (resp.
positive) integers. For $n \in \N_0$, we let $J^0_n = \{0, \ldots , n\}$ and
$J_n = J^0_n \cap \N$. For a prime $p$, the symbol $I_p$ will denote the set of all
positive integers prime to $p$. If $F$ is a functor (e.g., a cohomology group)
from any category of schemes to
another category, and $X = \Spec(A)$ is an affine scheme, we shall usually write
$F(X)$ as $F(A)$.

\section{Kato's filtration for discrete valued fields}\label{sec:KDRW}
In this section, we recall Kato's ramification filtration for the
{\'e}tale cohomology of ${\Q}/{\Z}(n)$ on henselian discrete valuation fields
(which we shall refer to as hdvf) of positive characteristic. We establish
several facts concerning this filtration, most of which follow from \cite{Kato-89}.
We fix once and for all a prime $p > 0$ and let $\F_p$ denote the prime field of
characteristic $p$. We begin by recalling the basic objects.

\subsection{de Rham-Witt forms, Milnor $K$-theory and and the Brauer group}
\label{sec:DRW-Br}
Given an  $\F_p$-scheme $X$, we let
$W_*\Omega^\bullet_X: = \{W_m\Omega^q_X\}_{q \ge 0, m \ge 1}$
denote the $p$-typical pro-de Rham-Witt complex of $X$ (cf. \cite{Illusie}).
This is a projective system (indexed by positive integers $m$) of Zariski sheaves of
commutative differential graded algebras on $X$
equipped with the operators $F \colon W_m\Omega^q_X \to W_{m-1}\Omega^q_X$
(called Frobenius), $V \colon W_m\Omega^q_X \to W_{m+1}\Omega^q_X$
(called Verschiebung) and $R \colon W_m\Omega^q_X \to W_{m-1}\Omega^q_X$
(called restriction) satisfying a set of properties (cf.
\cite[\S~2.4]{KP-Comp}). We let $d \colon  W_m\Omega^q_X \to W_{m}\Omega^{q+1}_X$
denote the differential of the complex $W_m\Omega^\bullet_X$ and let
$[.]_m  \colon \sO_X \to W_m\sO_X$ denote the Teichm{\"u}llar map.

We let $Z_1W_m\Omega^q_X = F(W_{m+1}\Omega^q_X) = \Ker(W_m\Omega^q_X
\xrightarrow{F^{m-1}d} \Omega^{q+1}_X)$ and let $C \colon Z_1W_m\Omega^q_X \to
W_m\Omega^q_X$ denote the Cartier homomorphism. Note that $Z_1W_m\Omega^q_X$
is a quasi-coherent sheaf of $W_{m+1}\sO_X$-modules via $F \colon W_{m+1}\sO_X \to
W_m\sO_X$ (cf. \cite[Cor.~5.10]{KM-1}).
For all the standard properties of the de Rham-Witt complex that we shall use in this
paper, we refer the reader to \cite{Illusie}. 

For a local $\F_p$-algebra $R$ with residue field $k$, we let
$K^M_*(R)$ denote the (improved) Milnor $K$-theory of $R$, introduced in
\cite{Kerz-JAG}. This coincides with the classically defined (e.g., in \cite{Kato-86})
Milnor  $K$-theory of $R$ if $k$ is infinite. We let
$m^r_R \colon K^M_r(R) \surj K^M_r(k)$ denote the canonical restriction map.
If $R$ is regular with quotient field $K$, then $K^M_r(R)$ is a subgroup
of $K^M_r(K)$ under the canonical restriction by \cite{Kerz-JAG} and we shall treat
elements of $K^M_r(R)$ also as elements of $K^M_r(K)$. 
We let $X \mapsto \sK^M_{r, X}$ denote the $r$-th Milnor $K$-theory sheaf on the big
{\'e}tale (or Zariski) site of $\F_p$-schemes. 

For $X \in \Sch_{\F_p}$ and $q, m \ge 1$, we let $W_m\Omega^q_{X, \log}$ be the image of
the map of {\'e}tale sheaves $\dlog \colon (\sO^{\times}_X)^q \to W_m\Omega^q_X$, given
by $\dlog(\{x_1, \ldots , x_q\}) = \dlog[x_1]_m \wedge \cdots
\wedge \dlog[x_q]_m$ (cf. \cite[Chap.~I, \S~5.7]{Illusie}).
This map uniquely factors through the composite quotient map
$(\sO^{\times}_X)^q \surj \sK^M_{q, X} \surj  {\sK^M_{q,X}}/{p^m}$
(cf. \cite[\S~1.2]{Morrow-ENS}, \cite[Lem.~3.2.8]{Zhao}).
We let $\dlog \colon {\sK^M_{q,X}}/{p^m} \to W_m\Omega^q_{X, \log}$ denote the induced map.
This latter map is bijective if $X$ is regular (cf. \cite[Thm.~5.1]{Morrow-ENS}).
We let $ W_m\Omega^q_{X, \log} = 0$ if either $q < 0$ or $m \le 0$.

For any $\F_p$-scheme $X$, we let $\Br(X)$ denote the cohomological
Brauer group $H^2_\et(X, \sO^\times_X)$ of $X$. There is a natural injective map
$\Br_{\rm Az}(X) \inj \Br(X)_\tor$ which is an isomorphism if $X$ admits an ample
invertible sheaf (e.g., $X$ is quasi-projective over an affine scheme), where
$\Br_{\rm Az}(X)$ is the Azumaya Brauer group of $X$ (cf. \cite[Thm.~4.2.1]{CTS}).
If $X$ is regular, then $\Br(X)$ is a torsion group and the restriction map
$\Br(X) \to \Br(U)$ is injective for any open $U \subset X$
(cf. \cite[Lem.~3.5.3, Thm.~3.5.7]{CTS}). In particular, the natural map
$\Br_{\rm Az}(X) \to \Br(X)$ is an isomorphism.
These facts about the Brauer group will be used without any reference throughout
this paper.

\subsection{The cohomology group $H^*(X)$}\label{sec:p-adic-coh}
We let ${\Zars}_{\F_p}$ (resp. ${\Et}_{\F_p}$) denote the big Zariski (resp. {\'e}tale)
site of the category of $\F_p$-schemes.
For $X \in \Sch_{\F_p}$, we let $\sD_{\et}(X)$ denote the bounded derived category of
{\'e}tale sheaves on $X$.
Given any $n \in I_p$ and $q \in \Z$, we let ${\Z}/n(q)$ denote $q$-th power of
the classical {\'e}tale sheaf $\mu_n$ of $n$-th roots of unity on
${\Et}_{\F_p}$ (cf. \cite[p.~163]{Milne-EC}).

For $q \in \Z$ and $n = p^mr \in \N$ with $r \in I_p$, we let
${\Z}/n(q) := {\Z}/r(q) \oplus W_m\Omega^q_{(-,\log)}[-q]$
and consider it as a complex of sheaves on ${\Et}_{\F_p}$.
For $X \in \Sch_{\F_p}$, we let $H^i_\et(X, \Z/n(j)) := \H^i_\et(X, \Z/n(j))$.
If $r' \in \N$ such that $r \mid r'$ and $r' \in I_p$, we have a natural map
$\Z/r(q) \to \Z/r'(q)$, induced by the maps ${\Z}/r \cong \frac{1}{r}\Z/\Z \inj
\frac{1}{r'}\Z/\Z \cong \Z/r'$ and $\mu_r \inj \mu_{r'}$. On the other hand, the
map $\ov{p} \colon W_m\Omega^q_{(-,\log)} \to W_{m+1}\Omega^q_{(-,\log)}$
(cf. \cite[Lem.~6.14]{KM-1}) induces a
canonical map $\Z/p^m(q) \to \Z/p^{m+1}(q)$. As a result, we have canonical maps
$\Z/n(q) \to \Z/n'(q)$, whenever $n \mid n'$. In particular,
$\{\Z/n(q)\}_{n \ge 1}$ form an ind-object in the category of complexes of
sheaves on ${\Et}_{\F_p}$ for any given $q \ge 0$. 
We let $H^i_\et(X, {\Q}/{\Z}(j)) = \varinjlim_n H^i_\et(X, \Z/n(j))$.
The cohomology groups of interest in this paper will be of the following kind.

\begin{defn}\label{defn:Coh-sheaf}
  For $X \in \Sch_{\F_p}$ and $q \in \N$, we let $H^q(X) =  H^q_\et(X, {\Q}/{\Z}(q-1))$
  so that $H^q(X)\{p\} = H^q_\et(X, {\Q_p}/{\Z_p}(q-1))$.
For $n \in \N$, we let $H^q_n(X) = H^q_\et(X, {\Z}/{n}(q-1))$.
We shall often use the notation $H^q_\infty(X)$ for $H^q(X)$.
\end{defn}

By \cite[Lem.~7.1]{KM-1}, the canonical map $H^q_{p^m}(R) \to H^q(R)$ is injective and
identifies $H^q_{p^m}(R)$ as $H^q(R)[p^m]$ if $R$ is a regular local $\F_p$-algebra.
This holds for the map $H^q_n(R) \to H^q(R)$ for every $n \ge 1$ if $q \le 2$.
There is a cup product pairing
\begin{equation}\label{eqn:Pairing-1}
    \{ \ , \ \} \colon H^q_n(R) \otimes K_r^M(R) \to H^{q+r}_n(R);
\end{equation}
\[
\chi \otimes \{a_1,\ldots,a_r\} \mapsto \{\chi, a_1,\ldots,a_r\} =
\chi \cup \beta^r_n(R)(\{a_1,\ldots,a_r\})
\]
for every $n \in \N$, where $\beta^q_n({R}) \colon K^M_q(R) \to H^q_\et(R, {\Z}/{n}(q))$
is the Norm-residue (also called the Bloch-Kato) map.

If $R$ is a regular $\F_p$-algebra, the exact sequence
\begin{equation}\label{eqn:Milnor-0}
  0 \to W_m\Omega^q_{{(-)}, \log} \to Z_1W_m\Omega^q_{(-)} \xrightarrow{1-C}
    W_m\Omega^q_{(-)} \to 0
\end{equation}
of sheaves on $\Spec(R)_\et$, and quasi-coherence of $Z_1W_m\Omega^q_{(-)}$
yield an exact sequence
\begin{equation}\label{eqn:Milnor-0.1}
  0 \to W_m\Omega^q_{R, \log} \to Z_1W_m\Omega^q_R \xrightarrow{1-C} W_m\Omega^q_R
  \xrightarrow{\delta^q_m} H^{q+1}_{p^m}(R) \to 0.
  \end{equation}

If $R$ is furthermore local, we
also have a commutative diagram (cf. \cite[\S~1.3]{Kato-89})
\begin{equation}\label{eqn:Milnor-1}
  \xymatrix@C1pc{
    {K^M_q(R)}/{p^m} \otimes W_m\Omega^{q'}_{R} \ar[r]^-{\wedge}
    \ar[d]_-{\id \otimes \delta^{q}_m} \ar[dr]^-{\gamma^{q+q'}_m} &
    W_m\Omega^{q+q'}_{R} \ar[d]^-{\delta^{q+q'}_{m}} \\
     {K^M_q(R)}/{p^m} \otimes H^{q'+1}_{p^m}(R) \ar[r]^-{\{ \ , \ \}} &
     H^{q+q'+1}_{p^m}(R),}
\end{equation}
where the top horizontal arrow is obtained by taking $\dlog \otimes \id$
and pre-composing it with the wedge product. We shall denote the composite map
$W_m\Omega^q_R \xrightarrow{\delta^q_m} H^{q+1}_{p^m}(R) \inj H^{q+1}(R)$
also by $\delta^q_m$.

\begin{exm}\label{exm:Brauer-coh}
  Let $X$ be a normal  $\F_p$-scheme.
  Using the Kummer sequence (cf. \cite[Chap.~II, Exm.~2.18(b)]{Milne-EC}),
and the exact sequence
\begin{equation}\label{eqn:Omega-1}
0 \to \sO^{\times}_X \xrightarrow{p^m}  \sO^{\times}_X \xrightarrow{\dlog}
W_m\Omega^1_{X, \log} \to 0
\end{equation}
of {\'e}tale sheaves on $X$, one gets that
${\Z}/n(1) \cong (\sO^{\times}_X \xrightarrow{n}  \sO^{\times}_X)[-1]$ in $\sD_\et(X)$
for every integer $n \in \N$. In particular, there is a canonical exact sequence
\begin{equation}\label{eqn:Kummer}
  0 \to {\Pic(X)}/n \to H^2_\et(X, {\Z}/n(1)) \to \Br(X)[n] \to 0.
  \end{equation}
It follows that there are canonical isomorphisms
$H^2_\et(R, {\Z}/n(1)) \xrightarrow{\cong} \Br(R)[n]$ for all integers $n \in \N$
and $H^2(R) \cong \Br(R)$ if $R$ is a normal $\F_p$-algebra which is either a
unique factorization domain (ufd) or a semi-local ring.
\end{exm}

\vskip .2cm

\subsection{The canonical lifting}\label{sec:Lifting}
Let $A$ be a henselian discrete valuation ring
(which we shall refer to as hdvr) containing $\F_p$ and 
let $\fm = (\pi)$ denote the maximal ideal of $A$. Let $\ff = {A}/{\fm}$ and
$K = Q(A)$. Let $n, q \in \N$. By \cite[\S~1.4]{Kato-89},
there is a homomorphism
\begin{equation}\label{eqn:Can-lift}
  l^q_n(A) \colon H^q_n(\ff) \to H^q_n(K),
\end{equation}
called the canonical lifting map. This is defined by letting
$l^q_n(A) = s^q_n(A) \circ t^q_n(A)^{-1}$, where
$s^q_n(A) \colon H^q_n(A) \to H^q_n(K)$ and
$t^q_n(A) \colon H^q_n(A) \to H^q_n(\ff)$ are the canonical restriction maps
and the latter is an isomorphism by \lemref{lem:Coh-iso}. It is clear that
$l^q_n(A)$ does not depend on the choice of generators for $\fm$. We
let $l^q(A) := {\varinjlim}_n l^q_n(A) \colon H^q(\ff) \to H^q(K)$.

It is clear from the definition of $l^q_n(A)$ that
for any $\chi \in H^q_n(\ff)$ and $x \in K^M_r(A)$, one has
\begin{equation}\label{eqn:Pairing-2}
  \begin{array}{lllll}
l^{q+r}_n(\{\chi,  m^r_A(x)\}) & = & 
l^{q+r}_n \circ t^{q+r}_n(\{\chi', x\}) & = & s^{q+r}_n(A)(\{\chi', x\}) \\
& = & \{s^{q}_n(A)(\chi'), x\} & = & \{l^q_n(\chi), x\}, \\
  \end{array}
  \end{equation}
where we let $\chi' = t^q_n(A)^{-1}(\chi)$.

Kato defined $l^q_n(A)$ in \cite[\S~1.4]{Kato-89} without assuming that $A$ is a
discrete valuation ring but ~\eqref{eqn:Pairing-2} shows that the
two definitions are identical in our case.
Since $\fm \subset (1-C)(K^p)$, one checks that lifting an element of
$\ff$ to $A$ and mapping it to $K$ defines a canonical map
$\wt{l}^1_p(A) \colon \frac{\ff}{(1-C)(\ff^p)} \to \frac{K}{(1-C)(K^p)}$
such that the diagram
\begin{equation}\label{eqn:delta-lift}
  \xymatrix@C1pc{
    \frac{\ff}{(1-C)(\ff^p)} \ar[rr]^-{\wt{l}^1_p(A)} \ar[d]_-{\delta^0_p} &&
    \frac{K}{(1-C)(K^p)} \ar[d]^-{\delta^0_p}
      \\
      H^1_p(\ff) \ar[rr]^-{l^1_p(A)} && H^1_p(K)}
\end{equation}
is commutative.

For $n, q \in \N$, we let
\begin{equation}\label{eqn:Can-lift-sum}
\lambda^q_n(\pi) \colon H^q_n(\ff)  \bigoplus  H^{q-1}_n(\ff) \to H^q_n(K)
\end{equation}
be defined by $\lambda^q_n(\pi)(\chi, \chi') = l^q_n(\chi) + \{l^{q-1}_n(\chi'), \pi\}$.
As this map depends on $\pi$, it is in general not functorial in $A$ but it is
functorial for unramified extensions of $A$.
We shall let $\wt{\lambda}^q_{p^m}(\pi)$ denote the composite map
\begin{equation}\label{eqn:Can-lift-sum-0}
W_m\Omega^{q-1}_\ff \bigoplus W_m\Omega^{q-2}_\ff \xrightarrow{\delta^{q-1}_m \oplus
  \delta^{q-2}_m}
H^q_{p^m}(\ff)  \bigoplus  H^{q-1}_{p^m}(\ff) \xrightarrow{\lambda^q_{p^m}(\pi)} H^q_{p^m}(K).
\end{equation}
Note that the first arrow here is surjective and the second arrow is injective
(cf. \cite[\S~1.9]{Kato-89}).
We let $\lambda^q(\pi) = {\varinjlim}_{n} \lambda^q_n(\pi) \colon
H^q(\ff)  \bigoplus  H^{q-1}(\ff) \to H^q(K)$.

\begin{defn}\label{defn:tame-ram}
  Let $A$ be as above. An element of $H^q(K)$ will be called unramified (resp.
  tamely ramified) if it lies in the image of $l^q(A)$ (resp. $\lambda^q(\pi)$).
  \end{defn}

Note that for $q =1$, the above definition agrees with the classically defined notion
of unramified and tamely ramified characters of {\'e}tale fundamental groups (cf.
\cite[Thm.~3]{Kato-Gal-coh}, \cite[(2.9)]{GKR}).

\begin{lem}\label{lem:Coh-iso}
Let $A$ be an hdvr containing $\F_p$ and 
let $\fm = (\pi)$ denote the maximal ideal of $A$. Let $\ff = {A}/{\fm}$ and
$K = Q(A)$. Let $n, q \in \N$. Then we have the following.
\begin{enumerate}
\item
  $\Ker(\Omega^q_A \surj \Omega^q_\ff)$ is contained in $(1-C)(Z_1\Omega^q_A)$.
  \item
The restriction maps $t^q_n(A) \colon H^q_n(A) \to H^q_n(\ff)$ and
$s^q_n(A) \colon H^q_n(A) \to H^q_n(K)$ are respectively, 
bijective and injective.
\end{enumerate}
\end{lem}
\begin{proof}
To prove item (1),  we first note (cf. \cite[Lem.~9.8(2)]{KM-1}) that
  $\Ker(\Omega^q_A \surj \Omega^q_\ff) = \pi\Omega^q_A(\log(\pi))$, where
  $\Omega^q_A(\log(\pi))$ is the module of $q$-th differential forms with
  log poles along $V((\pi))$ (cf. \S~\ref{sec:Fil-DRW}).
  Next, we note that $\pi\Omega^q_A(\log(\pi))$ is generated as an abelian group
by elements of the form
$a\pi\dlog (x_1) \wedge \cdots \wedge \dlog (x_q)$ with $a \in A$ and
$x_i \in A^\times \bigcup \{\pi\}$. Applying the
Hensel lemma to the polynomial $X^p-X -a\pi \in A[X]$, we see that there exists
$b \in (\pi)$ such that $b^p-b = a\pi$ which means that $(1-C)(b^p) = a\pi$
(cf. \cite[Claim~1.2.2]{JSZ}). But this yields
$(1-C)(b^p\dlog (x_1) \wedge \cdots \wedge \dlog (x_q)) =
a\pi\dlog (x_1) \wedge \cdots \wedge \dlog (x_q)$.
Since $b \in (\pi)$, one checks that
$b^p\dlog (x_1) \wedge \cdots \wedge \dlog (x_q) \in Z_1\Omega^q_K \bigcap \Omega^q_A
= Z_1\Omega^q_A$. This concludes the proof of (1).

We now prove (2).
  If $n \in I_p$, this is a special case of the well-known rigidity
  theorem and Bloch-Ogus injectivity theorem in {\'e}tale cohomology.
  We can therefore assume $n = p^m$.
We now use the sheaf exact sequence
  \begin{equation}\label{eqn:Coh-iso-0}
  0 \to \Omega^{q-1}_{X, \log} \xrightarrow{(\un{p})^m} W_{m+1}\Omega^{q-1}_{X, \log}
  \xrightarrow{R} W_{m}\Omega^{q-1}_{X, \log} \to 0
  \end{equation}
  for any regular $\F_p$-scheme $X$ (cf. \cite[Lem.~3.8]{GK-Duality})
  and the surjection
  $R \colon W_{m+1}\Omega^q_{S, \log} \surj W_{m}\Omega^q_{S, \log}$ to get a short exact
  sequence
  \begin{equation}\label{eqn:Coh-iso-1}
  0 \to H^q_{p}(S) \xrightarrow{(\un{p})^{m}} H^q_{p^{m+1}}(S) 
  \xrightarrow{R} H^q_{p^m}(S) \to 0
  \end{equation}
  for $S \in \{A, K, \ff\}$, where the surjectivity of $R$ follows from
  the proper base change theorem for $A$ and from the fact that the
  cohomological dimension for $p$-torsion sheaves on a characteristic $p$ field is
  one.
Comparing this exact sequence for $A$ and $\ff$ and using induction on $m$,
it suffices to show (for proving bijectivity of $t^q_n(A)$)
that $H^q_{p}(A) \to H^q_{p}(\ff)$ is an isomorphism.

To that end, we look at the commutative diagram of exact sequences
\begin{equation}\label{eqn:Coh-iso-2}
  \xymatrix@C1pc{
    Z_1\Omega^q_A \ar[rr]^-{1-C} \ar@{->>}[d] && \Omega^q_A \ar[rr]^-{\delta^1_1}
    \ar@{->>}[d] && H^q_{p}(A)  \ar[d]^-{t^q_p(A)} \ar[r] & 0 \\
   Z_1\Omega^q_\ff \ar[rr]^-{1-C} && \Omega^q_\ff \ar[rr]^-{\delta^1_1}
   && H^q_{p}(\ff)  \ar[r] & 0,}
\end{equation}
where the vertical arrows are the restriction maps.

The surjectivity of $t^q_p(A)$ follows directly from the above diagram.
Furthermore, a diagram chase reduces the proof of injectivity to showing that
$\Ker(\Omega^q_A \surj \Omega^q_\ff)$ is contained in
$(1-C)(Z_1\Omega^q_A)$. But this follows from item (1).
The injectivity of $s^q_n(A)$ follows
by combining the bijectivity of $t^q_n(A)$ with \cite[Thm.~3(1)]{Kato-Gal-coh}.
\end{proof}

\subsection{Kato filtration for hdvf}\label{sec:Kato-fil}
We now recall Kato's ramification filtration on the cohomology of an hdvf from
\cite[\S~2]{Kato-89}. We let $A$ be an hdvr containing $\F_p$ and 
let $\fm = (\pi)$ denote the maximal ideal of $A$. Let $\ff = {A}/{\fm}$ and
$K = Q(A)$.
For any $A$-algebra $R$, let $R^h$ denote the henselization of $R$ with respect to the
ideal $\fm R$. We fix integers $n, q \in \N_0$ and $m \in \N$. We let
$\Spec({R}/{\fm R}) \xrightarrow{\iota} \Spec(R) \xleftarrow{j} \Spec(R[\pi^{-1}])$
be the inclusions and let
\begin{equation}\label{eqn:Kato-0}
V^q_n(R) = \H^q_\et({R}/{\fm R}, \iota^*{\bf R}j_*({\Z}/{n}(q-1)))
\xrightarrow{\cong} V^q_n(R^h) \xrightarrow{\cong} H^q_\et(R^h[\pi^{-1}], {\Z}/{n}(q-1))
\end{equation}
\[
\mbox{and} \ V^q(R) = {\varinjlim}_{n} V^q_n(R) \cong V^q(R^h) \cong H^q(R^h[\pi^{-1}]).
\]
The following definition is due to Kato \cite[Defn.~2.1]{Kato-89}. 

\begin{defn}\label{defn:Kato-1}
For $n \in \N_0$ and $q \in \N$, we let $\Fil_nH^{q}(K)$ be the set of all elements
$\chi \in H^{q}(K)$ such that $\{\chi, 1 + \pi^{n+1}T\} :=
\sigma^q_{n}(\chi \otimes (1 + \pi^{n+1}T)) = 0$ in $V^{q+1}(A[T])$
under the composite map $\sigma^q_{n}$:
\[
H^{q}(K)\otimes ((A[T])^h[\pi^{-1}])^{\times}  \to
H^{q}(K)\otimes H^1_\et((A[T])^h[\pi^{-1}], {\Q}/{\Z}(1)) \xrightarrow{\cup}
V^{q+1}(A[T]).
\]
For $m \in \N$, we let $\Fil_nH^{q}_{p^m}(K) =
H^{q}_{p^m}(K) \bigcap \Fil_nH^{q}(K)$.
For $q \le 2$, we let $\Fil_nH^{q}_m(K) = H^{q}_{m}(K) \bigcap \Fil_nH^{q}(K)$.
\end{defn}

By \cite[Prop.~6.3]{Kato-89} and Example~\ref{exm:Brauer-coh}, the above definition of
$\Fil_nH^{q}(K)$ coincides with the one given in \S~\ref{sec:BN-thm} for $q = 2$.

The following result of Kato summarizes the key
properties of $\Fil_nH^{q}(K)$ and the map $\lambda^q(\pi)$ that we shall frequently
use.

\begin{prop}\label{prop:Kato-basic}
  For $m, q \in \N$ and $n \in \N_0$, we have the following.
  \begin{enumerate}
  \item
    The restriction $H^q(A) \to H^q(K)$ induces an inclusion
    $H^q(A) \subset \Fil_0H^{q}(K)$.
  \item
    $\Fil_nH^{q}(K) = H^q(K)\{p'\} \bigoplus \Fil_nH^{q}(K)\{p\}$.
  \item
   $\lambda^q_m(\pi)$ is injective. 
  \item
    $\lambda^q(\pi)$ induces isomorphisms
    \[
    \begin{array}{lll}
      H^q(\ff)  \bigoplus  H^{q-1}(\ff) & \xrightarrow{\cong}  & \Fil_0H^{q}(K)
      \xrightarrow{\cong}  \Ker(H^q(K) \to H^q(K_{nr})) \\
    & \xrightarrow{\cong} & H^q(K)\{p'\} \bigoplus
      \Ker(H^q(K)\{p\} \to H^q(K_{nr})\{p\}) \\
      & \xrightarrow{\cong} & \Ker(H^q(K) \to H^q(K_{tr})),
    \end{array}
    \]
    where $K_{nr}$ (resp. $K_{tr}$) is the maximal unramified (resp. tamely ramified)
    extension of $K$. 
  \item
    For $q \le 2$, $\lambda^q_m(\pi)$ induces an isomorphism
    $\lambda^q_m(\pi) \colon H^q_m(\ff) \bigoplus H^{q-1}_m(\ff) \xrightarrow{\cong}
    \Fil_0H^{q}_m(K)$.
  \item
    If  $[\ff : \ff^p] = p^c$ and $q > c+1$, then
    $\Fil_0H^{q}(K)  = H^{q}(K)$.
  \item
    If there is a chain of henselian discrete valuation fields
    $\ff = \ff_0, \ldots , \ff_c$ such that
    $\ff_{i+1}$ is the residue field of $\ff_i$ for $0 \le i \le c-1$ and
    $\ff_c$ is a perfect field, then $\Fil_0H^{q}(K)  = H^{q}(K)$ for
    $q > c+1$.
  \item
    $H^{q}(K) = \bigcup_{n \ge 0} \Fil_nH^{q}(K)$. For
    $q \le 2$, $H^{q}_m(K) = \bigcup_{n \ge 0} \Fil_nH^{q}_m(K)$.
  \item
    The product map of ~\eqref{eqn:Pairing-1} for $K$ induces 
    $\Fil_nH^{q}(K) \otimes K^M_r(K) \to \Fil_nH^{q+r}(K)$ such that
    the map $\Fil_nH^{q}(K)\{p\} \otimes K^M_r(K) \to
    \Fil_nH^{q+r}(K)\{p\}$ is surjective. 
  \end{enumerate}
  \end{prop}
\begin{proof}
Item (1) follows from \lemref{lem:Coh-iso} and
\cite[Thm.~3.2]{Kato-89} (see also \cite[Thm.~7.15]{KM-1}), (2) follows
from Cor.~2.5 of op. cit., (3) follows from \cite[Thm.~3]{Kato-Gal-coh},
(5) and (6) and the first three isomorphisms of (4) follow from Prop.~6.1 of op. cit.,
(8) follows from Lem.~2.2 of op. cit. and (7) follows by combining Lem.~7.2 of op. cit.
and (6). We now prove the last isomorphism of (4).

We first note that
the absolute Galois group $\Gal(K_{tr})$ is a pro-$p$ group and therefore
$H^q(K_{tr})\{p'\} = 0$ by a standard Galois cohomology argument
(cf. \cite[Prop.~4.4.3]{Lei-Fu}). On the other hand, if $L/{K_{nr}}$
is a tamely ramified finite extension, then $[L:K_{nr}]$ must be prime to $p$
which implies using the restriction-corestriction maps that $H^q(K_{nr})\{p\} \to
H^q(L)\{p\}$ is injective. Taking the limit, we get
$H^q(K_{nr}) \to H^q(K_{tr})$ is injective. This finishes the proof of (4).

We now prove (9). That the image of
$\Fil_nH^{q}(K)\{p'\} \otimes K^M_r(K)$ lies in $\Fil_nH^{q+r}(K)\{p'\}$
follows directly from ~\eqref{eqn:Pairing-1} and item (2).
The $p$-part of the claim for $q =1$ follows from \cite[Thm.~3.2(2)]{Kato-89}.
For $q \ge 2$, we consider the diagram
\[
\xymatrix@C.8pc{
  \Fil_nH^{1}(K)\{p\} \otimes K^M_{q-1}(K) \otimes K^M_{r}(K)
  \ar@{->>}[r] \ar@{->>}[d]_-{\{\ , \ \} \otimes id} &
  \Fil_nH^{1}(K)\{p\} \otimes K^M_{q+r-1}(K) \ar@{->>}[r] &
  \Fil_nH^{q+r}(K)\{p\} \ar@{^{(}->}[d] \\
  \Fil_nH^{q}(K)\{p\} \otimes  K^M_{r}(K) \ar[r] & H^{q}(K)\{p\} \otimes  K^M_{r}(K)
  \ar[r] & H^{q+r}(K).}
\]

It is easy to see that this diagram is commutative and the left horizontal arrow
on the top row is surjective. The right horizontal arrow on the top row and
the left vertical arrow are surjective by \cite[Thm.~3.2(2)]{Kato-89}.
A diagram chase then proves the claim and completes the proof of the proposition.
\end{proof}

Item (2) of \propref{prop:Kato-basic} allows one to define $\Fil_nH^{q}_m(K)$
for every pair of integers $m \in \N, n \in \N_0$ as follows.
Write $m = p^tr$, where $t \in \N_0$ and $r \in I_p$. We let
\begin{equation}\label{eqn:Fil-prime-to-p}
  \Fil_nH^{q}_m(K) = H^q_r(K) \bigoplus \Fil_nH^{q}_{p^t}(K).
\end{equation}
This coincides with the earlier definition (cf. Definition~\ref{defn:Kato-1})
when $q \le 2$.
Also, it easily follows with this definition that $\Fil_nH^{q}(K) =
{\varinjlim}_m \Fil_nH^{q}_m(K)$ for every $n \in \N_0$.
In view of item (8) of \propref{prop:Kato-basic}, we define the Swan conductor
$\Sw(\chi)$ of any $\chi \in H^q(K)$ to be the smallest integer $n \ge 0$ such that
$\chi \in \Fil_n H^q(K)$.

\vskip .2cm

Using \lemref{lem:Coh-iso} and \propref{prop:Kato-basic}, we
get the following.

\begin{cor}\label{cor:Kato-basic-0}
For $A$ as in \propref{prop:Kato-basic} and $m \in \N$, there is a split exact sequence
\begin{equation}\label{eqn:Fil-0-seq}
  0 \to H^q_m(A) \to \Fil_0H^{q}_m(K)
  \xrightarrow{\kappa^q_m(A)} H^{q-1}_m(\ff) \to 0,
\end{equation}
and $\kappa^q_m(A)$ has the property that $\kappa^q_m(A) (\{l^{q-1}_m(A)(\chi), \pi\}) =
\chi$ for all $\chi \in H^{q-1}_m(\ff)$.
In this exact sequence, we can replace $\Fil_0H^{q}_m(K)$
 by $H^{q}_m(K)$ if $[\ff : \ff^p] = p^c$ and $q > c+1$.
\end{cor}

Since ~\eqref{eqn:Fil-0-seq} is compatible with the
change in values of $m \in \N$, we can pass to the limit which yields a split exact
sequence
\begin{equation}\label{eqn:Fil-0-seq-0}
  0 \to H^q(A) \to \Fil_0H^{q}(K) \xrightarrow{\kappa^q(A)} H^{q-1}(\ff) \to 0.
\end{equation}

\begin{exm}\label{exm:Kato-basic-2}
  Suppose that $A$ in \propref{prop:Kato-basic} is such that $\ff$ is finite.
  We then have $\Fil_0H^{2}(K) = \Br(K)$ and  $\Br(A)\cong \Br(\ff) = 0$.
  We therefore get an isomorphism $\kappa^2(A) \colon \Br(K) \xrightarrow{\cong}
  H^1(\ff)$.
 Using the isomorphism $H^2(A) \bigoplus H^1(\ff) \xrightarrow{\cong} H^2(K)$ and
  comparing the description of $\kappa^2(A)$ given in \corref{cor:Kato-basic-0}
  with the identity \cite[(1.28)]{CTS}, one deduces that $\kappa^2(A)$ coincides
  with the classical Witt-residue map $r_W \colon \Br(K) \to H^1(\ff)$ 
  (cf. Defn.~1.4.11 of op. cit.). The evaluation at the Frobenius element yields
  an isomorphism ${\rm ev}_F(\ff) \colon  H^1(\ff) \xrightarrow{\cong} {\Q}/{\Z}$
  such that the composition
  ${\rm ev}_F(\ff) \circ r_W = {\rm ev}_F(\ff) \circ \kappa^2(A)$
  is the classical Hasse invariant map
  $\inv \colon \Br(K) \xrightarrow{\cong} {\Q}/{\Z}$.

  If $K'/K$ is a finite field extension with the ring of integers $A'$ and
  the residue field $\ff'$, then the diagram
  \begin{equation}\label{eqn:Kato-basic-0-0}
    \xymatrix@C1pc{
      \Br(K') \ar[rr]^-{\kappa^2(A')} \ar[d]_-{f_*} & & H^1(\ff')
      \ar[rr]^-{{\rm ev}_F(\ff')} \ar[d]^-{f'_*} && {\Q}/{\Z} \ar@{=}[d] \\
      \Br(K) \ar[rr]^-{\kappa^2(A)} && H^1(\ff)
      \ar[rr]^-{{\rm ev}_F(\ff)} && {\Q}/{\Z}}
  \end{equation}
  is commutative (cf. \cite[Lem.~9.8]{GKR}), where $f \colon \Spec(K') \to \Spec(K)$
  and $f' \colon \Spec(\ff') \to \Spec(\ff)$ are the projection maps.
\end{exm}

Let $A$ be as above and let ${\rm sp}^r(\pi) \colon K^M_{r}(K) \to K^M_r(\ff)$ denote
the specialization map. Recall from \cite[Chap.~III,  Thm.~7.3]{Weibel-K}
that this map is uniquely defined by the property that
${\rm sp^r}(\pi)(\{u_1\pi^{i_1}, \ldots , u_r\pi^{i_r}\}) = m^r_A(\{u_1, \ldots , u_r\})$
if $u_j \in A^\times$ and $i_j \in \Z$ for $j \in J_r$. This depends on the
choice of $\pi$ as a generator of $\fm$.
Let $\partial^r_A \colon K^M_r(K) \to K^M_{r-1}(\ff)$ denote the tame symbol map
which is uniquely defined by the property that
$\partial^r_A(\{\pi, u_1, \ldots , u_{r-1}\}) =
m^{r-1}_A(\{{u_1}, \ldots , u_{r-1}\})$
and $\partial^r_A(\{u_1, \ldots , u_{r}\}) = 0$ for $u_i \in A^\times$
(cf. \cite[\S~1.1]{Kato-hasse}). This differs from the classical tame symbol
for $K_2$ by a sign.

\begin{lem}\label{lem:Kato-basic-1}
For $n, q \in \N$ and $r \in \N_0$, the diagram
  \begin{equation}\label{eqn:Kato-basic-1-0}
    \xymatrix@C.8pc{
      H^q_n(\ff) \otimes K^M_r(K) \ar[d]_-{\phi^q_r} \ar[r]^-{\psi^q_r} &
      \left(H^q_n(\ff) \otimes K^M_r(\ff)\right) \bigoplus
      \left(H^q_n(\ff) \otimes K^M_{r-1}(\ff)\right)
      \ar[d]^-{\{\ , \ \} \oplus (-1)^{r-1}\{\ , \ \}} \\
      H^{q+r}_n(K) & H^{q+r}_n(\ff) \bigoplus  H^{q+r-1}_n(\ff)
      \ar[l]_-{\lambda^{q+r}_n(\pi)} }
    \end{equation}
 is commutative if we let $\phi^q_r(\chi \otimes b) =
 \{l^q_n(\chi), b\}$ and $\psi^q_r(\chi \otimes b) =
 \left(\chi \otimes {\rm sp}^r(\pi)(b), \chi \otimes {\partial}^r_A(b)\right)$. 
\end{lem}
\begin{proof}
Since $K^M_r(K)$ is generated by elements of the form $\{u_1,\ldots,u_{r-1}, \pi\}$
  and $\{u_1,\ldots, u_{r}\}$ with $u_i \in A^\times$, it suffices to check the
  commutativity for these two types of generators in $K^M_r(K)$.
  We let $\chi \in H^q_{n}(\ff)$. 

If $b = \{u_1,\ldots,u_{r-1}, \pi\}$, then ${\rm sp}^r(\pi)(b) =0$ and we get
  $\lambda^{q+r}_n(\pi)(\{\chi, {\rm sp}^r(\pi)(b)\} +
  (-1)^{r-1}\{\chi, \partial^r_A(b)\}) =
  (-1)^{r-1} \lambda^{q+r}_n(\pi)((0,\{\chi, \partial^r_A(b)\}))
  = \{l^{q+r-1}_n(\{\chi, \ov{u_1}, \ldots , \ov{u_{r-1}}\}), \pi\}$.
But the latter element is same as 
$\{l^q_n(\chi), u_1, \ldots , u_{r-1}, \pi\} = \{l^q_n(\chi), b\}$ by
~\eqref{eqn:Pairing-2}.
If $b = \{u_1,\ldots, u_{r}\}$, we get
 $\lambda^{q+r}_n(\pi)(\{\chi, {\rm sp}^r(\pi)(b)\} + \{\chi, \partial^r_A(b)\}) =
\lambda^{q+r}_n(\pi)(\{\chi, {\rm sp}^r(\pi)(b)\}) \ =$ \\
$l^{q+r-1}_n(\{\chi, \ov{u_1}, \ldots , \ov{u_{r}}\})$.
But the latter element is same as $\{l^q_n(\chi), u_1, \ldots , u_{r}\} =
  \{l^q_n(\chi), b\}$ by ~\eqref{eqn:Pairing-2}. This completes the proof.
\end{proof}

\begin{cor}\label{cor:Kato-basic-2}
  For $A$ as in \propref{prop:Kato-basic}, there is a commutative diagram
  \begin{equation}\label{eqn:Kato-basic-2-0}
    \xymatrix@C1pc{
      H^q(\ff) \otimes K^M_r(K) \ar[rr]^-{\id \otimes \partial^r_A}
      \ar[d]_-{\phi^q_r} && H^q(\ff) \otimes K^M_{r-1}(\ff)
      \ar[d]^-{(-1)^{r-1}\{\ , \ \}} \\
      \Fil_0H^{q+r}(K) \ar[rr]^-{\kappa^{q+r}(A)} && H^{q+r-1}(\ff).}
  \end{equation}
  If $q \le 2$, the same holds if we replace $H^*(-)$ by $H^*_n(-)$ for any $n \in \N$.
 \end{cor}
\begin{proof}
  The image of $\phi^q_r$ lies in $\Fil_0H^{q+r}(K)$ by
  \propref{prop:Kato-basic}(9). The commutativity of ~\eqref{eqn:Kato-basic-2-0}
  now follows directly from \lemref{lem:Kato-basic-1}.
  \end{proof}

\section{Filtrations of the Brauer group}\label{sec:Kat-Ev*}
In this section, we introduce the main objects appearing in \thmref{thm:Main-1}, namely,
the Kato and evaluation filtrations of the Brauer group of regular schemes over a local
field. The main result of this section is \propref{prop:Kato-Br-0} which is the first
key step in the proof of \thmref{thm:Main-1}.
We begin with the necessary definitions.

If $X$ is a Noetherian scheme and $Z \subset X$ is a closed subset, we let $\Div_Z(X)$
denote the free abelian group of Weil divisors on $X$ whose supports are contained in
$Z$. We let $\Div^{\eff}_Z(X)$ be the subgroup of $\Div_Z(X)$ consisting of
effective divisors. We let $\Div(X) = \Div_X(X)$ and
$\Div^{\eff}(X) = \Div^{\eff}_X(X)$.
If $X$ is a connected Noetherian regular scheme and $E \subset X$ is a
simple normal crossing divisor (which we shall refer to as snc divisor) on $X$,
we shall say that $(X,E)$ is an snc-pair.
Recall from \cite[\S~2.1]{KM-1} that a morphism $f \colon (X', E') \to (X,E)$ of
snc-pairs is a morphism of schemes $f \colon X' \to X$ such that $f^*(E)
\in \Div_{E'}(X')$. In particular, $f(U') \subset U$, where
$U' = X' \setminus E'$.
If $X$ is a  Noetherian integral scheme with function field $K$ and $y \in X^{(1)}$, 
we shall let $K_y$ (resp. $\wh{K}_y$) denote the field of fraction of
$\sO^h_{X,y}$ (resp. $\wh{\sO_{X,y}}$).

\subsection{Kato filtration for schemes}\label{sec:KFil-sc}
Let $(X,E)$ be an snc-pair such that $X$ is an $\F_p$-scheme. Let
$E_1, \ldots , E_r$ denote the irreducible components
of $E$ and let $j \colon U \inj X$ be the inclusion of the complement of $E$ in $X$.
Let $K$ denote the function field of $X$. Let $y_i$ denote the generic
point of $E_i$ and let $K_i$ (resp. $\wh{K}_i$) be the quotient field of
$\sO^h_{X, y_i}$ (resp. $\wh{\sO_{X, y_i}}$). Let
$\ov{h}_i \colon \Spec({\sO^h_{X, y_i}}) \to X$ and
${h}_i \colon \Spec(K_i) \to U$ denote the canonical maps.
If $\chi \in H^q(K)$ for some $q \ge 1$ and $Y \subset X$ is an integral divisor with
generic point $y$, we let $\Sw_Y(\chi)$ denote the Swan conductor of the image of
$\chi$ under the restriction map $H^q(K) \to H^q(K_y)$.

We fix integers $m, q \in \N$ and let $D = \sum_i n_iE_i \in \Div^{\eff}(X)$.
For $i \in J_r$, let $j_i \colon H^{q}_{m}(U)
\to H^{q}_{m}(K_i)$ be the canonical restriction map. 
The Kato filtrations of $H^{q}_{m}(U)$ and $H^{q}(U)$ (as defined by Kato
\cite[\S~7,8]{Kato-89}) are defined as follows.

\begin{defn}(cf. \cite[\S~7,8]{Kato-89}, \cite[Defn.~2.7]{Kerz-Saito-ANT})
  \label{defn:Log-fil-D} 
  We  let
 \[
  \Fil_DH^{q}_{m}(U) :=
  \Ker \left(H^{q}_{m}(U) \xrightarrow{\bigoplus j_i} \
  \stackrel{r}{\underset{i =1}\bigoplus} 
  \frac{H^{q}_{m}(K_i)}{\Fil_{n_i} H^{q}_m(K_i)}\right);
  \]
  \[
  \Fil_DH^{q}(U) := {\varinjlim}_m \Fil_DH^{q}_{m}(U).
\]
\end{defn}

\begin{remk}\label{remk:Log-fil-D-0}
  (1) \ We note that this definition of  $\Fil_DH^{q}_{m}(U)$ is not affected if
  we replace any $K_i$ by $\wh{K}_i$ (cf. \cite[Cor.~7.16]{KM-1}).

  (2) \ Another point to note concerning the above definition is that
  $\Fil_DH^{q}_{m}(U)$ depends both on the support of $D$ and the multiplicities of
  the components of $D$. For instance, if we let $1 \le r' < r, \ 
  E' = \stackrel{r'}{\underset{i =1}\sum} E_i, \
  D = \stackrel{r'}{\underset{i =1}\sum} n_iE_i$ (i.e., $n_i = 0 \ \forall \ i > r'$)
  and $U' = X \setminus E'$, then
      $\Fil_{D}H^{q}_{m}(U') \neq \Fil_DH^{q}_{m}(U)$.

  (3) \ Finally, note that in \cite{KM-1}, $\Fil_DH^{q}_{m}(U)$ is denoted by
  $\Fil^\log_DH^{q}_{m}(U)$.
\end{remk}

If $D = 0$ (i.e., $n_i = 0$ for every $i$), we shall write
$\Fil_DH^{q}_{m}(U)$ (resp. $\Fil_DH^{q}(U)$) as
$\Fil_0H^{q}_{m}(U)$ (resp. $\Fil_0H^{q}(U)$).
It is clear that the definition of $\Fil_DH^{q}_{m}(U)$ given above agrees with 
Definition~\ref{defn:Kato-1} when $X = \Spec(A)$ and $D = \Spec({A}/{(\pi^n)})$,
where $A$ is an hdvr containing $\F_p$.

\begin{lem}\label{lem:Log-fil-D-1} 
For $q \in \N$ and $1 \le m \le \infty$, we have the following.
  \begin{enumerate}
  \item
$j^*(H^q_m(X)) + H^{q}_m(U)\{p'\} \subset \Fil_0H^{q}_m(U)$.
\item
  If $0 \le D \le D'$, then $\Fil_DH^{q}_{m}(U) \subset \Fil_{D'}H^{q}_{m}(U)$.
\item
  $H^{q}_{m}(U) = {\underset{D \in \Div^\eff_E(X)}\varinjlim} \Fil_DH^{q}_{m}(U)$.
\item
  $\Fil_DH^{q}_m(U) = H^{q}_m(U)\{p'\} \bigoplus \Fil_DH^{q}_m(U)\{p\}$.
\item
  $\Fil_DH^{q}(U) =
  \Ker \left(H^{q}(U) \xrightarrow{\bigoplus j_i} \
  \stackrel{r}{\underset{i =1}\bigoplus} 
  \frac{H^{q}(K_i)}{\Fil_{n_i} H^{q}(K_i)}\right)$.
\end{enumerate}
\end{lem}
\begin{proof}
  Items (1) and (3) follow from
  \cite[Lem.~2.2, Cor.~2.5]{Kato-89} and \propref{prop:Kato-basic}(1).
  Item (4) follows from (1) while (2) is clear from definition.
  Item (5) is also clear because
  \[
    {\varinjlim}_m \Ker \left(H^{q}_{m}(U) \xrightarrow{\bigoplus j_i} \
  \stackrel{r}{\underset{i =1}\bigoplus} 
  \frac{H^{q}_{m}(K_i)}{\Fil_{n_i} H^{q}_m(K_i)}\right) \xrightarrow{\cong}
  \Ker \left(H^{q}(U) \xrightarrow{\bigoplus j_i} \
  \stackrel{r}{\underset{i =1}\bigoplus} 
  \frac{H^{q}(K_i)}{\Fil_{n_i} H^{q}(K_i)}\right).
  \]
\end{proof}

If we let $q =2$ in  Definition~\ref{defn:Log-fil-D}, then
the canonical isomorphism $H^2(K_i) \xrightarrow{\cong} \Br(K_i)$ shows that
the map $j_i \colon H^{2}(U) \to H^2(K_i)$ factors through $\Br(U)$
(cf. ~\eqref{eqn:Kummer}). We may therefore make the following definition.

\begin{defn}\label{defn:Log-Br}
  We let
  \[
  \Fil_D\Br(U):=
  \Ker \left(\Br(U) \xrightarrow{\bigoplus j_i} \
  \stackrel{r}{\underset{i =1}\bigoplus} 
  \frac{H^2(K_i)}{\Fil_{n_i} H^{2}(K_i)}\right).
  \]
\end{defn}
If $D = 0$, we shall write $\Fil_D\Br(U)$ as $\Fil_0 \Br(U)$.
One has inclusions
$\Br(X) \subset \Fil_0 \Br(U) \subset \Fil_D\Br(U) \subset \Fil_{D'}\Br(U)$ if
$D' \ge D \ge 0$.
The family $\{\Fil_D\Br(U)\}_{D \in \Div^\eff_E(X)}$ will be referred to as the
`Kato filtration' of $\Br(U)$.

It follows immediately (cf. the paragraph following Definition~\ref{defn:Coh-sheaf})
that
\begin{equation}\label{eqn:Fil-Br-m}
 \Fil_D\Br(U)[m]  = \Ker \left(\Br(U)[m] \xrightarrow{\bigoplus j_i} \
  \stackrel{r}{\underset{i =1}\bigoplus} 
  \frac{H^2_m(K_i)}{\Fil_{n_i} H^{2}_m(K_i)}\right)
\end{equation}
for every $m \in \N$.
Furthermore, an argument similar to that of \lemref{lem:Log-fil-D-1}
gives the following.
\begin{cor}\label{cor:Kato-exhaust}
  Let $D \in \Div^\eff_E(X)$ and let $m = p^tr$ with
  $t \in \N_0$ and $r \in I_p$. We then have the following.
  \begin{enumerate}
    \item
      $\Fil_D\Br(U)[m] = \Br(U)[r] \bigoplus \Fil_D\Br(U)[p^t]$.
    \item
      $\Fil_D\Br(U) = \Br(U)\{p'\} \bigoplus \Fil_D\Br(U)\{p\}$.
      \item
        $\Br(U)[m] = {\underset{D \in \Div^\eff_E(X)}\bigcup} \Fil_D\Br(U)[m]$.
        \item
          $\Br(U) = {\underset{D \in \Div^\eff_E(X)}\bigcup} \Fil_D\Br(U)$.
          \end{enumerate}
   \end{cor}

The following result about the Kato filtration on
cohomology and the Brauer group follows immediately by combining
Definitions~\ref{defn:Log-fil-D} and ~\ref{defn:Log-Br} with ~\eqref{eqn:Kummer}.

\begin{lem}\label{lem:Br-fil-0}
  Let $(X,E)$ be as above. Then for any
  $D \in \Div_E^\eff(X)$ and $m \in \N$, there is a short exact sequence
  \[
  0 \to {\Pic(U)}/{m} \to \Fil_D H^2_{m}(U) \xrightarrow{\kappa_U}
  \Fil_D \Br(U)[m] \to 0
  \]
  which is natural in the snc-pair $(X,E)$. In particular, the canonical map
  $\frac{\Fil_D H^2_{m}(U)}{\Fil_{D'} H^2_{m}(U)}
  \to \frac{\Fil_D \Br(U)[m]}{\Fil_{D'} \Br(U)[m]}$ is an isomorphism for
    every $0 \le D' \le D$.
\end{lem}

The above lemma says that if we are given an element $\omega \in \Br(U)[m]$, then
we can always lift it to a class $\chi \in H^2_{m}(U)$. Moreover, any such lift
$\chi$ has the property that for any $D \in \Div_E^\eff (X)$, one has
$\chi \in \Fil_D H^2_{m}(U)$ if and only if $\omega \in \Fil_D \Br(U)[m]$.
These facts will be implicitly used every time we have to prove (or disprove) that
certain Brauer class on $U$ lies in a given subgroup of the form $\Fil_D \Br(U)[m]$.

\subsection{A description of $\Fil_0 \Br(U)$}\label{sec:Fil-0**}
In this subsection, our goal is to prove \propref{prop:Kato-Br-0}. This is one  of
the key results that we shall use in the proofs of the main results.
In particular, it will prove a part of \thmref{thm:Main-1}.
The setting is the following.

Let $f \colon (X',E') \to (X,E)$ be a morphism of snc-pairs of $\F_p$-schemes.
We let $U' = X' \setminus E'$ and let $h \colon U' \to U$ be the restriction of
$f$. If $U' = f^{-1}(U)$, we let $g \colon E' \to E$ denote the restriction of
$f$ to $E'$.
We let $X^o = X \setminus E_\sing$ and $X'^o = X' \setminus E'_\sing$. We let
$\{E_1, \ldots , E_r\}$ (resp. $\{E'_1, \ldots , E'_s\}$) be the set of irreducible
components of $E$ (resp. $E'$) and let $E^o_i = E_i \bigcap X^o$ (resp.
$E'^o_i = E'_i \bigcap X'^o$). If $f^{-1}(E_\sing) \subset E'_\sing$ (this includes the
case when $f^{-1}(E_\sing) = \emptyset$), then
$f \colon (X'^o, E'^o) \to (X^o, E^o)$ is a morphism of snc-pairs,
where we let $E^o = E_\reg$ and $E'^o = E'_\reg$.

For $i \in J_r$, we let $f^*(E_i) = \stackrel{s}{\underset{j =1}\sum} n_{i j} E'_j
\in \Div(X')$. If $U' = f^{-1}(U)$, we define the map
\begin{equation}\label{eqn:Kato-Br-0*}
\wt{g}^* \colon \stackrel{r}{\underset{i =1}\bigoplus} H^1_\et(E_i) \to \
\stackrel{s}{\underset{j =1}\bigoplus} H^1_\et(E'_j); \ \
(a_i) \mapsto ((n_{ij}g^*_{ij}(a_i)),
\end{equation}
which takes $a \in H^1_\et(E_i)$ to
$(n_{i 1} g^*_{i 1}(a), \ldots , n_{i s} g^*_{i s}(a))$, where $g_{ij} \colon
E'_j \to E_i$ is the restriction of $f$ if $E'_j \subset f^{-1}(E_i)$.
If $E'_j \nsubset f^{-1}(E_i)$, we let $g^*_{i j}(a) = 0$. If
$f^{-1}(E_\sing) \subset E'_\sing$, then this map restricts to a similar map
$\wt{g}^* \colon \stackrel{r}{\underset{i =1}\bigoplus} H^1_\et(E^o_i) \to
\stackrel{s}{\underset{j =1}\bigoplus} H^1_\et(E'^o_j)$.

For each $i \in J_r$, we let $y_i$ be the generic point of $E_i$ and
  let $u_i \colon \Spec(k(y_i)) \inj E_i, \ \wt{u}_i \colon \Spec(k(y_i)) \inj X$
  and $v_i \colon E_i \inj X$ be the canonical inclusions. For each $j \in J_s$,
  we let $y'_j$ be the generic point of $E'_j$ and
  let $u'_j \colon \Spec(k(y'_j)) \inj E'_j, \ \wt{u}'_j \colon \Spec(k(y'_j)) \inj X'$
  and $v'_j \colon E'_j \inj X'$ be the canonical inclusions. 
  We then have $\wt{u}_i = v_i \circ u_i$ and $\wt{u}'_j = v'_j \circ u'_j$ for each
  $i \in J_r$ and $j \in J_s$. We let $j \colon U \inj X$ and $j' \colon U' \inj X'$
  be the inclusion maps. Let $\iota \colon \{\eta_X\} \inj X$ and
    $\iota' \colon \{\eta_X\} \to U$ denote the inclusions of the generic point of
      $X$.

\begin{lem}\label{lem:K-Br-0}
    For $i \in J_r$ and $q \ge 2$, we have a canonical isomorphism
    $H^{q-1}_\et(E_i, {\Q}/{\Z}) \xrightarrow{\cong} H^q_\et(X, \wt{u}_{i *} \Z)$.
  \end{lem}
  \begin{proof}
    It is clear that $H^q_\et(X, \wt{u}_{i *} \Z) \cong H^q_\et(X, v_{i *}(u_{i *} \Z))
    \cong H^q_\et(X, v_{i *}\Z) \cong H^q_\et(E_i, \Z)$.
    The lemma now follows because $H^q_\et(E_i, \Q) \cong H^q_\et(E_i, u_{i *}\Q)
    \cong H^q_\et(k(y_i), \Q) = 0$ for $q \in \N$, where the second isomorphism
    is easily checked using the Leray spectral sequence.
\end{proof}

The definition of $F$-finite $\F_p$-schemes is recalled in \S~\ref{sec:Fil-DRW}.

  \begin{lem}\label{lem:K-Br-1}
    Assume that $X$ is $F$-finite.
    Then the canonical map $H^2_\et(X, j_* \sO^\times_{U}) \to
    \H^2_\et(X, {\bf R}j_*\sO^\times_U) \cong \Br(U)$ induces an isomorphism
    \[
    H^2_\et(X, j_* \sO^\times_{U})[p^m] \xrightarrow{\cong} \Fil_0 \Br(U)[p^m]
    \]
    for every $m \in \N$.
  \end{lem}
  \begin{proof}
Since $\sO^\times_U\{p\} = 0 = {\bf R}^1j_* \sO^\times_{U}$,
    there is a short exact sequence of sheaves (cf. the proof of \lemref{lem:Coh-iso})
  \[
0 \to j_*\sO^\times_U \xrightarrow{p^m} j_*\sO^\times_U \to j_* W_m\Omega^1_{U, \log} \to 0
\]
on $X_\et$.
Taking the cohomology, we get a short exact cohomology sequence
\begin{equation}\label{eqn:K-Br-1-0}
  0 \to {H^1_\et(X, j_* \sO^\times_U)}/{p^m} \to H^1_\et(X, j_* W_m\Omega^1_{U, \log}) \to
    H^2_\et(X, j_* \sO^\times_U)[p^m] \to 0.
\end{equation}
The canonical map $j_* \sF \to {\bf R}j_*\sF$ for any sheaf $\sF$ on $U_\et$
induces a homomorphism $H^1_\et(X, j_* W_m\Omega^1_{U, \log}) \to H^2_{p^m}(U)$.
By \thmref{thm:H^1-fil} and \cite[Lem.~11.2(4)]{KM-1} (see its proof), this map
factors through an isomorphism
$H^1_\et(X, j_* W_m\Omega^1_{U, \log})  \xrightarrow{\cong} \Fil_0 H^2_{p^m}(U)$.

We now look at the commutative  diagram of exact sequences
\begin{equation}\label{eqn:K-Br-1-1}
  \xymatrix@C1pc{
    0 \ar[r] & {H^1_\et(X, j_* \sO^\times_U)}/{p^m} \ar[r] \ar[d] &
    H^1_\et(X, j_* W_m\Omega^1_{U, \log}) \ar[r] \ar[d]^-{\cong} &
    H^2_\et(X, j_* \sO^\times_U)[p^m] \ar[r] \ar[d] & 0 \\
     0 \ar[r] & {\Pic(U)}/{p^m} \ar[r] &
     \Fil_0 H^2_{p^m}(U) \ar[r] & \Fil_0 \Br(U)[p^m] \ar[r] & 0,}
  \end{equation}
where the bottom sequence is exact by \lemref{lem:Br-fil-0} and the vertical
arrows are induced by the canonical map $j_* \sF \to {\bf R}j_*\sF$.

It is an easy exercise using the Leray spectral sequence that the
canonical map $H^1_\et(X, j_* \sO^\times_U) \to \Pic(U)$ is an isomorphism.
The exactness of the rows and a diagram chase now imply that the canonical map
$H^2_\et(X, j_* \sO^\times_U)[p^m] \to  \Br(U)$ factors through
$H^2_\et(X, j_* \sO^\times_U)[p^m] \to \Fil_0 \Br(U)[p^m]$ which is an isomorphism.
This completes the proof.
\end{proof}

 \begin{lem}\label{lem:K-Br-3}
Let $\alpha \colon U \inj X^o$ be the inclusion. Then the canonical map
  $H^2_\et(X^o, \alpha_* \sO^\times_{U}) \to
    \H^2_\et(X^o, {\bf R}\alpha_*\sO^\times_U) \cong \Br(U)$ induces an isomorphism
    \[
    H^2_\et(X^o, \alpha_* \sO^\times_{U})[m] \xrightarrow{\cong} \Fil_0 \Br(U)[m]
    \]
    for every $m \in I_p$.
  \end{lem}
  \begin{proof}
 By \corref{cor:Kato-exhaust}(1), we can replace $\Fil_0 \Br(U)[m]$ by
    $\Br(U)[m]$. Meanwhile, the Leray spectral sequence for
    $\alpha \colon U \inj X^o$ gives rise to a canonical exact sequence
    \begin{equation}\label{eqn:K-Br-3-0}
      0 \to H^2_\et(X^o, \alpha_* \sO^\times_{U}) \to \Br(U) \to
      H^0_\et(X^o, {\bf R}^2\alpha_* \sO^\times_U).
    \end{equation}
    It suffices therefore to show that
    $H^0_\et(X^o, {\bf R}^2\alpha_* \sO^\times_U)\{p'\} = 0$.
    For this, it is enough to show that $({\bf R}^2\alpha_* \sO^\times_U)_x\{p'\} =0$
    for every $x \in E^o$.

 We now fix $x \in E^o$ and let $A = \sO^{sh}_{X,x}$. Since the irreducible
    components of $E^o$ are pairwise disjoint, there is a unique $i \in J_r$ such that
    $x \in E^o_i$. We let $B = A[\pi^{-1}]$, where $(\pi) \subset A$
    is the ideal that defines $E^o_i$ locally at $x$. We need to show that
    $\Br(B)[m] = 0$ if $m \in I_p$. Since $\Pic(B) = 0$, this is equivalent to
    showing that $H^2_m(B) = 0$. But this follows by the
    exact sequence
    \[
    H^2_m(A) \to H^2_m(B) \to H^1_m({A}/{(\pi)}),
    \]
    whose end terms are zero because $A$ and ${A}/{(\pi)}$ are strictly henselian
    local rings. 
 \end{proof}

\begin{lem}\label{lem:K-Br-2}
    On $X_\et$, there is a canonical short exact sequence of sheaves
    \begin{equation}\label{eqn:K-Ev-0-2}
0 \to \sO^\times_{X} \to j_* \sO^\times_{U} \to \
  \stackrel{r}{\underset{i =1}\bigoplus} \wt{u}_{i *} \Z \to 0.
\end{equation}
\end{lem}
  \begin{proof}
For any $y \in X^{(1)}$ and
  $w \in U^{(1)}$, we let
  $\iota_y \colon \Spec(k(y)) \to X$ and $\iota'_w \colon \Spec(k(w)) \to U$
  denote the inclusion maps.
We then have short exact sequences of {\'e}tale sheaves
  (cf. \cite[(3.8)]{CTS})
  \begin{equation}\label{eqn:K-Ev-0-1}
    0 \to \sO^\times_{X} \to \iota_* \sO^\times_{\eta_X} \to
    {\underset{y \in {X}^{(1)}}\bigoplus} \iota_{y *} \Z \to 0 \ \ {\rm and} \ \
 0 \to \sO^\times_{U} \to \iota'_* \sO^\times_{\eta_X} \to
    {\underset{y \in U^{(1)}}\bigoplus} \iota'_{y *} \Z \to 0
    \end{equation}
  on $X$ and $U$, respectively.

Since ${\bf R}^1j_* \sO^\times_{U} = 0$, the second sequence in ~\eqref{eqn:K-Ev-0-1}
  yields a short exact sequence
 \begin{equation}\label{eqn:K-Ev-0-2.5} 
0 \to j_* \sO^\times_{U} \to \iota_* \sO^\times_{\eta_X} \to
    {\underset{y \in U^{(1)}}\bigoplus} \iota_{y *} \Z \to 0
 \end{equation}
By comparing this with the first exact sequence of
 ~\eqref{eqn:K-Ev-0-1}, we get ~\eqref{eqn:K-Ev-0-2}.
\end{proof}

  \begin{prop}\label{prop:Kato-Br-0}
 Let $f \colon (X',E') \to (X,E)$ be as above such that $U' = f^{-1}(U)$ and
  $f^{-1}(E_\sing) \subset E'_\sing$.  Assume that $X$ and $X'$ are $F$-finite.
  Then there exists an exact sequence
  \begin{equation}\label{eqn:K-Ev-0-0}
  0 \to \Br(X) \to \Fil_0 \Br(U) \xrightarrow{\partial_U} H^1_\et(E^o)
  \end{equation}
  such that $\partial_U(\Fil_0 \Br(U)\{p\}) \subset \ \stackrel{r}{\underset{i =1}
    \bigoplus} H^1_\et(E_i) \subset H^1_\et(E^o)$. Furthermore, there is a 
  commutative diagram
  \begin{equation}\label{eqn:Kato-Br-0-1}
    \xymatrix@C1.7pc{
      \Fil_0\Br(U)\{p\} \ar[d]_-{f^*} \ar[r]^-{\partial_U} &
      \stackrel{r}{\underset{i =1}\bigoplus} H^1_\et(E_i)\{p\}
      \ar[d]^-{\wt{g}^*} \ar@{^{(}->}[r] &
      \stackrel{r}{\underset{j =1} \bigoplus} H^1_\et(E^o_i) \ar[d]^-{\wt{g}^*} &
      \Fil_0\Br(U)\{p'\} \ar[d]^-{f^*} \ar[l]_-{\partial_U} & \\
      \Fil_0\Br(U')\{p\} \ar[r]^-{\partial_{U'}} & \stackrel{s}{\underset{j =1}
        \bigoplus} H^1_\et(E'_j)\{p\} \ar@{^{(}->}[r] &
      \stackrel{s}{\underset{j =1} \bigoplus} H^1_\et(E'^o_j)
     & \Fil_0\Br(U')\{p'\} \ar[l]_-{\partial_{U'}}.}
    \end{equation}
   \end{prop}   
\begin{proof}
To prove the $p$-part of our claim,
  we look at the commutative diagram of exact sequences of
  {\'e}tale sheaves (on $X$)
  
\begin{equation}\label{eqn:Kato-Br-0-2}
    \xymatrix@C1.7pc{
      0 \ar[r] & \sO^\times_X \ar[r]^-{j^*} \ar[d]_-{f^*} &
      j_* \sO^\times_U \ar[r]^-{\divf} \ar[d]^-{h^*} &
      \stackrel{r}{\underset{i =1}\bigoplus} \wt{u}_{i *} \Z \ar[d]^-{\wt{g}^*}
      \ar[r] & 0 \\
      0 \ar[r] & f_* \sO^\times_{X'} \ar[r]^-{j'^*} \ar[d] &
      f_*(j'_* \sO^\times_{U'}) \ar[r]^-{\divf} \ar[d] &
      \stackrel{s}{\underset{j =1}\bigoplus} f_* \wt{u}_{i *} \Z
      \ar[d] \\
         & {\bf R}f_* \sO^\times_{X'} \ar[r] & {\bf R}f_* (j'_* \sO^\times_{U'})
         \ar[r]^-{\divf} &
         \stackrel{s}{\underset{j =1}\bigoplus} {\bf R}f_* \wt{u}_{j *} \Z.}
\end{equation}

In this diagram, the top two rows are exact by \lemref{lem:K-Br-2}
and the bottom row is an exact triangle in $\sD_\et(X)$.
The vertical arrows in the lower level are
the canonical maps and $\wt{g}^* \colon \wt{u}_{i *} \Z \to
\stackrel{s}{\underset{j =1}\bigoplus} f_* \wt{u}_{i *} \Z$
is defined by letting $\wt{g}^*(1) = (n_{i1}, \ldots , n_{is})$ if one has $f^*(E_i)
= \stackrel{s}{\underset{j =1}\sum} n_{ij} E'_j \in \Div(X')$.
The arrow `$\divf$' is the map which sends a rational function to its divisor
of zeros and poles.

By considering the ({\'e}tale) cohomology and using \lemref{lem:K-Br-0},
we get a commutative diagram of exact sequences of abelian groups
\begin{equation}\label{eqn:Kato-Br-0-3}
    \xymatrix@C1.7pc{
      0 \ar[r] & \Br(X) \ar[r]^-{j^*} \ar[d]_-{f^*} & H^2_\et(X, j_* \sO^\times_U)
      \ar[r]^-{\partial_U} \ar[d]^-{h^*} &
      \stackrel{r}{\underset{i =1}\bigoplus} H^1_\et(E_i, {\Q}/{\Z})
      \ar[d]^-{\wt{g}^*} \\
      0 \ar[r] & \Br(X') \ar[r]^-{j'^*} & H^2_\et(X', j'_*\sO^\times_{U'})
      \ar[r]^-{\partial_{U'}} &  
      \stackrel{s}{\underset{j =1}\bigoplus} H^1_\et(E'_j, {\Q}/{\Z}),}
    \end{equation}
in which $j^*$ is injective because its composition with the canonical map
$H^2_\et(X, j_* \sO^\times_U) \to \H^2_\et(X, {\bf R}j_*\sO^\times_U) \cong \Br(U)$
is injective by \cite[Thm.~3.5.7]{CTS}. The same is true for $j'^*$.

Since the terms on the two ends of the exact sequences of ~\eqref{eqn:Kato-Br-0-3}
are torsion groups, it follows that $H^2_\et(X, j_* \sO^\times_U)$ and
$H^2_\et(X', j'_*\sO^\times_{U'})$ are torsion groups.
The $p$-primary part of ~\eqref{eqn:K-Ev-0-0} and the
commutativity of the left square in ~\eqref{eqn:Kato-Br-0-1} now follow
by using ~\eqref{eqn:Kato-Br-0-3} and applying \lemref{lem:K-Br-1} to $X$ and
$X'$.

To prove the prime-to-$p$ part of ~\eqref{eqn:K-Ev-0-0} and the
commutativity of the right square in ~\eqref{eqn:Kato-Br-0-1}, we repeat the
above argument mutatis-mutandis to the morphism of snc-pairs
$f \colon (X'^o, E'^o) \to (X^o, E^o)$ (this map is defined by our
assumption), and use \lemref{lem:K-Br-3} 
and \cite[Thm.~3.7.6]{CTS}, the latter saying that the restriction map
$\Br(X) \to \Br(X^o)$ is an isomorphism.
The commutativity of the middle square in ~\eqref{eqn:Kato-Br-0-1}
is clear. The proof of the proposition is now complete.
\end{proof}

\subsection{Evaluation filtration of the Brauer group}\label{sec:Ev-Br}
We shall now define the evaluation filtration on the Brauer group in a special
situation. This filtration has much more explicit definition than the Kato
filtration.
Let $k$ be an hdvf of characteristic $p$ with finite residue field $\F$. We let
$\sO_k$ be the ring of integers of $k$ with maximal ideal $\fm = (\pi)$. We assume
$\sO_k$ to be an excellent ring. We let $S = \Spec(\sO_k)$ and $\eta = \Spec(k)$.

Let $\sX$ be a finite type,  separated, connected, regular, faithfully flat
$\sO_k$-scheme. We write $X=\sX \times_S \Spec(k)$ and
$\sX_s=\sX\times_S\Spec(\F)$. We let $i \colon \sX_s \inj \sX$ and
$j \colon X \inj \sX$ denote the inclusions. We assume that $Y = (\sX_s)_\red$
is an snc divisor on $\sX$ with irreducible components $Y_1,\ldots,Y_r$. We then
have $Y_\sing = \bigcup_{i \neq j} (Y_i \bigcap Y_j)$.  
We let $\sX^o = \sX \setminus Y_\sing$ and $Y^o = Y_\reg = Y^o_1 \amalg \cdots
\amalg Y^o_r$, where $Y^o_i = Y_i \setminus Y_\sing$. For a subscheme
$\sZ \subset \sX$, we let $\sZ_\eta = \sZ \times_S \Spec(k)$ and
$\sZ_s = \sZ\times_S\Spec(\F)$. Note that it is not assumed in our set-up
that $\sX_s$ is a reduced $\F$-scheme or $X$ is a smooth $k$-scheme.

For a finite field extension (which we shall abbreviate as ffe)
${k'}/k$ with ramification index $e({k'}/k)$,
ring of integers $\sO_{k'}$, maximal ideal $\fm' = (\pi')$ and residue field
$\F'$, we let $\sZ({\sO_{k'}}/I) = \Hom_{\Sch_{\sO_k}}(\Spec({\sO_{k'}}/I), \sZ)$ for any
ideal $I \subset \sO_{k'}$ and subscheme $\sZ \subset \sX$.
For a subscheme $Z \subset X$, we let $Z(k') = \Hom_{\Sch_k}(\Spec(k'), Z)$.
Note that the canonical restriction map $\sZ(\sO_{k'}) \to \sZ_\eta(k')$ is an
injection for any subscheme $\sZ \subset \sX$. This is a bijection if $\sZ$ is proper
over $S$.
For $P \in X_{(0)}$, we shall let $\ov{P}$ denote the scheme-theoretic closure of
$P$ in $\sX$. The latter is an integral closed subscheme of $\sX$. Moreover,
we have the following.

\begin{lem}\label{lem:qs-fin}
  $\ov{P}$ is finite over $S$ if and only if $\ov{P} \bigcap Y \neq \emptyset$.
  In either case, $\ov{P}$ is the spectrum of a 1-dimensional henselian local domain
whose residue field is finite over $\F$.
\end{lem}
\begin{proof}
The second statement is clear because $\sO_k$ is an hdvr. The `only if'
part of the first statement is also clear and the `if' part follows from
\cite[Chap.~I, Thm.~4.2]{Milne-EC}.
\end{proof}

If $P \in X_{(0)}$ is such that $\ov{P}$ is finite over $S$, then the
$\ov{P} \bigcap \sX_s$ is a 0-dimensional closed subscheme $\sX_s$
whose support consists of the unique closed point of $\ov{P}$. We shall denote
this point by $P_0$ and call it the reduction (or specialization) of $P$. We shall
often say that $P$ specializes to $P_0$. Letting $X_{\fin}$ denote the
set of points $P \in X_{(0)}$ such that $\ov{P}$ is finite over $S$, we get a
specialization map $\esp \colon X_{\fin} \to Y_{(0)}$.
For $P \in X_\fin$, we let
$e(P) = {\rm length}({{\sO_{\sX, P_0}}/{(I(\ov{P}) + I(Y))}})$ denote the intersection
multiplicity of $\ov{P} \bigcap Y$ at $P_0$.
We let $e_i(P) = {\rm length}({{\sO_{\sX, P_0}}/{(I(\ov{P}) + I(Y_i))}})$ denote the
intersection multiplicity of $\ov{P} \bigcap Y_i$ at $P_0$.
We let $\sX_s = \stackrel{r}{\underset{i =1}\sum} n_i Y_i \in \Div(\sX)$.

\begin{lem}\label{lem:Length-0}
  We have $[k(P):k] = \left(\stackrel{r}{\underset{i =1}\sum} n_i e_i(P)\right)
  [k(P_0):\F]$.
\end{lem}
\begin{proof}
  We let $A = {\sO_{\sX, P_0}}/{I(\ov{P})}$ and let
  $\pi_i = \stackrel{i}{\underset{j =1}\prod} f^{n_j}_j$, where
  $f_j$ defines $Y_i$ at $P_0$ and $\pi = \stackrel{r}{\underset{j =1}\prod} f^{n_j}_j$
  in $A$. Using the exact sequence of $A$-modules
  \[
  0 \to \frac{A}{(f^{n_{i+1}}_{i+1})} \to \frac{A}{(\pi_{i+1})} \to \frac{A}{(\pi_{i})}
  \to 0
  \]
  and the additivity of the length function, we get
  \[
  {\rm length}({{\sO_{\sX, P_0}}/{(I(\ov{P}) + I(\sX_s))}}) =
  \stackrel{r}{\underset{i =1}\sum} {\rm length}(A/{(f^{n_i}_i)}) =
  \stackrel{r}{\underset{i =1}\sum} n_i e_i(P),
  \]
where the second equality is an easy exercise (cf. \cite[Lem.~A.1.3]{Fulton}).
As $[k(P):k] = \left({\rm length}({{\sO_{\sX, P_0}}/{(I(\ov{P}) + I(\sX_s))}})\right)
[k(P_0):\F]$, the claim follows. 
\end{proof}

We let  
\begin{equation}\label{eqn:qs-fin-0}
X^o_{\fin}:=\{P \in X_{\fin}| \text{$\ov{P}$ is regular and $P_0  \in Y^o$}\} 
\end{equation}
\[=\{P \in X_{(0)}| \text{$\ov{P}$ is regular and }\emptyset \neq |\ov{P} \cap Y|
\subset Y^o \}  \text{ (by \lemref{lem:qs-fin})};
\]
\[
X^o_{\tr} := \{P \in X^o_\fin| \text{$\ov{P}$ intersects $Y^o$ transversely}\};
\]
\[
X_{\ur} := \{P \in X^o_\fin| \text{$\ov{P}$ is {\'e}tale over $S$}\}.
\]

It follows from \lemref{lem:Length-0} that if $P \in X^o_{\fin}$, then there
is exactly one $i \in J_r$ such that $P_0 \in Y_i$ and we have
\begin{equation}\label{eqn:mult}
[k(P):k] = n_i e_i(P) [k(P_0):\F] = n_i e(P) [k(P_0):\F].
\end{equation}

We let $r \in \N$. We let $rY \subset \sX$ be the closed subscheme locally given by
$\Spec({\sO_{\sX,x}}/{(t^r)})$, where $t$ is a local equation for $Y$ at a point
$x \in \sX$. For a finite unramified field extension 
(which we shall abbreviate as fufe) ${k'}/k$ with residue field $\F'$ and
$P\in (X_{k'})^o_{\fin}$, we let
\begin{equation}\label{eqn:qs-fin-1}
  \begin{array}{lll}
    B_{k'}(P,r) & := & \{Q\in (X_{k'})^o_{\fin}| \ov{P}\times_\sX rY =
    \ov{Q}\times_\sX rY\} \\
    & = & \{Q\in (X_{k'})^o_{\fin}| \ov{P}\times_\sX rY^o = \ov{Q}\times_\sX rY^o\}.
    \end{array}
\end{equation}
We let
\[
B_{k'}(P,0) = \{Q\in (X_{k'})^o_{\fin}|[k(Q) : k] = [k(P) : k]\}
\]
\[ \hspace*{2cm} = \{Q\in (X_{k'})^o_{\fin}|[k(Q) : k'] = [k(P) : k']\}.
\]
We shall usually write $B_{k}(P,r)$ simply as $B(P,r)$ for all $r \ge 0$.

For an element $\sA \in \Br(X)$ and an ffe ${k'}/k$, we let
\begin{equation}\label{eqn:Ev-map}
  \ev^{k'}_\sA \colon (X_{k'})_{(0)} \to {\Q}/{\Z};
  \ \  {\ev^{k'}_\sA}(P) = \inv_{k(P)}(\sA|_{k'}),
\end{equation}
where $\sA|_{k'}$ is the pull-back of $\sA$ under
the map $\Spec(k(P)) \to X$ and $\inv_{k(P)} \colon \Br(k(P)) \to {\Q}/{\Z}$ is the
Hasse invariant map for $k(P)$ (cf. Example~\ref{exm:Kato-basic-2}).
We shall usually write $\ev^{k}_\sA$ as $\ev_\sA$.

We can now define the evaluation filtration on $\Br(X)$ as follows.

\begin{defn}\label{defn:EV-filt-def}
  For $n \ge -1$, we let
  \[
  {\Ev}_n \Br(X) = \{\sA \in \Br(X)| \ev_\sA \ \text{is constant on} \ B(P,n+1) \
  \text{ for all } P \in X^o_{\fin} \};
  \]
  \[
  {\Ev}_{-2}\Br(X) = \{\sA \in \Br(X)|  \ev_\sA \ \text{is zero on}\ X^o_\fin\}.
  \]
\end{defn}

In other words, a Brauer class $\sA \in \Br(X)$ lies in ${\Ev}_n \Br(X)$
for $n \ge 0$ if and only if it takes a constant value in the
$\pi$-adic ball of radius $n+1$
around every point in $X^o_{\fin}$. Such classes can often be explicitly
described.

It is clear from ~\eqref{eqn:mult} and the above definition that
each ${\Ev}_n \Br(X)$ is a subgroup of $\Br(X)$ and
${\Ev}_n \Br(X) \subset {\Ev}_{n+1} \Br(X)$ for every $n \ge -2$.
It is not however clear that $\{{\Ev}_n \Br(X)\}_{n \ge -2}$
is an exhaustive filtration on $\Br(X)$.
We shall show nonetheless that this is actually true.

When $\sX$ is smooth over $S$ and has geometrically integral fibers, then
the evaluation filtration was originally defined by  Bright--Newton \cite{Bright-Newton}.
Although it is not clear a priori that their filtration agrees with the one
defined above, we shall show later in this paper that the two filtrations are
actually the same under the smoothness condition.

\vskip .2cm

One of the main objectives of this paper is to show that a non-log version of the
Kato filtration on $\Br(X)$ coincides with the evaluation filtration.

\section{Generalized refined Swan conductor}\label{sec:RSCKC}
Several of the key ingredients used in the proofs of the main theorems of this paper
were established in \cite{KM-1}. These results concern the
cohomological description of the Kato filtration and the refined Swan conductors
on the graded pieces of this filtration. In this section, we recall these results and
derive several applications needed later in the paper.

\subsection{The filtration on de Rham-Witt complex}\label{sec:Fil-DRW}
Recall that an $\F_p$-scheme $X$ is called $F$-finite if the
absolute Frobenius $\psi \colon X \to X$ is a finite morphism.
We shall say that a commutative ring $A$ is $F$-finite if it is an $\F_p$-algebra
such that $\Spec(A)$ is $F$-finite. In the rest of this paper, we shall assume all
$\F_p$-schemes to be $F$-finite. We refer to
\cite[Prop.~2.2]{KM-1} for several useful properties of $F$-finite rings and schemes
(and their modules of absolute K{\"a}hler differentials) that we shall use.
A particular case of this result that will be implicitly used throughout this paper
is that if $X$ is an essentially of finite type scheme over an $\F_p$-algebra
which is a quasi-excellent hdvr with perfect residue field,
then $X$ is $F$-finite. We shall also use the fact that
$W_m\Omega^q_X$ is a coherent sheaf of $W_m\sO_X$-modules if $X$ is Noetherian and
$F$-finite.

We let  $X$ be a Noetherian regular $\F_p$-scheme and let $E \subset X$ be an snc
divisor with irreducible components $E_1, \ldots , E_r$.  Let $y_i$ denote the
generic point of $E_i$ for $i \in J_r$. Let
$j \colon U \inj X$ be the inclusion of the complement of $E$ in $X$.
For $m \ge 1$, we let
$W_m\Omega^\bullet_X(\log E) \subset j_*W_m\Omega^\bullet_U$ be the de Rham-Witt complex
with log poles along $E$ (cf. \cite[\S~2.1]{KM-1}). This coincides with
the de Rham-Witt complex $W_m\Omega^\bullet_{X_\log}$ of the log scheme $X_\log$,
where the latter is defined by the canonical log
structure on $X$ given by $E$ (cf. \cite[\S~4]{Hyodo-Kato},
\cite[\S~1.2]{Geisser-Hesselholt-JAMS}, \cite[Thm.~2.8]{KM-1}).

For $\n=(n_1, \ldots , n_r)\in \Z^r$ and $t \in \N$, we define $\n/t:=
(\lfloor n_1/t \rfloor, {\ldots}, \lfloor n_r/t \rfloor)$, where
$\lfloor \cdot \rfloor$ (resp. $\lceil  \cdot \rceil$) is the greatest (resp. least)
integer function. We shall say that $\n \ge n'$ (equivalently, $\n' \le \n$)
in $\Z^r$ if $n_i \ge n'_i$ for
every $i$. We shall say that $\n \ge 0$ (resp.  $\n > 0$) if $n_i \ge 0$ for every $i$
(resp. $n_i > 0$ for some $i$). We shall say that $\n < 0$ if $n_i < 0$ for each $i$.
For any $\n \in \Z^r$, we let
$E_\n = n_1 E_1 + \cdots + n_r E_r \in \Div_E(X)$. For $t \in \N$, we shall let
$E_\n/t$ denote the Weil divisor $E_{\n/t}$. 
For any $\F_p$-scheme $Y$, we let $W_m(Y) = W_m(\sO(Y))$.

\begin{defn}\label{defn:Log-fil-3}
  Let $D = \stackrel{r}{\underset{i =1}\sum} n_iE_i \in \Div_E(X)$ with
  $n_i \in \Z$. For a morphism $u \colon \wt{X} \to X$ of schemes, we let
  $\wt{U} = U \times_X \wt{X}$ and $\wt{D} = D\times_X \wt{X}$. We let
      \[
      \Fil_D W_m(\wt{U}) = \left\{\underline{a}=(a_{m-1},\ldots, a_0)
      \in W_m(\wt{U})| \ 
        a_i^{p^i} \in H^0_\et(\wt{X}, \sO_{\wt{X}}(\wt{D})) \ \forall \ i\right\}.
      \]
      Let $\Fil_D W_m \sO_U$ be the presheaf on $\Sch_X$ such that
      $\Fil_D W_m \sO_U(\wt{X}) = \Fil_D W_m(\wt{U})$. One easily checks that
      $\Fil_D\sO_U=\sO_X(D)$ and $\Fil_D W_m \sO_U$ actually restricts
      to a sheaf on $X_\et$. For $q \ge 0$, we define
  \[
  \Fil_D W_m \Omega^q_U =  \Fil_D W_m \sO_U \cdot W_m\Omega^q_X(\log E) +
  d(\Fil_D W_m \sO_U) \cdot W_m\Omega^{q-1}_X(\log E),
  \]
considered as a subsheaf of $j_*W_m\Omega^q_U$.
\end{defn}

{\bf {Notation.}}
In the rest of this paper, we shall use $\Fil_\n W_m \Omega^q_U$ as an alternative
notation for $\Fil_{E_\n} W_m \Omega^q_U$, where $\n = (n_1, \ldots , n_r)$.
Similarly, we shall write $\Fil_{E_\n} H^{q}_{m}(U)$ alternatively
as $\Fil_\n H^{q}_{m}(U)$ for any $1 \le m \le \infty$ whenever it is convenient.
We shall write $\Fil_\n W_1 \Omega^q_U$ as $\Fil_\n \Omega^q_U$. 
  For an hdvf $K$ containing $\F_p$ with uniformizer $\pi$, we shall write
  $\Fil_\n W_m\Omega^q_U$ as $\Fil_n W_m\Omega^q_K$ if $X=\Spec(\sO_{K})$,
  $E=V((\pi))$ and $E_\n = nE$ with $n \in \Z$.
  We shall write $\Fil_\n H^{q}_{m}(K)$ as $\Fil_n H^{q}_{m}(K)$.

\begin{lem}\label{lem:Fid-DRW-basic}
The sheaf $\Fil_\n W_m \Omega^q_U$ satisfies the following properties.
\begin{enumerate}
\item
    $d(\Fil_\n W_m \Omega^q_U) \subset \Fil_\n W_m \Omega^{q+1}_U$.
 \item
    $F(\Fil_\n W_{m+1} \Omega^q_U) \subset \Fil_\n W_m \Omega^{q}_U$.
  \item
    $V(\Fil_\n W_{m} \Omega^q_U) \subset \Fil_\n W_{m+1} \Omega^{q}_U$.
  \item
    $R(\Fil_\n W_{m+1} \Omega^q_U) = \Fil_{\n/p} W_m \Omega^{q}_U$.
  \item    $\varinjlim_n \Fil_{nE} W_m \Omega^q_U \xrightarrow{\cong}
    j_*W_m\Omega^{q}_U$.
  \item
  $\Fil_\n W_{m} \Omega^q_U \subset \Fil_{\n'} W_{m} \Omega^q_U$ if
$\n \le \n'$ and $\Fil_\n W_{m} \Omega^q_U \subset  W_{m} \Omega^q_X$ if
    $\n < 0$.
  \item
  $\Fil_0 W_m\Omega^q_U = W_m\Omega^q_X(\log E)$, where $0 = (0, \ldots , 0) \in \Z^r$.
  \item
    $\Fil_\n \Omega^q_U = \Omega^q_X(\log E)(E_\n)$.
    \item
 $\Fil_{-E}W_m\Omega^q_U = W_m\Omega^q_X$ if $q$ is the rank of $\Omega^1_X$ as a
    locally free sheaf on $X$.
      \end{enumerate}
\end{lem}
\begin{proof}
See \cite[Lem.~3.6, 3.10, Thm.~6.19]{KM-1}.
\end{proof}

We let $Z_1 \Fil_\n W_m\Omega^q_U =  \Fil_\n W_m\Omega^q_U \bigcap
j_* Z_1W_m\Omega^q_U$, where 
$Z_1W_m\Omega^q_U := F(W_{m+1}\Omega^q_U) =
\Ker(F^{m-1}d \colon W_m\Omega^q_U \to \Omega^{q+1}_U)$.
Let $C \colon Z_1W_m\Omega^q_U \to W_m\Omega^q_U$ denote the classical Cartier
operator (cf. \cite[\S~10]{GK-Duality}).

\begin{thm}\label{thm:V-R-C}
  There exists a short exact sequence
  \[
 0 \to V^{m-1}(\Fil_\n \Omega^q_U) + dV^{m-1}(\Fil_\n \Omega^{q-1}_U) \to
    \Fil_\n W_m\Omega^q_U \xrightarrow{R}  \Fil_{\n/p} W_{m-1}\Omega^q_U \to 0.
    \] 
The Cartier operator $C$ restricts to a map
$C \colon Z_1 \Fil_\n W_m\Omega^q_U \to  \Fil_{\n/p} W_m\Omega^q_U$
satisfying the following.
\[
    \ov{F} \colon \Fil_{{\n}/p} W_{m}\Omega^q_U \xrightarrow{\cong}
    \frac{Z_1\Fil_{\n} W_m\Omega^q_U}{dV^{m-1}(\Fil_{\n} \Omega^{q-1}_U)}; \ \
    \ov{C} \colon  \frac{Z_1\Fil_{\n} W_m\Omega^q_U}{dV^{m-1}(\Fil_{\n}\Omega^{q-1}_U)}
     \xrightarrow{\cong} \Fil_{{\n}/p}W_{m}\Omega^q_U.
    \]
\end{thm}
\begin{proof}
See \cite[Thm.~6.19]{KM-1}.
\end{proof}

\begin{defn}\label{defn:1-c complex} 
  If $q \in \N_0$, and $\n \in \Z^r$ is such that $n_i + 1 \ge 0$ for every $i$,
  we define
\[
W_m\sF^{q,\bullet}_{\n} =
\left(Z_1\Fil_{\n} W_m\Omega^q_U \xrightarrow{1-C} \Fil_{\n} W_m\Omega^q_U\right)
\]
as an object of $\sD_\et(X)$.
\end{defn}

The following result is a generalization of the classical $V$-$R$ sequence
for $W_m\Omega^q_{X, \log}$.
\begin{thm}\label{thm:V-R exact}
  There exists a distinguished triangle
    
    \[
     W_1\sF^{q, \bullet}_{\n} \xrightarrow{V^{m-1}} W_m\sF^{q,\bullet}_{\n}
     \xrightarrow{R} W_{m-1}\sF^{q,\bullet}_{\n/p} \xrightarrow{+}  W_1\sF^{q, \bullet}_{\n}
                 [1]
    \]
    in $\sD_\et(X)$.
    There is a short exact sequence 
\begin{equation}\label{eqn:fil_0-log}
    0 \to j_* W_m\Omega^q_{U,\log} \to Z_1W_m\Omega^q_X(\log E)
    \xrightarrow{1 - C}
    W_m\Omega^q_X(\log E) \to 0.
\end{equation}
\end{thm}
\begin{proof}
  See \cite[Thm.~6.19, Lem.~11.2]{KM-1}.
  \end{proof}

The main result that will be crucial for us is the following cohomological
description of Kato's filtration on $H^q(U)$. 

\begin{thm}\label{thm:H^1-fil}
  For any $q \in \N_0, m \in \N$ and $\n \in \N^r_0$, there exists a functorial
  isomorphism  
\[
\Fil_\n H^{q+1}_{p^m}(U) \xrightarrow{\cong} \H^1_\et(X, W_m\sF^{q,\bullet}_\n).
\]
In particular, if $n \in \N, X=\Spec(\sO_K)$ and $E=V((\pi))$, where $\sO_K$ is an
$F$-finite hdvr with quotient field $K$ and uniformizer $\pi$, we have
an exact sequence
\[
0\to W_m\Omega^{q}_{K,\log}\to Z_1\Fil_{n}W_m\Omega^{q}_K \xrightarrow{1-C}
\Fil_{n}W_m\Omega^{q}_K \xrightarrow{\delta^q_m} \Fil_nH^{q+1}_{p^m}(K)\to 0.
\]
 \end{thm}
\begin{proof}
  See \cite[Thm.~8.7]{KM-1}.
  \end{proof}

The above results show that the complex  $W_m\sF^{q,\bullet}_{\n}$ generalizes 
the classically known sheaves
$W_m\Omega^q_{X,\log}$ and $j_* W_m\Omega^q_{U,\log}$. The latter sheaves make it
possible to study the unramified and tamely ramified cohomology classes in
$H^{q+1}(U)$ using sheaf-theoretic methods, whereas the former provides a framework
for studying cohomology classes with arbitrary ramification.

An immediate consequence of \thmref{thm:H^1-fil} is the following.

\begin {cor}\label{cor:Fil-functorial}
Let $f \colon (X',E') \to (X,E)$ be a morphism of snc-pairs and let
$U' = X'\setminus E'$. Let $\n \in \N^r_0$ and let $\n' = (n'_1, \ldots , n'_s)
\in \N^r_0$, where $E'_1, \ldots , E'_s$ are the irreducible components of $E'$
and $f^*(\sum_i n_i E_i) = \sum_j n'_j E'_j$. Then $f^* \colon
H^q(U) \to H^q(U')$ restricts to a homomorphism
$f^* \colon \Fil_\n H^{q}_{m}(U) \to \Fil_{\un{m}} H^{q}_{m}(U')$ for every $m,q \in \N$
and $\un{m} \ge \n'$.
\end{cor}
\begin{proof}
 See \cite[Cor.~8.8]{KM-1}
\end{proof}

\subsection{The refined Swan conductor for schemes}\label{sec:RSC}
If $A$ is an $F$-finite hdvr with maximal ideal $(\pi)$,
quotient field $K$ and residue field $\ff$, then Kato constructed a refined Swan
conductor map
\begin{equation}\label{eqn:Kato-SC}
  {\rsw}^q_{K,n} \colon \gr_nH^{q}(K) \to \pi^{-n} \Omega^{q}_A(\log \pi)\otimes_A \ff
  \xrightarrow{\cong}  \pi^{-n}
  \left(\Omega^q_\ff \bigoplus \Omega^{q-1}_\ff \dlog(\pi)\right)
  \end{equation}
in \cite[\S~5]{Kato-89} for $n, q \in \N$ and showed that it is injective if 
one lets $\gr_nH^{q}(K) = \frac{\Fil_nH^{q}(K)}{\Fil_{n-1}H^{q}(K)}$.

We now let $(X,E)$ be an snc-pair as in ~\S~\ref{sec:Fil-DRW}. Let
$n = (n_1, \ldots , n_r), \n' = (n'_1, \ldots , n'_r) \in \N^r_0$ be such that
${\n}/p \le \n' \le \n$. We let $F = E_\n - E_{\n'} \in \Div^{\eff}_E(X)$.
One of the main results of \cite{KM-1} is the following
generalization of Kato's refined Swan conductor.
Recall that a strict morphism of ind-abelian (resp. pro-abelian) groups
$f \colon \{A_n\} \to \{B_n\}$ is the one which preserves the level, i.e.,
$f = (f_n)$, where $f_n \colon A_n \to B_n$ for every $n$.

\begin{thm}\label{thm:RSW-gen}
For $m, q \in \N$, we have the following.
  
$(1)$ There exists a strict morphism of ind-abelian groups 
  \begin{equation}\label{eqn:RSW-gen-0}
  \left\{ \frac{\Fil_{E_\n} H^{q}_{p^m}(U)}{\Fil_{E_{\n'}} H^{q}_{p^m}(U)}\right\}_{m \ge 1}
      \xrightarrow{\Rsw^{m,q}_{X|(\n,\n')}} 
H^0_\zar\left(F, \Omega^{q}_X(\log E)(E_\n) \otimes_{\sO_X} \sO_F\right)
  \end{equation}
  which is a levelwise monomorphism.

  $(2)$ Let $f \colon (X',E') \to (X,E)$ be a morphism of snc-pairs with
  $U' = X'\setminus E'$. Then letting $F = E_\n - {E_\n}/p$ and
  $F' = f^*(E_\n) - {f^*(E_\n)}/p$,
  we have a commutative diagram of ind-abelian groups
\begin{equation}\label{eqn:RSW-gen-1-0}
  \xymatrix@C1pc{
    \left\{ \frac{\Fil_{E_\n} H^{q}_{p^m}(U)}{\Fil_{{E_\n}/p} H^{q}_{p^m}(U)}\right\}_{m \ge 1}
    \ar[r] \ar[d]_-{f^*} &
    H^0_\zar\left(F, \Omega^{q}_X(\log E)(E_\n) \otimes_{\sO_X} \sO_{F}\right)
    \ar[d]^-{f^*} \\
    \left\{ \frac{\Fil_{f^*(E_\n)} H^{q}_{p^m}(U')}
            {\Fil_{{f^*(E_\n)}/p} H^{q}_{p^m}(U')}\right\}_{m \ge 1} \ar[r] &
H^0_\zar\left(F', \Omega^{q}_{X'}(\log E')(f^*(E_\n)) \otimes_{\sO_{X'}} \sO_{F'}\right),}
\end{equation}
where the horizontal arrows are the refined Swan conductor maps of
~\eqref{eqn:RSW-gen-0}.

$(3)$ If $f^*(E)$ is reduced, then we can replace ${E_\n}/p$ (resp. ${f^*(E_\n)}/p$) in
~\eqref{eqn:RSW-gen-1-0} by any $D'$ (resp. $f^*(D')$) such that
$E_\n \ge D' \ge {E_\n}/p \ge 0$.

$(4)$ If $X=\Spec(\sO_K)$, $E= \Spec(\ff)$, $E_\n = nE$ and $E_{\n'} = (n-1)E$
with $n \in \N$, then
$\varinjlim_{m}\Rsw^{q,m}_{X|(\n,\n')}$ coincides with the map $ {\rsw}^q_{K,n}$ in \eqref{eqn:Kato-SC}. 
\end{thm}
\begin{proof}
  See \cite[Thm.~9.4]{KM-1}.
\end{proof}

When $X = \Spec(A)$ is affine, the refined Swan conductor can be
described by the following commutative diagram.

\begin{equation}\label{eqn:RSW-gen-1-0-0}
  \xymatrix@C1pc{
    & \Fil_\n W_m\Omega^{q-1}_{U} \ar@{->>}[r]^-{\delta^{q-1}_m}
    \ar[d]^-{(-1)^{q}F^{m-1}d}
    \ar@{->>}[drr]^-{\can} \ar[dl] & \Fil_\n H^{q}_{p^m}(U) 
    \ar@{->>}[r] & \frac{\Fil_\n H^{q}_{p^m}(U)}{\Fil_{\n'} H^{q}_{p^m}(U)}
    \ar[d]^-{\alpha^{m,q-1}_{X|(\n, \n')}} \\
   H^0_\zar(F, \Omega^q_X(\log E)(E_\n) \otimes_X \sO_F) &
\frac{\Fil_\n \Omega^{q}_{U}}{\Fil_{\n'}\Omega^{q}_{U}} \ar[l]_-{\cong} & &
\frac{\Fil_\n W_m\Omega^{q-1}_{U}}{\Fil_{\n'} W_m\Omega^{q-1}_{U} +
  Z_1 \Fil_\n W_m\Omega^{q-1}_{U}} \ar[ll]^-{(-1)^{q}F^{m-1}d}.}
  \end{equation}
In this diagram, the right diagonal arrow is the canonical surjection,
$\alpha^{m,q-1}_{X|(\n, \n')}$ is an injective map induced by the isomorphism given in
\thmref{thm:H^1-fil} (cf. \cite[(9.7)]{KM-1}). In particular,
$\alpha^{m,q-1}_{X|(\n, \n')}$ is a bijection. We have
$\Rsw^{m,q}_{X|(\n,\n')} = (-1)^{q}F^{m-1}d \circ \alpha^{m,q-1}_{X|(\n, \n')}$.
This diagram is compatible with the change in values of $m$.

We let $\Rsw^{q}_{X|(\n,\n')} := {\varinjlim}_m \Rsw^{m,q}_{X|(\n,\n')}$.
We shall refer to $\Rsw^{m,q}_{X|(\n,\n')}$ and $\Rsw^{q}_{X|(\n,\n')}$
as the generalized refined Swan conductor (or simply the refined
Swan conductor) maps.
For an element $\chi \in \Fil_{E_\n} H^{q}_{p^m}(U)$, the notation
$\Rsw^{m,q}_{X|(\n,\n')}(\chi)$ will mean the value of the generalized
refined Swan conductor on the image of $\chi$ in the quotient
$\frac{\Fil_{E_\n} H^{q}_{p^m}(U)}{\Fil_{E_{\n'}} H^{q}_{p^m}(U)}$.

\vskip.2cm

A special case of \thmref{thm:RSW-gen} that will be useful to us is the
following.

\begin{cor}\label{cor:RSW-gen-1}
  Assume in \thmref{thm:RSW-gen} that $0 \le n'_i \le n_i -1$ for some
  $1 \le i \le r$. 
Let $f \colon (X',E') \to (X,E)$ be a morphism of snc-pairs
and let $U' = X'\setminus E'$. Let $F = E_{\n} - E_{\n'}$.
Let $D_\n = f^*(E_\n), \ E^*_i = f^*(E_i)_\red, \  
D'_\n \in \{D_\n - E^*_i\} \bigcup \{D_\n - E'_j| E'_j \subset E^*_i\}$ and
$F' = D_\n - D'_\n$.
Then one has a commutative diagram of ind-abelian groups
\begin{equation}\label{eqn:RSW-gen-1-1}
  \xymatrix@C1pc{
    \left\{\frac{\Fil_{E_\n} H^{q}_{p^m}(U)}{\Fil_{E_{\n'}} H^{q}_{p^m}(U)}\right\}_{m \ge 1}
    \ar[r] \ar[d]_-{f^*} &
    H^0_\zar\left(F, \Omega^{q}_X(\log E)(E_\n) \otimes_{\sO_X} \sO_{F}\right)
    \ar[d]^-{f^*} \\
    \left\{ \frac{\Fil_{D_\n} H^{q}_{p^m}(U')}
            {\Fil_{D'_\n} H^{q}_{p^m}(U')}\right\}_{m \ge 1} \ar[r] &
H^0_\zar\left(F', \Omega^{q}_{X'}(\log E')(D_\n) \otimes_{\sO_{X'}} \sO_{F'}\right),}
\end{equation}
where the horizontal arrows are the refined Swan conductor maps of
\eqref{eqn:RSW-gen-0}.
\end{cor}
\begin{proof}
Since the diagram
  \begin{equation}\label{eqn:RSW-gen-1-5}
    \xymatrix@C1pc{
\left\{ \frac{\Fil_{D_\n} H^{q}_{p^m}(U')}
            {\Fil_{D_\n - E^*_i} H^{q}_{p^m}(U')}\right\}_{m \ge 1} \ar[r] \ar@{->>}[d] &
H^0_\zar\left(E^*_i, \Omega^{q}_{X'}(\log E')(D_\n) \otimes_{\sO_{X'}} \sO_{E^*_i}\right)
\ar[d] \\
\left\{ \frac{\Fil_{D_\n} H^{q}_{p^m}(U')}
            {\Fil_{D_\n - E'_j} H^{q}_{p^m}(U')}\right\}_{m \ge 1} \ar[r] &
 H^0_\zar\left(E'_j, \Omega^{q}_{X'}(\log E')(D_\n) \otimes_{\sO_{X'}} \sO_{E'_j}\right)}
  \end{equation}
  (where the vertical arrows are the canonical restrictions)
  clearly commutes for every $E'_j$ in the support of $E^*_i$,
  it suffices to prove the corollary when $D'_\n = D_\n - E^*_i$.
  Moreover, it suffices to prove this commutativity  for each $m \ge 1$.

We now let $F_0 = E_\n - {E_\n}/p$ and $F^*_0 = D_\n - {D_\n}/p$. There are
  inequalities
  \begin{equation}\label{eqn:RSW-gen-1-2}
  0 \le f^*({E_\n}/p) \le {D_\n}/p \le D_\n, \ 0 \le f^*(E_{\n'}) \le D'_\n \le D_\n, \
  0 \le {D'_{\n}}/p \le D'_{\n} \le D_\n.
  \end{equation}
  We let
  \[
  \sF_{\n/p} = \Omega^{q}_X(\log E)(E_\n) \otimes_{\sO_X} \sO_{F_0}, \
  \sF_{\n'} = \Omega^{q}_X(\log E)(E_\n) \otimes_{\sO_X} \sO_{F},
  \]
  \[
  \sG_{\n/p} = \Omega^{q}_{X'}(\log E')(D_\n) \otimes_{\sO_{X'}} \sO_{F^*_0} \
  \ {\rm and} \ \
  \sG_\n =  \Omega^{q}_{X'}(\log E')(D_n) \otimes_{\sO_{X'}} \sO_{F'}.
  \]
  By ~\eqref{eqn:RSW-gen-1-2}, there are canonical maps
\begin{equation}\label{eqn:RSW-gen-1-3}
  \sG_\n \leftarrow f^*(\sF_{\n'}) \twoheadleftarrow f^*(\sF_{\n/p}) \to \sG_{\n/p} \surj
  \sG_\n.
 \end{equation}  
 
We now look at the diagram
  \begin{equation}\label{eqn:RSW-gen-1-4}
\xymatrix@C1pc{
  \frac{\Fil_{E_\n} H^{q}_{p^m}(U)}{\Fil_{{E_\n}/p} H^{q}_{p^m}(U)} \ar[rr] \ar[dd]
  \ar@{->>}[dr] &
  & H^0_\zar(X, \sF_{\n/p}) \ar[dd]
  \ar[dr] & \\ 
  & \frac{\Fil_{E_\n} H^{q}_{p^m}(U)}{\Fil_{E_{\n'}} H^{q}_{p^m}(U)} \ar[rr]
    \ar[dd] & &  H^0_\zar(X, \sF_{\n'})
    \ar[dd] \\
 \frac{\Fil_{D_\n} H^{q}_{p^m}(U')}{\Fil_{{D_\n}/p} H^{q}_{p^m}(U')}   
\ar[dr] \ar[rr] & & 
  H^0_\zar(X',\sG_{\n/p})
  \ar[dr] & \\
  & \frac{\Fil_{D_\n} H^{q}_{p^m}(U')}{\Fil_{D'_\n} H^{q}_{p^m}(U')} \ar[rr] & &
  H^0_\zar(X', \sG_\n)}
  \end{equation}
  in which the horizontal arrows are the generalized refined Swan conductor maps and
  the vertical arrows are the canonical maps induced by ~\eqref{eqn:RSW-gen-1-2} and
   ~\eqref{eqn:RSW-gen-1-3}. One checks that all arrows are defined.
The left and the right faces of this cube commute by the functoriality.
  The back, top and bottom faces commute by \thmref{thm:RSW-gen}. A diagram
  chase shows that the front face also commutes, as desired.
\end{proof}

The following result (cf. \cite[Lem.~9.9]{KM-1}) provides a simplified description of
$\frac{\Fil_{E_\n} \Omega^{q}_U}{\Fil_{E_{n'}} \Omega^{q}_U} \cong
\Omega^{q}_X(\log E)(E_\n) \otimes_{\sO_X} \sO_{F}$ when
$E_{\n'} = E_\n - E_i$ in \corref{cor:RSW-gen-1}. We let
$F_i =  \sum\limits_{\substack{j \in J_r \setminus \{i\}}} (E_i \bigcap E_j) \in
\Div^{\eff}(E_i)$ and
$D_i = E_\n \times_X E_i$. Let $u_i \colon E_i \inj X$ be the inclusion.

\begin{lem}\label{lem:RSW-spl}
There exists an exact sequence
  \[
  0 \to \Omega^{q}_{E_i}(\log F_i)({D}_i) \xrightarrow{\nu'_i}
  \Omega^{q}_X(\log E)(E_\n)\otimes_{\sO_X} \sO_{E_i} \xrightarrow{\nu_i}
  \Omega^{q-1}_{E_i}(\log F_i)({D}_i) \to 0,
  \]
  where $\nu_i$ is induced by the Poincar{\'e} residue map
  $\Omega^{q}_X(\log E)\xrightarrow{\res} \Omega^{q-1}_{E_i}(\log F_i)$.
\end{lem}

\subsection{Matsuda filtration on higher cohomology groups}\label{sec:MAF}
In \cite{Matsuda}, Matsuda introduced a non-logarithmic version of Kato's filtration
on $H^1(K)$ for an $F$-finite hdvf $K$ of characteristic $p$. However, his methods
were not sufficiently general to construct analogous filtrations on $H^q(K)$ for
$q \ge 2$. Using the generalized refined Swan conductor, we obtain such an
extension. This subsection is devoted to its construction and basic properties. The
resulting higher-degree analogue of Matsuda's filtration will play a central role
in the proof of \thmref{thm:Main-1}. We begin by introducing the following notion.

We let $K$ be an $F$-finite hdvf of characteristic $p$ with the ring of integers
$\sO_K$ which is an excellent ring. Let $\pi$ be a uniformizer of $K$ and
let $\ff$ denote its residue field. Let $m, q \in \N$ and $\chi \in H^{q}_{p^m}(K)$.
Suppose that $\Sw(\chi) = n \ge 1$.
\begin{defn}\label{defn:Type-I-local}
We shall say that $\chi$ has type I if the image of $\chi$
is not zero under the composite map (cf. ~\eqref{eqn:Kato-SC})
  \[
  \Fil_n H^q_{p^m}(K) \surj \frac{\Fil_n H^q_{p^m}(K)}{\Fil_{n-1} H^q_{p^m}(K)}
  \xrightarrow{\rsw_{K,n}^{m,q}} \pi^{-n}\Omega^q_{\sO_K}(\log \pi)\otimes_{\sO_K}\ff
  \xrightarrow{\nu}
  \pi^{-n} \Omega^{q-1}_\ff \dlog(\pi).
  \]
Otherwise, we shall say that $\chi$ has type II.
\end{defn}

Let us now assume that $(X, E)$ is an snc-pair as in \S~\ref{sec:Fil-DRW}. Let
$q \in \N$ and $\chi \in H^{q}_{p^m}(U)$. Let $i \in J_r$ be such that
$\Sw_{E_i}(\chi) = n_i \ge 1$.
\begin{defn}\label{defn:Type-I-global}
  We shall say that $\chi$ has type I (resp. type II) at $E_i$ if the restriction of
  $\chi$ to $\Spec(K_{y_i})$ has type I (resp. type II), where $y_i$ is
  the generic point of $E_i$ and $K_{y_i}$ is the henselization of the function field
  of $X$ at $y_i$.
\end{defn}

The following result follows by comparing the exact sequence of
\lemref{lem:RSW-spl} for $(X,E)$ and $(\Spec(\sO_{K_{y_i}}), \Spec(\ff))$
via the canonical map of snc-pairs $(\Spec(\sO_{K_{y_i}}), \Spec(\ff)) \to
(X,E)$ and using that the map $H^0_\zar(E_i, \Omega^{q-1}_{E_i}(\log F_i)({D}_i)) \to
\Omega^{q-1}_{k(y_i)}$ is injective.

\begin{cor}\label{cor:Type-I-global-0}
  We let $i \in J_r$.
  Let $D = \stackrel{r}{\underset{j =1}\sum} n_j E_j \in \Div^{\eff}_E(X)$
  and let $\chi \in \Fil_D H^{q}_{p^m}(U)$ be such that $\Sw_{E_i}(\chi) = n_i \ge 1$.
Then $\chi$ has type I (resp. type II) at $E_i$ if and only if
  $\nu^0_i \circ \Rsw_{X|(D, D - E_i)}^{m,q}(\chi) \neq 0$ (resp. $= 0$), where
  $\nu^0_i \colon H^0_\zar(X, \Omega^{q}_X(\log E)(E_\n)\otimes_{\sO_X} \sO_{E_i}) \to
  H^0_\zar(E_i, \Omega^{q-1}_{E_i}(\log F_i)({D}_i))$ is the map induced by $\nu_i$
  (cf. \lemref{lem:RSW-spl}) on the spaces of global sections.
\end{cor}

Let $(X, E)$ be an snc-pair as in \S~\ref{sec:Fil-DRW} and let
$D = \sum_i n_i E_i \in \Div^{\eff}_E(X)$. \lemref{lem:RSW-spl}
allows us to define the following filtration on $H^q_{p^m}(U)$ for $m \ge 1$ and
for $ H^q(U)$.

\begin{defn}\label{defn:Matsuda-0}
  We let $\Fil'_D H^q_{p^m}(U)$ be the subgroup of $\Fil_{D+E} H^q_{p^m}(U)$
  consisting of elements which die under the composite diagonal map
  \[
  \xymatrix@C1pc{
    \Fil_{D+E} H^q_{p^m}(U) \ar[r] \ar[drr] &
    \frac{\Fil_{D+E} H^q(U)}{\Fil_{D+E-E_i} H^q(U)} \ar[r] &
    H^0_\zar(E_i, \Omega^{q}_X(\log E)(D+ E)\otimes_{\sO_X} \sO_{E_i})
    \ar[d]^-{\nu^0_i} \\
    & & H^0_\zar(E_i, \Omega^{q-1}_{E_i}(\log F_i)({(D+E)}_i))}
\]
for every $i \in J_r$, where the right horizontal arrow on the top is the
generalized refined Swan conductor map ${\Rsw^{m,q}_{X|(D+E,D+E-E_i)}}$.

If $m = p^tr$ is an integer with $t \in \N_0$ and $r \in I_p$, we let
\[
\Fil'_D H^q_m(U) = \Fil'_D H^q_{p^m}(U) \bigoplus H^q_{r}(U);
\]
\[
\Fil'_D H^q(U) = {\varinjlim}_{m} \Fil'_D H^q_m(U) \cong
\left({\varinjlim}_{m} \Fil'_D H^q_{p^m}(U) \right)
\bigoplus H^q(U)\{p'\}.
\]
\end{defn}
We shall refer to $\Fil'_D H^q_m(U)$ (resp. $\Fil'_D H^q(U)$) as the {\sl Matsuda
filtration} of $H^q_m(U)$ (resp. $H^q(U)$).

To define the Matsuda filtration on $\Br(U)$, we note that according
to \lemref{lem:Br-fil-0}, all refined Swan conductor maps
$\Rsw^{m,2}_{X|(D,D')}$ (where $D \ge D' \ge D/p \ge 0$) are actually
defined on the quotients $\frac{\Fil_D \Br(U)[m]}{\Fil_{D'} \Br(U)[m]}$.
We can therefore define $\Fil'_D \Br(U)[m]$ in the same manner as 
$\Fil'_D H^2_{m}(U)$ above.
It is then immediate that there is a functorial exact sequence
\begin{equation}\label{eqn:Br-fil-Matsuda}
  0 \to {\Pic(U)}/m \to \Fil'_D H^2_{m}(U) \to \Fil'_D \Br(U)[m] \to 0
\end{equation}
for every $D \in \Div^{\eff}_E(X)$ and $m \in \N$.

\vskip .2cm

The following is clear from the definition of $\Fil'_D H^q_{p^m}(U)$ using the
observation that $H^0_\zar(E_i, \Omega^{q-1}_{E_i}(\log F_i)({(D+E)}_i)) \to
\Omega^{q-1}_{k(y_i)}$
is injective for every $i$ and $q$, where $y_i$ is the generic point of $E_i$.

\begin{cor}\label{cor:Matsuda-1}
Let $(X, E)$ be an snc-pair as in \S~\ref{sec:Fil-DRW} and let
$D = \sum_i n_i E_i \in \Div^{\eff}_E(X)$. Then
an element $\chi \in  H^q_{p^m}(U)$ lies in $\Fil'_D H^q_{p^m}(U)$ if and only if
its image under the map $H^q_{p^m}(U) \to H^q_{p^m}(K_{y_i})$ lies in
$\Fil'_{n_i} H^q_{p^m}(K_{y_i})$ for every $i$, where $K_{y_i}$ is the henselization of
the function field of $X$ at $y_i$.
\end{cor}

Applying this result for $q = 1$, we get the following.
  
\begin{cor}\label{cor:Matsuda-4}
  If $K$ is an $F$-finite hdvf with residue field $\ff$ and $n \in \N_0$, then
  $\Fil'_n H^1(K)$ coincides with the ramification filtration on
  $H^1(K)$ defined by Matsuda.
\end{cor}
\begin{proof}
  It is shown in the proof of \cite[Thm.~6.3(2)]{GK-JA} that if $\Fil^{\ms}_n H^1(K)$
  denotes Matsuda's filtration, then
  $\Fil^{\ms}_n H^1(K)$ is the subgroup of $\Fil_{n+1} H^1(K)$ consisting of characters
  whose refined Swan conductor dies in $H^0_\zar(\ff, \ff (\pi^{-n} \dlog(\pi))) \cong
  \ff$, where $\pi$ is a uniformizer of $K$.
 We now now apply \corref{cor:Matsuda-1} to conclude the proof.
\end{proof}

\vskip .2cm

In this paper, we shall use the following results about the functorial property of
the Matsuda filtration and its relation to the Kato filtration.

\begin{cor}\label{cor:Matsuda-2}
One has inclusion $\Fil_D H^{q+1}_{p^m}(X) \subset \Fil'_D H^{q+1}_{p^m}(X) \subset \Fil_{D+E} H^{q+1}_{p^m}(X)$ for any
$D \in \Div^{\eff}_E(X)$, and the equality $\Fil_D H^{q+1}_{p^m}(X) = \Fil'_D H^{q+1}_{p^m}(X)$ holds if $p\nmid (n_i+1)$ for all $i$, where $D=\sum_i n_i E_i$. In
particular, $\Fil'_0 H^{q+1}_{p^m}(X) = \Fil_0 H^{q+1}_{p^m}(X)$.
  If each $E_i$ is a finite type scheme over a perfect
  field and has dimension $q$, then $\Fil'_D H^{q+1}_{p^m}(X) = \Fil_D H^{q+1}_{p^m}(X)$.
  \end{cor}
\begin{proof}
 The equality $\Fil_D H^{q+1}_{p^m}(X) = \Fil'_D H^{q+1}_{p^m}(X)$, when $p \nmid (n_i+1)$ for all $i$, follows easily from \corref{cor:Matsuda-1} and \cite[Thm.~3.2(3)]{Kato-89}. Everything else is clear
  from \thmref{thm:RSW-gen}, \lemref{lem:RSW-spl} and \corref{cor:Matsuda-1}.
  \end{proof}

\begin{prop}\label{prop:Matsuda-3}
  Let $f \colon (X',E') \to (X,E)$ be a morphism of snc-pairs and let
  $\wt{E} = f^*(E)_{\red}$. Let $U' = X' \setminus E'$.
  Assume that
  $f^*(E_i)$ is irreducible for each $i$ and $f^*(E_i)_\red \neq f^*(E_j)_\red$ if
  $i \neq j$. Let $D \in \Div^{\eff}_E(X)$ and
  $D' = f^*(D+E) - \wt{E}$. Then $f^* \colon H^q_{p^m}(U) \to H^q_{p^m}(U')$
  induces $f^* \colon \Fil'_D H^q_{p^m}(U) \to \Fil'_{D'} H^q_{p^m}(U')$.
  In particular, $f^*$ takes $\Fil'_D H^q_{p^m}(U)$ into
  $\Fil'_{f^*(D)} H^q_{p^m}(U')$ if $f^*(E)$ is reduced.
\end{prop}
\begin{proof}
We let $\chi \in \Fil'_D H^q_{p^m}(U) \subset \Fil_{D+E} H^q_{p^m}(U)$.
  Let $E'_1, \ldots , E'_s$ be the irreducible components of $E' \subset X'$.
  We then note that $D', D' + E' - E'_j \in \Div^{\eff}_{E'}(X')$ for every
  $j \in J_s$.
  We fix one such $j$. If $E'_j$ is not a component of $\wt{E}$, then
  $E' - \wt{E} - E'_j \ge 0$ and hence $D' +E' -E'_j \ge f^*(D +E)$ which
  implies that
  $f^*(\chi) \in \Fil_{D' +E' -E'_j} H^q(U')$. In particular, $f^*(\chi)$ already
  dies in $\frac{\Fil_{D'+E'} H^q(U')}{\Fil_{D'+E'-E'_j} H^q(U')}$.

We assume now that $E'_j$ is a component of $\wt{E}$.
  In this case, there is unique $E_i$ such that $f^*(E_i)_\red = E'_j$.
  We let $n_j$ be the multiplicity of $E'_j$ in $\wt{E}$ and look at the diagram
\begin{equation}\label{eqn:Matsuda-3-0}
  \xymatrix@C.8pc{
    \frac{\Fil_{D+E} H^q(U)}{\Fil_{D+E-E_i} H^q(U)} \ar[d]_-{f^*} \ar[r] &
    H^0_\zar(E_i, \Omega^{q}_X(\log E)(M)\otimes_{\sO_X} \sO_{E_i})  \ar[r]^-{\nu^0_i}
    \ar[d]^-{f^*}  
    & H^0_\zar(E_i, \Omega^{q-1}_{E_i}(\log F_i)(M|_{E_i})) \ar[d]^-{n_j f^*} \\
    \frac{\Fil_{D'+E'} H^q(U')}{\Fil_{D'+E'-E'_j} H^q(U')} \ar[r] &
    H^0_\zar(E'_j, \Omega^{q}_{X'}(\log E')(M')\otimes_{\sO_{X'}} \sO_{E'_j})
    \ar[r]^-{\nu'^0_j} &
    H^0_\zar(E'_j, \Omega^{q-1}_{E'_j}(\log F'_i)(M'|_{E'_j})),}
\end{equation}
where $M = D+E$ and $M' = D'+ E'$.

The left square in ~\eqref{eqn:Matsuda-3-0} commutes by
\corref{cor:RSW-gen-1} (this uses that $0 \le f^*(D+E) -E'_j \le D'+E'-E'_j$).
To prove the corollary, it suffices therefore to show that the right square is
commutative. But this can be easily checked by computing the residues
locally because all the underlying sheaves are locally free. We leave out the
details.
\end{proof}

\begin{cor}\label{cor:Matsuda-5}
  The assertions of Corollaries~\ref{cor:Matsuda-1}, 
  ~\ref{cor:Matsuda-2} and Proposition~\ref{prop:Matsuda-3} remain true if
  we replace $H^2_m(U)$ by $\Br(U)[m]$ in their statements.
\end{cor}
\begin{proof}
  Using the results for  $H^2_m(U)$ proven above, the corollary is a direct
  consequence of definitions of the Kato and Matsuda filtrations on $H^2_m(U)$ and
  $\Br(U)[m]$ once we use Lemma~\ref{lem:Br-fil-0} and
  ~\eqref{eqn:Br-fil-Matsuda}.
  \end{proof}

\section{Kato Complex and applications}\label{sec:Kato-complex}
For an excellent Noetherian scheme, Kato introduced in \cite[\S~1]{Kato-hasse}
(see also \cite[\S~7.5]{Kato-89}) a Bloch--Ogus type complex, now commonly referred
to as the \emph{Kato complex}. In this section, we recall the special case of this
complex that is relevant to our work and establish several applications that will
play a crucial role in the proof of \thmref{thm:Main-1}.

\subsection{The Kato complex}\label{sec:Kato-complex-0}
Let $A$ be an excellent regular Noetherian local domain of dimension two containing
$\F_p$. Let $K$ (resp. $\ff$) denote the quotient (resp. residue) field of $A$.
Let $x$ denote the closed  point of $X=\Spec(A)$. We assume
that $[\ff : \ff^p]=p^{c-1}$ for some integer $c \ge 1$. This is equivalent to the
fact that $\dim_\ff \Omega^1_\ff=c-1$ (cf. \cite[Thm.~26.5]{Matsumura}).
In this special case, the Kato complex $K(A, c)$ has the form:
\begin{equation}\label{eqn:Kato com}
  H^{c+2}(K) \xrightarrow{\partial=\bigoplus \partial_z}
  \bigoplus\limits_{z \in X^{(1)}} H^{c+1}(k(z))
  \xrightarrow{\partial'=\bigoplus \partial'_z} H^c(\ff).
\end{equation}

One checks using \lemref{lem:Rank-module-diff} that $[k(z) : k(z)^p] = p^c$ for
every $z \in X^{(1)}$ and $[K : K^p] = p^{c+1}$. 
The boundary maps $\partial$ and $\partial'$ are defined as follows.
Let $z \in  X^{(1)}$ and let $K_z = Q(\sO_{X,z}^h)$. Then $K_z$ is an hdvf
with the residue field $k(z)$. We let $\pi$ be its uniformizer.
As $[K_z: K^p_z] = p^{c+1}$, \propref{prop:Kato-basic} implies that
$H^{c+2}(K_z) = \Fil_0H^{c+2}(K_z)$ and we have an isomorphism
$\lambda^{c+2}(\pi) \colon  H^{c+2}(k(z)) \bigoplus H^{c+1}(k(z)) \xrightarrow{\cong}
H^{c+2}(K_z)$.
This  allows one to define $\partial_z$ as the composite 
\[
H^{c+2}(K) \xrightarrow{} H^{c+2}(K_z)
\xrightarrow{\pr_2 \circ (\lambda^{c+2}(\pi))^{-1}} H^{c+1}(k(z)),
\]
where the first map is the restriction map in Galois cohomology.

To define $\partial'_z$, we let $Z$ denote the normalization of $\ov{\{z\}}$. Then
$k(z)= Q(\sO_{Z,x'})$ for each $x'\in Z_{(0)}$. Using that
$[k(x'):k(x')^p]=p^{c-1}$ (cf. \cite[Thm.~25.1, 26.10]{Matsumura}), we can define
$\partial^z_{x'} \colon H^{c+1}(k(z))\to H^{c}(k(x'))$ exactly
as we defined $\partial_z$ above. We then let
$\partial_{z}'= \sum\limits_{x' \in Z_{(0)}}\Cor_{k(x')/{\ff}} \circ \partial^z_{x'}$,
where $\Cor_{k(x')/{\ff}}$ is the corestriction map in Galois cohomology.

\vskip .3cm

Let $A$ be as above and let $\fm=(\pi,t)$ be the maximal ideal of $A$.
Let $A_1 = A/(\pi), \ E = \Spec(A_1)$ and let $y$ be the generic
point of $E$. We let $A_y = \sO_{X,y}^h, \ K_{y} = Q(A_y)$ and $F = Q(A^h_1)$.
We fix an integer $q \in J_c$ and an element
$\chi \in H^{q+1}(A[\pi^{-1}])\{p\}$ such that $\Sw_{E}(\chi) = n \ge 1$.
We let $n' = n-1$ and denote $\Rsw_{X,(n,n')}^{q+1}(\chi)$ by
$\Rsw_{y,n}(\chi)$. We let $\nu^0 \colon (\pi^{-n}) \Omega^{q+1}_A(\log(\pi))
\cong  \Omega^{q+1}_A(\log(\pi)) \to \Omega^q_{A_1} \cong
\frac{(\pi^{-n})}{(\pi^{-n+1})} \Omega^q_{A_1}$
be the  map $\nu_i$ of \lemref{lem:RSW-spl} after taking global sections.

\subsection{Application I: Non-vanishing of elements in $H^*(K_y)$}\label{sec:Non-V}
In this subsection, our goal is to show using the Kato complex
that certain elements in
$H^i(K)$ for $i = q+2, q+3$ do not die when we pass to $H^i(K_y)$.

Our standing assumption throughout this subsection is that $H^c_p(\ff)\neq 0$.
In \S~\ref{sec:non-neg}, we shall apply the results of this section when
$\ff$ is a finite field and $c = 1$ in which case this assumption is satisfied.

We shall consider two situations.
Suppose first that $\nu^0 \circ \Rsw_{y,n}(\chi) \neq 0$ in $\Omega^q_{A_1}$
(in particular, $\chi$ has type I at $E$) and let
$\beta := \nu^0 \circ \Rsw_{y,n}(\chi)$. 
For any $i \ge 0$ and $a \in \Omega^i_A$, we let $a'$ (resp. $\ov{a}$) denote the
image of $a$ in $\Omega^i_{A_1}$ (resp. $\Omega^i_\ff$).
We let $u_1, \ldots , u_{c-1} \in A^\times$ be such that
  $\{d(t'), \dlog(u'_1), \ldots , \dlog(u'_{c-1})\}$ is a basis of the
free $A_1$-module $\Omega^1_{A_1}$. For any $I = \{i_1, \ldots , i_r\}
\subset J_{c-1}$ of cardinality $r$,
we let $u_I = \dlog(u_{i_1}) \wedge \cdots \wedge \dlog(u_{i_r}) \in \Omega^{r}_A$. We
let $u'_I$ (resp. $\ov{u}_I$) be the image of $u_I$ in $\Omega^r_{A_1}$
(resp. $\Omega^r_\ff$).
  We can then write
  \begin{equation}\label{eqn:w-1-1}
    \beta = \left({\underset{|I| = q-1}\sum} \beta'_I d(t') \wedge u'_I\right)
    + \left({\underset{|J| = q}\sum} \beta'_J u'_J \right) \in \Omega^q_{A_1} \ \
    {\rm with} \ \ \beta_I, \beta_J \in A.
    \end{equation}

\begin{lem}\label{lem:w-1}
  Assume that there exists $I = \{i_1, \ldots , i_{q-1}\} \subset J_{c-1}$
  such that $\beta_I \in A^\times$. Then there exists
  $a, u_1,\dots, u_{c-q} \in A^\times$ such that 
  $$\partial'_y \circ \partial_y(\{ \chi, 1+\pi^n at^{-1}, u_1,\ldots, u_{c-q}\})\neq
  0 \text{ in } H^c(\ff).$$
  
  In particular, the image of $\{\chi, 1+\pi^n a t^{-1}\} \in H^{q+2}(K)$ in
  $H^{q+2}(K_{y})$ is non-zero
  under the restriction map $H^{q+2}(K) \to H^{q+2}(K_{y})$.
\end{lem}
\begin{proof}
We first note that there exists $I' = \{j_1, \ldots , j_{c-q}\}
  \subset J_{c-1}$ such that
  $\beta \wedge u'_{I'} = \beta'_I d(t') \wedge u'_I \wedge u'_{I'} \in \Omega^c_{A_1}$.
  Moreover, $\ov{\beta}_I \ov{u}_I \wedge \ov{u}_{I'}$ is a free generator
  of the $\ff$-vector space $\Omega^{c-1}_\ff$.
Using the surjective maps $\Omega^c_{A_1}(\log(t')) \stackrel{\res}{\surj}
  \Omega^{c-1}_\ff \stackrel{\delta^{c-1}_1(\ff)}{\surj} H^c_p(\ff)$ and our assumption
  $H^c_p(\ff)\neq 0$, we can thus find an element $a \in A^\times$ such that
  $\delta^{c-1}_1(\ff) \circ \res(a'\beta'_I {t'}^{-1} d(t')
  \wedge u'_I \wedge u'_{I'}) \neq 0$.
Equivalently,
  \begin{equation}\label{eqn:w-1-0}
    \delta^{c-1}_1(\ff)(\ov{a} \ov{\beta}_I  \ov{u}_I \wedge \ov{u}_{I'}) \neq 0.
  \end{equation}
  We let $\gamma_I = \{u_{i_1}, \ldots , u_{i_{q-1}}\} \in K^M_{q-1}(A)$,
  $\gamma_{I'} = \{u_{j_1}, \ldots , u_{j_{c-q}}\} \in K^M_{c-q}(A)$ and
  $\theta = \{\chi, 1+\pi^n a t^{-1}, \gamma_{I'}\} \in H^{c+2}(K_y)$.

By \cite[Thm.~5.1]{Kato-89}, we can write 
$\{\chi, 1+\pi^n a t^{-1}, \gamma_{I'}\} {=}  \{\lambda^{q+2}(\pi)(a' t'^{-1} \alpha,
    a' t'^{-1}\beta), \gamma_{I'}\}$ for some  $\alpha \in \Omega^{q+1}_{A_1}$. 
  Since $H^{c+2}(K_y) = \Fil_0 H^{c+2}(K_y)$ (note $[K_y: K^p_y] = p^{c+1}$), the map
  $\partial_y \colon H^{c+2}(K_y) \to H^{c+1}(k(y))$ is defined and we get
\begin{equation}\label{eqn:w-1-3}
    \begin{array}{lll}
      \partial_y(\theta) & = & \partial_y(\{\chi, 1+\pi^n a t^{-1}, \gamma_{I'}\}) \\
      & {=} &  \partial_y(\{\lambda^{q+2}(\pi)(a' t'^{-1} \alpha,
    a' t'^{-1}\beta), \gamma_{I'}\}) \\  
 & {=}^1 &  \partial_y \circ \lambda^{c+2}(\pi)(a' t'^{-1} \alpha \wedge u'_{I'},
   (-1)^{c-q} a' t'^{-1}\beta \wedge u'_{I'}\}) \\ 
    & = & (-1)^{c-q} \delta^{c}_1(k(y))(a' t'^{-1}\beta \wedge u'_{I'}) \\
    & = & (-1)^{c-q} \delta^{c}_1(k(y))(a' \beta'_{I}t'^{-1} d(t') \wedge u'_I \wedge u'_{I'}) \\
    & = & (-1)^{c-q} \delta^{c}_1(k(y))(a' \beta'_{I} \dlog(t') \wedge u'_I \wedge u'_{I'}) \\
    & {=}^2 & (-1)^{c-q} \{\delta^{0}_1(k(y))(a' \beta'_{I}), t', \gamma'_I, \gamma'_{I'}\}
     \end{array}
\end{equation}
where ${=}^1$ holds by \lemref{lem:Kato-basic-1} since $\gamma_{I'} \in K^M_{q-c}(A)$, and ${=}^2$ holds by
~\eqref{eqn:Milnor-1}.

If we consider the images of the two sides of ~\eqref{eqn:w-1-3} under the
restriction $\iota_y \colon H^{c+1}(k(y)) \to H^{c+1}(F)$, we get
\begin{equation}\label{eqn:w-1-4}
    \begin{array}{lll}
      \iota_y \circ \partial_y(\theta) 
   & {=} & (-1)^{c-q} \iota_y(\{\delta^{0}_1(k(y))(a' \beta'_{I}), t', \gamma'_I, \gamma'_{I'}\})
      \\
& {=} & (-1)^{c-q} \{\delta^{0}_1(F)(a' \beta'_{I}), t', \gamma'_I, \gamma'_{I'}\})
      \\
      & {=}^3 & (-1)^{c-q} \{l^1(A^h_1) \circ \delta^{0}_1(\ff)(\ov{a} \ov{\beta}_{I}),
      t', \gamma'_I, \gamma'_{I'}\}
      \\      
& {=} & (-1)^{c-q} (-1)^{c-1}\{\lambda^2(t')(0, \ov{a} \ov{\beta}_{I}), \gamma'_I, \gamma'_{I'}\} \\
    & {=}^4 & (-1)^{q+1} \lambda^{c+1}(t')(0, \ov{a} \ov{\beta}_{I} \ov{u}_I \wedge \ov{u}_{I'}), \\
    \end{array}
\end{equation}    
where ${=}^3$ holds by ~\eqref{eqn:delta-lift} and ${=}^4$ holds by 
~\eqref{eqn:Milnor-1}.

In particular,
$\partial'_y \circ  \partial_y(\theta)   =
(-1)^{q+1} \delta^{c-1}_1(\ff)(\ov{a} \ov{\beta}_{I} \ov{u}_I \wedge \ov{u}_{I'})$.
     But the last element is non-zero by ~\eqref{eqn:w-1-0}.
 This implies that $\theta \neq 0$, and hence $\{\chi, 1+\pi^n a t^{-1}\} \neq 0$.
 This completes the proof.      
\end{proof}

\begin{cor}\label{cor:w-3}
  Let $A$ be as above and $q = c$. Suppose that
  $\beta := \nu^0 \circ \Rsw_{y,n}(\chi) \in \Omega^c_{A_1}$ is such that
  the image of $\beta$ in $\Omega^c_{A_1}\otimes_A \ff$ is non-zero.
Then there exists $a \in A^\times$ such that the image of 
$\partial'_y \circ \partial_y(\{\chi, 1+\pi^n a t^{-1}\})\neq 0 \in H^{c}(\ff)$.
Moreover, the map
\[
\tau_\chi \colon A^\times \to H^c_p(\ff); \ \mbox{given by} \ \tau_\chi(a) =
\partial'_y \circ \partial_y(\{\chi, 1+\pi^n a t^{-1}\}),
\]
is surjective.
\end{cor}
\begin{proof}
When $q = c$ in \lemref{lem:w-1}, then each $\beta'_J$
in ~\eqref{eqn:w-1-1} can be chosen to be zero which implies
$\beta = \beta'_I d(t') \wedge u'_I$
with $I = J_{c-1}$ and $\beta_I \in A^\times$. This proves the first part. For the
surjectivity of $\tau_\chi$, note that as $\ov{\beta}_I \ov{u}_I$ is a free
generator of the $\ff$-vector space $\Omega^{c-1}_\ff$, the elements
$\ov{a}\ov{\beta}_I\ov{u}_I$, as $a$
ranges over $A^\times$, exhaust all of $\Omega^{q-1}_\ff$.
The surjectivity of $\delta^{c-1}_1$ therefore implies that the elements
$\partial'_y \circ \partial_y(\{\chi, 1+\pi^n a t^{-1}\}) =
\delta^{c-1}_1(\ff)(\ov{a}\ov{\beta}_I\ov{u}_I)$, as $a$
ranges over $A^\times$, exhaust all of $H^c_p(\ff)$, thereby yielding the desired
result.
\end{proof}

\vskip.2cm

Suppose next that $\Rsw_{y,n}(\chi) \in \Omega^{q+1}_{A_1}$ (in particular, $\chi$ has
type II at $E$) and let $\alpha = \Rsw_{y,n}(\chi)$.
We can then write
  \begin{equation}\label{eqn:w-2-1}
    \alpha = \left({\underset{|I| = q}\sum} \alpha'_I d(t') \wedge u'_I\right)
    + \left({\underset{|J| = q+1}\sum} \alpha'_J u'_J \right) \in
    \Omega^{q+1}_{A_1} \ \
    {\rm with} \ \ \alpha_I, \alpha_J \in A.
    \end{equation}

\begin{lem}\label{lem:w-2}
    Assume that $q \in J^0_{c-1}$ and 
    there exists $I = \{i_1, \ldots , i_{q}\} \subset J_{c-1}$
    such that $\alpha_I \in A^\times$. Then there exist
    $a, u_1,\ldots, u_{c-q-1} \in A^\times$ such that
    $$\partial'_y \circ
    \partial_y(\{ \chi, 1+\pi^n at^{-1}, \pi, u_1,\ldots, u_{c-q-1}\})\neq 0
    \text{ in } H^c(\ff).$$ 
    In particular, the image of $\{\chi, 1+\pi^n a t^{-1}, \pi\} \in H^{q+3}(K)$ in
    $H^{q+3}(K_{y})$ is
  non-zero under the restriction map $H^{q+3}(K) \to H^{q+3}(K_{y})$. 
\end{lem}
  \begin{proof}
The proof is very similar to that of \lemref{lem:w-1} and we shall mostly
    explain the modifications.
    As before, there exists $I' = \{j_1, \ldots , j_{c-q-1}\} \subset J_{c-1}$
    such that
    $\alpha \wedge u'_{I'} = \alpha'_I d(t') \wedge u'_I \wedge u'_{I'} \in
    \Omega^c_{A_1}$.
    Moreover,  we can find an element $a \in A^\times$ such that
  \begin{equation}\label{eqn:w-2-0}
    \delta^{c-1}_1(\ff)(\ov{a} \ \ov{\alpha}_I  \ov{u}_I \wedge \ov{u}_{I'}) \neq 0
    \ \mbox{in} \ H^c(\ff).
  \end{equation}

We let $\tau_I = \{u_{i_1}, \ldots , u_{i_{q}}\} \in K^M_{q}(A)$,
  $\tau_{I'} = \{u_{j_1}, \ldots , u_{j_{c-q-1}}\} \in K^M_{c-q-1}(A)$ and
  $\psi = \{\chi, 1+\pi^n a t^{-1}, \pi, \tau_{I'}\} \in H^{c+2}(K_y)$.
We then get 
\begin{equation}\label{eqn:w-2-3}
    \begin{array}{lll}
      \partial_y(\psi) & = & \partial_y(\{\chi, 1+\pi^n a t^{-1},\pi,  \tau_{I'}\}) \\
      & {=} &  \partial_y(\{\lambda^{q+2}(\pi)(a't'^{-1} \alpha, 0), \pi,\tau_{I'} \})
      \\
& {=} &  \partial_y(\{\lambda^{q+3}(\pi)(0,at^{-1} \alpha), \tau_{I'}\}) \\
& {=} & (-1)^{c-q-1} \partial_y(\lambda^{c+2}(\pi)(0, at^{-1} \alpha \wedge u_{I'})) \\
      & {=} &
      (-1)^{c-q-1} \{\delta^{0}_1(k(y))(a' \alpha'_{I}), t', \tau'_I, \tau'_{I'}\}.
    \end{array}
    \end{equation}

If we consider the images of the two sides of ~\eqref{eqn:w-2-3} under
$\iota_y \colon H^{c+1}(k(y)) \to H^{c+1}(F)$, we get
\begin{equation}\label{eqn:w-2-4}
    \begin{array}{lll}
      \iota_y \circ \partial_y(\psi) 
   & {=} & (-1)^{c-q-1}\iota_y(\{\delta^{0}_1(k(y))(a' \alpha'_{I}), t', \tau'_I, \tau'_{I'}\})
      \\
& {=} & (-1)^{c-q-1}\{\delta^{0}_1(F)(a' \alpha'_{I}), t', \tau'_I, \tau'_{I'}\}
      \\
      & {=} & (-1)^{c-q-1} \{l^1(A^h_1) \circ \delta^{0}_1(\ff)(\ov{a} \ \ov{\alpha}_{I}),
      t', \tau'_I, \tau'_{I'}\}
      \\      
& {=} &(-1)^{c-q-1} \{\lambda^2(t')(0, \ov{a} \ \ov{\alpha}_{I}), \tau'_I, \tau'_{I'}\} \\
    & {=} & (-1)^{q}\lambda^{c+1}(t')(0, \ov{a} \ \ov{\alpha}_{I} \ov{u}_I \wedge \ov{u}_{I'}). \\
    \end{array}
\end{equation}       

In particular,
$\partial'_y \circ  \partial_y(\psi) =
(-1)^q \delta^{c-1}_1(\ff)(\ov{a} \ \ov{\alpha}_{I} \ov{u}_I \wedge \ov{u}_{I'})$. 
But the last element is non-zero by ~\eqref{eqn:w-2-0}.
This implies that $\psi \neq 0$ (equivalently,
$\{\chi, 1+\pi^n a t^{-1}, \pi, \tau_{I'}\} \neq 0$), and hence
$\{\chi, 1+\pi^n a t^{-1}, \pi\} \neq 0$.
 This completes the proof.
\end{proof}

\subsection{Application II: Swan conductor when $n = p =2$}\label{sec::Non-V-0}
In this subsection, we apply the Kato complex to estimate the drop in the
value of the Swan conductor of $\chi$ when we restrict it to certain divisors on
$X$. We first prove some general results.

\begin{lem}\label{lem:p-2-1}
  Let $K$ be an $F$-finite hdvf of characteristic $p =2$
  with the ring of integers $\sO_K$ and the uniformizer $\pi$. Let
  $a \in \Fil_2K=\pi^{-2} \sO_K$ and $u,v \in \sO^\times_K$. Then
  $\delta_1^2\left(a \dlog(1+\pi^2uv)\wedge
  \dlog(\pi) - a \dlog(1+\pi u)\wedge \dlog(1+\pi v)\right)$ is an
  unramified element (cf. Definition~\ref{defn:tame-ram}) of $H^3_p(K)$.
\end{lem}
\begin{proof}
 We let $\alpha = a \dlog(1+\pi^2uv)\wedge \dlog(\pi)$, 
  $\beta = a \dlog(1+\pi u)\wedge \dlog(1+\pi v)$ and
  $\gamma = \frac{a\pi^2 d(uv)\wedge \dlog(\pi)}{(1+\pi u)(1+\pi v)}$.
  One checks that
  $\alpha = \frac{a\pi^2 d(uv) \wedge \dlog(\pi)}{1 + \pi^2uv}$.
On the other hand,
\begin{equation}\label{eqn:p-2-1-0}
  \begin{array}{lll}
          \beta & = &
          a \dlog(1+\pi u)\wedge \dlog(1+\pi v) \\
          &=&\frac{a}{(1+\pi u)(1+\pi v)}
          ((\pi du+u d\pi )\wedge (\pi dv+ vd\pi ))\\
          & {=}^1 &
          \frac{1}{(1+\pi u)(1+\pi v)} (a \pi^2 dudv + a \pi^2(vdu +u dv)
          \wedge \dlog \pi)\\
          &= & \frac{a\pi^2 uv  \dlog(u) \wedge \dlog(v)}{(1+\pi u)(1+\pi v)}+
          \frac{a\pi^2 d(uv)\wedge \dlog(\pi)}{(1+\pi u)(1+\pi v)} \\
          & = & \frac{a\pi^2 uv  \dlog(u) \wedge \dlog(v)}{(1+\pi u)(1+\pi v)} +
          \gamma,
           \end{array}
\end{equation}
where ${=}^1$ holds because $p = 2$. We let $\gamma' =
\frac{a\pi^2 uv  \dlog(u) \wedge \dlog(v)}{(1+\pi u)(1+\pi v)}$.

We let $\theta_1 = a\pi^2 d(uv) \wedge \dlog(\pi)$ and
        $\theta_2 = a\pi^2 uv  \dlog(u) \wedge \dlog(v)$.
        Since $(1+ \pi^2 uv)^{-1} = 1 + \pi^2 x$ and $((1+ \pi u)(1+\pi v))^{-1} =
        1 + \pi y$ for some $x, y \in \sO_K$, we get
        $\alpha = \theta_1 + \pi^2x \theta_1$. Note that $a' := a\pi^2 \in \sO_K$ by our assumption. In particular,
        $\delta^2_1(\alpha) = \delta^2_1(\theta_1) + \delta^2_1(\pi^2x \theta_1) =
        \delta^2_1(\theta_1)$ 
        and $\delta^2_1(\gamma) = \delta^2_1(\theta_1) +
        \delta^2_1(\pi y \theta_1) = \delta^2_1(\theta_1)$ by 
        \lemref{lem:delta-fil-0} since
        $\pi^2x \theta_1, \pi y \theta_1 \in \Fil_{-1}\Omega^2_K$.
        We similarly get $\delta^2_1(\gamma') =
        \delta^2_1(\theta_2) + \delta^2_1(\pi y \theta_2) =
        \delta^2_1(\theta_2)$.
        It remains therefore to show that $\delta^2_1(\theta_2)$ is unramified.

        To that end, note again that $a' := a\pi^2 \in \sO_K$.
        If $a' \in \sO^\times_K$, then $\theta_2 \in \Omega^2_{\sO_K}$ and hence
        $\delta^2_1(\theta_2)$ is unramified. If $a' \in (\pi)$, then
        $\theta_2 \in \pi\Omega^2_{\sO_K} \subset \Fil_{-1}\Omega^2_K$.
        In particular, $\delta^2_1(\theta_2) = 0$. This completes the proof.
\end{proof}

\begin{lem}\label{lem:p-2-2}
  Let $K$ be as in \lemref{lem:p-2-1}. Let $a \in \Fil_1 K =\pi^{-1} \sO_K$ and
  $u, v \in \sO^\times_K$. Then
  $\delta_1^2\left(a \dlog(1+\pi u)\wedge \dlog(1+\pi v)\right) = 0
  = \delta_1^3\left(a \dlog(\pi) \wedge \dlog(1+\pi u)\wedge \dlog(1+\pi v)\right)$.
  \end{lem}
\begin{proof}
Letting $b = (1 + \pi u)(1 + \pi v) \in \sO_K$, we have
\begin{equation}\label{eqn:p-2-2-0}
  \begin{array}{lll}
    a \dlog(1+\pi u)\wedge \dlog(1+\pi v) & = &
    ab^{-1} (ud\pi + \pi du)\wedge (vd\pi + \pi dv) \\
    & = & ab^{-1} \pi^2 uv (\dlog(u) \wedge \dlog(v) - \dlog(\pi) \wedge \dlog(u) \\
      & & + \ \dlog(\pi) \wedge \dlog(v)).
  \end{array}
  \end{equation}
\begin{equation}\label{eqn:p-2-2-1}
\begin{array}{lll}
    a \dlog(\pi) \wedge \dlog(1+\pi u)\wedge \dlog(1+\pi v) & = &
    a \dlog(1+\pi u)\wedge \dlog(1+\pi v) \wedge  \dlog(\pi) \\
    & = & ab^{-1} (ud\pi + \pi du)\wedge (vd\pi + \pi dv)  \wedge \dlog(\pi) \\
    & = & ab^{-1} \pi^2 uv \dlog(u) \wedge \dlog(v) \wedge \dlog(\pi).
  \end{array}
\end{equation}

Since $a \in \Fil_1 K$, the right side side of ~\eqref{eqn:p-2-2-0}
(resp. ~\eqref{eqn:p-2-2-1}) lies in $\Fil_{-1}\Omega^2_K$
(resp. $\Fil_{-1}\Omega^3_K$).
In particular, it dies under $\delta^2_1$ (resp. $\delta^3_1$)  by
\lemref{lem:delta-fil-0}. This concludes the proof.
\end{proof}

\begin{lem}\label{lem:delta-fil-0}
  For $K$ as in \lemref{lem:p-2-1}, we have $\Fil_{-1}\Omega^q_K \subset
  (1-C)(Z_1\Omega^q_{\sO_K})$ for every $q \ge 0$.
\end{lem}
\begin{proof}
  Combine Lemmas~\ref{lem:Coh-iso} and ~\ref{lem:Fid-DRW-basic}(8), and
  \cite[Lem.~9.8(2)]{KM-1}.
\end{proof}

\begin{lem}\label{lem:w-4}
  Suppose we are in the situation of \lemref{lem:w-2}. Assume additionally
  that  $n = p =2$ and $\chi \in H^{q+1}_2(A[\pi^{-1}])$. Then we have
    \begin{equation}\label{eqn:w-3-0}
      \{\chi, 1+\pi^2 a t^{-1}, \pi, u_1,\ldots,u_{c-q-1}\} =
      \{\chi, 1+\pi a , 1+\pi t^{-1}, u_1,\ldots,u_{c-q-1}\}
\end{equation}
    when we pass from $H^{c+2}_2(K)$ to $H^{c+2}_2({K_y})$, where
    $a, u_1,\ldots,u_{c-q-1}$ are as in
    \lemref{lem:w-2}.
\end{lem}
\begin{proof}
It suffices to show that
  $\{\chi_y, 1+\pi^2 a t^{-1}, \pi\} - \{\chi_y, 1+\pi a , 1+\pi t^{-1}\}$ is an
  unramified element of $H^{q+3}_2(K_y)$. Indeed, if $\alpha \in H^{q+3}_2(K_y)$
  is unramified, then so is $\{\alpha, \tau\} \in H^{c+2}_2(K_y)$
  for any $\tau \in K^M_{c-q-1}(A)$  (cf. \lemref{lem:Kato-basic-1}). Hence,
  $\{\alpha, \tau\}$ must be zero since $H^{c+2}_2(k(y))$
  is a quotient of $\Omega^{c+1}_{k(y)}$ which is zero (cf. \cite[Lem.~7.2]{Kato-89})

We now note that there are canonical maps
  \begin{equation}\label{eqn:w-3-1}
    \Fil_2 K_y \otimes K^M_q(K_y) \to \Fil_2 \Omega^{q}_{K_y}
    \xrightarrow{\delta^{q}_1} \Fil_2 H^{q+1}_2(K_y),
  \end{equation}
  where the left arrow takes $u \otimes v$ to $u\dlog(v)$. The composite map
  is surjective by \thmref{thm:H^1-fil} and \cite[Cor.~7.3]{KM-1}. 
  Recalling the generators of $K^M_q(K_y)$ (cf. proof of \lemref{lem:Kato-basic-1})
  and noting that $\chi_y \in \Fil_2 H^{q+1}_2(K_y)$, we can therefore write 
  \begin{equation}\label{eqn:w-3-2}
    \chi_y = \delta^{q}_1 \left(\sum_i a_i \dlog(\alpha_i) + \sum_j b_j \dlog(\beta_j)
    \wedge \dlog(\pi)\right),
  \end{equation}
  where $a_i, b_j \in \Fil_2 K_y, \ \alpha_i \in K^M_{q}(A_y)$ and $\beta_j \in
  K^M_{q-1}(A_y)$ are Milnor symbols.
  We let $\gamma_1 = \sum_i a_i \dlog(\alpha_i)$ and
  $\gamma_2 = \sum_j b_j \dlog(\beta_j) \wedge \dlog(\pi)$.

  \vskip.2cm
  
{\bf{Claim:}} $\delta^q_1 \left(\gamma_2\right)
\in \Fil_1 H^{q+1}_2(K_y)$. 

To prove the claim, we let
$\omega = \sum_i a_i \dlog(\alpha_i) + \sum_j b_j \dlog(\beta_j) \wedge \dlog(\pi)$.
Since $a_i, b_j \in  \Fil_2 K_y$, we can write $a_i = \pi^{-2}a'_i$ and
$b_j = \pi^{-2}b'_j$ for some $a'_i, b'_j \in A_y$. In particular,
$d(a_i) = \pi^{-2}da'_i$ and $d(b_j) = \pi^{-2}db'_j$ as $p =2$.
We thus get $d\omega = \pi^{-2}\sum_i da'_i \dlog(\alpha_i) +
\pi^{-2}\sum_j db'_j \dlog(\beta_j) \wedge \dlog(\pi) =
\pi^{-2}\omega'$, say.

We now note using ~\eqref{eqn:RSW-gen-1-0-0} and our assumption on
$\chi$ (that it has type II at $E$) that $(-1)^{q+1}\nu^0 \circ d(\omega) = 0$.
Equivalently, $d\left(\sum_j \ov{b'_j} \dlog(\ov{\beta_j} \right) =
\sum_j d\ov{b_j} \dlog(\ov{\beta_j}) = \res(\omega') = 0$ in
$\Omega^{q-1}_{k(y)}$, where $\ov{b'_j}$ (resp. $\ov{\beta_j}$) is the image of $b'_j$
(resp. $\beta_j$) in $k(y)$ (resp. $K^M_{q-1}(k(y))$).
In other words, $\sum_j \ov{b'}_j \dlog(\ov{\beta_j}) \in Z_1\Omega^{q-1}_{k(y)}$.
Using the surjection $Z_1\Omega^{q-1}_{A_y} \surj Z_1\Omega^{q-1}_{k(y)}$ (recall that
a closed form is the image of an element under the Frobenius of the de Rham-Witt
complex), we can find a lift $\alpha$ of $\sum_j \ov{b'_j} \dlog(\ov{\beta_j})$
in $Z_1\Omega^{q-1}_{A_y}$. In particular, $\alpha - \sum_j b'_j \dlog(\beta_j)
\in \Ker(Z_1\Omega^{q-1}_{A_y} \surj Z_1\Omega^{q-1}_{k(y)}) \subset
\Fil_{-1} \Omega^{q-1}_{K_y}$ by \lemref{lem:Fid-DRW-basic} and
\cite[Lem.~9.8(2)]{KM-1}. We let this element be $\alpha'$.

As $\alpha' \in \Fil_{-1} \Omega^{q-1}_{K_y}$, we see that
$\pi^{-2} \alpha' \wedge \dlog(\pi) \in \Fil_{1} \Omega^{q}_{K_y}$
(cf. \cite[Defn.~3.8, Lem.~3.10]{KM-1}).
In particular,
$\delta^q_1(\pi^{-2} \alpha' \wedge \dlog(\pi)) \in \Fil_1 H^{q+1}(K_y)$.
To prove the claim, it remains therefore to show that
$\delta^q_1(\pi^{-2}\alpha \wedge \dlog(\pi)) \in \Fil_1 H^{q+1}(K_y)$.
To that end, we note that $\pi^{-2}\alpha \wedge \dlog(\pi) \in Z_1\Omega^q_{K_y}$
since $\alpha \in Z_1\Omega^{q-1}_{K_y}$ and $p =2$. This implies that
$\delta^q_1(\pi^{-2}\alpha \wedge \dlog(\pi)) =
\delta^q_1 \circ C(\pi^{-2}\alpha \wedge \dlog(\pi))
= \delta^q_1(\pi^{-1}C(\alpha) \wedge \dlog(\pi))$, where the second
equality follows from \cite[Lem.~1.1]{Milne-ENS}.
We are done by \thmref{thm:H^1-fil} because
$\pi^{-1}C(\alpha) \wedge \dlog(\pi) \in \Fil_{1} \Omega^{q}_{K_y}$.

Using the claim and \cite[Prop.~6.3(iii)]{Kato-89}, we get
$\{\delta^q_1(\gamma_2), 1+\pi^2 a t^{-1}\} = 0$ (note that the images of
$a$ and $t$ in $A_y$ are units). In particular,
$\{\delta^q_1(\gamma_2), 1+\pi^2 a t^{-1}, \pi\} {=} 
\{\{\delta^q_1(\gamma_2), 1+\pi^2 a t^{-1}\}, \pi\} = 0$,
where the first equality follows from the last equality at the bottom of
\cite[\S~(1.3)]{Kato-89} (cf. \eqref{eqn:Milnor-1}).

To compute $\{\delta^q_1(\gamma_2), 1+\pi a, 1+ \pi t^{-1}\}$, we first use our
claim and \thmref{thm:H^1-fil} to write 
$\delta^q_1(\gamma_2) = \delta^q_1\left(\sum_i x_i \dlog(y_i) +
\sum_j x'_j \dlog y'_j \wedge \dlog(\pi)\right)$,
where $x_i, x'_j \in \Fil_1 K_y, \ y_i \in K^M_{q}(A_y)$ and $y'_j \in K^M_{q-1}(A_y)$
(cf. \lemref{lem:Fid-DRW-basic}(8)).
By \lemref{lem:p-2-2}, we have
$\{\delta^q_1(x_i \dlog(y_i)), 1+\pi a, 1+ \pi t^{-1}\} =
\{\delta^2_1\left(x_i \dlog( 1+\pi a)\wedge \dlog(1+ \pi t^{-1})\right), y_i\}
= 0$ for every $i$. Similarly, we have 
$\{\delta^q_1\left(x'_j \dlog y'_j \wedge \dlog(\pi)\right),
1+\pi a, 1+ \pi t^{-1}\} = 0$ for every $j$.
It follows that $\{\delta^q_1(\gamma_2), 1+\pi a, 1+ \pi t^{-1}\} =0$.
It remains therefore to show that
$\{\delta^q_1(\gamma_1), 1+\pi^2 a t^{-1}, \pi\} -
\{\delta^{q}_1(\gamma_1), 1+\pi a, 1+ \pi t^{-1}\}$
is an unramified element of $H^3(K_y)$.

To that end, we compute
\[
\begin{array}{lll}
\{\delta^q_1(\gamma_1), 1+\pi^2 a t^{-1}, \pi\} & = &
\{\delta^{q}_1 \left(\sum_i a_i \dlog(\alpha_i)\right),
1+\pi^2 a t^{-1}, \pi\} \\
& {=}^1 &
\sum_i \{\delta^2_1\left(a_i \dlog( 1+\pi^2 a t^{-1})\wedge \dlog(\pi)\right),
\alpha_i\},
\end{array}
\]
where ${=}^1$ follows again from \cite[\S~(1.3)]{Kato-89}
(cf. \eqref{eqn:Milnor-1}).
A similar computation yields
$\{\delta^q_1(\gamma_1),  1+\pi a, 1+ \pi t^{-1}\}
= \sum_i \{\delta^2_1\left(a_i\dlog(1+\pi a) \wedge \dlog(1+ \pi t^{-1})\right),
\alpha_i\}$. 
Since $a_i \in \Fil_2 K$, \lemref{lem:p-2-1} meanwhile says that
$$\delta^2_1\left(a_i \dlog( 1+\pi^2 a t^{-1})\wedge \dlog(\pi)\right) -
\delta^2_1\left(a_i \dlog(1+\pi a)\wedge \dlog(1+ \pi t^{-1})\right)$$
is unramified in $H^3(K_y)$ for every $i$.
This implies by \lemref{lem:Kato-basic-1} that for every $i$, the element
$\{\delta^2_1\left(a_i \dlog( 1+\pi^2 a t^{-1})\wedge \dlog(\pi)\right), \alpha_i\}
- \{\delta^2_1\left(a_i\dlog(1+\pi a) \wedge \dlog(1+ \pi t^{-1})\right),
\alpha_i\}$ is unramified in $H^{q+3}(K_y)$.
Taking the sum, we conclude that
$\{\delta^q_1(\gamma_1), 1+\pi^2 a t^{-1}, \pi\} -
\{\delta^{q}_1(\gamma_1), 1+\pi a, 1+ \pi t^{-1}\}$ is unramified in
$H^{q+3}(K_y)$. This completes the proof.
\end{proof}

The final result of this section is the following lemma about the Swan conductor.

\begin{lem}\label{lem:spl-loc-1}
  Assume that $q \in J_{c-1}$ and $\chi \in H^{q+1}_p(A[\pi^{-1}])$.
  Assume also the following.
  \begin{enumerate}
    \item
      $n = p =2$.
    \item $\chi$ has type II at $E$. 
    \item
      $\Rsw_{y,n}(\chi)$ does not die under the quotient map
      $\Omega^{q+1}_{A_1} \surj \Omega^{q+1}_{A_1} \otimes_{A_1} \ff$.
    \item
      $H^c_p(\ff)\neq 0$.
      \end{enumerate}
 Then there exists $a \in A^\times$ such that  
    $\{\chi_y, 1+\pi a , 1+\pi t^{-1}\} \neq 0$ in $H^{q+3}_p(K_y)$. Moreover,
    there exists $X' \in  \{\Spec(A/(t)),\Spec(A/(t+\pi))\}$ such that
    $1 \le \Sw_{x}(\chi|_{U'})\le 2$, where $U'=X'\setminus \{x\}$.
\end{lem}
\begin{proof}
We let $B = {A}/{(t)}$ and $B' = {A}/{(t+\pi)}$. 
  We let $\alpha = \Rsw_{y,2}(\chi) \in \Omega^{q+1}_{A_1}$ be as in ~\eqref{eqn:w-2-1}.
  We let
  \[
  \wt{\alpha} = \left({\underset{|I| = q}\sum} \alpha_I d(t) \wedge u_I\right)
    + \left({\underset{|J| = q+1}\sum} \alpha_J u_J \right) \in \Omega^{q+1}_{A} \ \
    {\rm with} \ \ \alpha_I, \alpha_J \in A
    \]
    be a lift of $\alpha$ in $ \Omega^{q+1}_{A}$.

We first consider the case when $\alpha_{J} \in A^\times$ for some $J$.
This implies that $q \le c-2$. Since
$\{\dlog(\ov{u}_1), \ldots , \dlog(\ov{u}_{c-1})\}$
is an $\ff$-basis of $\Omega^1_\ff$, this also implies that $\wt{\alpha}$
does not die in $\Omega^{q+1}_\ff$ (cf. \cite[Lem.~7.2]{Kato-89}).
We let $\chi_B$ be the image of $\chi$ in $H^{q+1}_p(B[\pi^{-1}])$.
We consider the diagram
\begin{equation}\label{eqn:spl-loc-1-0}
  \xymatrix@C2pc{
    \frac{\Fil_2 H^{q+1}_p(A[\pi^{-1}])}{\Fil_1 H^{q+1}_p(A[\pi^{-1}])}
    \ar[r]^-{\Rsw_{y,2}} \ar[d] &
    \frac{\Fil_2 \Omega^{q+1}_{A[\pi^{-1}]}}{\Fil_1 \Omega^{q+1}_{A[\pi^{-1}]}}
    \ar[d] & \Omega^{q+1}_{A_1} \ar[d] \ar[l]_-{\nu'} \\
\frac{\Fil_2 H^{q+1}_p(B[\pi^{-1}])}{\Fil_1 H^{q+1}_p(B[\pi^{-1}])}
    \ar[r]^-{\Rsw_{y,2}}  &
    \frac{\Fil_2 \Omega^{q+1}_{B[\pi^{-1}]}}{\Fil_1 \Omega^{q+1}_{B[\pi^{-1}]}}
     & \Omega^{q+1}_{\ff} \ar[l]^-{\nu'}}
\end{equation}
whose vertical arrows are the restriction maps and $\nu'$ is from
\lemref{lem:RSW-spl}. The left square of this diagram is commutative by
\corref{cor:RSW-gen-1} and the right square is commutative because $\nu'$ is
functorial. Since $\Rsw_{y,2}(\chi) \in  \Omega^{q+1}_{A_1}$ whose image in
$\Omega^{q+1}_{k(x)}$ is non-zero, it follows that $\Rsw_{y,2}(\chi_B) \neq 0$.
This shows that $\Sw_{\ff}(\chi_B) = 2$.

Under the condition (3) of the lemma, the only other case to consider is when
$\alpha_I \in  A^\times$ for some $I$.
In this case, \lemref{lem:w-2} says that there exist $a, u_1,\ldots,u_{c-q-1} \in
A^\times$ such that for
  $\tau = \{u_{1}, \ldots , u_{c-q-1}\} \in K^M_{c-q-1}(A)$
  and $\psi = \{\chi, 1+\pi^2 a t^{-1}, \pi, \tau\} \in H^{c+2}_p(K)$,
  one has that  $\partial'_y \circ \partial_y(\psi) \neq 0$ in $H^{c}_p(\ff)$.
  Letting $\psi' = \{\chi, 1+\pi a, 1+ \pi t^{-1},  \tau_{I'}\}  \in H^{c+2}_p(K)$,
  it follows from \lemref{lem:w-4} that $\psi|_{k(y)} = \psi'|_{k(y)}$ in
  $H^{c+2}_p(K_y)$. In particular, $\partial'_y \circ \partial_y(\psi') \neq 0$.

We now consider the Kato complex
  \begin{equation}\label{eqn:spl-loc-1-1}
  H^{c+2}(K) \xrightarrow{\partial=\bigoplus \partial_z}
  \bigoplus\limits_{z \in X^{(1)}} H^{c+1}(k(z))
  \xrightarrow{\partial'=\bigoplus \partial'_z} H^c(\ff).
\end{equation}
We note that $\partial_z(\psi') = 0$ if $\psi'$ is unramified
  at $z \in X^{(1)}$, and  $\psi'$ is possibly ramified only
  along the prime ideals $(\pi), \fp = (t)$ and $\fp' = (t + \pi)$
(cf. \lemref{lem:Coh-iso}).
This implies that $\partial(\psi') = \partial_y(\psi') + \partial_\fp(\psi') +
  \partial_{\fp'}(\psi)$. If $\partial_\fp(\psi') = 0 = \partial_{\fp'}(\psi')$,
  we get $0 = \partial' \circ \partial (\psi') =
  \partial'_y \circ \partial_y(\psi')$.
  This leads to a contradiction as $\partial'_y \circ \partial_y(\psi') \neq 0$.
  It follows that one of $\partial_\fp(\psi')$ and $\partial_{\fp'}(\psi)$ is not zero.
But this implies that either $\{\chi_\fp, 1+\pi a\} \neq 0$ in
  $H^{q+2}_p(K_\fp)$ or $\{\chi_{\fp'},  1+\pi a\} \neq 0$ in  $H^{q+2}_p(K_{\fp'})$.
We conclude from this that either $\Sw_\ff(\chi_\fp) \ge 1$ or
$\Sw_\ff(\chi_{\fp'}) \ge 1$ (cf. \cite[Thm~6.3]{Kato-89}). In other words,
$\Sw_\ff(\chi|_{B})\ge 1$ or $\Sw_\ff(\chi|_{B'})\ge 1$. Since
\corref{cor:Fil-functorial} implies that
$\Sw_\ff(\chi|_{B}), \Sw_\ff(\chi|_{B'})\le 2$, the proof is complete.
\end{proof}

\section{Specialization of Swan conductor I}\label{sec:Esp}
In this section and the next, we apply \thmref{thm:RSW-gen} to establish results
bounding the drop in the Swan conductor of a cohomology class on a scheme upon
restriction to subschemes satisfying certain conditions. The results of this
section hold for all cohomology groups $H^q(U)$ with $q \ge 1$ (where $U$ is as in
\S~\ref{sec:Fil-DRW}) and are of independent interest in the study of the
ramification theory. In this paper, however, we use these results only in the case
$q = 2$.

\subsection{The setting and basic results}\label{sec:Set-up}
Our setting and notation for Sections~\ref{sec:Esp} and~\ref{sec:Esp-0}
are as follows.
We let $R$ be an $F$-finite hdvr containing $\F_p$. Let $k$ and $\F$ denote the
fraction field and residue field of $R$, respectively, and let $\fm=(\pi)$ be the
maximal ideal of $R$. We write $p_\F$ for the $p$-rank of $\F$, recalling that the
$p$-rank of a field $L$ of characteristic $p$ is the integer $n$ such that
$[L:L^p]=p^n$. Finally, set
$S = \Spec(R), \ \eta = \Spec(k)$ and $s = \Spec(\F)$.

If $\iota \colon Z \inj X$ is a closed immersion of schemes and $\sF$
is a sheaf on $Z_\et$ (or $Z_\zar$), we shall abuse notations and
usually write $\iota_*(\sF)$ simply as $\sF$ when it is clear that the sheaves are
being considered on $X$ (and not on $Z$).
For any line bundle $\sL$ on $X$, we let $|\sL|$ denote
  the complete linear system of effective divisors on $X$ associated to $\sL$. 
 We begin by collecting some basic results.

\begin{lem}\label{lem:Rank-module-diff}
  Let $E$ be a finite type, connected, regular, faithfully flat $A$-scheme of
  dimension $q$. Then $\Omega^1_E$ is a locally free sheaf on $E$ of rank $q+p_\F$
  (resp. $q+p_\F+1$) if $A \in \{R, \F\}$ (resp. $A = K$).
  \end{lem}
\begin{proof}
  Using \cite[Prop.~2.2]{KM-1}, the lemma is easily reduced to showing the folklore
  statement that if $L/L'$ is a finite extension of $F$-finite fields containing
  $\F_p$, then $L'$ and $L$ have the same $p$-rank. Equivalently,
  $\dim_{L'} (\Omega^1_{L'}) = \dim_L(\Omega^1_L)$. But the latter equality follows
  from \cite[Thm.~26.10]{Matsumura} using the first fundamental exact sequence
  for K{\"a}hler differentials (cf. Thm.~ 25.1(1) of op. cit.).
\end{proof}

\begin{lem}\label{lem:Kahler-seq}
  Let $(X,E)$ be an snc-pair of $\F_p$-schemes and let $u \colon X' \inj X$ be the
  inclusion of a regular divisor such that $E' := X' \bigcap E$ is an snc divisor
  on $X'$ with components $E'_1, \ldots , E'_r$, where $E'_i = X' \bigcap E_i$.
  Then we have the following exact sequences for $q \ge 0$.
  \begin{enumerate}
    \item
  $0 \to \Omega^q_X(\log E) \to \Omega^q_X(\log (E + X')) \xrightarrow{\res}
    \Omega^{q-1}_{X'}(\log E') \to 0$;
    \item
 $0 \to \Omega^q_X(\log (E + X'))(-X') \to \Omega^q_X(\log E) \xrightarrow{\theta}
      \Omega^q_{X'}(\log E') \to 0$,
  \end{enumerate}
  where $\theta$ is the canonical pull-back map induced by $u$.
\end{lem}
\begin{proof}
See \cite[Lem.~9.8]{KM-1}.
\end{proof}

\begin{lem}\label{lem:Coh-vanishing}
  Let $k$ be an $F$-finite field  and let $E$ be a projective
  $k$-scheme which is connected and regular of dimension $d_E \ge 2$. Let $F, H, F'
  \in \Div(E)$ such that $F$ is an snc divisor and $H$ is ample. Let $q \in \N_0$.
  Then there is $N \gg 0$ such that for every $m \ge N$ and
  $Z \in |\sO_E(mH)|(k)$, one has 
  $H^i_\zar(Z, \Omega^q_Z(\log F_Z)((F' - nZ)|_Z) = 0$ for all $i \le d_E -q-2$
  and $n \ge 1$ provided $Z$ is regular and $F_Z := Z \times_E F$ is an snc divisor
  on $Z$.
\end{lem}
\begin{proof}
By \cite[Lem.~5.1]{GK-Jussieu}, we choose $N \gg 0$ such that
  \begin{equation}\label{eqn:Coh-vanishing-0}
    H^i_\zar(E,
    \Omega^q_E(\log F)(F'- mH)) = 0 \ \ \mbox{for  all} \ q \ge 0, i \le d_E -1 \ \
    \mbox{and} \ m \ge N.
    \end{equation}
  We now fix $m \ge N$ and let $Z \in |\sO_E(mH)|(k)$ such that $Z$ is regular
  (note that $Z$ is connected because $E$ is and $d_E \ge 2$) and
  $F_Z$ is an snc divisor on $Z$. We shall prove by induction on
  $q \ge 0$ that
  $H^i_\zar(Z, \Omega^q_Z(\log F_Z)((F' - nZ)|_Z)) = 0$ for all $n \ge 1$
  and $i \le d_E-q-2$.

The long cohomology sequence associated to the exact sequence
  \[
  0 \to \sO_E(F'-(n+1)mH) \to \sO_E(F'- nmH) \to \sO_Z(F' - nZ) \to 0
  \]
  implies our claim for $q = 0$. For $q \ge 1$, we let $\iota \colon Z \inj E$ be
  the inclusion and use the long cohomology
  sequence associated to the exact sequence
  \[
  0 \to \Omega^q_E(\log F)(F' -(n+1)mH) \to  \Omega^q_E(\log F)(F' - nmH) \to
  \iota^* \Omega^q_E(\log F)(F' - nmH) \to 0
  \]
  and ~\eqref{eqn:Coh-vanishing-0} to get
  $H^i_\zar(Z,  \iota^* \Omega^q_E(\log F)(F' - nmH)) = 0$ for all $i \le d_E-2$
  and $n \ge 1$.
 We combine this with the long cohomology sequence associated to the exact sequence
  \[
  0 \to \Omega^{q-1}_Z(\log F_Z)(F' - (n+1)mH) \to \iota^* \Omega^q_E(\log F)(F' - nmH)
  \to \Omega^q_Z(\log F_Z)(F' - nmH) \to 0
  \]
and induction on $q$ to conclude the proof of the lemma. 
\end{proof}

\begin{cor}\label{cor:Coh-vanishing-1}
  Let $k$ be an $F$-finite field. Let
  $\{(E_1, F_1, D_1), \ldots , (E_r, F_r, D_r)\}$ be a collection of triples, where
  $E_j$ is a connected and regular closed subscheme of $\P^N_k$ of dimension at least
  two, $F_j \subset E_j$ is an snc divisor and $D_j \in \Div(E_j)$.
  Then there is $m_0 \gg 0$ such that for every
  $m \ge m_0$ and $H \in |\sO_{\P^N_k}(m))|(k)$, one has 
  \[
  H^i_\zar(E^H_j, \Omega^q_{E^H_j}(\log F^H_j)((D_j - nE^H_j)|_{E^H_j}) = 0
  \]
  for all $q \in \N_0$,  $n \in \N$, $j \in J_r$ and $i \le \dim(E_j) -q-2$ provided
  $E^H_j :=  E_j \bigcap H$ is regular and ${F^H_j} := F_j \bigcap E^H_j$ is an snc
  divisor on $E^H_j$.
\end{cor}
\begin{proof}
Let $W$ be set of ordered pairs $(j,q)$ such that $j \in J_r$ and
  $q \in J^0_{N+ p_k}$, where $p_k$ is the $p$-rank of $k$.
Then for a given pair $a = (j,q) \in W$, \lemref{lem:Coh-vanishing} says that
  there exists an integer $m(a) \gg 0$ such that for every
  $m \ge m(a)$ and $H \in |\sO_{\P^N_k}(m)|(k)$, one has 
  $H^i_\zar(E^H_j, \Omega^q_{E^H_j}(\log F^H_j)(\iota^*_j(D_j - nE^H_j)) = 0$ for all
  $n \ge 1$ and $i \le \dim(E_j) -q-2$ provided $E^H_j$ is
  regular and ${F^H_i}$ is an snc divisor on $E^H_j$.
 Letting $m =
  {\rm max}\{m(a)| a \in W\}$ and noting that $\Omega^q_{E^H_j}(\log F^H_j) = 0$
  for $q > N+ p_k$, we arrive at the desired conclusion.
\end{proof}

We let $\sX$ be a separated and finite type $S$-scheme of dimension $d_\sX \ge 1$.
Let $f \colon \sX \to S$ denote the structure map.
Given an irreducible subscheme $Z \subset \sX$, we  shall let $Z_\eta$ (resp. $Z_s$)
denote the generic (resp. special) fiber of $Z$.  We shall say that $Z$ is
vertical if the support of $\ov{Z}$ is contained in the special fiber of $f$. Else,
we shall say that $Z$ is horizontal (equivalently, $R$-flat). 

\begin{defn}\label{defn:QSS}
  We shall say that $\sX$ is a quasi-semistable $R$-scheme if every connected
  component $\sX'$ of $\sX$ is regular and faithfully flat over $S$ and its reduced
  special fiber $(\sX'_s)_\red$ is an snc divisor on $\sX'$. We shall say that
  $\sX$ is semistable if it is quasi-semistable and its generic fiber $\sX_\eta$
  is smooth over $\Spec(k)$. 
  \end{defn}

By a morphism of quasi-semistable schemes
  $\phi \colon (T', S') \to (T, S)$, we shall mean a commutative diagram
  \begin{equation}\label{eqn:SS-map}
    \xymatrix@C1.2pc{
      T' \ar[r]^-{\phi} \ar[d]_-{f'} & T \ar[d]^-{f} \\
      S' \ar[r] & S,}
  \end{equation}
  where the bottom horizontal arrow is the morphism between spectra induced by
  a dominant local ring homomorphism of henselian discrete valuation rings and the
  vertical arrows are
  the structure maps of quasi-semistable schemes over them. Note that
  $\phi$ preserves the generic and special fibers, and moreover,
  $\phi \colon (T', Y') \to (T, Y)$ is a morphism of snc-pairs.

For the rest of Sections~\ref{sec:Esp} and ~\ref{sec:Esp-0}, we fix
a quasi-semistable $R$-scheme $\sX$ which is connected of dimension
$d_\sX \ge 2$ and let $f \colon \sX \to S$ be the structure map.
We let $Y=(\sX_s)_\red$ and let $Y_1, \ldots , Y_r$ denote the
irreducible components of $Y$. We let $Y^o$ denote the regular locus of
$Y$ and let $Y^o_i = Y^o \bigcap Y_i$. We let $j \colon \sX_\eta \inj \sX$ and
$j' \colon \sX_s \inj \sX$ be the inclusions. We shall write $X = \sX_\eta$. 
We shall say that $\sX$ is a relative curve over $S$ if $d_\sX = 2$.
By a {\sl subscheme} of $\sX$, we shall mean a locally closed subscheme of $\sX$.

\begin{defn}\label{defn:good-sub}
  $(1)$ Let $\sX' \subset \sX$ be a subscheme. Consider the following conditions.
\begin{listabc}
\item
  $\sX'$ is connected, regular and faithfully flat over $S$.
\item
  $\sX' \times_\sX Y$ is an snc divisor in $\sX'$.
  \item
    $\sX' \times_{\sX} Y_i$ is regular and connected for every $i \in J_r$.
  \item
    $\sX'_\eta$ is projective and geometrically integral over $k$ if $\sX_\eta$ is
    projective and geometrically integral over $k$.
  \item
$\sX'_\eta$ is smooth over $k$ if $\sX_\eta$ is smooth over $k$.
\end{listabc}

We shall say that $\sX'$ is a `quasi-admissible' subscheme of $\sX$ if it
satisfies the first three of the above conditions. We shall say that $\sX'$ is
an `admissible' subscheme of $\sX$ if it satisfies all five conditions above.

If $\sX'$ is closed (note that it is always separated) in $\sX$ and 
quasi-admissible (resp. admissible), we shall say that $\sX'$ is a quasi-admissible
(resp.  admissible) closed subscheme. If $\sX'$ is a quasi-admissible (resp.
admissible) subscheme of $\sX$ with $d_{\sX'} = 2$, we shall
say that $\sX'$ is a quasi-admissible (resp. admissible) relative curve 
inside $\sX$.

$(2)$ If $\sX$ is quasi-projective over $S$ with a given embedding
  $\sX \inj \P^N_R$, we shall say that a hypersurface
$H \in |\sO_{\P^N_R}(n)|(R)$ is quasi-admissible (resp. admissible)
relative to $\sX$ if $\sX' := \sX \bigcap H$
is a quasi-admissible (resp. admissible) subscheme of $\sX$. In this case, we shall
say that $\sX'$ is a 
  quasi-admissible (resp. admissible) hypersurface section of $\sX$ (of degree $n$).

  $(3)$ If $\sX' \subset \sX$ is a quasi-admissible subscheme, we shall let
  $X' = \sX'_\eta$,
  $\sX'_s = \sX' \bigcap \sX_s, \ Y' = \sX' \bigcap Y, \ Y'_i =  \sX' \bigcap Y_i$
  and $Y'^o = Y'_\reg$.
\end{defn}

\begin{lem}\label{lem:good-sub-0}
  If $\sX' \subset \sX$ is a quasi-admissible subscheme of $\sX$, then the following
  hold.
  \begin{enumerate}
  \item
    $Y'_i \neq Y'_j$ if $i \neq j$.
  \item
    $Y'_\sing = Y_\sing \bigcap \sX'$.
  \item
    $Y'^o = Y^o \times_{\sX} \sX'$.
  \item
    $X'^o_{\fin} \subset X^o_\fin$.
  \item
    $\ov{P} \times_{\sX} rY \cong \ov{P} \times_{\sX'} rY'$ for every
    $P \in X'^o_{\fin}$ and $r \ge 1$.
  \item
   The inclusion $(\sX', Y') \inj (\sX, Y)$ is a morphism of snc-pairs. 
    \end{enumerate}
\end{lem}
\begin{proof}
  If $Y'_i = Y'_j$ for some $i \neq j$, then one checks that $\sX' \times_\sX Y$
  is a non-reduced closed subscheme of $\sX'$ which contradicts the
  condition (b) of the above definition. This proves item (1).
If $\sX' \subset \sX$ is a quasi-admissible subscheme, one checks that
  $Y'_\sing = {\underset{i \neq j}\bigcup} (Y'_i \bigcap Y'_j) =
  {\underset{i \neq j}\bigcup}
  \left((Y_i \bigcap \sX') \bigcap (Y_j \bigcap \sX')\right) =
  \left({\underset{i \neq j}\bigcup} (Y_i \bigcap Y_j)\right) \bigcap \sX' =
  Y_\sing \bigcap \sX'$. This proves (2), and (3) follows from (2) and
  the definition of a quasi-admissible subscheme.

  Item (4) follows from (3) because the closure of $P$ in $\sX'$ is finite over $S$
  and hence is also a (regular) closed
  subscheme of $\sX$ for every $P \in X'^o_{(0)}$. To prove (5), it suffices to show
  that $\sX' \times_{\sX} rY^o \cong rY'^o$. But this follows from (3) because
  $Y^o$ is locally closed in $\sX$. Finally, item (6) is clear from the definition of
  a quasi-admissible subscheme.
\end{proof}

\subsection{Specialization of Swan conductor: case~I}\label{sec:SP-1}
Let $\sX$ be the above quasi-semistable $R$-scheme.
We let $m \in \N, q \in \N_0$ and $\chi \in H^{q+1}_{p^m}(X)$.
If $Z \subset X$ is a subscheme,
then $\chi|_Z$ will denote the restriction of $\chi$ to $H^{q+1}_{p^m}(Z)$.
For $j \in J_r$, we let $n_j = \Sw_{Y_j}(\chi) \ge 0$.
We let $D_\chi = \stackrel{r}{\underset{j =1}\sum} n_j Y_j$.  This is sometimes
referred to as the `Swan divisor' of $\chi$.

We now fix $i \in J_r$ and assume that $n_i \in \N$.
If $D = \stackrel{r}{\underset{j =1}\sum} m_j Y_j \in \Div_Y(\sX)$, then
we shall say that $D$ is $(\chi, i)$-admissible if $m_j \ge n_j$ for
every $j$ and $m_i = n_i$. We let $D^{\<i\>}_n = D - nY_i$ for any $n \in \N_0$. Write
$F_i = \sum\limits_{\substack{j \in J_r \setminus \{i\}}}  (Y_j \cap Y_i) \in \Div(Y_i)$ and
$Y^o_i = Y_i \setminus F_i$.

Our first result concerning the behavior of Swan conductor under restriction to
subschemes is the following.

\begin{prop}\label{prop:SC-change}
 Assume that one of the following conditions holds.
 \begin{enumerate}
 \item
   $\chi$ has type I at $Y_i$ and $q \le d_\sX -2$.
 \item
   $\chi$ has type II at $Y_i$ and $q \le d_\sX -3$.
 \end{enumerate}
 Then we have the following.

 $(A)$ There exists a dense open $U_i \subset Y_i$ such that for every closed
 point $x \in U_i$, we can find a quasi-admissible subscheme $\sX' \subset \sX$ of
 dimension $d_\sX -1$ containing $x$ which satisfies the following.
\begin{list}{(\alph{elno-abc})}{\usecounter{elno-abc}}
\item
  $\Sw_{Y'_i}(\chi|_{X'}) = n_i$.
\item $\chi|_{X'}$ has the same type at $Y'_i$ as $\chi$ has at
  $Y_i$. 
\end{list}

$(B)$ If $\sX$ is semistable and there exists a closed embedding
$\sX \inj \P^N_R$ of $R$-schemes, then for all $n \gg 0$, we can find an admissible
hypersurface section $\sX' \subset \sX$ of degree $n$ which satisfies conditions
(a) and (b) of $(A)$.
\end{prop}
\begin{proof}
We can assume without loss of generality that $i =1$.
We let $D = \stackrel{r}{\underset{j =1}\sum} m_j Y_j$ be any
$(\chi, 1)$-admissible divisor so that $\chi \in \Fil_D H^{q+1}_{p^m}(X)$.
Let $D_1 = D|_{Y_1}$.

Suppose that $\chi$ has type I at $Y_1$ and $q \le d_\sX -2$.
If we let $\alpha$ denote the image of $\chi$ under the composition
$\nu^0_1 \circ \Rsw^{m, q+1}_{\sX|(D, D^{\<1\>}_1)}$, then our assumption is that
$\alpha \neq 0$ in $H^0_\zar(Y_1, \Omega^q_{Y_1}(\log F_1)(D_1))$.
Since $Y_1$ is a finite type integral $\F$-scheme of positive dimension and
$\Omega^q_{Y_1}(\log F_1)(D_1)$ is a finite rank locally free sheaf on $Y_1$,
an easy exercise shows that there is an affine  dense open 
$\sW \subset \sX \setminus F_1$ such that $\sW \bigcap Y_1$ is a principal
divisor on $\sW$ and the image of
$\alpha$ in $\Omega^q_{Y_1}(\log F_1)(D_1) \otimes_{\sO_{Y_1}} k(x)$
is not zero for all points $x \in \sW \bigcap Y_1$.
We let $\sW = \Spec(A)$ and $U_1 = \sW \bigcap Y_1 = \Spec({A}/{(\pi_1)})$.

We fix a closed point $x \in U_1$.
We assume $\sW$ to be small enough so that $A$ has elements
  $\pi_1, t_2, t_3, \ldots , t_{d_\sX}, u_1, \ldots , u_{p_\F}$ such that
$\fm_x = (\pi_1, t_2, \ldots , t_{d_\sX})$ is the ideal with quotient $k(x)$
and $\Omega^1_{k(x)} = \ \stackrel{p_\F}{\underset{j =1}\oplus}
  k(x)d\ov{u}_j$, where $\ov{a} = a$ mod $\fm_x$ for any $a \in A$.
We then have $\Omega^1_{A} = A d\pi_1 \bigoplus
  \left(\stackrel{d_{\sX}}{\underset{j =2}\oplus}
  Adt_j\right) \bigoplus \left(\stackrel{p_\F}{\underset{j =1}\oplus}
  Adu_j\right)$ (cf. \lemref{lem:Rank-module-diff}).
  We let $A' = {A}/{(t_2)}, \ A_1 = A/{(\pi_1)}, \ A'_1 = {A}/{(\pi_1, t_2)}$,
  $M_1 = \stackrel{d_\sX}{\underset{j =3}\oplus} Adt_j$,
  $M_2 = \stackrel{p_\F}{\underset{j =1}\oplus} Adu_j$ and
  $M = M_1 \bigoplus M_2$.

We then get
\begin{equation}\label{eqn:SC-change-0}
    \begin{array}{lll}
\Omega^{q}_A & = & {\bigwedge}^{q}_A(M) \bigoplus
  \left({\bigwedge}^{q-1}_A(M) \otimes_A Ad\pi_1\right)
  \bigoplus \left({\bigwedge}^{q-1}_A(M) \otimes_A Adt_2 \right) \\
  & & \hspace*{1.6cm}
  \bigoplus \left( {\bigwedge}^{q-2}_A(M)\otimes_A A(d\pi_1 \wedge dt_2) \right); \\ 
\Omega^{q}_{A_1} & = &\left({\bigwedge}^{q}_A(M) \otimes_A A_1\right) \bigoplus
\left({\bigwedge}^{q-1}_A(M) \otimes_A Adt_2\right) \otimes_A A_1. \\
\Omega^{q}_{A'_1} & = & {\bigwedge}^{q}_A(M) \otimes_A A'_1.
\end{array}
\end{equation}

If we let $\alpha^o$ denote the restriction of $\alpha$ to $U_1$, we
  can write (after ignoring the twist)
  \begin{equation}\label{eqn:SC-change-1}
  \alpha^o = \sum_I a_I u_I  + \sum_{J} c_{J} v_{J} \wedge dt_2,
  \end{equation}
  where $a_I, c_{J} \in A_1$ and  $u_I$
  (resp. $v_J$) runs through the wedge products of $q$ (resp. $q-1$) distinct 
  elements of the set
  $\sigma(M) = \{d t_i, du_j| 3 \le i \le d_\sX, 1 \le j \le p_\F\}$.
  By our assumption, either $a_I \notin {\fm_x}/{(\pi_1)}$ for some
  $I$ or $c_J \notin {\fm_x}/{(\pi_1)}$ for some $J$.

In the first case, we take $\sX' = \Spec(A')_\reg$ and
  let $u \colon \sX' \inj \sX$ be the inclusion. It easily follows from
  ~\eqref{eqn:SC-change-0} that $\alpha$ does not die under the map
\begin{equation}\label{eqn:SC-1}
    u^* \colon  H^0_\zar(Y_1, \Omega^q_{Y_1} (\log F_1)(D_1)) \to
  H^0_\zar(Y'_1, \Omega^q_{Y'_1} (\log F'_1)(D'_1))
  \end{equation}
if we let $D'_1 = D_1 \times_{Y_1} Y'_1$ and $F'_1 = F_1 \times_{Y_1} Y'_1$.
We conclude from \corref{cor:RSW-gen-1} that $\Sw_{Y'_i}(\chi|_{X'}) = n_1$, and
from \corref{cor:Type-I-global-0} that $\chi|_{X'}$ has type I at $Y'_1$.
Since $\sX'$ is regular and so is $\sX' \bigcap Y = \sX' \bigcap Y_1 = \Spec(A'_1)$,
it follows that $\sX'$ is a quasi-admissible subscheme of $\sX$.
Note that $\sX'$ is faithfully flat over
$S$ since it is connected, regular, contains $x$ and is not contained in $\sX_s$.

Suppose on the other hand that $c_J \notin {\fm_x}/{(\pi_1)}$ for some $J$.
In this case, not all $t_j$ in the set $\{t_3, \ldots , t_{d_\sX}\}$ can appear
when we write $v_J$ as a wedge product of $q-1$ distinct elements of $\sigma(M)$
since $q-1 \le d_\sX-3$ by our assumption. We assume without loss of generality
that $t_3$ does not appear in the wedge product presentation of $v_J$.
We let $\sX' = \Spec({A}/{(t_3)})_\reg$. Then an identical argument as in the
first case shows that $\sX'$ is a quasi-admissible subscheme of $\sX$ and
$\Sw_{Y'_i}(\chi|_{X'}) = n_1$. If $\chi$ has type II at $Y_1$ and $q \le d_\sX -3$,
we repeat the above argument verbatim with $q$ replaced by $q+1$ and note that
$q = (q+1)-1 \le d_\sX -3$. This completes the proof of $(A)$.

We assume now that $\sX$ is semistable and there is a closed  embedding
$\sX \inj \P^N_R$ of $R$-schemes.
It follows in this case from \thmref{thm:Bertini-dvr} that there exists an admissible
hypersurface section $\sX' \subset \sX$ of degree $n$ for all $n \gg 0$.
To prove  part {\sl ($a$)} for $\sX'$,
we assume first that  $\chi$ has type I at $Y_1$ and $q \le d_\sX -2$. 
By \lemref{lem:RSW-spl} and \corref{cor:Type-I-global-0}, it suffices to show that
the map 
  \begin{equation}\label{eqn:SC-2}
      H^0_\zar(Y_1, \Omega^q_{Y_1} (\log F_1)(D_1)) \to
  H^0_\zar(Y'_1, \Omega^q_{Y'_1} (\log F'_1)(D'_1))
  \end{equation}
  is injective if we let
  $D'_1 = D_1 \times_{Y_1} Y'_1$ and $F'_1 = F_1 \times_{Y_1} Y'_1$.
  By the second exact sequence of \lemref{lem:Kahler-seq}, this is equivalent
  to showing that $H^0_\zar(Y_1, \sF) = 0$, where
  $\sF = \Omega^q_{Y_1}(\log (F_1 + Y'_1))(D_1 - Y'_1)$.

To prove the vanishing of $H^0_\zar(Y_1, \sF)$, we note that $Y_1$ is regular and
  projective of dimension $d_\sX-1$ over $\F$. We now use the first exact sequence of
  \lemref{lem:Kahler-seq} to get an exact sequence
\begin{equation}\label{eqn:SC-change-00}
  0 \to \Omega^q_{Y_1}(\log F_1)(D_1 - Y'_1) \to
  \Omega^q_{Y_1}(\log (F_1 + Y'_1))(D_1 - Y'_1) \xrightarrow{\res}
  \Omega^{q-1}_{Y'_1}(\log F'_1)(D_1 - Y'_1) \to 0.
\end{equation}
Note that $F_1 + Y'_1$ is an snc divisor on $Y_1$ as $\sX'$ is quasi-admissible
(cf. \lemref{lem:good-sub-0}(1)).  Since $n \gg 0$, we have that
$H^i_\zar(Y_1, \Omega^q_{Y_1}(\log F_1)(D_1 - Y'_1)) =0$
for all $i \le d_\sX-2$ by ~\eqref{eqn:Coh-vanishing-0}, and
$H^i_\zar(Y'_1, \Omega^{q-1}_{Y'_1}(\log F'_1)(D_1 - Y'_1)) = 0$ for $i \le d_\sX-q-2$ by
\lemref{lem:Coh-vanishing}.
It follows from ~\eqref{eqn:SC-change-00} that
$H^i_\zar(Y_1, \Omega^q_{Y_1}(\log (F_1 + Y'_1))(D_1 - Y'_1)) = 0$
for $i \le d_\sX -q-2$. Since $d_\sX -q -2 \ge 0$, we get the desired vanishing,
and hence ~\eqref{eqn:SC-2}.

If $\chi$ has type II at $Y_1$ and $q \le d_\sX -3$, then an identical proof shows
that the map
\[
H^0_\zar(Y_1, \Omega^{q+1}_{Y_1} (\log F_1)(D_1)) \to
H^0_\zar(Y'_1, \Omega^{q+1}_{Y'_1} (\log F'_1)(D'_1))
\]
is injective and this implies that $\Sw_{Y'_1}(\chi|_{X'}) = n_1$.
To conclude the proof of the lemma, we note that the proof above shows that
$\chi|_{X'}$ has the same type at $Y'_1$ as $\chi$ has at $Y_1$.
\end{proof}

\subsection{Specialization of Swan conductor: case~II}\label{sec:SP-2}
We continue with the setting and notation of \S~\ref{sec:SP-1}.
We now study the drop in the Swan conductor in the case where 
$\chi \in H^{q+1}_{p^m}(X)$ has type II at $Y_i$ and $q = d_\sX -2 \ge 0$.
If $n_i=1$, then \corref{cor:RSW-gen-1} implies that
$0\le \Sw_{Y'_i}(\chi|_{X'})\le 1$ for any quasi-admissible subscheme $\sX'$ of $\sX$.
In what follows, the following four conditions will collectively be referred to as
\textbf{conditions {($\star$)}}.

\begin{enumerate}
  \item
      $\chi$ is of type II at $Y_i$.
      \item
    $q = d_\sX -2 \ge 0$.
    \item
      $n_i \ge 2$.
    \item
      Either $p \neq 2$ or $n_i \ge 3$.
  \end{enumerate}

We let $D = \stackrel{r}{\underset{j =1}\sum} m_j Y_j$ be any
$(\chi, i)$-admissible divisor so that $\chi \in \Fil_D H^{q+1}_{p^m}(X)$.
Assume that  the conditions {$(\star)$} hold.
We then have the inequalities
$D \ge D^{\<i\>}_1 \ge D^{\<i\>}_2 \ge D/p \ge 0$. We therefore have an injective map
\begin{equation}\label{eqn:Type-2-1*}
\Rsw^{m,q+1}_{\sX|(D, D^{\<i\>}_s)} \colon
\frac{\Fil_D H^{q+1}_{p^m}(X)}{\Fil_{D^{\<i\>}_s} H^{q+1}_{p^m}(X)} \to
H^0_\zar(Y_i, \Omega^{q+1}_\sX(\log Y)(D) \otimes_{\sO_{\sX}} \sO_{sY_i})
\end{equation}
for $s \in \{1,2\}$ by \thmref{thm:RSW-gen}.
In particular, the image of $\chi$ under this map is not zero.
We let $\wt{\chi} = \Rsw^{m,q+1}_{\sX|(D, D^{\<i\>}_2)}(\chi) \neq 0$.
Let $Y_o= Y - Y_i \in \Div(\sX)$.

\begin{lem}\label{lem:Type-2-0}
The sequence
  \begin{equation}\label{eqn:Type-2-1}
    0 \to \frac{\Omega^{q+1}_\sX(\log Y_o)(D)}{\Fil_{D^{\<i\>}_2}\Omega^{q+1}_{X}}
    \xrightarrow{\delta_i}
  \frac{\Fil_{D}\Omega^{q+1}_{X}}{\Fil_{D^{\<i\>}_2}\Omega^{q+1}_{X}} \xrightarrow{\psi_i}
  \Omega^q_{Y_i}(\log F_i)(D_i) \to 0
  \end{equation}
  is exact, where $\psi_i$ is the composite map
  $\frac{\Fil_{D}\Omega^{q+1}_{X}}{\Fil_{D^{\<i\>}_2}\Omega^{q+1}_{X}}
  \stackrel{\eta_i}{\surj}
  \frac{\Fil_{D}\Omega^{q+1}_{X}}{\Fil_{D^{\<i\>}_1}\Omega^{q+1}_{X}} \xrightarrow{\nu_i}
  \Omega^q_{Y_i}(\log F_i)(D_i)$.
\end{lem}
\begin{proof}
See \cite[Cor.~9.10]{KM-1}.
\end{proof}

Since $\chi$ has type II at $Y_i$, it follows from \lemref{lem:Type-2-0} that
$\wt{\chi}$ is a non-zero element of
$H^0_\zar\left(\sX, \frac{\Omega^{q+1}_\sX(\log Y_o)(D)}{\Fil_{D^{\<i\>}_2}\Omega^{q+1}_X}\right)$. 
We now look at the commutative diagram
\begin{equation}\label{eqn:Type-2-2}
  \xymatrix@C.7pc{
    \frac{\Omega^{q+1}_\sX(\log Y_o)(D)}{\Fil_{D^{\<i\>}_2}\Omega^{q+1}_X} 
    \ar@{->>}[r]^-{\phi_i} \ar[dr]_-{\theta_i} &
    \frac{\Omega^{q+1}_\sX(\log Y_o)(D)}{\Omega^{q+1}_\sX(\log Y_o)(D^{\<i\>}_1)} 
    \ar[r]^-{\cong} & \Omega^{q+1}_\sX(\log Y_o)(D) \otimes_{\sO_\sX} \sO_{Y_i}
    \ar@{->>}[dl]^-{\gamma_i} \\
    & \Omega^{q+1}_{Y_i}(\log F_i)(D_i), &}
\end{equation}
where $\theta_i$ and $\gamma_i$ are the canonical maps on the log differential
forms induced by the inclusion $Y_i \inj \sX$ and
$\phi_i$ is the canonical surjection induced by the inclusions
$\Fil_{D^{\<i\>}_2}\Omega^{q+1}_X  \subset \Omega^{q+1}_\sX(\log Y_o)(D^{\<i\>}_1)
\subset \Omega^{q+1}_\sX(\log Y_o)(D)$. 

\begin{lem}\label{lem:Type-2-3}
  We have $(\gamma^*_1 \circ \phi^*_1)
  \left(\wt{\chi}\right) = \theta^*_i(\wt{\chi}) =
  \Rsw^{m,q+1}_{\sX|(D, D^{\<i\>}_1)}(\chi) \neq 0$. In particular,
  $\phi^*_i(\wt{\chi}) \neq 0$.
\end{lem}
\begin{proof}
By ~\eqref{eqn:Type-2-1*} and ~\eqref{eqn:Type-2-2}, we only have
to show that $\theta^*_i(\wt{\chi}) = \Rsw^{m,q+1}_{\sX|(D, D^{\<i\>}_1)}(\chi)$.
To that end, one first checks easily that the diagram of complexes of sheaves
  (where $\nu'_i$ is from \lemref{lem:RSW-spl})
  \begin{equation}\label{eqn:Type-2-4}
  \xymatrix@C2pc{
    \frac{\Omega^{q+1}_\sX(\log Y_o)(D)}{\Fil_{D^{\<i\>}_2}\Omega^{q+1}_X}
    \ar@{^{(}->}[r]^-{\delta_i} \ar[d]_-{\theta_i} &
    \frac{\Fil_{D}\Omega^{q+1}_X}{\Fil_{D^{\<i\>}_2}\Omega^{q+1}_X}
    \ar@{->>}[d]^-{\eta_i} &
    \frac{W_m\sF^{q, \bullet}_{\sX|D}}{W_m\sF^{q, \bullet}_{\sX|D^{\<i\>}_2}}[1] \ar[d]
    \ar[l]_{F^{m-1}d} \\
    \Omega^{q+1}_{Y_i}(\log F_i)(D_i) \ar[r]^-{\nu'_i} &
    \frac{\Fil_{D}\Omega^{q+1}_X}{\Fil_{D^{\<i\>}_1}\Omega^{q+1}_X} &
    \frac{W_m\sF^{q, \bullet}_{\sX|D}}{W_m\sF^{q, \bullet}_{\sX|D^{\<i\>}_1}}[1]
    \ar[l]_-{F^{m-1}d}}
  \end{equation}
  commutes (cf. the proof of \cite[Lem.~9.9]{KM-1}). In particular, the diagram
  resulting from applying $\H^0_\et(-)$) also commutes.
  The commutativity of the right square shows that
  $\eta_i \circ \delta_i(\wt{\chi})$ is the refined Swan
  conductor $\Rsw^{m,q+1}_{\sX|(D, D^{\<i\>}_1)}(\chi)$. Since $\chi$ has type II with
  Swan conductor $n_i$ at $Y_i$, this refined Swan conductor actually lies
  in $H^0_\zar(Y_i,\Omega^{q+1}_{Y_i}(\log F_i)(D_i))$. This proves the desired
  identity.
 \end{proof}

Suppose now that there exists a quasi-admissible subscheme $u \colon \sX' \inj \sX$.
Let $Y'_o = Y' - Y'_i$, $D'=D\times_\sX \sX'$ (and similarly for $D'^{\<i\>}_j$).
We consider the diagram
\begin{equation}\label{eqn:Type-2-5}
  \xymatrix@C.7pc{
    \frac{\Fil_D H^{q+1}_{p^m}(X)}{\Fil_{D^{\<i\>}_2} H^{q+1}_{p^m}(X)} \ar[r] \ar[d] &
    \frac{\Fil_{D'} H^{q+1}_{p^m}(X')}{\Fil_{D'^{\<i\>}_2} H^{q+1}_{p^m}(X')} \ar[d] \\
    H^0_\zar\left(\sX, \frac{\Fil_{D}\Omega^{q+1}_X}{\Fil_{D^{\<i\>}_2}\Omega^{q+1}_X}\right)
    \ar[r] &
    H^0_\zar\left(\sX', \frac{\Fil_{D'}\Omega^{q+1}_{X'}}
    {\Fil_{D'^{\<i\>}_2}\Omega^{q+1}_{X'}}\right) \\
 H^0_\zar\left(\sX, \frac{\Omega^{q+1}_\sX(\log Y_o)(D)}{\Fil_{D^{\<i\>}_2}\Omega^{q+1}_X}\right)
 \ar[r] \ar[d]_-{\phi^*_i} \ar@{^{(}->}[u]^-{\delta_i} &
 H^0_\zar\left(\sX', \frac{\Omega^{q+1}_{\sX'}(\log Y'_o)(D')}
 {\Fil_{D'^{\<i\>}_2}\Omega^{q+1}_{X'}}\right)
 \ar[d]^-{\phi'^*_i} \ar@{^{(}->}[u]_-{\delta'_i} \\
H^0_\zar\left(\sX,
\frac{\Omega^{q+1}_\sX(\log Y_o)(D)}{\Omega^{q+1}_\sX(\log Y_o)(D^{\<i\>}_1)}\right)
\ar[r] \ar[d]_-{\cong} & 
H^0_\zar\left(\sX',
\frac{\Omega^{q+1}_{\sX'}(\log Y'_o)(D')}{\Omega^{q+1}_{\sX'}(\log Y'_o)(D'^{\<i\>}_1)}\right)
\ar[d]^-{\cong} \\
H^0_\zar\left(\sX, \Omega^{q+1}_\sX(\log Y_o)(D)\otimes_{\sO_\sX} \sO_{Y_i}\right)
\ar[r]^-{u^*} &
H^0_\zar\left(\sX', \Omega^{q+1}_{\sX'}(\log Y'_o)(D')\otimes_{\sO_{\sX'}} \sO_{Y'_i}\right).}
\end{equation}

In the above diagram, the horizontal arrows are induced by the pull-back via the
inclusion $u \colon \sX' \inj \sX$. The vertical arrows in the top square
are the refined Swan conductor maps.
The top square commutes by \thmref{thm:RSW-gen}(3) because $\sX$ is a
quasi-admissible subscheme of $\sX$. It is clear that
the other squares commute. Using \lemref{lem:Type-2-3} and a diagram chase in
~\eqref{eqn:Type-2-5}, we conclude the following.

\begin{cor}\label{cor:Type-2-6}
  Assume that the conditions {$(\star)$} hold and $u \colon \sX' \inj \sX$ is the
  inclusion of a quasi-admissible subscheme such that
  $u^* \circ \phi^*_i(\wt{\chi}) \neq 0$ under the map
  \[
  u^* \colon
  H^0_\zar\left(\sX, \Omega^{q+1}_\sX(\log Y_o)(D)\otimes_{\sO_\sX} \sO_{Y_i}\right)
  \to
H^0_\zar\left(\sX', \Omega^{q+1}_{\sX'}(\log Y'_o)(D')\otimes_{\sO_{\sX'}} \sO_{Y'_i}\right)
  \]
  for some $(\chi, i)$-admissible divisor $D$.
Then $n_i - 1 \le \Sw_{Y'_i}(\chi|_{X'}) \le n_i$.
\end{cor}

\vskip.2cm

In order to find $u \colon \sX' \inj \sX$ for which
$u^* \circ \phi^*_i(\wt{\chi}) \neq 0$, we shall use the following local
result.

Assume that the conditions {$(\star)$} hold. Fix $i \in J_r$ and a closed
  point $x \in Y^o_i$. Let $A = \sO_{\sX,x}$ and let
  $\fm = (\pi_i, t_2, \ldots , t_{d_\sX})$ be the maximal ideal of $A$, where
  $\pi_i$ defines $Y_i$ locally at $x$. Let $t^\dagger_2 = t_2 - b\pi_i$, where
  $b \in A^\times$. Let $\sX_x = \Spec(A), \sX^1_x = \Spec({A}/{(t_2)})$ and
  $\sX^2_x = \Spec({A}/{(t^\dagger_2)})$. Let $Y^1_i = \Spec({A}/{(\pi_i, t_2)})$ and
  $Y^2_i = \Spec({A}/{(\pi_i, t^\dagger_2)})$. Let
  $D = \stackrel{r}{\underset{j =1}\sum} m_j Y_j$ be a
   $(\chi, i)$-admissible divisor.

\begin{lem}\label{lem:Type-2-key}
    Let $\beta \in
    H^0_\zar\left(\sX, \Omega^{q+1}_\sX(\log Y_o)(D)\otimes_{\sO_\sX} \sO_{Y_i}\right)$
    be such that the image of $\gamma^*_i(\beta)$ is not zero in
 $\Omega^{q+1}_{Y_i}(\log F_i)(D_i) \otimes_{\sO_{Y_i}} k(x)$.
 Then there exists $j \in J_2$ such that the image of $\beta$ is not zero in 
 $H^0_\zar(\sX^j_x, \Omega^{q+1}_{\sX^j_x} \otimes_{\sO_{\sX^j_x}} \sO_{Y^j_i})$.
 \end{lem}
\begin{proof}
We can assume without loss of generality that $i = 1$. 
We let $A_1 = {A}/{(\pi_1)}, A' = A/{(t_2)}, A^\dagger = A/{(t^\dagger_2)},
\ A'_1 = {A}/{(\pi_1, t_2)}$ and $A^\dagger_1 = {A}/{(\pi_1, t^\dagger_2)}$. We let
$M = M_1 \bigoplus M_2$ be as in \propref{prop:SC-change}. We then get

\begin{equation}\label{eqn:Type-2-8-0}
    \begin{array}{lll}
  \Omega^{q+1}_A & = & {\bigwedge}^{q+1}_A(M) \bigoplus
  \left({\bigwedge}^{q}_A(M) \otimes_A Ad\pi_1 \right)
  \bigoplus \left({\bigwedge}^{q}_A(M) \otimes_A Adt_2 \right) \\
  & & \hspace*{1.6cm}
  \bigoplus \left({\bigwedge}^{q-1}_A(M) \otimes_A A(d\pi_1 \wedge dt_2)\right); \\ 
  \Omega^{q+1}_{A'} & = & \left({\bigwedge}^{q+1}_A(M) \otimes_A A'\right) \bigoplus
  \left({\bigwedge}^{q}_A(M) \otimes_A Ad\pi_1\right)  \otimes_A A'; \\
  \Omega^{q+1}_{A_1} & = &\left({\bigwedge}^{q+1}_A(M) \otimes_A A_1\right) \bigoplus
  \left({\bigwedge}^{q}_A(M) \otimes_A Adt_2 \right) \otimes_A A_1; \\
  \Omega^{q+1}_{A'_1} & = & {\bigwedge}^{q+1}_A(M) \otimes_A A'_1.
\end{array}
\end{equation}

We write (after ignoring the twist)
\begin{equation}\label{eqn:Type-2-8-1}
{\beta|_{A}} =
  \sum_I a_I u_I + \sum_J b_J v_J \wedge d\pi_1 +  
  \sum_{J} c_{J} v_{J} \wedge dt_2 + \sum_{L} e_L w_L \wedge d\pi_1 \wedge dt_2,
  \end{equation}
where $a_I, b_J, c_{J},  e_L \in A$, and $u_I$, $v_J$ and $w_L$
run through the wedge products of $q+1$, $q$, and $q-1$ distinct 
  elements, respectively, of the set
  $\sigma(M) = \{d t_i, du_j| 3 \le i \le d_\sX, 1 \le j \le p_\F\}$
  (note that $F_1 \bigcap \sX_x = \emptyset$).
In particular,
  \begin{equation}\label{eqn:Type-2-8-2}
  {\beta}|_{A'} = \sum_I a_I u_I + \sum_J b_J v_J \wedge d\pi_1
  \  \mbox{mod}  \ (t_2) \ \ {\rm and} \ \
  {\beta}|_{A_1} =  \sum_I a_I u_I + \sum_{J} c_{J} v_{J}
\wedge dt_2 \ \mbox{mod} \ (\pi_1).
  \end{equation}

Assume first that one of the following holds:\\
\hspace*{.5cm}(1) there exists an index $I$ such that
$a_I$ does not die in $A'_1$;\\
\hspace*{.5cm}(2) there exists an index $J$ such that
$b_J$ does not die in $A'_1$. \\
In this case, the expression of ${\beta}|_{A'}$ shows that
$u^*({\beta}|_{A}) \neq 0$ in
$H^0_\zar(\sX^1_x, \Omega^{q+1}_{\sX^1_x} \otimes_{\sO_{\sX^1_x}} \sO_{Y^1_1})$ if
$u \colon \sX^1_x \inj \sX_x$ is the inclusion.

Suppose now that $a_I$ as well as $b_J$ dies in $A'_1$ for every
index $I$ and $J$. Since the image of ${\beta}|_{A_1}$
in $\Omega^{q+1}_{A_1} \otimes_{A_1} k(x)$ is not zero, we see from
~\eqref{eqn:Type-2-8-2} that there is an index $J$ such that $c_J$ does not die in
$k(x)$. Since $(t_2, \pi_1) = (t^\dagger_2, \pi_1)$ and
$\fm_x = (\pi_1, t^\dagger_2, t_3, \ldots , t_{d_\sX})$, we get
\[
  {\beta}|_{A^\dagger} = \sum_I a_I u_I + \sum_J (b_J +b c_J)v_J \wedge d\pi_1
+ \sum_J \pi_1 c_J v_J \wedge db +
\sum_{L} \pi_1 e_L w_L \wedge d\pi_1 \wedge db
\ \mbox{mod} \ (t^\dagger_2).
\]

Letting  $\beta'$ denote the image of
${\beta}|_{A^\dagger}$ in $\Omega^{q+1}_{A^\dagger} \otimes_{A^\dagger} k(x)$,
we get $\beta' = \sum_J b c_J v_J \wedge d\pi_1 \ \mbox{mod} \ \fm_x$.
Since $\{v_J \wedge d\pi_1\}_{J}$ is part of a basis of the $k(x)$-vector space
  $\Omega^{q+1}_{A^\dagger} \otimes_{A^\dagger} k(x)$ and $bc_J \notin \fm_x$, it
  follows that ${\beta}|_{A^\dagger} \neq 0$.
 Equivalently, $v^*({{\beta}|_{A}}) \neq 0$ in
$H^0_\zar(\sX^2_x, \Omega^{q+1}_{\sX^2_x} \otimes_{\sO_{\sX^2_x}} \sO_{Y^2_1})$ if
$v \colon \sX^2_x \inj \sX_x$ is the inclusion. This completes the proof.
\end{proof}

Let $D$ be the $(\chi,i)$-admissible divisor as above and let
$\beta$ denote the element $\phi^*_1(\wt{\chi})$ 
of $H^0_\zar(\sX, \Omega^{q+1}_{\sX}(\log Y_o)(D) \otimes_{\sX} \sO_{Y_i})$.
Combining  \corref{cor:Type-2-6} and \lemref{lem:Type-2-key}, we
get the following result which will a key role in handling type II cohomology
classes.

\begin{prop}\label{prop:Type-2-key-0}
        Let $x \in Y^o_i$ be a closed point and suppose that the image of
        $\gamma^*_1(\beta)$ is not zero in
        $\Omega^{q+1}_{Y_i}(\log F_i)(D_i)) \otimes_{\sO_{Y_i}} k(x)$.
        Then there exists $j \in J_2$ and a quasi-admissible subscheme
        $\sX' \subset \sX$ of dimension $d_\sX -1$
        such that $\sX^j_x = \sX' \times _{\sX} \Spec(A)$ and
        $n_i -1 \le \Sw_{Y'_i}(\chi|_{X'}) \le n_i$.
        \end{prop}
      \begin{proof}
 By \lemref{lem:Type-2-key}, the image of $\beta$
is not zero in $H^0_\zar(\sX^j_x, \Omega^{q+1}_{\sX^j_x} \otimes_{\sO_{\sX^j_x}} \sO_{Y^j_i})$
for some $j \in J_2$. We can clearly choose a quasi-admissible subscheme
$\sX' \subset \sX$ with the property that $\sX' \times_\sX \Spec(A) = \sX^j_x$
and $(\sX'_s)_\red \subset Y^o_i$. We then get from \corref{cor:Type-2-6} and the
functoriality of the refined Swan conductor (cf. \corref{cor:RSW-gen-1}) that
$n_i -1 \le \Sw_{Y'_i}(\chi|_{X'}) \le n_i$.
        \end{proof}

\begin{cor}\label{cor:Type-2-8}
  Assume that the conditions {$(\star)$} hold.
  Then there exists a dense open $U_i \subset Y_i$ such that for every closed
 point $x \in U_i$, we can find a quasi-admissible subscheme $\sX' \subset \sX$ of
 dimension $d_\sX -1$ containing $x$ such that
  $n_i - 1 \le \Sw_{Y'_i}(\chi|_{X'}) \le n_i$.
\end{cor}
\begin{proof}
  We can assume that $i =1$.
  We let $\beta = \phi^*_1 \left(\wt{\chi}\right)$ in
  $H^0_\zar\left(\sX, \Omega^{q+1}_\sX(\log Y_o)(D)\otimes_{\sO_\sX} \sO_{Y_1}\right)$.
  Then \lemref{lem:Type-2-3} implies that $\gamma^*_1(\beta)$ is not zero in
  $H^0_\zar(Y_1, \Omega^{q+1}_{Y_1}(\log F_1)(D_1))$.
  As in the proof of \propref{prop:SC-change}, there is an affine  dense open 
$\sW \subset \sX \setminus F_1$ such that the image of $\gamma^*_1(\beta)$ in
  $\Omega^{q+1}_{Y_1}(\log F_1)(D_1) \otimes_{\sO_{Y_1}} k(x)
  \cong \Omega^{q+1}_{Y_1} \otimes_{\sO_{Y_1}} k(x)$
is not zero for all points $x \in \sW \bigcap Y_1$.
The desired assertion now follows from \propref{prop:Type-2-key-0}.
 \end{proof}

\section{Specialization of Swan conductor II}\label{sec:Esp-0}
This section is a continuation of \S~\ref{sec:Esp}, in which we investigate the
specialization of the Swan conductor in the
remaining cases. In particular, we establish stronger results under the additional
assumption that $f \colon \sX \to S$ is projective. We keep the setting and
notation of \S~\ref{sec:Esp}.

\subsection{Specialization of Swan conductor: case III}\label{sec:SP-3}
In this subsection, we assume that $\sX$ is projective and semistable
over $S$. We fix a closed embedding $\sX \inj \P^N_R$.
Our goal is to prove a stronger version of \corref{cor:Type-2-8}
where we ask $\sX'$ to be an admissible hypersurface section of $\sX$ and not
just a quasi-admissible subscheme. We begin with some general results.

Let $\Sigma = \{x_1, \ldots , x_r\}$ be a set of closed points on $Y$ such that
$x_i \in Y^o_i$ for each $i \in J_r$. We let $B_\Sigma = \sO_{\sX, \Sigma}$. We let
$Y_i$ be locally defined at the points of $\Sigma$ by $\pi_i \in B_\Sigma$, and let
$\pi_0 = \stackrel{r}{\underset{i =1}\prod} \pi_i$. Recall here that $B_\Sigma$ is
a ufd (cf. \cite[Exc.~20.5]{Matsumura}).

\begin{lem}\label{lem:Type-2-10}
    For every $n \in \N$, there exists $m \ge n$ and
    $L \in |\sO_{\P^N_R}(m)|(R)$ such
    that a local defining equation of $L \bigcap \sX$ in $\Spec(B_\Sigma)$ is
    $u\pi_0 = 0$ for
    some $u \in B^\times_\Sigma$.
  \end{lem}
\begin{proof}
  We first apply \cite[Thm.~3.4, Lem.~4.1, 4.2, 6.2]{GK-JLMS} to find for all
  $n \gg 0$, an element $H \in |\sO_{\P^N_R}(n)|(R)$ such that $H \bigcap \Sigma =
    \emptyset$. We pick any $m \gg 0$ and $f \in H^0_\zar(\P^N_R, \sO_{\P^N_R}(m))$ such
    that $H = Z(f)$ has the property that $H \bigcap \Sigma = \emptyset$ and let
    $\Spec(\wt{B}) = \P^N_R \setminus H$.
    We let $\wt{B}_\Sigma = \sO_{\P^N_R, \Sigma}$ and let $\wt{\pi}_0$ be a pre-image of
    $\pi_0$ under the canonical surjection $\wt{B}_\Sigma \surj B_\Sigma$.
    We can then write $\wt{\pi}_0 = as^{-1}$, where $a, s \in \wt{B}$ and
    $s \in (\wt{B}_\Sigma)^\times$.

Using \cite[Lem.~II.5.14]{Hartshorne-AG} (with $\sF = \sO_{\P^N_R}$), we get that for
    all $n \gg 0$, there exists $g \in H^0_\zar(\P^N_R, \sO_{\P^N_R}(nm))$ whose
    restriction to $\Spec(\wt{B})$ is $af^n$. In particular, $a = gf^{-n}$
    and hence $\wt{\pi}_0 = gf^{-n}s^{-1}$. Letting $L = Z(g)$, it follows that a
    local defining equation of $L$ in $\Spec(\wt{B}_\Sigma)$ is $s\wt{\pi}_0 = 0$.
    The lemma now follows if we let $u$ be the image of $s$ under the canonical
    surjection $\wt{B} \surj B_\Sigma$.
\end{proof}

\begin{cor}\label{cor:Type-2-10-0}
  If $L$ is as obtained in \lemref{lem:Type-2-10} and $i \in J_r$, then
  $L \bigcap \sX$ is locally defined at $x_i$ by $u_i\pi_i$ for some
    $u_i \in \sO^\times_{\sX,x_i}$. Furthermore, there is a scheme-theoretic inclusion
    $Y \subset L \bigcap \sX$.
  \end{cor}
  \begin{proof}
    The first claim is clear. To prove the second claim, one first notes easily that
    each $Y_i$ is a closed subscheme of $L \bigcap \sX$ since $Y_i$ is integral.
    From this it follows that if $y \in Y_{(0)}$ is any point at which 
    $L \bigcap \sX$ is locally defined by $f$ and $Y$ is defined at $y$ by
    $t_{i_1} \cdots t_{i_q}$ (where $q \le r$), then $t_{i_j} \mid f$ for each $i_j$.
    But this implies that $t_{i_1} \cdots t_{i_q} \mid f$ as $\sO_{\sX,y}$ is a
    ufd. This proves the second claim.
\end{proof}

Before we state the next lemma, we recall from \cite[Thm.~10.5.1]{CTS} and
  \cite[Prop.~2.1]{Conrad-adic} that for every separated and finite type
  $k$-scheme $X$, the set $X(k)$ is endowed with the adic topology which is induced
  by the adic topology of $k$, given by its absolute value. This assignment defines a
  functor from the category of finite type and separated $k$-schemes with regular
  morphisms to the category of Hausdorff topological spaces with continuous maps.
  This functor takes open (resp. closed) immersions to open (resp. closed)
  inclusions of topological spaces. If $k$ is locally compact, then $X(k)$ is
  locally compact and is compact if $X$ is projective over $k$.

If we let $|a| = 2^{-v(a)}$ (where $v$ is the normalized
    valuation of $k$) denote the absolute value on $k$, then
    the adic topology of $\A^q_k(k) = k^q$ is induced by latter's normed linear space
    structure with respect to the sup-norm  given by
    $\left\|{\un{a}}\right\| = {\underset{1 \le i \le q}{\rm max}} \{|a_i|\}$.
     We let 
     \[
     D^{q}_k(1) = \{\un{a} \in k^q|
    \left\|{\un{a}}\right\| \le 1\} \text{ and } U^{q}_k(1) =  \{\un{a} \in k^q|
    \left\|{\un{a}}\right\| < 1\}
    \]
    denote the closed and open unit disks,
    respectively, in $k^q$ in the adic topology. Let
    \[
    S^{q}_k(1) = \{\un{a} \in k^q| \left\|{\un{a}}\right\| = 1\}
    \]
    denote the
    unit sphere. We note that $D^q_k(1), U^q_k(1)$ (hence $S^q_k(1)$) are
    both open and closed in $k^q$ with respect to the adic topology,
    and they are all contained in $R^q$.

\begin{lem}\label{lem:Type-2-9}
  Let $x \in Y^o$ be any closed point and let $Y$ be locally defined at $x$
  by $\pi_0 \in \sO_{\sX,x}$. Let $n \in \N$. Then
      there exists $m \ge n$ and admissible
    hypersurfaces $H_1, H_2 \in |\sO_{\P^N_R}(m)|(R)$ containing $x$
    such that letting  $\sX_i = \sX \bigcap H_i$,
    one has that if $\sX_1$ is locally defined by $t \in \sO_{\sX,x}$ at $x$,
    then $\sX_2$ is locally defined by $t+ u\pi_0$ at $x$ for some
    $u \in \sO^\times_{\sX,x}$.
\end{lem}
\begin{proof}
  We may assume that $x \in Y^o_1$ and let $\Sigma = \{x=x_1,x_2,\ldots, x_r\}$.
  We choose a finite set of closed points $\Sigma' \subset Y$
 such that $\Sigma' \bigcap Y_i  \neq \emptyset$ for every $i$
 and $\Sigma \bigcap \Sigma' = \emptyset$ (this is always possible).

 It is shown in the first paragraph of the proof of
  \thmref{thm:Bertini-dvr} that
    there exists $d \gg 0$ such that for every $m \ge d$, we can find a
    Zariski dense open subscheme $U_{m,k}$ of $|\sO_{\P^N_k}(m)|$ such
    that for every $\un{a} := (a_0, \cdots , a_{m'-1}) \in U_{m,k}(k)$
    (where $m'$ is the rank of the free $R$-module $H^0_\zar(\P^N_R, \sO_{\P^N_R}(m))$),
    one has that $H_{\un{a}} \bigcap X$ is a smooth and connected $k$-scheme and
    the restriction map
    $H^0_\zar(X, \sO_X) \to H^0_\zar(H_{\un{a}} \bigcap X, \sO_{H_{\un{a}} \bigcap X})$
    is an isomorphism, where $H_{\un{a}}$ is the hypersurface corresponding to the
    point $\un{a}$. By \thmref{thm:Bertini-dvr} (where we take $T = \Sigma', \
    Z_1 = \Spec(k(x))$ and $Z_2 = \emptyset$), we can take $d$ large enough so
    that for all $m \ge d$, there exists an admissible hypersurface
    $T \in |\sO_{\P^N_R}(m)|(R)$ which contains $x$. We now fix 
    $n \in \N$ and choose $m \ge {\rm max}\{n,d\}$ for which the conclusion of
    \lemref{lem:Type-2-10} is also satisfied.

We let 
\[
\phi \colon \A^{m'}_k \setminus \{0\} \to \P^{m'-1}_k
\]
    denote the canonical quotient map by the $\G_m$-action (through scalar
    multiplication). Then $U'_{m,k}:= \phi^{-1}(U_{m,k})$  is Zariski dense in
    $\A^{m'}_k \setminus \{0\}$. It follows from \cite[Thm.~10.5.1]{CTS}
    (this requires semi-stability of $\sX$)
    that $U'_{m,k}(k)$ is dense open in $k^{m'} \setminus \{0\}$
    (hence in $k^{m'}$) in the adic topology. 

We let $L$ be as in \lemref{lem:Type-2-10} and choose
    an admissible hypersurface $T  \in |\sO_{\P^N_R}(m)|(R)$ containing $x$. Write
    $L = Z(F)$ and  $T = Z(G)$ with
    $F = \stackrel{m'-1}{\underset{i =0}\sum} c_i F_i$
and $G = \stackrel{m'-1}{\underset{i =0}\sum} d_i F_i$, where
$c_i, d_i \in R$ and not all $d_i$ die in $\F$. Here, $\{F_i\}_{i=0}^{m'-1}$ is a 
basis of the free $R$ module $H^0_\zar(\P^N_R, \sO_{\P^N_R}(m))$. For
$\un{a} = (a_0, \ldots , a_{m'-1}) \in R^{m'}$, we let
${\un{a}}_\F$ be the canonical image of $\un{a}$ in $\F^{m'}$.
We shall prove the lemma by first assuming that
$\un{c} = (c_0, \ldots , c_{m'-1}) \in S^{m'}_k(1)$. Equivalently, $c_i \in R^\times$ for
some $i$.

We let $\psi_{\un{c}} \colon D^{m'}_k(1) \to D^{m'}_k(1)$ denote the map
$\psi_{\un{c}}(\un{a}) = \un{a} + \un{c} = (a_0 + c_0, \ldots , a_{m'-1} + c_{m'-1})$.
Using the fact that $| \cdot |$ is an ultrametric on $k$, one easily checks that
$\psi_{\un{c}}$ is well-defined and is a homeomorphism whose inverse is
$\psi_{-\un{c}}$. If we let $V_{m,k} = U'_{m,k} \bigcap D^{m'}_k(1)$, then
$V_{m,k}$ is an adic dense open subset of $D^{m'}_k(1)$. In particular,
$\psi_{\un{c}}(V_{m,k})$ is also dense open in $D^{m'}_k(1)$.
As $S^{m'}_k(1) \subset D^{m'}_k(1)$ open, it follows that
$V_{m,k} \bigcap \psi_{\un{c}}(V_{m,k}) \bigcap S^{m'}_k(1)$ is dense open in
$S^{m'}_k(1)$.

We let ${\rm sp} \colon D^{m'}_k(1) =  R^{m'} = \A^{m'}_R(R) \to \A^{m'}_\F(\F) =
\F^{m'}$ denote the restriction map
which takes an $R$-point of $\A^{m'}_k$ to $\A^{m'}_\F$ by reducing the
coordinates modulo $(\pi)$. As it is easy to see that this map is
continuous with respect to the discrete topology of $\F^{m'}$
(cf. \cite[Rem.~2.4]{Ochiai-Shimomoto}), ${\rm sp}^{-1}({\un{d}}_\F)$ is open
in $R^{m'}$ whence in $k^{m'}$, where
$\un{d} = (d_0, \ldots , d_{m'-1}) \in R^{m'}$ is the point
corresponding to the coefficients of $G$. It follows that
\[
Q={\rm sp}^{-1}({\un{d}}_\F) \bigcap \psi_{\un{c}}(V_{m,k}) \bigcap \left(V_{m,k} \bigcap
S^{m'}_k(1)\right) \subset R^{m'}
\]
is open dense in ${\rm sp}^{-1}({\un{d}}_\F) \bigcap S^{m'}_k(1)$. Since
$\un{d} \in {\rm sp}^{-1}({\un{d}}_\F) \bigcap S^{m'}_k(1)$ by our assumption, we
conclude that $Q$ is non-empty. We choose a point $\un{a} \in Q$. We then have the
following.

\begin{enumerate}
  \item $H_{\un{a}} \in |\sO_{\P^N_R}(m)|(R)$.
  \item
    $H_{\un{a}} \bigcap X$ is a smooth $k$-scheme which is connected and the
    map $H^0_\zar(X, \sO_X) \to
    H^0_\zar(H_{\un{a}} \bigcap X, \sO_{H_{\un{a}} \bigcap X})$ is an isomorphism
    (since $\un{a} \in V_{m,k}$).
  \item
    The closure of $H_{\un{a}} \bigcap X$ in $\sX$ is $H_{\un{a}} \bigcap \sX$ (since
    $\un{a} \in S^{m'}_k(1)$).
  \item
    $H_{\un{a}} \bigcap \sX_s = T \bigcap \sX_s$
    (since ${\rm sp}(\un{a}) = \un{d}_\F$).
\end{enumerate}

The property (4) implies that $x \in H_{\un{a}}$ and the
closed subscheme $H_{\un{a}} \bigcap Y_i$
(resp. $H_{\un{a}} \bigcap Y$) of $\sX$ coincides with $T \bigcap Y_i$
(resp. $T \bigcap Y$) for each $i \in J_r$. In particular, the admissibility of $T$
implies that
$H_{\un{a}} \bigcap Y$ is an snc scheme (cf. \S~\ref{sec:NCS}) and
each $H_{\un{a}} \bigcap Y_i$ is a regular closed subscheme of $\sX$ of codimension
two. Given this, an argument identical to that of the proof of item (1) of
\thmref{thm:Bertini-dvr} implies that $H_{\un{a}} \bigcap \sX$ must be regular.
Furthermore, \lemref{lem:SNC-0}(1) implies that $H_{\un{a}} \bigcap Y =
(H_{\un{a}} \bigcap \sX_s)_\red$ is an snc divisor on $H_{\un{a}} \bigcap \sX$.
The property (2) above implies that $H_{\un{a}} \bigcap X$ is geometrically integral
over $k$ if $X$ is so. It follows that $H_{\un{a}}$ is admissible and
$x \in H_{\un{a}}$.

Since $\un{a} \in \psi_{\un{c}}(V_{m,k})$, we can write $\un{a} = \un{a}' + \un{c}$
  for some $\un{a}' \in V_{m,k} \subset D^{m'}_k(1)$. We let $G' = \sum_i a'_i F_i$ so
  that $H_{\un{a}'} = Z(G')$. Then $H_{\un{a}'} \bigcap X$ is a smooth $k$-scheme. 
  We claim that $\un{a}' \in S^{m'}_k(1)$. To prove this, suppose $\un{a}' \notin
  S^{m'}_k(1)$.
The latter condition implies that $\pi \mid a'_i$ for all $i$.
  In particular, $\un{a}' + \un{c} \in S^{m'}_k(1)$ (because $\un{c} \in S^{m'}_k(1)$)
  and ${\rm sp}(\un{a}' + \un{c}) = {\rm sp}(\un{c}) = \un{c}_\F$.
On the other hand, $\un{a} = \un{a}' + \un{c}$ implies that
  ${\rm sp}(\un{a}' + \un{c}) = {\rm sp}(\un{a}) = \un{d}_\F$. But this a
contradiction since we can not have ${\rm sp}(\un{c}) = {\rm sp}(\un{d})$
because $L \bigcap \sX$ contains $Y$ while $T \bigcap Y$ has codimension one
in $T \bigcap \sX$ at $x$.
   This proves the claim, from which we conclude that the
  closure of $H_{\un{a}'} \bigcap X$ in $\sX$ is $H_{\un{a}'} \bigcap \sX$.

We next claim that $H_{\un{a}'} \bigcap Y = H_{\un{a}} \bigcap Y$. To prove this, we
  choose a closed point $y \in Y$ and let $\alpha = \pi_{i_1} \cdots \pi_{i_s} \in
  \sO_{\sX,y}$ define $Y$ at $y$, where $1 \le i_1 < \cdots  < i_s \le r$.
 We let $\beta, \beta' \in \sO_{\sX,y}$ be elements which define
  $H_{\un{a}} \bigcap \sX$ and $H_{\un{a}'} \bigcap \sX$,
  respectively, at $y$. We can then write $\beta = \beta' + \gamma$, where $\gamma
  \in \sO_{\sX,y}$ defines $L$ at $y$. Letting $W = H_{\un{a}} \bigcap Y$ and
  $W' = H_{\un{a}'} \bigcap Y$, we then get
  \begin{equation}\label{eqn:Type-2-9-0}
    \sO_{W,y} = \frac{\sO_{\sX,y}}{(\alpha, \beta'+ \gamma)}; \ \
    \sO_{W',y} = \frac{\sO_{\sX,y}}{(\alpha, \beta')}.
    \end{equation}
On the other hand, \corref{cor:Type-2-10-0} implies that $Y \subset L \bigcap \sX$,
  and this in turn implies that $\alpha \mid
  \gamma$. In particular, $(\alpha, \beta'+ \gamma) = (\alpha, \beta')$. 
  This implies the claim, as one easily checks.
  
From the above claims and property (4) of $H_{\un{a}}$, we get the following. \\
   \hspace*{.5cm} (1') \ $H_{\un{a'}}  \in |\sO_{\P^N_R}(m)|(R)$. \\
   \hspace*{.5cm} (2') \ $H_{\un{a'}} \bigcap X$ is a smooth $k$-scheme which is
   connected and the map $H^0_\zar(X, \sO_X) \to
    H^0_\zar(H_{\un{a'}} \bigcap X, \sO_{H_{\un{a'}} \bigcap X})$ is an isomorphism
    (since $\un{a'} \in V_{m,k}$). \\
  \hspace*{.5cm} (3') \ The closure of $H_{\un{a'}} \bigcap X$ in $\sX$ is
  $H_{\un{a'}} \bigcap \sX$. \\
  \hspace*{.5cm} (4') \ $H_{\un{a'}} \bigcap Y = T \bigcap Y$.
  In particular, $(H_{\un{a'}} \bigcap \sX)_\red =  H_{\un{a'}} \bigcap Y$. \\
These properties imply (as we showed for $H_{\un{a}}$)
  that $H_{\un{a}'}$ is admissible and $x \in H_{\un{a'}}$. 
 Moreover, it follows from the choice of $L$ (cf. \lemref{lem:Type-2-10}) that if
  $t \in \sO_{\sX,x}$ defines
  $H_{\un{a}'} \bigcap \sX$ at $x$, then $t + u\pi_1$ defines $H_{\un{a}} \bigcap \sX$
  at $x$ for some $u \in \sO^\times_{\sX,x}$. Letting $H_1 = H_{\un{a}'}$ and
  $H_2 = H_{\un{a}}$, we see that the assertion of the lemma is satisfied.

It remains to prove the lemma when $\un{c} \notin S^{m'}_k(1)$. Equivalently,
$\pi \mid c_i$ for all $i$. In this case, we let
$\psi_{\un{c}} \colon S^{m'}_k(1) \to S^{m'}_k(1)$ be the map which takes $\un{a}$ to
$\un{a} + \un{c}$. One checks as before that this is well-defined
and a topological homeomorphism. Letting $V'_{m,k} = U'_{m,k} \bigcap S^{m'}_k(1)$ and
repeating the argument of the previous case, we get again that
$Q' := {\rm sp}^{-1}(\un{d}_\F) \bigcap V'_{m,k} \bigcap \psi_{\un{c}}(V'_{m,k})$ is open
dense in ${\rm sp}^{-1}(\un{d}_\F) \bigcap S^{m'}_k(1)$. In particular, it is
non-empty. We now let  $\un{a} \in Q'$ be any point and write
$\un{a} = \un{a}' + \un{c}$ with $\un{a}' \in V'_{m,k}$.  Noting that
$a' \in S^{m'}_k(1)$, an argument identical to that in the previous case shows that
$H_1 = H_{\un{a}'}$ and $H_2 = H_{\un{a}}$ satisfy the assertion of the lemma.
This completes the proof.
\end{proof}

Continuing with the notation of \S~\ref{sec:Esp}, let $\chi \in H^{q+1}_{p^m}(X)$.
Let $i \in J_r$ and let $D = \stackrel{r}{\underset{j=1}\sum} m_j Y_j$ be a
$(\chi,i)$-admissible divisor. Assume that the conditions {$(\star)$} hold.
Let $\beta$ denote the element $\phi^*_1(\wt{\chi})$ 
of $H^0_\zar(Y_i, \Omega^{q+1}_{\sX}(\log Y_o)(D) \otimes_{\sX} \sO_{Y_i})$.

\begin{cor}\label{cor:Type-2-9-key}
  Let $x \in Y^o_i$ be a closed point and let $n \in \N$. Let
  $m \ge n$ and $H_1, H_2 \in |\sO_{\P^N_R}(m)|(R)$ be admissible hypersurfaces
  containing $x$ as in \lemref{lem:Type-2-9}.
Assume that the image of $\gamma^*_i(\beta)$ is not zero in
 $\Omega^{q+1}_{Y_i}(\log F_i)(D_i) \otimes_{\sO_{Y_i}} k(x)$.
Then there exists $j \in J_2$ such that letting $\sX_j = \sX \bigcap H_j$,
one has $n_i -1 \le \Sw_{Y'_i}(\chi|_{X_j}) \le n_i$.
\end{cor}
\begin{proof}
Suppose that $\sX_1$ is locally defined at $x$ by an element $t \in \sO_{\sX,x}$.
Then \lemref{lem:Type-2-9} says that $\sX_2$ is locally defined at $x$ by 
  $t + u \pi_i$ for some $u \in \sO^\times_{\sX,x}$. Since $\sX_1$ is
  admissible, we can write the maximal ideal of $\fm$ of $\sO_{\sX,x}$ as
  $\fm = (\pi_i, t, t_3, \ldots , t_{d_\sX})$. As argued in the proof of
  \propref{prop:Type-2-key-0}, this implies using \lemref{lem:Type-2-key} 
that $n_i -1 \le \Sw_{Y'_i}(\chi|_{X_j}) \le n_i$ for some $j \in J_2$.
\end{proof}

The following is case III of the specialization of Swan conductor.

\begin{cor}\label{cor:Type-2-11}
  Let $i \in J_r$ and assume that the conditions {$(\star)$} hold. Assume
  also that $\sX$ is semistable and there is a closed embedding $\sX \inj \P^N_R$
  of $R$-schemes.
  Then there exists a dense open $U_i \subset Y_i$ such that for every closed
  point $x \in U_i$ and every $n \in \N$, there exists $m \ge n$ and an admissible 
  hypersurface section $\sX'$ of $\sX$ of degree $m$ containing $x$ such that
 $n_i - 1 \le \Sw_{Y'_i}(\chi|_{X'}) \le n_i$.
\end{cor}
\begin{proof}
  This follows from \corref{cor:Type-2-9-key} using an argument similar to that
  of \corref{cor:Type-2-8}.
  \end{proof}

\subsection{Specialization of Swan conductor: case~IV}\label{sec:SP-3.1}
We now consider the specialization of Swan conductor in the last remaining
case, namely, when part (4) of conditions ($\star$) is not satisfied.
In this subsection, we shall assume the
following, and refer to them  collectively as  \textbf{conditions {($\star'$)}}.
\begin{enumerate}
  \item
      $\chi \in \Fil_D H^{q+1}_{p^m}(X)$ is of type II at $Y_i$.
  \item
    $q = d_\sX -2 \ge 0$.
  \item
       $p = 2 = n_i$.
  \end{enumerate}

\vskip.2cm

 Assume that the conditions {$(\star')$} hold. Fix $i \in J_r$ and a closed
  point $x \in Y^o_i$. Let $A = \sO_{\sX,x}$ and let
  $\fm = (\pi_i, t_2, \ldots , t_{d_\sX})$ be the maximal ideal of $A$, where
  $\pi_i$ defines $Y_i$ locally at $x$. Let $t^\dagger_2 = t_2 - b\pi_i$, where
  $b \in A^\times$. Let $A_1 = {A}/{(t_2)}, \ A_2 = {A}/{(t^\dagger_2)}$ and
  $\ov{A} = {A}/{(\pi_i)}$. We let $B = A[\pi^{-1}_i], \ B_1 = {A_1}[\pi^{-1}_i]$ and
  $B_2 = {A_2}[\pi^{-1}_i]$. We let $\sX_x=\Spec(A)$, $X_x=\Spec(B)$ and
  $\chi_x = \chi|_{X_x}$.

 Let $D = \stackrel{r}{\underset{j =1}\sum} m_j Y_j$ be a $(\chi, i)$-admissible
  divisor. Let $\ov{\chi}$ denote the image of $\chi$ in
  $\frac{\Fil_D H^{q+1}_{p^m}(X)}{\Fil_{D^{\<i\>}_1} H^{q+1}_{p^m}(X)}$ and let
  $\wt{\chi} = \Rsw^{m, q+1}_{\sX|(D,D^{\<i\>}_1)}(\ov{\chi})
  \in  H^0_\zar\left(\sX, \Omega^{q+1}_\sX(\log Y)(D)\otimes_{\sO_\sX} \sO_{Y_i}\right)$.
  Since $\chi$ has type II at $Y_i$, we know that
  $\wt{\chi} \in H^0_\zar(Y_i, \Omega^{q+1}_{Y_i}(\log F_i)(D_i))$.
Let $R_* \colon \Fil_D H^{q+1}_{p^m}(X) \to \Fil_{D/2} H^{q+1}_{p^{m-1}}(X)$ be the 
projection map induced by the map $R$ in \thmref{thm:V-R exact} after
taking hypercohomology and using \thmref{thm:H^1-fil}. We have a similar map for
$D^{\<i\>}_1$. Taking the quotients, we get
\begin{equation}\label{eqn:Rest**}
 R_* \colon \frac{\Fil_D H^{q+1}_{p^m}(X)}{\Fil_{D^{\<i\>}_1} H^{q+1}_{p^m}(X)} \to
\frac{\Fil_{D/2} H^{q+1}_{p^{m-1}}(X)}{\Fil_{{D^{\<i\>}_1}/2} H^{q+1}_{p^{m-1}}(X)}.
\end{equation}

\begin{lem}\label{lem:Type-12-key}
 Assume that $R_*(\ov{\chi}) = 0$ but the image of $\wt{\chi}$ in
 $\Omega^{q+1}_{Y_i}(\log F_i)(D_i) \otimes_{\sO_{Y_i}} k(x)$ is not zero.
 Then there exists $j \in J_2$ such the image of $\ov{\chi}$ is not zero in
$\frac{\Fil_2 H^{q+1}_{p^m}(B_j)}{\Fil_0 H^{q+1}_{p^m}(B_j)}$.
\end{lem}
\begin{proof}
We can assume without loss of generality that $i = 1$. We let $\pi'_1 = (-b)\pi_1$.
We consider the commutative diagram
  \begin{equation}\label{eqn:Type-2-12-1}
    \xymatrix@C1.2pc{
      0 \ar[r] & \Fil_i H^{q+1}_{p}(B) \ar[r]^-{V^{m-1}} \ar[d] &
      \Fil_i H^{q+1}_{p^m}(B) \ar[r]^-{R_*} \ar[d] &
      \Fil_0 H^{q+1}_{p^{m-1}}(B) \ar[r] \ar[d] & 0 \\
       0 \ar[r] & \Fil_2 H^{q+1}_{p}(B) \ar[r]^-{V^{m-1}} &
      \Fil_2 H^{q+1}_{p^m}(B) \ar[r]^-{R_*}  &
      \Fil_1 H^{q+1}_{p^{m-1}}(B) \ar[r] & 0}
  \end{equation}
  for $i \in \{0,1\}$, obtained by applying hypercohomology to the exact triangle of
  \thmref{thm:V-R exact} and using
  \thmref{thm:H^1-fil}, where the vertical arrows are the canonical inclusions.
  We consider a similar diagram for $B_j$ for $j \in J_2$.

The rows of the above diagram are left exact by \thmref{thm:H^1-fil} and
  \cite[Lem.~7.1]{KM-1} using the observation that
  $V \colon \Fil_l H^{q+1}_{p^s}(B) \to \Fil_l H^{q+1}_{p^{s+1}}(B)$ is the
  canonical inclusion of $\Fil_l H^{q+1}_{p^{s+1}}(B)[p^s]$ into
  $\Fil_l H^{q+1}_{p^{s+1}}(B)$ for every $l \in \N_0, s \in \N$
  (cf. diagram at the end of \cite[\S~(1.3)]{Kato-89}). On the other hand, the maps
  $R_*$ in \eqref{eqn:Type-2-12-1} are surjective, since
  $\H^t_\et(X, W_m\sF^{2,\bullet}_\n)=0$ for any affine scheme $X$ and for all
  $t \ge 2$ (see \cite[Lem.~7.2(2)]{KM-1}, whose proof works for any affine scheme).
  
Taking the quotients of vertical inclusions for $B$ and comparing them with
  the corresponding quotients for $B_j$,
  we get a commutative diagram of exact sequences
  \begin{equation}\label{eqn:Type-2-12-2}
      \xymatrix@C1.7pc{
        0 \to \frac{\Fil_2 H^{q+1}_{p}(B)}{\Fil_i H^{q+1}_{p}(B)}
        \ar[r]^-{V^{m-1}} \ar[d] &
       \frac{\Fil_2 H^{q+1}_{p^m}(B)}{\Fil_i H^{q+1}_{p^m}(B)} 
\ar[r]^-{R_*} \ar[d] & \frac{\Fil_1 H^{q+1}_{p^{m-1}}(B)}{\Fil_0 H^{q+1}_{p^{m-1}}(B)} 
\ar[r] \ar[d] & 0 \\
0 \to \frac{\Fil_2 H^{q+1}_{p}(B_j)}{\Fil_i H^{q+1}_{p}(B_j)}
        \ar[r]^-{V^{m-1}} &
       \frac{\Fil_2 H^{q+1}_{p^m}(B_j)}{\Fil_i H^{q+1}_{p^m}(B_j)} 
\ar[r]^-{R_*} & \frac{\Fil_1 H^{q+1}_{p^{m-1}}(B_j)}{\Fil_0 H^{q+1}_{p^{m-1}}(B_j)} 
\ar[r] & 0}
      \end{equation}
  for $i \in J^0_1$ and $j \in J_2$.

Since $R_*(\ov{\chi}_x) =0$ in the top row, we can find an element
  $\chi_1 \in \Fil_2 H^{q+1}_{p}(B)$ such that letting
$\ov{\chi}_1$ denote its image in
$\frac{\Fil_2 H^{q+1}_{p}(B)}{\Fil_0 H^{q+1}_{p}(B)}$,
we have $V^{m-1}(\ov{\chi}_1) = \ov{\chi}_x$ in
$\frac{\Fil_2 H^{q+1}_{p^m}(B)}{\Fil_0 H^{q+1}_{p^m}(B)}$. Moreover, the commutative
diagram 
  \begin{equation}
      \xymatrix@C1.7pc{
        \frac{\Fil_2 H^{q+1}_{p}(B)}{\Fil_0 H^{q+1}_{p}(B)}
        \ar[r]^-{V^{m-1}}\ar@{->>}[d] &
        \frac{\Fil_2 H^{q+1}_{p^m}(B)}{\Fil_0 H^{q+1}_{p^m}(B)} \ar@{->>}[d]\\
        \frac{\Fil_2 H^{q+1}_{p}(B)}{\Fil_1 H^{q+1}_{p}(B)}\ar[r]^-{V^{m-1}} &
        \frac{\Fil_2 H^{q+1}_{p^m}(B)}{\Fil_1 H^{q+1}_{p^m}(B)}
      }
  \end{equation}
  implies that the equality $V^{m-1}(\ov{\chi}_1) = \ov{\chi}_x$ also holds in
  $\frac{\Fil_2 H^{q+1}_{p^m}(B)}{\Fil_1 H^{q+1}_{p^m}(B)}$, where vertical maps are
  canonical quotient maps, and the bar denotes the image of an element under the
  canonical quotient maps. In particular,
  $\Rsw^{1,q}_{A|(2, 1)}(\ov{\chi}_1) = \Rsw^{m,q}_{A|(2, 1)}(\ov{\chi}_x)$ by
 \thmref{thm:RSW-gen} (which says that the refined Swan conductor does not
 see the value of $m$).

 If we let $\wt{\chi}_1 = \Rsw^{1,q}_{A|(2, 1)}(\ov{\chi}_1)$, it follows that
$\wt{\chi}_1 = \wt{\chi}|_{\sX_x}$ so that
$\wt{\chi}_1 \in \Omega^{q+1}_{\ov A}$ does not
  die when we pass to $\Omega^{q+1}_{\ov{A}} \otimes_{\ov{A}} k(x)$. We let
  $\chi^j_1$ denote the image of $\chi_1$ under the canonical restriction
  $\Fil_2 H^{q+1}_{p}(B) \to \Fil_2 H^{q+1}_{p}(B_j)$ and let
 $\chi^j_x$ be the image of $\chi_x$ in $\Fil_2 H^{q+1}_{p^m}(B_j)$.

To prove the lemma, it is enough to show by using ~\eqref{eqn:Type-2-12-2} for
  $i =0$ that $\chi^j_1$ does not die in
  $\frac{\Fil_2 H^{q+1}_{p}(B_j)}{\Fil_0 H^{q+1}_{p}(B_j)}$ for some $j \in J_2$.
Equivalently, the Swan conductor of $\chi^j_1$ at
$\Spec({A_j}/{(\pi'_1)})$ is either one or two for some $j \in J_2$. 
  To prove the latter claim, we can replace $A$ by $A_{\fp}$ and $B$ by
  $B_{\fp} = A_\fp[1/{\pi'_1}]$, where
  $\fp = (\pi'_1, t_2) = (\pi'_1, t^\dagger_2) \subset A$.
  We let $L = Q(A/{\fp})$ denote the residue field of $A_\fp$ and
  let $[L: L^p] = p^{c-1}$. It follows then that
  $c-1 = d_\sX -2 + p_\F = q + p_\F \ge q$. We would like now to apply
  \lemref{lem:spl-loc-1} but we can't because it is not guaranteed that
  $H^c_p(L) \neq 0$. We remedy this by using a trick due to Kato, as follows.

By \cite[Lem.~7.7]{Kato-89} and \cite[Lem.~1]{Kato80-3}
  (see the proof of \cite[Thm.~7.1]{Kato-89}), there is a complete regular
  local ring $A^\dagger_\fp$ having residue field $L^\dagger$ together with a flat
  local ring homomorphism $\psi \colon A_\fp \to A^\dagger_\fp$ such that the
  following hold.
  \begin{enumerate}
  \item
    $A^\dagger_\fp \otimes_{A_\fp} L \to L^\dagger$ is an isomorphism.
  \item
    Any $p$-basis of $L$ is also a $p$-basis of $L^\dagger$.
  \item
    $H^{c}_p(L^\dagger) \neq 0$.
  \end{enumerate}
  
The above properties imply that $A^\dagger_\fp$ is $F$-finite (cf.
  \cite[Prop.~2.2]{KM-1}) and $\fm^{\dagger} = (\pi'_1, t_2)A^\dagger_\fp$ is the maximal
  ideal of $A^\dagger_\fp$. If we let $A' \in \{A_1, A_2\}, \ B' \in \{B_1, B_2\}$
  and
  \[
  C_\fp = {A_\fp}/{(\pi'_1)}, \ {C}^\dagger_\fp = {A^\dagger_\fp}/{(\pi'_1)} \cong
 C_\fp \otimes_{A_\fp} A^\dagger_\fp, \
  A'^\dagger_\fp  = A' \otimes_{A_\fp} A^\dagger_\fp, \
  B'^\dagger_\fp = B'_{\fp} \otimes_{A_\fp} A^\dagger_\fp,
  \]
  then, for every $i \ge 0$, $\Omega^i_{C^\dagger_\fp}$ is a free $C^\dagger_\fp$-module
  with the same basis as that of the free $C_\fp$-module $\Omega^i_{C_\fp}$.
In particular,
    $\Omega^i_{C_\fp} \to \Omega^i_{C^\dagger_\fp}$ is injective for every $i \ge 0$.
  We let $B^\dagger_\fp = A^\dagger_\fp[1/{\pi'_1}]$.

We now look at the diagram
  
\begin{equation}\label{eqn:Type-2-12-3}
      \xymatrix@C1.2pc{
\frac{\Fil_2 H^{q+1}_{p}(B)}{\Fil_1 H^{q+1}_{p}(B)}
\ar[r] \ar[dr] & \frac{\Fil_2 H^{q+1}_{p}(B_\fp)}{\Fil_1 H^{q+1}_{p}(B_\fp)}
\ar[r]^-{\psi^*} \ar[d] &
\frac{\Fil_2 H^{q+1}_{p}(B^\dagger_\fp)}{\Fil_1 H^{q+1}_{p}(B^\dagger_\fp)}
  \ar[d] \\
&  \Omega^{q+1}_{C_\fp} \bigoplus \Omega^q_{C_\fp} \ar[r]^-{\psi^*} &
  \Omega^{q+1}_{C^\dagger_\fp} \bigoplus \Omega^q_{C^\dagger_\fp},}
  \end{equation}
where the vertical arrows are the refined Swan conductors,
the top horizontal arrow on the left is the canonical map induced by the
localization and the diagonal arrow is the composition.
This diagram is commutative by \thmref{thm:RSW-gen}.

Since $\wt{\chi}_1$ does not die in
$\Omega^{q+1}_{\ov{A}} \otimes_{\ov{A}} k(x)$, it is not zero in $\Omega^{q+1}_{\ov{A}}$.
In particular, it does not die in $\Omega^{q+1}_{\ov{A}_{\fp}}=\Omega^{q+1}_{C_\fp}$. It
follows that $\ov{\chi}_1 \in \frac{\Fil_2 H^{q+1}_{p}(B)}{\Fil_1 H^{q+1}_{p}(B)}$
does not die under the diagonal arrow on the left in ~\eqref{eqn:Type-2-12-3}.
We let $\chi^\dagger_1$ denote the image of $\chi_1$ under the map
$\Fil_2 H^{q+1}_{p}(B) \to \Fil_2 H^{q+1}_{p}(B^\dagger_\fp)$ and let
$\ov{\chi}^\dagger_1$ be the image of $\ov{\chi}_1$ under top composite map
in ~\eqref{eqn:Type-2-12-3}. We let $\omega^\dagger :=
\psi^* \circ \Rsw^{1,q}_{A|(2,1)}(\ov{\chi}_1)
= \Rsw^{1,q}_{A^\dagger_\fp|(2,1)}(\ov{\chi}^\dagger_1)$ denote the refined Swan conductor
of $\chi^\dagger_1$ in $\Omega^{q+1}_{C^\dagger_\fp}$.

{\bf{Claim:}}  $\omega^\dagger$ does not die in
$\Omega^{q+1}_{C^\dagger_\fp} \otimes_{C_\fp} L^\dagger$.

To prove this claim, we observe that there are
canonical quotient maps $\Omega^{q+1}_{\ov{A}} \surj
\Omega^{q+1}_{\ov{A}} \otimes_{\ov{A}} A/{\fp} \surj
\Omega^{q+1}_{\ov{A}} \otimes_{\ov{A}} k(x)$
and $\wt{\chi}_1$ does not die under the composite quotient map.
It follows that it does not die in $\Omega^{q+1}_{\ov{A}} \otimes_{\ov{A}} A/{\fp}$.
Since the latter is a free $A/{\fp}$-module, it injects inside
$\Omega^{q+1}_{\ov{A}} \otimes_{\ov{A}} L$. But the latter is easily seen to be
isomorphic to $\Omega^{q+1}_{C_\fp} \otimes_{C_\fp} L$ as an $L$-vector space. It
follows that $\wt{\chi}_1$ does not die in $\Omega^{q+1}_{C_\fp} \otimes_{C_\fp} L$.
In particular, it does not die in
$\left(\Omega^{q+1}_{C_\fp} \otimes_{C_\fp} L\right) \otimes_{L} L^\dagger$.
The latter is easily seen to be isomorphic to
$\left(\Omega^{q+1}_{C_\fp} \otimes_{C_\fp} C^\dagger_\fp\right) \otimes_{C^\dagger_\fp}
L^\dagger$.

Since the canonical map
$\Omega^{q+1}_{C_\fp} \otimes_{C_\fp} C^\dagger_\fp \to \Omega^{q+1}_{C^\dagger_\fp}$ is
an isomorphism of free $C^\dagger_\fp$-modules as observed before,
we get an isomorphism
$\left(\Omega^{q+1}_{C_\fp} \otimes_{C_\fp} C^\dagger_\fp\right) \otimes_{C^\dagger_\fp}
L^\dagger \cong \Omega^{q+1}_{C^\dagger_\fp} \otimes_{C^\dagger_\fp} L^\dagger$.
It follows that $\wt{\chi}_1$ does not die in
$\Omega^{q+1}_{C^\dagger_\fp} \otimes_{C^\dagger_\fp} L^\dagger$.
We now use the commutative diagram
\begin{equation}\label{eqn:Type-2-12-4}
      \xymatrix@C1.2pc{
\Omega^{q+1}_{A_1} \ar[r]^-{\psi^*} \ar[d] & \Omega^{q+1}_{C^\dagger_\fp} \ar[d] \\ 
\Omega^{q+1}_{C_\fp} \otimes_{C_\fp} L \ar[r]^-{\psi^*} &
\Omega^{q+1}_{C^\dagger_\fp} \otimes_{C^\dagger_\fp} L^\dagger}
\end{equation}
to conclude that $\omega^\dagger$ does not die in
$\Omega^{q+1}_{C^\dagger_\fp} \otimes_{C_\fp} L^\dagger$.
This proves the claim.

Using the claim, diagram ~\eqref{eqn:Type-2-12-3} and the property (3) above, we
see that  $\chi^\dagger_1$ is an element of $\Fil_2 H^{q+1}_{p}(B^\dagger_\fp)$ which
has Swan conductor two and type II at $C^\dagger_\fp$ and for which all assumptions of
\lemref{lem:spl-loc-1} are satisfied. Applying this lemma (with
$t = t_2, t + \pi'_1 = t^\dagger_2$), we conclude that there exists
$A' \in \{A_1, A_2\}$ such that
$1 \le \Sw_{C'^\dagger_\fp}(\chi^\dagger_1|_{B'^\dagger_\fp}) \le 2$,
where $C'^\dagger_\fp = C^\dagger_\fp \otimes_{A_\fp} A'$.
Equivalently, the image of ${\chi}'_1$ under the composite map
$\Fil_2 H^{q+1}_{p}(B') \surj \frac{\Fil_2 H^{q+1}_{p}(B')}{\Fil_0 H^{q+1}_{p}(B')} \to
\frac{\Fil_2 H^{q+1}_{p}(B'_\fp)}{\Fil_0 H^{q+1}_{p}(B'_\fp)} \xrightarrow{\psi^*}
\frac{\Fil_2 H^{q+1}_{p}(B'^\dagger_\fp)}{\Fil_0 H^{q+1}_{p}(B'^\dagger_\fp)}$
is not zero for some $B' \in \{B_1, B_2\}$ and $\chi'_1 \in \{\chi^1_1, \chi^2_1\}$.
It follows that ${\chi}^j_1$ does not die in
$\frac{\Fil_2 H^{q+1}_{p}(B_j)}{\Fil_0 H^{q+1}_{p}(B_j)}$ for some $j \in J_2$, as
desired. This completes the proof.
\end{proof}

\begin{remk}\label{remk:Type-2-12-5}
  In the proof of \lemref{lem:Type-12-key}, the condition $H^c_p(L) \neq 0$ is
  automatically satisfied when $\F$ is a finite field and $q = 1$. Indeed, in this
  case, $L$ is the function field of a smooth projective curve over $\F$, and the
  Hasse principle for the Brauer group of global fields implies that
  $H^2_{p^m}(L) \neq 0$ for every $m \in \N$. Consequently, the proof of
  \lemref{lem:Type-12-key} becomes considerably shorter in this setting. Moreover,
  this special case is sufficient for establishing the main results of the present
  paper.
\end{remk}

\begin{cor}\label{cor:Type-2-12-6}
  Assume in \lemref{lem:Type-12-key} that $\sX$ is semistable and there is a
  closed embedding of $R$-schemes $\sX \inj \P^N_R$. 
  Let $x \in Y^o_i$ be a closed point as in the lemma, and let $n \in \N$. Let
  $m \ge n$ and let $H_1, H_2 \in |\sO_{\P^N_R}(m)|(R)$ be admissible hypersurfaces
  containing $x$ as in \lemref{lem:Type-2-9}.
Then there exists $j \in J_2$ such that letting $\sX_j = \sX \bigcap H_j$,
one has $1 \le \Sw_{Y'_i}(\chi|_{\sX_j}) \le 2$.
\end{cor}
\begin{proof}
Suppose that $\sX_1$ is locally defined at $x$ by an element $t \in \sO_{\sX,x}$.
  Then \lemref{lem:Type-2-9} says that $\sX_2$ is locally defined at $x$ by 
  $t + u \pi_i$ for some $u \in \sO^\times_{\sX,x}$. Since $\sX_1$ is
  admissible, we can write the maximal ideal of $\sO_{\sX,x}$ as
  $\fm = (\pi_i, t, t_3, \ldots , t_{d_\sX})$. As argued in the proof of
  \propref{prop:Type-2-key-0}, this implies using \lemref{lem:Type-12-key} 
that $1 \le \Sw_{Y'_i}(\chi|_{\sX_j}) \le 2$ for some $j \in J_2$.
\end{proof}

\begin{prop}\label{prop:Type-2-12}
  Assume that the conditions {$(\star')$} hold.
  Then there exists a dense open $U_i \subset Y_i$ such that for every closed
 point $x \in U_i$, we can find a quasi-admissible subscheme $\sX' \subset \sX$ of
 dimension $d_\sX -1$ containing $x$ such that
 $1 \le \Sw_{Y'_i}(\chi|_{X'}) \le 2$.
 If $\sX$ is semistable and
 there is a closed embedding $\sX \inj \P^N_R$ of $R$-schemes,
  then we can find, for every closed
  point $x \in U_i$ and every $n \in \N$, an integer $n' \ge n$ and an admissible
  hypersurface section $\sX'$ of $\sX$ of degree $n'$ containing $x$ such that
 $1 \le \Sw_{Y'_i}(\chi|_{X'}) \le 2$.
\end{prop}
\begin{proof}
Under conditions {$(\star')$}, the inequalities $D\ge D^{\<i\>}_2\ge D/2$ do not hold.
Consequently, the argument used under conditions {$(\star)$} does not work
and we have to proceed differently. We argue as follows.

We can assume without loss of generality that $i =1$.
We let $y$ be the generic point of $Y_1$.
We let $\chi_{y}$ be the image of $\chi$ in $H^{q+1}_{p^m}(K_{y})$.
By our assumption, the image of $\chi$ is non-zero in
$\frac{\Fil_D H^{q+1}_{p^m}(X)}{\Fil_{D^{\<1\>}_1} H^{q+1}_{p^m}(X)}$.
We denote this image by $\ov{\chi}$. 
If $R_*(\ov{\chi}) = 0$ (cf. ~\eqref{eqn:Rest**}),
we apply \lemref{lem:Type-12-key},
\corref{cor:Type-2-12-6} and repeat the
proofs of Corollaries~\ref{cor:Type-2-8} and ~\ref{cor:Type-2-11} to conclude.
We assume in the rest of proof that $R_*(\ov{\chi}) \neq 0$.
 
We note that under the condition $R_*(\ov{\chi}) \neq 0$, we must have $m \ge 2$. 
Since $D/2 = Y_1 + \lfloor{{n_2}/2}\rfloor Y_2 + \cdots 
+ \lfloor{{n_r}/2}\rfloor Y_r$ and ${D^{\<1\>}_1}/2 = 0.
Y_1 + \lfloor{{n_2}/2}\rfloor Y_2 +
  \cdots + \lfloor{{n_r}/2}\rfloor Y_r$, it follows that
  $R_*(\ov{\chi})$ remains non-zero when we pass from $X$ to $\Spec(K_{y})$.
Equivalently, $R_*({\chi}_{y})$ does not die in
  $\frac{\Fil_1 H^{q+1}_{p^{m-1}}(K_{y})}{\Fil_{0} H^{q+1}_{p^{m-1}}(K_{y})}$.
It follows that $\Sw_{Y_1}(R_*(\chi)) = 1$. In particular,
  $R_*(\chi)$ has type I at $Y_1$ by \cite[Thm.~3.2(3)]{Kato-89}.
In view of condition (2) of {$(\star')$} therefore, \propref{prop:SC-change} is
  applicable to $R_*(\chi)$. We now choose
  $U_1 \subset Y_1$ and $\sX' \subset \sX$ as in that proposition and let
  $E = u^*(D), \ E' = u^*(D^{\<1\>}_1) = E - Y'_1$, where $u \colon \sX' \inj \sX$
  is the inclusion. 

We now look at the diagram
  \begin{equation}\label{eqn:Type-2-12-0}
    \xymatrix@C1.5pc{
      \frac{\Fil_D H^{q+1}_{p^m}(X)}{\Fil_{D^{\<1\>}_1} H^{q+1}_{p^m}(X)}
      \ar[d]_-{u^*} \ar[r]^-{R_*} &
   \frac{\Fil_{D/2} H^{q+1}_{p^{m-1}}(X)}{\Fil_{{D^{\<1\>}_1}/2} H^{q+1}_{p^{m-1}}(X)}
   \ar[d]^-{u^*} \\
   \frac{\Fil_E H^{q+1}_{p^m}(X')}{\Fil_{E'} H^{q+1}_{p^m}(X')} \ar[r]^-{R_*} &
   \frac{\Fil_{E/2} H^{q+1}_{p^{m-1}}(X')}{\Fil_{{E'}/2} H^{q+1}_{p^{m-1}}(X')}.}
    \end{equation}
  This diagram is commutative by \thmref{thm:H^1-fil} because $\sX'$
  is quasi-admissible. Since $R_*(\ov{\chi})$ maps to a non-zero element under the
  right vertical arrow, it follows that $\ov{\chi}$ maps to a non-zero element
  under the left vertical arrow. This implies that $\Sw_{Y'_1}(\chi|_{X'}) = 2$. 
\end{proof}

The following theorem brings together all of the specialization results established
in this and the preceding section.
  
\begin{thm}\label{thm:SC-change-main}
  Let $q \in J_{d_\sX - 2}, \ i \in J_r$ and $m \in \N$. Let
  $\chi \in H^{q+1}_{p^m}(X)$ with $\Sw_{Y_i}(\chi)=n_i\ge 1$.
  Then there exists a dense open $U_i \subset Y^o_i$ such that for every closed
 point $x \in U_i$, the following hold.
  \begin{enumerate}
  \item
 There exists a quasi-admissible subscheme $\sX' \subset \sX$ of
 dimension $d_\sX -1$ containing $x$ such that
 $n_i -1 \le \Sw_{Y'_i}(\chi|_{X'}) \le n_i$.
\item
  If $\sX$ is semistable and there is a closed embedding $\sX \inj \P^N_R$ of
  $R$-schemes, then for every $n \in \N$, there exists an integer $n' \ge n$ and an 
  admissible hypersurface section $\sX'$ of $\sX$ of degree $n'$ containing $x$ such
  that $n_i -1 \le \Sw_{Y'_i}(\chi|_{X'}) \le n_i$.
\item
  In any of the above cases, if $\chi$ has either type I at $Y_i$,
  or, has type II at $Y_i$ and $q \le d_\sX -3$, then $\Sw_{Y'_i}(\chi|_{X'})= n_i$ and
  $\chi|_{X'}$ has the same type at $Y'_i$ as $\chi$ has at $Y_i$.
  \end{enumerate}
\end{thm}

\section{Kato vs Evaluation filtration}\label{sec:KEF}
In this section, we prove \thmref{thm:Main-1}. We recall the setting of this
theorem. We let $k$ be an hdvf of characteristic $p$
with excellent ring of integers $\sO_k$ and finite residue field $\F$.
We let $\fm = (\pi)$ be the maximal ideal of $\sO_k$.
We let $S = \Spec(\sO_k), \eta = \Spec(k)$ and $s = \Spec(\F)$.

We fix a quasi-semistable $\sO_k$-scheme $\sX$. Thus,  $\sX$ is a finite type,
connected, regular, separated, faithfully flat $\sO_k$-scheme
such that the reduced special fiber of $\sX$
is an snc divisor on $\sX$. Let $K$ denote the function field of $\sX$.
We let $i \colon \sX_s \inj \sX$ and
$j \colon X = \sX_\eta \inj \sX$ denote the inclusions of special and
generic fibers, respectively. We let $f \colon X \to \Spec(k)$ denote the
structure map.
We let $Y_1,\ldots,Y_r$ be the irreducible components of $Y = (\sX_s)_\red$.
We let $\sX^o = \sX \setminus Y_\sing$ and $Y^o = Y_\reg = Y^o_1 \amalg \cdots
\amalg Y^o_r$, where $Y^o_i = Y_i \setminus Y_\sing$. Let $d_\sX = \dim(\sX) \ge 2$.
We let $\eta_X$ denote the generic point of $X$ (or $\sX$) and let
$y_i$ denote the generic point of $Y_i$.
Let $\iota' \colon \eta_X \inj \sX$ and $\iota \colon \eta_X \inj X$ denote the
inclusion maps. We let $\sX_s = \stackrel{r}{\underset{i=1}\sum} \ov{n}_i Y_i
\in \Div(\sX)$.

\vskip.2cm

For $n \in \N_0$ and $m, q \in \N$, we let $\Fil_n H^q_m(X) = \Fil_{nY}H^q_m(X)$
and $\Fil_n H^q(X) = \Fil_{nY}H^q(X)$. We similarly write
$\Fil_n \Br(X) = \Fil_{nY} \Br(X)$.
For $D' \le D \in \Div_Y(\sX)$ and $\chi \in \Fil_D H^q(X)$,
we shall denote the image of $\chi$ in
$\frac{\Fil_D H^q(X)}{\Fil_{D'} H^q(X)}$ by $\ov{\chi}$.
We shall adopt a similar notation for Brauer classes.
For any $\chi \in \Fil_{nY} \Br(X)[p^m]$ (or in $\Fil_{nY} H^2(X)$), we shall write
$\Rsw^{m,2}_{\sX|(nY, nY-Y_i)}(\ov{\chi})$ as $\Rsw^m_{\sX, n,i}(\ov{\chi})$.
Recall here from \S~\ref{sec:MAF} that all refined Swan conductor maps for
$q = 2$ are defined on Brauer groups because of the isomorphism
$\frac{\Fil_D H^2_m(X)}{\Fil_{D'} H^2_m(X)} \xrightarrow{\cong}
\frac{\Fil_D \Br(X)[m]}{\Fil_{D'} \Br(X)[m]}$ (cf. \lemref{lem:Br-fil-0}).

We refer the reader to \S~\ref{sec:Ev-Br} for the definition of the evaluation
filtration on $\Br(X)$.
We shall now compare the Matsuda filtration on $\Br(X)$ with
latter's evaluation filtration in several cases. We shall use the
following result at several places in our proofs.

\begin{lem}\label{lem:K-Ev-3}
  Let $u \colon \sX' \inj \sX$ be the inclusion of a quasi-admissible subscheme.
  Then $u^* \colon \Br(X) \to \Br(X')$ preserves the evaluation filtration.
\end{lem}
\begin{proof}
  For $n \ge 0$, this follows easily from \lemref{lem:good-sub-0} while
  for $n < 0$, this follows directly from the definition of ${\Ev}_{n} \Br(X)$.
\end{proof}

\subsection{The case $n \ge 0$}\label{sec:non-neg}
The first step in proving \thmref{thm:Main-1} is the following result
whose prime-to-$p$ part can be also deduced using the argument of
\cite[Prop.~5.1]{Bright}.

\begin{lem}\label{lem:K-Ev-1}
We have $\Fil'_0 \Br(X) = \Fil_0 \Br(X) \subset {\Ev}_0 \Br(X)$.
\end{lem}
\begin{proof}
  The equality $\Fil'_0 \Br(X) = \Fil_0 \Br(X)$ is already shown in
\corref{cor:Matsuda-5}. We show the inclusion part of the lemma.
For any $P \in X^o_\fin$ which  specializes to $P_0 \in Y^o$,
\propref{prop:Kato-Br-0} says that the diagram
\begin{equation}\label{eqn:K-Ev-1-0}
\xymatrix@C1.2pc{
  \Fil_0 \Br(X) \ar[r]^-{\partial_X} \ar[d]_-{\iota^*_P} &  H^1_\et(Y^o)
  \ar[d]^-{\wt{\iota}^*_{P_0}} & \\
  \Br(k(P)) \ar[r]^-{\partial_P} & H^1_\et(k(P_0)) \ar[r]^-{\cong} & H^1_\et(\F),}
\end{equation}
is commutative, where $\iota_P \colon \Spec(k(P)) \inj X$ and
$\iota_{P_0} \colon \Spec(k(P_0)) \inj Y^o$ are the inclusions, the lower
horizontal arrow on the right is the push-forward map associated to the map
$\Spec(k(P_0)) \to \Spec(\F)$, and the right vertical arrow was defined in
~\S~\ref{sec:Fil-0**}. This implies
$\ev_\chi(P)= \wt{\iota}^*_{P_0} \circ \partial_X(\chi)$. Now, for any $Q \in B(P,1)$,
we have $\ov P \times_\sX Y = \ov Q \times_\sX Y$. In particular, $P_0=Q_0$ and
$\wt{\iota}_{P_0}^*=\wt{\iota}_{Q_0}^*$. But this implies $\ev_\chi(P) = \ev_\chi(Q)$.
\end{proof}

To prove our results for the $p$-primary Brauer classes, we shall use the
following.

\begin{lem}\label{lem:K-Ev-2}
  Let $n \in \N_0$ and let $\chi \in \Br(X)\{p\}$ be such that
  $\chi \notin \Fil'_n \Br(X)$.
  Then we can find a quasi-admissible relative curve $\sX' \subset \sX$ such that
  $\chi|_{X'} \notin \Fil'_n \Br(X')$.
\end{lem}
\begin{proof}
We can assume that $\chi \in \Br(X)[p^m]$ for some $m \in \N$. We can also
  assume that $d_\sX \ge 3$ as there is nothing to prove otherwise.
  In particular, $q = 1 \le d_\sX -2$.
  Recall now that $\Fil'_n \Br(X) \subset \Fil'_{n+1} \Br(X)$ by definition and
  $\Fil_{n+1} \Br(X) \subset \Fil'_{n+1} \Br(X)$ by \corref{cor:Matsuda-2}.
Suppose first that $\chi \in \Fil_{n+1} \Br(X)[p^m]$. Then our hypothesis implies
that $\nu^0_i \circ \Rsw^{m}_{\sX, n+1, i}(\ov{\chi}) \neq 0$ for some
$i \in J_r$. In particular, $\chi$ has type I at $Y_i$.
We can thus apply \thmref{thm:SC-change-main} repeatedly to get a quasi-admissible
relative curve
$\sX' \subset \sX$ such that $\Sw_{Y'_i}(\chi|_{X'}) = n+1$ and $\chi|_{X'}$ has
type I at $Y'_i$. But this implies that $\chi|_{X'} \notin \Fil'_n \Br(X')$.

We assume now that $\chi \notin \Fil_{n+1} \Br(X)[p^m]$. This means that
$\Sw_{Y_i}(\chi) > n+1$ for some $i$. If $\chi$ has type I at $Y_i$, then
we can apply \thmref{thm:SC-change-main} to get a quasi-admissible relative curve
$\sX' \subset \sX$ such that $\Sw_{Y'_i}(\chi|_{X'}) > n+1$ and $\chi|_{X'}$ has
type I at $Y'_i$. But this implies that $\chi|_{X'} \notin \Fil'_n \Br(X')$.
If $\chi$ has type II at $Y_i$, we can
apply \thmref{thm:SC-change-main} to get a quasi-admissible relative curve
$\sX' \subset \sX$ such that $\Sw_{Y'_i}(\chi|_{X'}) > n$.
In particular, $\chi|_{X'} \notin \Fil_n \Br(X')$. Since $\sX'$ is a relative
curve, it follows from \corref{cor:Matsuda-5}
(more precisely,  \corref{cor:Matsuda-2}) that $\chi|_{X'} \notin \Fil'_n \Br(X')$.
This completes the proof.
\end{proof}

We can now prove the following.

\begin{lem}\label{lem:K-Ev-4}
We have ${\Ev}_n \Br(X) \subset \Fil'_n \Br(X)$ for every $n \in \N_0$.
\end{lem}
\begin{proof}
If $\chi \in {\Ev}_n \Br(X)$, we can write $\chi = \chi_1 + \chi_2$,
where $\chi_1 \in {\Ev}_n \Br(X)[p^m]$ for some $m \in \N$ and
$\chi_2 \in {\Ev}_n\Br(X)\{p'\}$. 
  In particular, $\chi_2 \in \Fil_0 \Br(X) \subset \Fil_n \Br(X) \subset
  \Fil'_n \Br(X)$ by Lemmas~\ref{lem:Log-fil-D-1} and ~\ref{lem:Br-fil-0},
  and \corref{cor:Matsuda-5}.
 We can thus assume that $\chi \in {\Ev}_n \Br(X)[p^m]$.
 By Lemmas~\ref{lem:K-Ev-3} and ~\ref{lem:K-Ev-2}, we can assume that $d_\sX = 2$.
 In particular,  $\Omega^2_{k(y_i)} = 0$ for every $i \in J_r$, and
 $\Fil_n \Br(X) = \Fil'_n \Br(X)$ (cf. \corref{cor:Matsuda-5}).

Now, suppose on the contrary that $\chi \notin \Fil_n \Br(X)$. We can then find
  $i \in J_r$ such that $\Sw_{Y_i}(\chi) > n$. We assume without loss of
  generality that $i = 1$. We let $\Sw_{Y_j}(\chi) = n_j \ge 0$ for $j \in J_r$ and
  let $D = \sum_j n_j Y_j$ so that $n_1 > n$. We let $D' = D- Y_1$.
  We let $\wt{\chi} = \Rsw^{m,2}_{\sX|(D,D')}(\ov{\chi}) \in
  H^0_\zar(Y_1, \Omega^1_{Y_1}(\log F_1)(D_1))$. Since $\wt{\chi} \neq 0$,
  we can find a closed point $x \in Y^o_1$ such that $\wt{\chi}$ does not die
  in $\Omega^1_{Y_1}(\log F_1)(D_1) \otimes_{\sO_{Y_1}} k(x) \cong
  \Omega^1_{Y^o_1} \otimes_{\sO_{Y^o_1}} k(x)$.

We let $A = \sO_{\sX,x}$ with maximal ideal $(\pi_1, t)$, where $(\pi_1)$ is
the ideal of $A$ defining $Y_1$ at $x$. We have $H^1_p(k(x)) \cong {\Z}/p$.
We let $B = A[\pi^{-1}_1]$.
  We let $\chi_x$ denote the image of $\chi$ in
  $\Br(B)[p^m]$ and let $\wt{\chi}_x$ denote the image of $\wt{\chi}$ in
  $\Omega^1_{A_1} \otimes_{A_1} k(x)$, where $A_1 = {A}/{(\pi_1)}$. It follows from
  \corref{cor:w-3} that there exists $a \in A^\times$ such that
  letting $\psi = \{\chi_x, 1 + \pi^{n_1}_1at^{-1}\} \in H^3(K)$ and $\fp = (\pi_1)$,
  we have that $\partial'_{\fp} \circ \partial_\fp (\psi) \neq  0$ in $H^1(k(x))$.
  We let
  $\fp_1 = (t)$ and $\fp_2 = (t + \pi^{n_1}_1a)$ be the
  prime ideals of $A$.

We now look at the Kato complex $K(A,1)$:
  \begin{equation}\label{eqn:K-Ev-4-0}
    H^3(K) \xrightarrow{\partial = \oplus \partial_{\fp'}}
    {\underset{{\rm ht}(\fp') =1}\bigoplus} H^2(k(\fp')) 
    \xrightarrow{\partial' = \sum_{\fp'} \partial'_{\fp'}} H^1(k(x)).
  \end{equation}
As in the proof of \lemref{lem:spl-loc-1}, one has
  $\partial(\psi) = \partial_{\fp}(\psi) + \partial_{\fp_1}(\psi) +
\partial_{\fp_2}(\psi)$.
Since $\chi_x \in \Br(B)$, it follows from the definition of the boundary maps in
the Kato complex (cf. \S~\ref{sec:Kato-complex}) that
$\partial'_{\fp_1} \circ \partial_{\fp_1}(\psi) = - \partial'_{\fp_1}(\chi_x|_{k(\fp_1)})$
and
$\partial'_{\fp_2} \circ \partial_{\fp_2}(\psi) = \partial'_{\fp_2}(\chi_x|_{k(\fp_2)})$.
Since $\partial' \circ \partial(\psi) = 0$  
and $\partial'_{\fp} \circ \partial_\fp (\psi) \neq  0$, it follows that
\begin{equation}\label{eqn:K-Ev-4-1}
\partial'_{\fp_2}(\chi_x|_{k(\fp_2)}) - \partial'_{\fp_1}(\chi_x|_{k(\fp_1)}) \neq 0.
\end{equation}

We let $P = \Spec(k(\fp_1))$ and $Q = \Spec(k(\fp_2))$
so that $\ov{P} = \Spec(A/{\fp_1})$ and $\ov{Q} = \Spec(A/{\fp_2})$.
Then it is clear that $P, Q \in X^o_\tr \subset X^o_\fin$ (cf. \lemref{lem:qs-fin}),
and ~\eqref{eqn:K-Ev-4-1} says that
${\rm ev}_{\chi}(P) \neq {\rm ev}_{\chi}(Q)$. Since
$(t, \pi^{n_1}_1) = (t, t + \pi^{n_1}_1a)$ in $A$, we get that
$Q \in B(P, n_1)$. It follows that $\chi \notin {\Ev}_{n_1-1} \Br(X)$.
On the other hand, ${\Ev}_{n} \Br(X) \subset {\Ev}_{n_1-1} \Br(X)$ as
$n \le n_1-1$. We conclude that $\chi \notin {\Ev}_{n} \Br(X)$ which contradicts
our assumption. This completes the proof.
\end{proof}

\begin{lem}\label{lem:K-Ev-5}
We have $\Fil'_n \Br(X) \subset {\Ev}_n \Br(X)$ for every $n \in \N_0$.
\end{lem}
\begin{proof}
We can assume $n \in \N$ in view of \lemref{lem:K-Ev-1}. 
If $\chi \in \Fil'_n \Br(X)$, we can write $\chi = \chi_1 + \chi_2$,
where $\chi_1 \in \Fil'_n\Br(X)[p^m]$ for some $m \ge 1$ and
$\chi_2 \in \Fil'_n\Br(X)\{p'\}$. In particular,
$\chi_2 \in \Fil_0 \Br(X) \subset  {\Ev}_0 \Br(X) \subset  {\Ev}_n \Br(X)$ by
Lemmas~\ref{lem:Log-fil-D-1}, ~\ref{lem:Br-fil-0} and ~\ref{lem:K-Ev-1}.
  We have thus assume that $\chi \in \Fil'_n \Br(X)[p^m]$.

We now let $P, Q \in X^o_\fin$ such that $Q \in B(P, n+1)$. We let $x = P_0 = Q_0$
  and $A = \sO_{\sX,x}$. We let $\fp, \fq \subset A$ be ideals such that
  $\ov{P} = \Spec(A/{\fp})$ and $\ov{Q} = \Spec(A/{\fq})$. Since $\ov{P}$ is
  regular, we can write $\fp = (t_1, \ldots , t_d)$, where
  $(t_1, \ldots , t_d, t_{d_\sX})$ is the maximal ideal of $A$. Since $x \in Y^o$,
  we can assume without loss of generality that $x \in Y^o_1$. Let $(\pi_1)$
  be the ideal of $A$ defining $Y_1$ at $x$. 

The assumption $Q \in B(P, n+1)$ means that $(\fp, \pi^{n+1}_1) = (\fq, \pi^{n+1}_1)$
  as ideals. This implies that for each $i \in J_d$, there is
  $t'_i \in \fq$ and $a_i \in A$ such that $t_i = t'_i + a_i \pi^{n+1}_1$.
  As $n+1 \ge 2$, it follows that
\begin{equation}\label{eqn:K-Ev-5-0}
  t'_i \in \fm \  {\rm and} \  t'_i = t_i \ \mbox{mod} \ \fm^2.
  \end{equation}
In particular, $\fm = (t'_1, \ldots , t'_d, t_{d_\sX})$. As
$(t'_1, \ldots , t'_d) \subset \fq$, this implies that
  $\fq = (t'_1, \ldots , t'_d)$. We now consider a sequence of ideals
  $\fp_0, \ldots , \fp_d$ in $A$ given by
  \begin{equation}\label{eqn:K-Ev-5-1}
\fp_0 = \fp \  {\rm and} \  \fp_i = (t'_1, \ldots , t'_i, t_{i+1}, \ldots , t_d)
\ {\rm for} \ i \in J_d. \ \mbox{In particular}, \ \fp_d = \fq. 
\end{equation}

It follows from the definition of $\fp_i$ and ~\eqref{eqn:K-Ev-5-0}
that $\fm = (\fp_i, t_{d_{\sX}})$ for each $i \in J^0_d$. In particular,
$A/{\fp_i}$ is regular. Moreover, $\pi_1 \notin \fp_i$ for any $i$. Indeed,
if $\pi_1 \in \fp_i$, then $\fp \subset \fp_i$, which implies that
$\fp_i = \fp$. This forces $\pi_1$ to lie in $\fp$ which is clearly not possible.
A similar argument tells us that
$(\fp_i, \pi^{n+1}_1) = (\fp, \pi^{n+1}_1)$. It follows that
  there exists a unique point $P_i \in X_{(0)}$ such that
  $\ov{P}_i = \Spec(A/{\fp_i})$ and the specialization $P_i$ is $x$.
  \lemref{lem:qs-fin} implies that $P_i \in B(P, n+1)$
  (note that it is possible that $P_i = P_j$ for $i \neq j$).
  To prove the lemma, it suffices to show that
  ${\rm ev}_{\chi}(P_i) = {\rm ev}_{\chi}(P_{i+1})$ for each $i \in J^0_{d-1}$.

We now define ideals
 \begin{equation}\label{eqn:K-Ev-5-2} 
   I_0 = (t_2, \ldots , t_d), \  I_i =
   (t'_1, \ldots , t'_i, t_{i+2}, \ldots , t_d)
   \ {\rm for} \ i \in J_{d-2} \ \mbox{and}
\end{equation}
   \[
   I_{d-1} = (t'_1, \ldots , t'_{d-1}).
   \]
   in $A$ so that $\fp_i = (I_i, t_{i+1})$ and $\fp_{i+1} = (I_i, t'_{i+1})$ for
 $i \in J^0_{d-1}$. In particular, each ${A}/{I_i}$ is a regular local ring of
 dimension two. We now fix $i \in J^0_{d-1}$ and let $R = A/{I_i}$.

For any $a \in A$, we let $\ov{a}$ denote its image in $R$, and for any ideal
 $I \subset A$, we let $\ov{I}$ be its extension in $R$.
We may write $\ov{\pi}_1 = us^{n_1}_1 \cdots s^{n_l}_l$ for some $u \in R^\times$ and
distinct prime elements $s_j \in R$ with $n_j \ge 1$ for each $j \in J_l$.
Observe that none of the $s_j$ is associate to either $\ov{t}_{i+1}$ or
$\ov{t'}_{i+1}$ in $R$. Indeed, if $s_j$ were associate to $\ov{t}_{i+1}$, then we
would have $\pi_1 \in (I_i, t_{i+1}) = \fp_i$ in $R$. Similarly, if $s_j$ were
associate to $\ov{t'}_{i+1}$, then we
would have $\pi_1 \in (I_i, t'_{i+1}) = \fp_{i+1}$. Both possibilities contradict
what we have already established.

For $j \in J_l$, we let $\fq_j = (s_j) \subset R$ and $R_j = R^h_{\fq_j}$.
The $\sO_k$-algebra homomorphisms $A \surj R \to R_j$ induce the maps
$f_j \colon \Spec(R_j) \to \Spec(A)$ which induce morphisms of snc-pairs
$f_j \colon (\Spec(R_j), \Spec(k(s_j))) \to (\Spec(A), \Spec(A/{(\pi_1)}))$.
If we let $B = A[{\pi}^{-1}_1]$ and let $\chi_x$ be the image of $\chi$ in $\Br(B)$,
then $\chi_x \in \Fil'_n \Br(B)[p^m]$ by \corref{cor:Matsuda-5} and
$Q(R_j) \cong B \otimes_A R_j$. 
Moreover, $f^*_j(\chi_x) \in \Fil'_{(n+1)n_j -1} \Br(Q(R_j))[p^m]$ by
\corref{cor:Matsuda-5}. Equivalently (in the notation of  ~\eqref{eqn:Kato-SC}),
\begin{equation}\label{eqn:K-Ev-5-3} 
\Rsw^m_{R_j, (n+1)n_j}(f^*_j(\chi_x)) =
(\alpha, 0) \in \Omega^2_{k(s_j)} \bigoplus \Omega^1_{k(s_j)}.
\end{equation}

We now look at the Kato complex $K(R,1)$:
  \begin{equation}\label{eqn:K-Ev-5-4}
    H^3(Q(R)) \xrightarrow{\partial = \oplus \partial_{\fp'}}
    {\underset{{\rm ht}(\fp') =1}\bigoplus} H^2(k(\fp')) 
    \xrightarrow{\partial' = \sum_{\fp'} \partial'_{\fp'}} H^1(k(x)).
  \end{equation}
 Letting $\theta = \{\chi_x|_{Q(R)}, \ov{t}_{i+1} (\ov{t'}_{i+1})^{-1}\} \in H^3(Q(R))$
  (note that this makes sense because $\chi_x|_{Q(R)} \in \Br(Q(R)) = H^2(Q(R))$),
  we get
\begin{equation}\label{eqn:K-Ev-5-5}
  \partial(\theta) = \stackrel{l}{\underset{j =1}\sum} \partial_{\fq_j}(\theta)
  + \partial_{\ov{\fp}_i}(\theta) + \partial_{\ov{\fp}_{i+1}}(\theta) =
  \stackrel{l}{\underset{j =1}\sum} \partial_{\fq_j}(\theta) +
  \chi_x|_{k(\fp_i)} - \chi_x|_{k(\fp_{i+1})},
  \end{equation}
where the second equality occurs because $(\ov{t}_{i+1}) = \ov{\fp}_i$ and
$(\ov{t'}_{i+1}) =  \ov{\fp}_{i+1}$ in $R$.

It follows that
\begin{equation}\label{eqn:K-Ev-5-6}
  \begin{array}{lll}
    0 = \partial' \circ \partial (\theta) & = &
 \stackrel{l}{\underset{j =1}\sum} \partial'_{\fq_j} \circ \partial_{\fq_j}(\theta) +
\partial'_{\ov{\fp}_i}(\chi_x|_{k(\fp_i)}) - \partial'_{\ov{\fp}_{i+1}}(\chi_x|_{k(\fp_{i+1})})
    \\
    & = &
    \stackrel{l}{\underset{j =1}\sum} \partial'_{\fq_j} \circ \partial_{\fq_j}(\theta)
    + {\rm ev}_{\chi}(P_i) -  {\rm ev}_{\chi}(P_{i+1}).
    \end{array}
\end{equation}
It remains therefore to show{\footnote{This is the step for which the condition
  $\chi \in \Fil'_n \Br(X)$ is crucial.}}
that $\partial'_{\fq_j} \circ \partial_{\fq_j}(\theta) = 0$ for each $j$.

To that end, we fix $j \in J_l$ and recall that $\partial_{\fq_j}(\theta)$ is the
image of $\theta$ under the composition of maps
$H^3(Q(R)) \xrightarrow{f^*_j} H^3(Q(R_j)) \xrightarrow{\partial_{\fq_j}} H^2(k(s_j))$
and $\partial_{\fq_j} = \pr_2 \circ (\lambda^3(s_j))^{-1}$
(cf. \propref{prop:Kato-basic} and \S~\ref{sec:Kato-complex}), where
$\pr_2 \colon H^3(k(s_j)) \bigoplus H^2(k(s_j)) \to H^2(k(s_j))$ is the projection.
Meanwhile, we know that
\[
\begin{array}{lll}
f^*_j(\theta)  & =  & f^*_j\left(\{\chi_x, \ov{t}_{i+1} (\ov{t'}_{i+1})^{-1}\}\right) 
=  \{f^*_j(\chi_x), 1 + s^{(n+1)n_j}_j {b}_{i+1} (\ov{t'}_{i+1})^{-1}\} \\
 & = & \wt{\lambda}^3(\ov{\pi}_1)\left(\delta^2_m({b}_{i+1} (\ov{t'}_{i+1})^{-1} \alpha),
\delta^1_m({b}_{i+1} (\ov{t'}_{i+1})^{-1} \beta)\right) 
\end{array}
\]
(cf. ~\eqref{eqn:Can-lift-sum-0}), where
$b_{i+1}$ is an associate of $\ov{a}_{i+1}$ in $R_j$ and
$(\alpha, \beta) = \Rsw^m_{R_j, (n+1)n_j}(f^*_j(\chi_x))$.
The second equality occurs by
\cite[Thm.~5.1]{Kato-89}. Since $\beta = 0$ by ~\eqref{eqn:K-Ev-5-3},
we get $\partial_{\fq_j}(\theta) = 0$, and this completes the proof.
\end{proof}

\vskip.2cm

We give an application of the proof technique used in \lemref{lem:K-Ev-5}.
We fix $P \in X^o_{\tr}$ and let
$P_0$ be the specialization of $P$ so that $P_0 \in Y^o_i$ for some unique
$i \in J_r$. We let $B_{\tr}(P,n) := X^o_{\tr} \bigcap  B(P,n)$.
We let $A = \sO_{\sX, P_0}$ and let $\fm = (\pi_i, t_1, \ldots , t_d)$
be the maximal ideal of $A$, where $\fp = I(\ov{P}) = (t_1, \ldots , t_d)$ and
$d = d_\sX -1$. We define an $A$-linear map $\Phi \colon A^d \to T_{P_0}(Y_i)$ by
$\Phi((\alpha_j)) = \stackrel{d}{\underset{j =1}\sum} \ov{\alpha}_j
\frac{\partial}{\partial {\ov{t}_j}}$, where $T_{P_0}(Y_i)$ is the tangent
space  of $Y_i$ at $P_0$.

Suppose next that $n \in \N$ and $Q \in B_{\tr}(P,n)$.
If we let $\fq = I(\ov{Q}) \subset A$, we have $(\fp, \pi^n_i) = (\fq, \pi^n_i)$.
This implies (cf. the proof of \lemref{lem:K-Ev-5}) that there exists
$t'_j \in \fq$ and $a_j \in A$ such that $t_j = t'_j + a_j\pi^n_i$ for $j \in J_r$.
It is then easy to check that $\fm = (t'_1, \ldots , t'_d, t_{d_{\sX}})$,
which implies that $\fq = (t'_1, \ldots , t'_d)$.
We thus get a surjective map of sets $\Psi \colon A^d \surj B_{\tr}(P,n)$, given by
$\Psi((\alpha_j)) = (t_1 + \alpha_1\pi^n_i, \ldots , t_d + \alpha_d\pi^n_i)$.

\begin{exm}\label{exm:Smooth-case*}
  Assume that $\sX$ is smooth over $S$. Suppose that $Q \in B_\tr(P,n)$ and
  write  $Q = \Psi((\alpha_j))$. Then $\Phi((\alpha_j))$ coincides with the tangent
  vector $[\stackrel{\longrightarrow}{PQ}]_n$, as defined in
  \cite[\S~7]{Bright-Newton}.
\end{exm}

Let $P \in X^o_{\tr}$ be as above. Let $n \in \N$ and
$\chi \in \Fil_n\Br(X)[p^m]$.  Let
$\beta = \nu^0_i \circ \Rsw^{m}_{\sX, n,i}(\ov{\chi})
\in H^0_\zar(Y_i, \Omega^1_{Y_i}(\log F_i)(nY|_{Y_i}))$.
Let $\beta_0$ be the image of $\beta$ in
$(\Omega^1_{Y_i})_{P_0} := \Omega^1_{Y_i} \otimes_{\sO_{Y_i}} k(P_0)$ under the canonical
quotient map.
Let $(\Omega^1_{Y_i})_{P_0} \times  T_{P_0}(Y_i) \to k(P_0)$; ($(x,y) \mapsto \<x,y\>$),
denote the duality map.
A review of the proof of \lemref{lem:K-Ev-5} yields the following
positive characteristic analogue of the first part of \cite[Thm.~B]{Bright-Newton}.

\begin{thm}\label{thm:K-Ev-5-extra}
  If $Q \in B_\tr(P,n)$ and we write $Q = \Psi((a_j))$, then 
  \[
    {\ev}_{\chi}(Q) = {\ev}_{\chi}(P) +
    {\Tr}_{{k(P_0)}/{\F}}\left(\<\beta_0, \Phi((a_j))\>\right).
    \]
    \end{thm}
\begin{proof}
We can assume without loss of generality that $i =1$. We let $A_1 = A/{(\pi_1)}$
and write $\beta = \stackrel{d}{\underset{i =1}\sum}c_i d\ov{t}_i
\in \Omega^1_{A_1}$.
We let $\fp_i$ for $i \in J^0_d$ and $I_i$ for $i \in J^0_{d-1}$ be the ideals of
$A$ chosen in the proof of \lemref{lem:K-Ev-5}. We let $B_i = {A}/{I_i}$ and
${A}_{1,i} := {A}/{((\pi_1) + I_i)}$.
We let $b_j, x_j, x'_j$ and $\wt{\pi}$ denote the images of $a_j, t_j, t'_j$ and
$\pi_1$, respectively, under the quotient map $A \surj B_i$.
We let ${\fq}_j = \fp_jB_i$.
For an element $\alpha \in B_i$ (resp. $A_1$, we let $\ov{\alpha}$ denote its
image in $A_{1,i}$ under the quotient map $B_i \surj A_{1,i}$
(resp. $A_1 \surj A_{1,i}$). We then get
\begin{equation}\label{eqn:K-Ev-5-extra-0}
  \frac{x_{i+1}}{x'_{i+1}} = 1 + {b}_{i+1} (x'_{i+1})^{-1} \wt{\pi}^n \
  \mbox{in} \ Q(B_i).
\end{equation}

We let $\chi_i = \chi|_{B'_i}$.
By the functoriality of the refined Swan conductor, we get
$\beta_{i+1} := \nu^0 \circ \Rsw^{m}_{B_i, n,1}(\ov{\chi}_i) =
\ov{c}_{i+1} d\ov{x}_{i+1} \in \Omega^1_{A_{1,i}} \cong A_{1,i} d\ov{x}_{i+1}$.
Let $\beta_{i+1,0}$ denote the image of $\beta_{i+1}$  in
$\Omega^1_{A_{1,i}} \otimes_{A_{1,i}} k(P_0)$. 

We let
$\psi_i = \{\chi_i, 1 + \wt{\pi}^n \ov{b}_{i+1} (x'_{i+1})^{-1}\} \in H^3(Q(B_i))$ and
consider the Kato complex
\begin{equation}\label{eqn:K-Ev-5-extra-1}
    H^3(Q(B_i)) \xrightarrow{\partial = \oplus \partial_{\fp'}}
    {\underset{{\rm ht}(\fp') =1}\bigoplus} H^2(k(\fp')) 
    \xrightarrow{\partial' = \sum_{\fp'} \partial'_{\fp'}} H^1(k(P_0)).
  \end{equation}
The same argument as in \lemref{lem:K-Ev-5} shows that there exists
$\alpha_{i+1} \in \Omega^2_{A_1}$ for $i \in J^0_{d-1}$ such that
\[
{\ev}_{\chi}(\fq_i) - {\ev}_{\chi}({\fq}_{i+1}) +
\Tr_{{k(P_0)}/{\F}}(\<\beta_{i+1,0}, \wt{a}_{i+1}\>)  =
\]
\[
\partial'_{\fq_i} \circ \partial_{\fq_i}(\psi_i) +
\partial'_{\fq_{i+1}}  \circ \partial_{\fq_{i+1}}(\psi_i) +
\partial'_{(\wt{\pi})} \circ \partial_{(\wt{\pi})}
\circ
\wt{\lambda}^1_p(\wt{\pi})\left(\alpha_{i+1}, \ov{c}_{i+1}\ov{b}_{i+1}
\dlog(\ov{x}_{i+1})\right) =
\]
\[
\partial'_{\fq_i} \circ \partial_{\fq_i}(\psi_i) +
\partial'_{\fq_{i+1}}  \circ \partial_{\fq_{i+1}}(\psi_i) + \partial'_{(\wt{\pi})} \circ
\partial_{(\wt{\pi})}(\psi_i) =  \partial' \circ \partial(\psi_i) = 0,
\]
where $\wt{a}_{i+1}$ is the image of $a_{i+1}$ in $k(P_0)$. In the first
equality, we used that $\ov{x}_{i+1} = \ov{x}'_{i+1}$.
Taking the sum over $i \in J^0_{d-1}$, we get
\[
{\ev}_{\chi}(Q) - {\ev}_{\chi}(P) = \stackrel{d-1}{\underset{i=0}\sum}  
\Tr_{{k(P_0)}/{\F}}(\<\beta_{i+1,0}, \wt{a}_{i+1}\>) =
   {\Tr}_{{k(P_0)}/{\F}}\left(\<\beta_0, \Phi((a_i))\>\right).
   \]
   This completes the proof of the theorem.
\end{proof}

\begin{cor}\label{cor:K-Ev-5-extra-2}
  Let $\chi$ be as in \thmref{thm:K-Ev-5-extra} and assume that $\beta_0 \neq 0$.
  Then the map ${\ev}^{P}_{\chi} \colon
      B_\tr(P,n) \to ({\Q}/{\Z})[p]$, given by
      ${\ev}^{P}_{\chi}(Q) =  {\ev}_{\chi}(Q) - {\ev}_{\chi}(P)$, is surjective.
      In particular, for any $\chi \in  \Br(X)[p] \setminus \Fil_0 \Br(X)$,
      the map ${\ev}_\chi \colon X^o_{\fin} \to ({\Q}/{\Z})[p]$
      is surjective.
\end{cor}
\begin{proof}
  Using the surjectivity of $\Psi$ and \thmref{thm:K-Ev-5-extra}, the first part of
  the lemma follows if one knows that the map $\psi \colon A^d \to \F$, given by
  $\psi((a_j)) =  {\Tr}_{{k(P_0)}/{\F}}\left(\<\beta_0, \Phi((a_j))\>\right)$,
  is surjective. But this is clear.
  To prove the second part, let $\alpha \in  ({\Q}/{\Z})[p]$ and choose
  $P \in X^o_{\tr}$. We then have $\alpha = (\alpha - {\ev}_\chi(P)) +
  {\ev}_\chi(P) = {\ev}^{P}_{\chi}(Q) + {\ev}_\chi(P) = {\ev}_\chi(Q)$ for
  some $Q \in  B_\tr(P,n)$ by the first part of the corollary.
\end{proof}

\subsection{The case $n < 0$}\label{sec:Neg}
To deal with the case $n < 0$ of \thmref{thm:Main-1}, we shall use the following.

\begin{lem}\label{lem:Sp-surj}
  The specialization map $\esp \colon X^o_\tr \to Y^o_{(0)}$ is
  surjective. If $\sX_s$ is reduced, then the composite map
  $X_{\rm ur} \inj X^o_\fin \xrightarrow{\esp} Y^o_{(0)}$
  is surjective.
\end{lem}
\begin{proof}
 Let $x$ be a closed point of $Y^o_i$ for
  some $i$ and let $A = \sO_{\sX, x}$ with maximal ideal
  $\fm_x = (\pi_i, t_1, \ldots , t_{d_{\sX-1}})$, where $(\pi_i)$ defines $Y^o_i$
  locally at $x$. We let $\fp = (t_1, \ldots , t_{d_{\sX-1}})$. Then it it clear that
  there exists a unique closed point $P \in X$ such that
  $\ov{P} = \Spec(A/{\fp})$. It follows then that $P \in X^o_\tr$ and $P_0 =
  \esp(P) = x$. If $\sX_s$ is reduced, then the map $\sO_k \to A/{\fp}$ is
  unramified. As it is clearly flat, it follows that $\ov{P}$ is {\'e}tale over
  $S$.
 \end{proof}

\begin{lem}\label{lem:K-Ev-6}
  We have ${\Ev}_{-2} \Br(X) = \Br(\sX)$ as subgroups of $\Br(X)$.
\end{lem}
\begin{proof}
  The inclusion $\Br(\sX) \subset {\Ev}_{-2} \Br(X)$ follows directly from the
  fact that for any $P \in X^o_{\fin}$, the pull-back map
  $\Br(\sX) \to \Br(k(P))$ has a factorization
  $\Br(\sX) \to \Br(\sO_{k(P)}) \to \Br(k(P))$ and $\Br(\sO_{k(P)}) =0$.
 To show the other inclusion, first note that ${\Ev}_{-2} \Br(X)$
  is a subgroup of ${\Ev}_{0} \Br(X)$ (cf. \S~\ref{sec:Ev-Br})
  which is a subgroup of $\Fil'_0 \Br(X)$ by \lemref{lem:K-Ev-4}. It follows that
  ${\Ev}_{-2} \Br(X) \subseteq \Fil'_0 \Br(X)$. 

  We now let $\chi \in {\Ev}_{-2} \Br(X)$ and let
  $\alpha = \partial_X(\chi)$ in ~\eqref{eqn:K-Ev-1-0}.
  By \propref{prop:Kato-Br-0}, we need to
  show that $\alpha = 0$. By \cite[Thm.~5.8.16]{Szamuely}, it suffices to show that
  $\alpha|_{k(x)} = 0$ for every $x \in Y^o_{(0)}$. To show the latter vanishing,
  we let $x \in Y^o_{(0)}$
  and apply \lemref{lem:Sp-surj} to get $P \in X^o_\tr$ such that $P_0 = x$.
  We then get $\alpha|_{k(x)} = \iota^*_{P_0}(\alpha) = \wt{\iota}^*_{P_0}(\alpha)
  = \partial_P(\chi|_{k(P)}) = {\ev}_{\chi}(P) = 0$, where the
  second equality occurs because $P \in X^o_\tr$, and the third equality occurs by
  \propref{prop:Kato-Br-0}. This completes the proof.
\end{proof}

We now consider the last remaining case $n = -1$. To deal with this case,
we let $g \colon Y^o \to \Spec(\F)$ be the structure map of $Y^o$ and let
$g_i$ be the restriction of $g$ to $Y^o_i$ for $i \in J_r$.
Recall that $\sX_s = \stackrel{r}{\underset{i=1}\sum} \ov{n}_i Y_i$.
We let $H^1_\et(Y^o)_c$ denote the image of the map
$\wt{g}^* \colon H^1_\et(\F) \to H^1_\et(Y^o)$, given by
$\wt{g}^*(a) = (\ov{n}_1 g^*_1(a), \ldots , \ov{n}_r g^*_r(a))$
(cf. ~\eqref{eqn:Kato-Br-0*}). For any point $x \in Y$, we let
$\psi_x \colon \Spec(k(x)) \to \Spec(\F)$ denote the projection map.
We let
\begin{equation}\label{eqn:Fil-neg}
  \Fil'_{-1} \Br(X) = \{\sA \in \Fil_0 \Br(X)| \partial_X(\sA) \in
  H^1_\et(Y^o)_c\}.
\end{equation}
Applying \propref{prop:Kato-Br-0} to the map of snc-pairs
($\sX, Y)  \to (S,s)$, and using the isomorphism
$\partial_k \colon \Br(k) = \Fil_0 \Br(k) \xrightarrow{\cong} H^1_\et(\F)$,
we see that
\begin{equation}\label{eqn:Fil-Br-spl}
  \Fil'_{-1} \Br(X) = \Br(\sX) + \Br_0(X),
\end{equation}
where recall that $\Br_0(X)$ is the image of the pull-back map
$f^* \colon \Br(k) \to \Br(X)$.

\begin{lem}\label{lem:K-Ev-7}
  We have $\Fil'_{-1} \Br(X) \subset {\Ev}_{-1} \Br(X)$. In particular,
   $f^*(\chi) \in {\Ev}_{-1} \Br(X)$ for any $\chi \in \Br(k)$.
\end{lem}
\begin{proof}
Suppose that $\chi \in \Fil'_{-1} \Br(X)$. Let $\alpha \in H^1_\et(\F)$ be such that
$\partial_{X}(\chi) = \wt{g}^*(\alpha)$.
  We let $P, Q \in X^o_\fin$ such that $Q \in B(P,0)$. Using ~\eqref{eqn:K-Ev-1-0},
  we get
  \begin{equation}\label{eqn:K-Ev-7-*-0}
    \begin{array}{lll}
  {\rm ev}_{\chi}(P) & = &
  \partial_P(\chi|_{k(P)}) = \wt{\iota}^*_{P_0}(\partial_X(\chi)) =
  \wt{\iota}^*_{P_0}(\wt{g}^*(\alpha)) =
  \ov{n}_i e(P) \iota^*_{P_0} (g^*_i(\alpha)) \\
  & = & \ov{n}_i e(P) \psi^*_{P_0}(\alpha) =
  \ov{n}_i e(P) [k(P_0): \F] \alpha = [k(P):k] \alpha,
    \end{array}
    \end{equation}
  where the last equality holds by ~\eqref{eqn:mult}.
   Similarly, we get ${\rm ev}_{\chi}(Q) = [k(Q), \F]\alpha$.
It follows that ${\rm ev}_{\chi}(P) = {\rm ev}_{\chi}(Q)$. This shows that
$\chi \in {\Ev}_{-1} \Br(X)$. The second part of the lemma is clear from the
first part and ~\eqref{eqn:Fil-Br-spl}.
\end{proof}

\begin{lem}\label{lem:K-Ev-8}
  We have $ {\Ev}_{-1} \Br(X) \subset \Fil'_{-1} \Br(X)$ if each $Y_i$ is
  geometrically integral and $\ov{n}_i = \ov{n}_j $ for all $i,j  \in J_r$.
\end{lem}
\begin{proof}
  In the proof of this lemma, we closely follow the strategy used in the proof of
  the second part of \cite[Lem.~9.5]{Bright-Newton}. However, additional arguments
  are required to handle the more general setting considered here. Let
  $\chi \in {\Ev}_{-1} \Br(X)$. We will show that $\chi \in \Fil'_{-1} \Br(X)$ by
  proceeding as follows.

For $i \in J_r$, we let $\iota_i \colon Y_i \inj Y$ denote the inclusion
 and let $\partial_i \colon \Fil_{0} \Br(X) \xrightarrow{\partial_X} H^1_\et(Y^o)
 \xrightarrow{\cong} \ \stackrel{r}{\underset{j =1}\bigoplus} H^1_\et(Y^o_j)
 \xrightarrow{\iota^*_i} H^1_\et(Y^o_i)$ be the composite map.
 We let $\wt{g}^*_i = \iota^* \circ \wt{g}^* \colon H^1_\et(\F) \to H^1_\et(Y^o_i)$
 be the composite map and let $H^1_\et(Y^o_i)_c$ denote the image of $\wt{g}^*_i$.
Note that $H^1_\et(Y^o_i)_c$ coincides with the image of the pull-back map
$g^*_i \colon H^1_\et(\F) \to H^1_\et(Y^o_i)$ because $\wt{g}^*_i = \ov{n}_i g^*_i$ and
$H^1_\et(\F)$ is divisible.

 \vskip .2cm

{\bf{Claim~1:}} $\partial_i(\chi) \in H^1_\et(Y^o_i)_c$ for each $i$.

  To prove the claim, suppose there exists $i \in J_r$ such that
  $\partial_i(\chi) \notin H^1_\et(Y^o_i)_c$. 
We let $\wt{Y}^o_i$ denote the
  base change of $Y^o_i$ to $\ov{\F}$ and let $\phi_i \colon \wt{Y}^o_i  \to Y^o_i$
  denote the projection. We let $\alpha = \partial_i(\chi)$ and
  $\wt{\alpha} = \phi^*_i(\alpha)$. 

{\bf{Step~1:}}
  Following the argument of \cite[Lem.~9.5]{Bright-Newton}, we use 
  the exact sequence (cf. \cite[Prop.~5.6.1]{Szamuely})
  \begin{equation}\label{eqn::K-Ev-8-0}
    H^1_\et(\F) \xrightarrow{\wt{g}^*_i} H^1_\et(Y^o_i) \xrightarrow{\phi^*_i}
    H^1_\et(\wt{Y}^o_i) \to 0
  \end{equation}
  and apply \propref{prop:Kato-Br-0} to the map of snc-pairs
  $(\sX,Y) \to (S,s)$. This allows us to find
  $\chi' \in \Br(k)$ such that letting $\alpha' := \partial_i(\chi - f^*(\chi'))$,
  we have that $\phi^*_i(\alpha') = \phi^*_i(\alpha)$, both sides of this
  equality have the same order as
  that of $\alpha'$, and the latter does not lie in $H^1_\et(Y^o_i)_c$.
  Furthermore, \lemref{lem:K-Ev-7} ensures that
  replacing $\chi$ by $\chi - f^*(\chi')$ 
  does not alter our assumption that $\chi \in {\Ev}_{-1} \Br(X)$.
  We can therefore replace $\alpha$ by $\alpha'$, allowing us to assume,
  without changing our hypotheses, that $\alpha$ and $\wt{\alpha}$ have the same
  order. We let $m_i$ be this order. We can assume $m_i \ge 2$, as 
  there is nothing to prove otherwise.

We let $h \colon T_i \to Y^o_i$ denote the ${\Z}/{m_i}$-torsor corresponding to
  $\alpha$ and let $\wt{h} \colon \wt{T}_i \to \wt{Y}^o_i$ denote the base change of
  $h$ to $\ov{\F}$. It follows from \cite[Lem.~5.15]{Bright} that $T_i$ is
  geometrically integral. It is clearly geometrically regular. 
  Since $|Y^o(\F)| < \infty$, there exists a dense open $U \subset Y^o$ such that
  $U(\F) = \emptyset$.  In particular, $T_{U_i} (\F) = \emptyset$ if we let
  $T_{U_i} = h^{-1}(U_i)$, where $U = U_1 \amalg \cdots \amalg U_r$. 
By shrinking $U$ further if necessary, we can assume that $U$ is
quasi-projective over $\F$. We let $\wt{T}_{U_i}$ denote the base change of
$T_{U_i}$ to $\ov{\F}$.

{\bf{Step~2:}}
We choose $a_i, b_i \in H^1_\et(\F, {\Z}/{m_i})$ such that $a_i \neq b_i$
(note that $m_i \ge 2$). We let $h_{a_i} \colon T_{a_i} \to U_i$ denote the
${\Z}/{m_i}$-torsor on $U_i$ corresponding to the restriction of the
element $\alpha - g^*_i(a_i) \in H^1_\et(Y^o_i, {\Z}/{m_i})$ to $U_i$.
We similarly define $h_{b_i} \colon T_{b_i} \to U_i$. If we let
$\wt{h}_{a_i} \colon \wt{T}_{a_i} \to \wt{U}_i$ denote the base change of $h_{a_i}$ to
$\ov{\F}$, we see that $\wt{T}_{a_i} \cong \wt{T}_{{U}_i}$. In particular,
$T_{a_i}$ is geometrically integral and geometrically regular. The same holds for
$T_{b_i}$.
Using \cite[Fact~5.13]{Bright}, we choose an ffe ${\F'}/{\F}$
such that $|\F'| > {\rm max}\{B(\wt{T}_{a_i}), B(\wt{T}_{b_i})| i \in J_r\}$
(see loc. cit. for the notation)
and $[\F': \F]$ is a prime number $\ell$ which does not divide any $m_i$, and
$T_{a_i}(\F') \neq \emptyset \neq T_{b_i}(\F')$ for every $i \in J_r$.
We let $\tau \colon \Spec(\F') \to \Spec(\F)$ denote the projection.

{\bf{Step~3:}}
We let $U'_i$ (resp. $h'_{a_i} \colon T'_{a_i} \to U'_i$) denote the base change of
$U_i$ (resp. $h_{a_i}$) to $\Spec(\F')$. We similarly define
$h'_{b_i} \colon T'_{b_i} \to U'_i$. We choose $x_{a_i} \in T'_{a_i}(\F')$ and let
$y_{a_i} = h'_{a_i}(x_{a_i})$. We let $\gamma_{a_i} \colon \Spec(\F') \to U'_i$
denote the corresponding map. We now look at the commutative diagrams
\begin{equation}\label{eqn::K-Ev-8-1}
  \xymatrix@C1pc{
    S'_{a_i} \ar[r] \ar[d]_-{s_{a_i}} & T'_{a_i} \ar[r] \ar[d]^-{h'_{a_i}} &
    T_{a_i} \ar[d]^-{h_{a_i}} & &
 S'_{b_i} \ar[r] \ar[d]_-{s_{b_i}} & T'_{b_i} \ar[r] \ar[d]^-{h'_{b_i}} &
    T_{b_i} \ar[d]^-{h_{b_i}} \\
    \Spec(\F') \ar[r]^-{\gamma_{a_i}} & U'_i \ar[r]^-{u_i} & U_i & &
\Spec(\F') \ar[r]^-{\gamma_{b_i}} & U'_i \ar[r]^-{u_i} & U_i,}
\end{equation}
in which all squares are Cartesian and the horizontal arrows in the right
squares in both diagrams are the base change maps.

Since $\gamma_{a_i}$ has a lift to $T'_{a_i}$,
the torsor $s_{a_i}$ has a section. In particular, it is trivial
(cf. \cite[Tag~0497]{SP}). If we let $\wt{\gamma}_{a_i} \colon  \Spec(\F') \to 
U'_i \to U_i$ denote the composite map, it follows that
$\wt{\gamma}^*_{a_i}((\alpha - g^*_i(a_i))|_{U_i}) = 0$. Equivalently,
$\wt{\gamma}^*_{a_i}(\alpha|_{U_i}) = \wt{\gamma}^*_{a_i}(g^*_i(a_i)|_{U_i}) =
\tau^*(a_i) = \ell a_i$. Similarly, $\wt{\gamma}^*_{b_i}(\alpha|_{U_i}) =
\tau^*(b_i) = \ell b_i$. Here, we have identified $H^1(\F)$ and $H^1(\F')$
with ${\Q}/{\Z}$ via the isomorphisms $H^1(\F') \xrightarrow{\tau_*} H^1(\F)
\xrightarrow{\inv_{\F}} {\Q}/{\Z}$.

If we let $z_{a_i} = u_i(y_{a_i})$ and $z_{b_i} = u_i(y_{b_i})$, then
$[k(z_{a_i}): \F]$ divides $[\F': \F] = \ell$. Since $U_i(\F) = \emptyset$, this
implies that $[k(z_{a_i}): \F] = \ell$. Similarly, $[k(z_{b_i}): \F] = \ell$.
By \lemref{lem:Sp-surj}, there are points $P, Q \in X^o_\tr$ such that
$P_0 = z_{a_i}$ and $Q_0 = z_{b_i}$. Since $e(P) = e(Q) = 1$,  the equality
$[k(z_{a_i}): \F] = [k(z_{b_i}): \F]$ implies that
$[k(P): k] = [k(Q): k]$ (cf. ~\eqref{eqn:mult}). In particular,
$P,Q \in X^o_\fin$ and $Q \in B(P,0)$. 

{\bf{Step~4:}}
We now look at the diagram
\begin{equation}\label{eqn:K-Ev-8-2}
  \xymatrix@C1pc{
    & H^1_\et(k(P_0), {\Z}/{m_i}) \ar[dr]^-{\id} & & \\
    H^1_\et(U_i, {\Z}/{m_i}) \ar[dr]_-{\iota^*_{Q_0}}
    \ar[ur]^-{\iota^*_{P_0}} & & H^1_\et(\F', {\Z}/{m_i}) \ar[r]^-{\tau_*}_-{\cong} &
    H^1_\et(\F, {\Z}/{m_i}) \\
    & H^1_\et(k(Q_0), {\Z}/{m_i}), \ar[ur]_-{\id} & & }
  \end{equation}
where $\iota_y \colon \Spec(k(y)) \inj U_i$ is the inclusion for a point
$y \in U_i$.
The composite arrow $H^1_\et(U_i, {\Z}/{m_i}) \to H^1_\et(\F, {\Z}/{m_i})$
through $\iota^*_{P_0}$ (resp.  $\iota^*_{Q_0}$) is $\tau_* \circ \wt{\gamma}^*_{a_i}$
(resp. $\tau_* \circ \wt{\gamma}^*_{b_i}$).

Using diagrams~\eqref{eqn:K-Ev-1-0} and ~\eqref{eqn:K-Ev-8-2}, we get
\[
\ev_\chi(P) = \partial_{P}(\chi|_{k(P)}) = \wt{\iota}^*_{P_0} \circ \partial_X(\chi)
   \ {=}^{\dagger} \ \iota^*_{P_0} \circ \partial_X(\chi) =
\iota^*_{P_0}(\alpha) = \iota^*_{P_0}(\alpha|_{U_i}) = \ell a_i; 
\]
\[
\ev_\chi(Q) = \partial_{Q}(\chi|_{k(Q)}) = \wt{\iota}^*_{Q_0} \circ \partial_X(\chi) \
   {=}^{\dagger} \ \iota^*_{Q_0} \circ \partial_X(\chi) =
\iota^*_{Q_0}(\alpha) = \iota^*_{Q_0}(\alpha|_{U_i}) = \ell b_i
\]
in $H^1_\et(\F)$, where the equalities ${=}^{\dagger}$ hold because $P, Q \in
X^o_\tr$.
Since $a_i, b_i$ are two distinct elements of $H^1_\et(\F, {\Z}/{m_i})$ and
$\ell \nmid m_i$, we get $\ell a_i \neq \ell b_i$ in $H^1_\et(\F, {\Z}/{m_i})$.
By \cite[Lem.~7.1]{KM-1}, we get $\ell a_i \neq \ell b_i$ in $H^1_\et(\F)$.
It follows that $\ev_\chi(P) \neq \ev_\chi(Q)$, contradicting
our assumption that $\chi \in \Ev_{-1} \Br(X)$. This proves the claim.

\vskip.2cm

{\bf{Claim~2:}} $\alpha_i = \alpha_j$ for all $i, j \in J_r$.

We let $U \subset Y^o$ be the open subset that we chose in Step~1 above.
As in Step~2 of the proof of Claim~1, we choose an ffe ${\F'}/{\F}$
such that $|\F'| > {\rm max}\{B(\ov{U}_i), m_i| i \in J_r\}$,
$[\F': \F]$ is a prime number $\ell$ which does not divide any $m_i$, and
$U_i(\F') \neq \emptyset$ for every $i \in J_r$.
We let $\tau \colon \Spec(\F') \to \Spec(\F)$ denote the projection.
We pick a point in $U_i(\F')$ and let
$\gamma_i \colon \Spec(\F') \to U_i$ be the corresponding map of $\F$-schemes. Let
$x_i \in U_i$ be the image of this map.
Then we have $[\F': \F] = [k(x_i): \F] = \ell$ by the same argument as in
Step~3 above (recall that $U(\F) = \emptyset$). By \lemref{lem:Sp-surj},
there is a point $P_i \in X^o_\tr$ such that
$\esp(P_i) = x_i$ for every $i$. This yields
\begin{equation}\label{eqn:K-Ev-8-3}
\ev_\chi(P_i) = \partial_{P_i}(\chi|_{k(P_i)}) = \iota^*_{x_i} \circ \partial_X(\chi) =
\iota^*_{x_i}(\alpha_i) = \tau^*(\alpha_i) = \ell \alpha_i \in H^1_\et(\F)
\end{equation}
for every $i$.

If we let $m = \stackrel{r}{\underset{i =1} \prod} m_i$, then a diagram analogous to
~\eqref{eqn:K-Ev-8-2} tells us that
$\ev_\chi(P_i) = \ell \alpha_i \in  H^1_\et(\F, {\Z}/{m})$ for every $i$.
We now note that $P_j \in B(P_i,0)$ for every $i,j$ because we have the
equalities
$[k(P_i):k] = \ov{n}_i[k(x_i): \F] = \ov{n}_j[k(x_j): \F] = [k(P_j):k]$. 
Since $\chi \in \Ev_{-1} \Br(X)$, this implies by ~\eqref{eqn:K-Ev-8-3} that
$\ell \alpha_i = \ell \alpha_j$ in $H^1_\et(\F, {\Z}/{m})$ for every $i,j$.
This forces $\alpha_i = \alpha_j$ for every $i,j$ and proves
Claim~2.

To finish the proof, we let $\alpha = \alpha_i$ and $n = \ov{n}_i$ for $i \in J_r$.
By Claim~1, there exists $\beta \in H^1_\et(\F)$ such that $n\beta = \alpha$.
It follows that
$\wt{g}^*(\beta) = (\ov{n}_1 g^*_1(\beta), \ldots , \ov{n}_r g^*_r(\beta))
= (\alpha_1, \ldots , \alpha_r) = \partial_X(\chi)$. Equivalently, $\chi \in
\Fil'_{-1} \Br(X)$.
The proof of the lemma is now complete.
\end{proof}

\vskip.2cm

\noindent
{\bf{Proof of \thmref{thm:Main-1}:}}
Combine \propref{prop:Kato-Br-0} with Lemmas~\ref{lem:K-Ev-4}, ~\ref{lem:K-Ev-5},
~\ref{lem:K-Ev-6}, ~\ref{lem:K-Ev-7} and ~\ref{lem:K-Ev-8}.
\qed

\vskip.2cm
We now use \thmref{thm:Main-1} to prove a refined version. It shows that
$\Fil'_* \Br(X)$ can be completely described by studying the evaluation of Brauer
classes only at those closed points of $X$ whose closures in $\sX$ meet $Y$
transversely. If $\sX_s$ is reduced, we can describe $\Fil'_* \Br(X)$ by further
restricting the evaluation to only unramified closed points of $X$.
We need to introduce some additional notation.

We let $n \in \N_0$. If $\sX_s$ is reduced and $P \in X_{\ur}$, we let
$B_{\ur}(P,n) = X_{\ur} \bigcap B(P,n )$.
For $n \ge -1$, we let
\begin{equation}\label{eqn:qs-fin-3}
{\Ev}^{\tr}_n \Br(X) = \left\{\sA \in \Br(X)| \ev_\sA \ \text{is constant on} \
B_{\tr}(P, n+1) \ \text{ for all } P \in X^o_{\tr} \right\};
\end{equation}
\[
{\Ev}^{\ur}_n \Br(X) = \left\{\sA \in \Br(X)| \ev_\sA \ \text{is constant on} \
B_{\ur}(P, n+1) \ \text{ for all } P \in X_{\ur} \right\}.
\]
We let
\[
{\Ev}^{\tr}_{-2} \Br(X) = \{\sA \in \Br(X)| {\ev}_{\chi} \ \text{is zero on} \
X^o_\tr\}
\]
and (if $\sX_s$ is reduced)
\[
  {\Ev}^{\ur}_{-2} \Br(X) = \{\sA \in \Br(X)| {\ev}_{\chi} \ \text{is zero on} \
  X^o_{\ur}\}.
  \]
These are all subgroups of $\Br(X)$.

\begin{thm}\label{thm:K-Ev-10}
  Let $n \ge -2$ be an integer. We have
  ${\Ev}^{\tr}_n \Br(X) = {\Ev}_n \Br(X) = \Fil'_n\Br(X)$.
  If $\sX_s$ is reduced, we also have
  ${\Ev}^{\ur}_n \Br(X) =  {\Ev}_n \Br(X) = \Fil'_n \Br(X)$.
\end{thm}
\begin{proof}
  It is clear that ${\Ev}_n \Br(X) \subset {\Ev}^{\tr}_n \Br(X)$,
  and ${\Ev}_n \Br(X) \subset {\Ev}^{\ur}_n \Br(X)$ if $\sX_s$ is reduced.
  On the other hand, a careful examination of the proofs of
  Lemmas~\ref{lem:K-Ev-4}, \ref{lem:K-Ev-6}, and \ref{lem:K-Ev-8} shows that we in
  fact established the inclusion $\Ev^{\tr}_n\Br(X)  \subset \Fil'_n\Br(X)$.
  Moreover, in view of \lemref{lem:Sp-surj}, those proofs also establish that
 ${\Ev}^{\ur}_n \Br(X) \subset  \Fil'_n \Br(X)$ whenever $\sX_s$ is reduced.
We thus get
  ${\Ev}_n \Br(X) \subset {\Ev}^{\tr}_n \Br(X) \subset  \Fil'_n \Br(X) \subset
{\Ev}_n \Br(X)$, where the last inclusion follows from Lemmas~\ref{lem:K-Ev-5},
~\ref{lem:K-Ev-6} and ~\ref{lem:K-Ev-7}.
  It follows that all these groups are the same.
  If $\sX_s$ is reduced, the same conclusion holds if we replace
  ${\Ev}^{\tr}_n \Br(X)$ by ${\Ev}^{\ur}_n \Br(X)$.
\end{proof}

Another corollary of \thmref{thm:Main-1} is the following property of the
evaluation filtration.

\begin{cor}\label{cor:K-Ev-11}
  Let $\sX$ be as in \thmref{thm:Main-1}. Then the evaluation filtration
  ${\Ev}_\bullet \Br(X)$ of $\Br(X)$ is exhaustive.
\end{cor}
\begin{proof}
  We know that ${\Fil}_\bullet \Br(X)$ is exhaustive and hence so is
  $\Fil'_\bullet \Br(X)$ because $\Fil_n \Br(X) \subset \Fil'_n \Br(X)$.
  We can now apply \thmref{thm:Main-1} to conclude.
\end{proof}

\begin{remk}\label{remk:K-Ev-5-extra}
  The reader can check from the proof of \lemref{lem:K-Ev-5} that
  the inclusion $\Fil'_n \Br(X) \subset {\Ev}_n \Br(X)$ for $n \ge 0$ holds even if
  $Y$ is not an snc divisor on $\sX$. That is, it holds for any $\sX$ which is
  regular, separated, finite type and faithfully flat over $S$.
  We expect that a suitable version of all of \thmref{thm:Main-1} holds true
  for any such $S$-scheme. 
  \end{remk}

\section{Relation with Bright--Newton filtration}\label{sec:BN**}
In this section, we show that the evaluation filtration ${\Ev}_\bullet \Br(X)$
introduced in this paper coincides with the filtration defined by
Bright--Newton \cite{Bright-Newton} when $\sX$ is smooth over $S$. In fact, we
establish this under the slightly weaker assumptions.
We keep the notation and assumptions of \S~\ref{sec:KEF}.
In particular, we continue to assume that
$\sX$ is a quasi-semistable scheme over $S$.

For an ffe ${k'}/k$, we let $e(k'/k)$ denote
the ramification index of $k'$ over $k$. We let $\fm' = (\pi')$ denote the maximal
ideal and $\F'$ the residue field of $\sO_{k'}$. We let
$S' = \Spec(\sO_{k'})$. We shall say that $k'$ is `good' if
$\sX' := \sX \times_{S} S'$ is regular. Note that if ${\sX}/S$ is smooth, then
every ffe ${k'}/k$ is good.

\begin{lem}\label{lem:good-ffe}
  Suppose that $\sX_s$ is reduced and let ${k'}/k$ be a good ffe.
  Then $\sX'$ is a quasi-semistable scheme over $S'$ with reduced special fiber.
  If the components of $\sX_s$ are
  furthermore geometrically integral, then so are the components of $\sX'_s$ and
  there is a one-to-one correspondence between the irreducible components of
  $\sX'_s$ and $\sX_s$.
\end{lem}
\begin{proof}
  For the first part, we only need to show that
  $\sX'_s$ is reduced and an snc divisor on $\sX'$.
  But this follows from \lemref{lem:SNC-0}(1) because
  $\sX'_s = \sX_s \times_{\Spec(\F)} \Spec(\F') \to \sX_s$ is a finite
  {\'e}tale map of $\F$-schemes under our assumption, and hence $\sX'_s$ is an snc
  $\F'$-scheme. The second part is clear.
\end{proof}

For $n \in \N$, an ffe ${k'}/k$ and $f \in \sX^o(\sO_{k'})$, we let
$B'(f, n)$ be the subset of $\sX^o(\sO_{k'})$ consisting of
$g \in \sX^o(\sO_{k'})$ such that $f$ and $g$ have the same image in
$\sX^o({\sO_{k'}}/{(\pi'^n)})$ under the restriction map
$\sX^o(\sO_{k'}) \to \sX^o({\sO_{k'}}/{(\pi'^n)})$.
Let $\pi'$ be a uniformizer of $k'$.
Given an element $\sA \in \Br(X)$ and an ffe ${k'}/k$, we get a function
$\ev_{\sA} \colon \sX^o(\sO_{k'}) \to {\Q}/{\Z}$ defined as follows.
If $f \colon \Spec(\sO_{k'}) \to \sX^o$ is a morphism of $S$-schemes,
it restricts to a morphism of $k$-schemes $f_\eta \colon \Spec(k') \to X$.
We let $\ev_{\sA}(f) = \inv_{k'}(f^*_\eta(\sA))$.

\begin{defn}\label{defn:BN-mod} 
  For $n \in \N_0$, we let
  \[
  {\Ev}^{\bn}_n \Br(X) = \{\sA \in \Br(X)| \ev_\sA \ \text{is constant on} \
  B'(f, e(k'/k)(n+1)) \hspace*{2cm}
  \]
  \[
  \hspace*{6cm}
  \text{ for all good ffe } {k'}/k \ \text{and all} \ f
  \in \sX^o(\sO_{k'})\};
  \]
  \[
 {\Ev}^{\bn}_{-1} \Br(X) = \{\sA \in \Br(X)| \ev_\sA \ \text{is constant on} \
 \sX^o(\sO_{k'})  \ \text{for all good ffe} \ {k'}/k \};
   \]
   \[
  {\Ev}^{\bn}_{-2} \Br(X) = \{\sA \in \Br(X)|  \ev_\sA \ \text{is zero on} \
  \sX^o(\sO_{k'})  \ \text{for all good ffe } \ {k'}/k \}.
    \]
\end{defn} 

It is clear from the above definition that if $\sX$ is smooth over $S$ with
geometrically integral special fiber (the case considered in \cite{Bright-Newton}),
then ${\Ev}^{\bn}_\bullet \Br(X)$
coincides with the filtration defined by Bright--Newton.

The following theorem explains the relation between the evaluation filtration
of this paper and the one considered in \cite{Bright-Newton}.

\begin{thm}\label{thm:K-Ev-11}
  Assume that $\sX_s$ is reduced and its irreducible components are geometrically
  integral. Then we have ${\Ev}^{\bn}_n \Br(X) = {\Ev}_n \Br(X)$ for every
  $n \ge -2$.  
\end{thm}
\begin{proof}
 We first assume $n \ge 0$ and show that ${\Ev}_n \Br(X) \subset
  {\Ev}^{\bn}_n \Br(X)$.  By \lemref{lem:K-Ev-4}, it suffices to show that
  $\Fil'_n \Br(X) \subset  {\Ev}^{\bn}_n \Br(X)$.
  We now fix an element $\chi \in \Fil'_n \Br(X)$.
  We let ${k'}/k$ be a good ffe, $f \in \sX^o(\sO_{k'})$ and
  $g \in  B'(f, e(k'/k)(n+1))$. 
We denote these maps (of $S$-schemes) by $f, g \colon S' \to \sX^o$. We let
$f_\eta, g_\eta \colon \Spec(k') \to X$ be the induced maps of $k$-schemes.
We let  $m = e(k'/k)(n+1)$ and let
$\iota_{n+1} \colon \Spec({\sO_{k'}}/(\pi'^{m})) \inj S'$ be the
inclusion. We are given that $f \circ \iota_{n+1} = g \circ \iota_{n+1}$.

We look at the commutative diagram
\begin{equation}\label{eqn:K-Ev-11-0}
  \xymatrix@C1.6pc{
    \Spec({\sO_{k'}}/(\pi'^{m})) \ar[rr]^-{\iota_{n+1}}
    \ar[drr]^-{f' \circ \iota_{n+1}}  \ar[ddrr]_-{\iota_{n+1}} & &
      S' \ar[dr]^-{f} \ar[d]_-{f'} & & \\
    & & \sX'^o \ar[r]^-{v'} \ar[d]_-{u'} & \sX^o \ar[d]^-{u} \\
    & & S' \ar[r]^-{v} & S,}
  \end{equation}
where $u$ and $u'$ are the projection maps, the lower right square is Cartesian and
$f'$ is the unique map such that $v' \circ f' = f$ and $u' \circ f' = \id_{S'}$. 

It follows from ~\eqref{eqn:K-Ev-11-0} that $f' \circ \iota_{n+1}$ is the unique map
such that $v' \circ f' \circ \iota_{n+1} = f \circ \iota_{n+1}$ and
$u' \circ f' \circ \iota_{n+1} = \iota_{n+1}$. We deduce from this uniqueness that
\begin{equation}\label{eqn:K-Ev-11-1}
  f' \circ \iota_{n+1} = g' \circ \iota_{n+1},
\end{equation}
where $g' \colon S' \to \sX'$ is the
section of $u'$ defined by $g \colon S' \to \sX$.
If we let $\ov{P} = f'(S')$ and $\ov{Q} = g'(S')$, we get two points
$P = \ov{P} \bigcap X'$ and $Q = \ov{Q} \bigcap X'$ in $X'^o_{\fin}$
such that $\ov{P} \times_{\sX'} mY'^o = \ov{Q} \times_{\sX'} mY'^o$.
In other words, $P, Q \in X'^o_{\fin}$ and $Q \in B(P, m)$.

We next note that the assumption $\chi \in  \Fil'_n \Br(X)$ implies by
\propref{prop:Matsuda-3} that $v^*_\eta(\chi) \in \Fil'_{m-1} \Br(X')$,
  where $v'_\eta \colon X' \to X$ is the restriction of $v'$ to the generic
  fibers. As $\Fil'_{m-1} \Br(X') \subset {\Ev}_{m-1} \Br(X')$ by
  Lemmas~\ref{lem:K-Ev-5} and ~\ref{lem:good-ffe},
  we get $\ev_{v^*_\eta(\chi)}(P) = \ev_{v^*_\eta(\chi)}(Q)$.
  But this is the same thing as saying that $\inv_{k'}(f^*_\eta(\chi)) =
  \inv_{k'}(g^*_\eta(\chi))$. Equivalently, $\ev_{\chi}(f) = \ev_{\chi}(g)$.
  This shows that $\chi \in {\Ev}^{\bn}_n \Br(X)$.

To show that ${\Ev}^{\bn}_n \Br(X) \subset {\Ev}_n \Br(X)$ when $n \ge 0$,
it suffices to show, using \thmref{thm:K-Ev-10}, that
  ${\Ev}^{\bn}_n \Br(X) \subset {\Ev}^{\ur}_n \Br(X)$.
  We let $\chi \in {\Ev}^{\bn}_n \Br(X)$ and let $P, Q \in X_{\ur}$ such that
  $Q \in B_{\ur}(P, n+1)$. If we let $P_0 = Q_0 = x$, then this condition
  implies that there is an $S$-isomorphism
  $\phi \colon \ov{P} \xrightarrow{\cong} \ov{Q}$
  which induces the identity map $\id_{k(x)} \colon \ov{P} \times_S \F
  \to \ov{Q} \times_S \F$ (cf. \cite[Chap.~I, Prop.~4.4]{Milne-EC}).
  In particular, it induces the identity map
  $\id \colon \ov{P} \times_{\sX} (n+1)Y^o \to \ov{Q} \times_{\sX} (n+1)Y^o$
  (cf. [Chap.~I, Thm.~3.23] of op. cit.).

We let $k' = k(P), S' = \Spec(\sO_{k'})$ and
  $\iota_n \colon S' \times_{\sX} nY^o \inj S'$.
  Since $P \in X_{\ur}$, we see that ${k'}/k$ is a good ffe.
  If we let $f \colon \ov{P} \inj \sX^o, \ f' \colon \ov{Q} \inj \sX^o$
  denote the inclusions and write $g = f' \circ \phi$, we get
  $f, g \in \sX^o(\sO_{k'})$ such that $f \circ \iota_{n+1} = g \circ \iota_{n+1}$.
  In particular, $g \in B'(f, n+1)$. Since $\chi \in {\Ev}^{\bn}_n \Br(X)$, it
  follows that $\ev_{\chi}(f) = \ev_{\chi}(g)$. But this is the same thing as saying
  that $\ev_{\chi}(P) = \ev_{\chi}(Q)$. This shows that
  $\chi \in {\Ev}^{\ur}_n \Br(X)$.

We next show ${\Ev}^{\bn}_{-1} \Br(X) = {\Ev}_{-1} \Br(X)$. To show the forward
inclusion, it suffices to show, using \thmref{thm:K-Ev-10}, that
${\Ev}^{\bn}_{-1} \Br(X) \subset {\Ev}^{\ur}_{-1} \Br(X)$.
We let $\chi \in {\Ev}^{\bn}_{-1} \Br(X)$ and let $P, Q \in X_{\ur}$ such that
  $Q \in B_{\ur}(P,0)$.  We need to show that $\ev_{\chi}(P) = \ev_{\chi}(Q)$.
To that end, we let $x = \esp(P)$ and $y = \esp(Q)$. We then we get 
an isomorphism of field extensions
$\phi_\F \colon {k(x)}/{\F} \xrightarrow{\cong} {k(y)}/{\F}$.
This implies, as shown above, that there is
  an isomorphism $\phi \colon \ov{P} \xrightarrow{\cong} \ov{Q}$ of $S$-schemes
  such that its reduction to the special fiber is $\phi_\F$.
We let $f \colon \ov{P} \inj \sX$ and $f' \colon \ov{Q} \inj \sX$ be the
  inclusions and let $g = f' \circ \phi$.
  Letting $k' = k(P), S' = \ov{P}$ and $\sX' = \sX \times _S S'$,
  and repeating the  argument in  $n \ge 0$ case above, we conclude that
$\ev_{\chi}(f) = \ev_\chi(g)$. Equivalently, $\ev_{\chi}(P) = \ev_{\chi}(Q)$.

To prove the backward inclusion, we let $\chi \in {\Ev}_{-1} \Br(X)$ and let
  $f, g \colon S' = \Spec(\sO_{k'}) \to \sX^o$ be morphisms of $S$-schemes, where
  ${k'}/k$ is a good ffe. We let $\sX' = \sX \times_S S'$ and let
$u \colon \sX' \to \sX$ be the base change map. Let
$v \colon \Spec(\F') \to \Spec(\F)$ and $w \colon Y' \to Y$ be the induced maps.
  By \lemref{lem:good-ffe}, $\sX'$ is a quasi-semistable $S'$-scheme.
  Moreover, $f$ and $g$ define unique points
  $P', Q' \in X'(k')$ such that $\esp(P'), \esp(Q') \in Y'^o(\F)$.
By \lemref{lem:K-Ev-8}, $\chi \in \Fil'_{-1} \Br(X)$.
In particular, $\chi \in \Fil_0 \Br(X)$ such that $\partial_X(\chi)$ lies in the
image of the map $g^* \colon H^1_\et(\F) \inj H^1_\et(Y^0)$, induced by the
projection $g \colon Y^o \to \Spec(\F)$. Note here that $\wt{g}^* = g^*$ because
$\sX_s$ is reduced.

We now look at the commutative diagram (cf. \propref{prop:Kato-Br-0})
\begin{equation}\label{eqn:K-Ev-11-2}
  \xymatrix@C1pc{
    \Fil_0 \Br(X) \ar[r]^-{\partial_X} \ar[d]_-{u^*} &
    H^1_\et(Y^o) \ar[d]^-{\wt w^*} &
    H^1_\et(\F) \ar[d]^-{\wt v^*} \ar[l]_-{\ {g}^*} \\
    \Fil_0 \Br(X') \ar[r]^-{\partial_{X'}} & H^1_\et(Y'^o) &
    H^1_\et(\F'), \ar[l]_-{\ {g'}^*}}
\end{equation}
where $g' \colon Y'^o \to \Spec(\F')$ is the structure map and
the vertical arrows are the base change maps. Note that $\wt w^*=e(k'/k) w^*$ and
$\wt v^*=e(k'/k) v^*$ by \lemref{lem:good-ffe} because $\sX_s$ is reduced. Since
$\wt{g'}^* = g'^*$ (again by \lemref{lem:good-ffe}), it follows from this diagram
that $u^*(\chi) \in \Fil_{-1} \Br(X') = {\Ev}_{-1} \Br(X')$.
In particular, $\ev_{u^*(\chi)}(P') = \ev_{u^*(\chi)}(Q')$. But this is same 
as saying that $\ev_{\chi}(f) = \ev_{\chi}(g)$. We conclude that
$\chi \in {\Ev}^{\bn}_{-1} \Br(X)$.

Finally, we show ${\Ev}^{\bn}_{-2} \Br(X) = {\Ev}_{-2} \Br(X)$.
To show the forward inclusion, it suffices to show, using
\thmref{thm:K-Ev-10}, that  ${\Ev}^{\bn}_{-2} \Br(X) \subset {\Ev}^{\ur}_{-2} \Br(X)$.
We now let $\chi \in {\Ev}^{\bn}_{-2} \Br(X)$ and let $P \in X^o_{\ur}$.
Let $\ov{P}$ denote the closure of $P$ in $\sX$ and let
$f \colon \ov{P} \inj \sX^o$ be the inclusion.
 Letting $k' = k(P)$, it follows that ${k'}/k$ is a good ffe
 and $\ov{P} = \Spec(\sO_{k'})$. Since $\chi  \in {\Ev}^{\bn}_{-2} \Br(X)$, we have
 that $\ev_{\chi}(f) = 0$. But this is equivalent to $\ev_{\chi}(P)$ being zero.
  
To show the backward inclusion, we let $\chi \in \Br(\sX)$ and let
  ${k'}/k$ be a good ffe. Let $f \colon S' = \Spec(\sO_{k'}) \to \sX^o$ be an
  $S$-morphism and let $\sX' = \sX \times_S S'$. Then $\sX'$ a quasi-semistable
  $S'$-scheme by \lemref{lem:good-ffe}. If we let $u \colon \sX' \to \sX$ denote the
  base change map, then $u^{-1}(\sX^o) = \sX'^o$, as one easily checks.
  If let $\chi' = u^*(\chi) \in \Br(\sX')$ and apply
  \propref{prop:Kato-Br-0} to $(\sX', Y') \to (\sX, Y)$,
  we get $\chi' \in \Fil_0 \Br(X')$ and
  $\partial_{X'}(\chi') = 0$. In particular, $\ev_{\chi'}(P) =
  \partial_{X'}(\chi')|_{k(x)} = 0$, where $P = \Spec(k') \in X'$ and $k(x) =
  \esp(P) \in Y'^o$. But this is same as saying that $\ev_{\chi}(f) = 0$.
  This shows that $\Br(\sX) \subset {\Ev}^{\bn}_{-2} \Br(X)$. The
  proof of the theorem is now complete.
\end{proof}

\section{The Brauer--Manin pairing in higher dimension}\label{sec:BMP}
In the rest of this paper, we give applications of \thmref{thm:Main-1}
to some open problems on Brauer groups and Chow groups of
zero-cycles on projective varieties over local fields.
 The goal of the next two sections is to prove \thmref{thm:Main-3} which
 extends the main result of \cite{Saito-Sato-ENS} to positive characteristics.
 We begin by recalling the Brauer--Manin pairing for smooth projective
 varieties over local fields.

 For any scheme $Z$, let $\CH_0(Z)$ denote the Chow group of zero-cycles on $Z$
 (see \cite[\S~1.6]{Fulton}). Now suppose that $Z$ is a complete $F$-scheme, where
 $F$ is a field. Let $A_0(Z)$ denote the kernel of the degree map
 $\deg \colon \CH_0(Z) \to \Z$.
 As recalled in \cite[Prop.~6.4.2]{CTS}, the evaluation map for Brauer classes
 (cf.\eqref{eqn:Ev-map}) induces in this setting a bilinear pairing
\begin{equation}\label{eqn:BMP-0}
 \< , \> \colon \CH_0(Z) \times \Br(Z) \to \Br(F),
\end{equation}
which is functorial with respect to proper morphisms of complete $F$-schemes. We
shall refer to this pairing as the Brauer--Manin pairing for $Z$.

For the next two sections, we shall work in the following setting.
We let $k$ be an hdvf of characteristic $p$ with
excellent ring of integers $\sO_k$, uniformizer $\pi$ and finite residue field $\F$.
We let $S =\Spec(\sO_k)$ and $\eta = \Spec(k)$. We let $\sX$ be a semistable
(cf. \S~\ref{sec:Set-up})
projective $S$-scheme of absolute dimension $d_\sX$ with generic fiber $X$ and
special fiber $\sX_s$. We let $Y = (\sX_s)_\red$ be the reduced special fiber.
We fix a closed embedding $\sX \inj \P^N_R$. Let $f \colon \sX \to S$ be the
structure map. We shall denote the map $X \to \Spec(k)$ also by $f$. We let
$g \colon Y \to \Spec(\F)$ denote the structure map for $Y$.
We keep the notation of \S~\ref{sec:KEF}.

Applying ~\eqref{eqn:BMP-0} to $X$ and using
Example~\ref{exm:Kato-basic-2},
we get a bilinear
pairing $\Br(X) \times \CH_0(X) \to {\Q}/{\Z}$. Using the functoriality of the
Brauer--Manin pairing with respect to the structure map $X \to \Spec(k)$, this
induces a bilinear pairing $ A_0(X) \times {\Br(X)}/{\Br_0(X)} \to {\Q}/{\Z}$
(recall that $\Br_0(X) = f^*(\Br(k))$).
If $x \in X_{(0)}$ and $\sA \in \Br(\sX)$, then
$\<\sA, [x]\> = \inv_{k(x)}(\sA|_{\sO_{k(x)}}) = 0$ since
$\Br(\sO_{k(x)}) = 0$. It follows that the Brauer--Manin pairing for $X$
canonically induces bilinear pairings
\begin{equation}\label{eqn:BMP-1}
   \CH_0(X) \times \frac{\Br(X)}{\Br(\sX)}  \xrightarrow{ \< , \>} {\Q}/{\Z}; \ \ \
  A_0(X) \times \frac{\Br(X)}{\Br(\sX) + \Br_0(X)} \xrightarrow{ \< , \>} {\Q}/{\Z}.
\end{equation}

It is known (cf. \cite{Parimala-Suresh}, \cite{CTSaito}) that these pairings in
general have non-trivial left and right kernels. 
Our main result is concerned with the study of their right kernels.
The characteristic zero version of this result was proven by Saito--Sato
\cite{Saito-Sato-ENS}.
In this section, we shall establish some intermediate steps which go into
the proof of this result. In particular, we shall prove an injectivity
result for ${\Br(X)}/{\Fil'_{-1} \Br(X)}$ when we pass to admissible
hypersurfaces.

\subsection{Restriction of ${\Fil'_0\Br(X)}/{\Fil'_{-1} \Br(X)}$
  to hypersurfaces}\label{sec:Prelim-0}
We assume in this subsection that $d_\sX \ge 3$. We begin with some lemmas.

\begin{lem}\label{lem:SC-change-2}
Let $\{(E_1, \sE_1), \ldots , (E_t, \sE_t)\}$ be a set of ordered pairs,
where each $(E_i, \sE_i)$ consists of a connected regular closed subscheme
$E_i$ of $\sX_s$ of dimension at least two and $\sE_i$ is
a coherent locally free sheaf on $E_i$. Then we can find an integer $n \gg 0$
  such that for all hypersurface sections $\sX' \subset \sX$ of degree at least
  $n$, we have $H^i_\zar(E_j, \sE_j(-\sX')) = 0$ for all $i \in J^0_1$ and
  $j \in J_t$.
\end{lem}
\begin{proof}
See the proof of \cite[Lem.~5.1]{GK-Jussieu}.
\end{proof}

\begin{lem}\label{lem:SC-change-3}
  Let $\{D_1, \ldots , D_t\}$ be a collection of
  effective divisors on $\sX$. Let $\Sigma$ be a finite set of closed points in
  $Y^o$. Then there exists an integer $n \gg 0$ such that for every $d \ge n$, 
  there are infinitely many admissible
  hypersurface sections $\sX' \subset \sX$ of degree $d$ containing $\Sigma$
  such that the restriction map
  $\alpha^q_{i,j} \colon H^0_\zar(Y_i, \Omega^q_{Y_i}(\log F_i)(D_j)) \to
  H^0_\zar(Y'_i, \Omega^q_{Y'_i}(\log F'_i)(D_j))$ is injective for all
  $i \in J_r$, $j \in J_t$ and $q \in J^0_1$ (resp. $q \in J^0_2$) if
  $d_\sX \ge 3$ (resp. $d_\sX \ge 4$).
\end{lem}
\begin{proof}
By the second exact sequence of \lemref{lem:Kahler-seq}, proving the injectivity of
  $\alpha^q_{i,j}$ is equivalent
  to showing that $H^0_\zar(Y_i, \sF^q_{i,j}) = 0$, where
  $\sF^q_{i,j} = \Omega^q_{Y_i}(\log (F_i + Y'_i))(D_{j} - Y'_i)$.
  To prove the latter statement,
  we  use the first exact sequence of
  \lemref{lem:Kahler-seq} which yields an exact sequence
\begin{equation}\label{eqn:SC-change-3-0}
  0 \to \Omega^q_{Y_i}(\log F_i)(D_{j} - Y'_i) \to \sF^q_{i,j} \xrightarrow{\res}
  \Omega^{q-1}_{Y'_i}(\log F'_i)(D'_{j} - Y'_i) \to 0.
\end{equation}

Let $\iota_i \colon Y_i \inj \sX$ be the inclusion for $i \in J_r$.
By applying \lemref{lem:SC-change-2} to the set of pairs
$\{(Y_i, \Omega^q_{Y_i}(\log F_i)(D_{j}))| \ i \in J_r, j \in J_t, q \in J^0_2\}$,
we get that there exists $n \gg 0$
such that for all hypersurface sections $\sX' \subset \sX$ of degree at least
$n$, we have $H^m_\zar(Y_i, \Omega^q_{Y_i}(\log F_i)(D_{j} - Y'_i) =0$ for all
$m \in J^0_1$, $i \in J_r$ and $j \in J_t$. By choosing
$n$ even larger and applying \corref{cor:Coh-vanishing-1} to the set of triples
$\{(Y_i, F_i, \sD_{i,j})| i \in J_r, j \in J_t\}$ (where
$\sD_{i,j}$ is the line bundle $\iota^*_i(\sL(D_j)$), we get that
for all hypersurface sections $\sX' \subset \sX$ of degree at least
$n$, we have $H^m_\zar(Y'_i, \Omega^{q-1}_{Y'_i}(\log F'_i)(D'_{j} - Y'_i)) =0$ for all
$m \in J^0_{ d_\sX-q-2}$, $i \in J_r$ and $j \in J_t$.

If either $q \in J^0_1$ and $d_\sX \ge 3$ or $q \in J^0_2$ and $d_\sX \ge 4$,
then $d_\sX-q-2 \ge 0$, and 
the long exact cohomology sequence corresponding to ~\eqref{eqn:SC-change-3-0}
shows that $H^0_\zar(Y_i, \sF^q_{i,j}) =0$
for all hypersurface sections $\sX' \subset \sX$ of degree at least
$n$, all $i \in J_r$, $j \in J_t$. 
We combine this with \thmref{thm:Bertini-dvr} to complete the
proof of the lemma.
\end{proof}

Since $\sX_s \in \Div(\sX)$, we can write 
$\sX_s = \stackrel{r}{\underset{i =1}\sum} m_i Y_i$ with $m_i \ge 1$ for each $i$.
We let $\Phi = \left\{\stackrel{r}{\underset{i =1}\sum} n_i Y_i| \  i \in J_r \
\mbox{and} \ n_i \in J_{m_i-1}\right\} \subset \Div(\sX)$. Let 
$H^1_\et(Y^N)_c$ denote the image of the map $\wt{g}^* \colon
H^1_\et(\F) \to H^1_\et(Y^N)$, where $Z^N$ denotes the normalization of a reduced
scheme $Z$. We define $H^1_\et(Y^o)_c$ similarly
(cf. ~\eqref{eqn:Kato-Br-0*} and \lemref{lem:K-Ev-7}).

\begin{cor}\label{cor:SC-change-4}
  Let $\Sigma$ be a finite set of closed points in $Y^o$. Then there
  exists an integer $n \gg 0$ such that for every $d \ge n$, there are infinitely
  many admissible hypersurface sections $\sX' \subset \sX$ of degree $d$ containing
  $\Sigma$ such that letting $Y'$ denote the reduced special fiber of $\sX'$,
  we have the following.
  \begin{enumerate}
  \item
  The restriction map
  \[
  H^0_\zar(Y_i, \Omega^q_{Y_i}\log(F_i)(D)) \to
  H^0_\zar(Y'_i, \Omega^q_{Y'_i}\log(F'_i)(D))
  \]
  is injective for every $D \in \Phi$, $i \in J_r$ and $q \in J^0_ 1$
  (resp. $q \in J^0_2$) if $d_\sX \ge 3$ (resp. $d_\sX \ge 4$).
 \item
  The restriction map
  \[
  \frac{H^1_\et(Y^N)}{H^1_\et(Y^N)_c} \{p\} \bigoplus
    \frac{H^1_\et(Y^o)}{H^1_\et(Y^o)_c}\{p'\} \to 
\frac{H^1_\et(Y')}{H^1_\et(Y'^N)_c} \{p\} \bigoplus
  \frac{H^1_\et(Y'^o)}{H^1_\et(Y'^o)_c}\{p'\}
  \]
  is injective.
  \end{enumerate}
  \end{cor}
  \begin{proof}
The first part of the corollary is a direct consequence of
    \lemref{lem:SC-change-3}. To prove the second part, we first note that the
    injectivity of the map $H^1_\et(Y^N) \to H^1_\et(Y'^N)$ is classical
    (cf. \cite[Exp.~X, Cor.~2.5]{SGA2}) for any admissible
    hypersurface section $\sX'$ of
    $\sX$. One can also deduce it  directly from \cite[Thm.~12.6, 12.8]{KM-1}.
    The injectivity of the map $H^1_\et(Y^o)\{p'\} \to H^1_\et(Y'^o)\{p'\}$ for any
    admissible hypersurface section $\sX'$ of $\sX$ follows from
    \cite[Thm.~12.8]{KM-1} (see also \cite[Thm.~1.1]{Esnault-Kindler}).

We now look at the diagram
    \begin{equation}\label{eqn:SC-change-4-0}
      \xymatrix@C1.2pc{
        H^1_\et(\F)\{p\} \ar[r]^-{\wt{g}^*} \ar[d]_-{\id} &
        H^1_\et(Y^N)\{p\} \ar[r] \ar[d] &
        \frac{H^1_\et(Y^N)}{H^1_\et(Y^N)_c} \{p\} \ar[d] \ar[r] & 0 \\
       H^1_\et(\F)\{p\} \ar[r]^-{\wt{g'}^*} &  H^1_\et(Y'^N)\{p\} \ar[r] &
       \frac{H^1_\et(Y'^N)}{H^1_\et(Y'^N)_c} \{p\}  \ar[r] & 0,}
    \end{equation}
    where $g' \colon Y' \to \Spec(\F)$ is the structure map and the middle vertical
    arrow is the pull-back map.
    The left (and hence the right) square commutes because $\sX'$ is
    admissible (cf. Definition~\ref{defn:good-sub}).
    Since the rows are exact, it follows that the right vertical arrow is injective.
    An identical argument shows that $\frac{H^1_\et(Y^o)}{H^1_\et(Y^o)_c}\{p'\} \to
    \frac{H^1_\et(Y'^o)}{H^1_\et(Y'^o)_c}\{p'\}$ is injective.
\end{proof}

\begin{lem}\label{lem:SC-change-5}
   Let $\sX' \subset \sX$ be any hypersurface section
    for which item (2) of \corref{cor:SC-change-4} holds. Then the
    restriction map
    \[
    \frac{\Fil'_0 \Br(X)}{\Fil'_{-1} \Br(X)} \to
    \frac{\Fil'_0 \Br(X')}{\Fil'_{-1} \Br(X')}
  \]
  is injective.
  \end{lem}
  \begin{proof}
We look at the diagram
     \begin{equation}\label{eqn:SC-change-5-0}
      \xymatrix@C1.2pc{
        \frac{\Fil'_0 \Br(X)}{\Fil'_{-1} \Br(X)}\{p\} \ar[r]^-{\partial_X} \ar[d] &
        \frac{H^1_\et(Y^N)}{H^1_\et(Y^N)_c} \{p\} \ar[d]  \\
       \frac{\Fil'_0 \Br(X')}{\Fil'_{-1} \Br(X')}\{p\} \ar[r]^-{\partial_{X'}}  &
    \frac{H^1_\et(Y'^N)}{H^1_\et(Y'^N)_c} \{p\}}
     \end{equation}
     which clearly commutes by the functoriality of $\Fil_\bullet \Br(X)$ and
     \propref{prop:Kato-Br-0} (recall from \corref{cor:Matsuda-2} that
     $\Fil'_0 \Br(X) = \Fil_0 \Br(X)$). The horizontal arrows in the diagram are
     injective by the definition of $\Fil'_{-1} \Br(X)$
     (cf. ~\eqref{eqn:Fil-neg}).
We now apply \corref{cor:SC-change-4} to conclude that the left vertical
     arrow in ~\eqref{eqn:SC-change-5-0} is injective.
 A similar argument shows that $\frac{\Fil'_0 \Br(X)}{\Fil'_{-1} \Br(X)}\{p'\}
     \to \frac{\Fil'_0 \Br(X')}{\Fil'_{-1} \Br(X')}\{p'\}$ is injective. This
     completes the proof.
\end{proof}

  In the proofs of \corref{cor:SC-change-4} and \lemref{lem:SC-change-5}
  (especially, in ~\eqref{eqn:SC-change-4-0} and ~\eqref{eqn:SC-change-5-0}), we
  implicitly used the following elementary group theory lemma.

  \begin{lem}\label{lem:Elem-lem}
    Let $\phi \colon A \to B$ be a homomorphism of torsion abelian groups.
    Then the sequences
    \[
    A\{p\} \xrightarrow{\phi} B\{p\} \to \coker(\phi)\{p\} \to 0; \ \
    A\{p'\} \xrightarrow{\phi} B\{p'\} \to \coker(\phi)\{p'\} \to 0
    \]
    are exact.
  \end{lem}
  \begin{proof}
    This is easily checked using that ${\rm Tor}^1({\Q_\ell}/{\Z_\ell}, -)$ is a right
      exact functor and ${\rm Tor}^2({\Q_\ell}/{\Z_\ell}, -) = 0$ in the category of
        torsion abelian groups for every prime $\ell$.
        \end{proof}

\subsection{Restriction of ${\Br(X)}/{\Fil'_{-1} \Br(X)}$ to hypersurfaces}
\label{sec:Prelim-1}
We assume that $d_\sX \ge 3$.
  As each $Y_i$ is projective over $\F$ which is finite,
  ${\underset{D \in \Phi, i \le r, q \le 2}\bigcup}
  H^0_\zar(Y_i, \Omega^q_{Y_i}(\log F_i)(D))$ is a finite set.
  We can therefore find a dense open $U_i \subset Y^o_i$ for each $i$ 
  such that the restriction map
  \begin{equation}\label{eqn:SC-change-5-1}
    H^0_\zar(Y_i, \Omega^q_{Y_i}(\log F_i)(D)) \to \Omega^q_{Y_i}(\log F_i)(D)
    \otimes_{\sO_{Y_i}} k(x)
  \end{equation}
  is injective for all $x \in U_i$, $D \in \Phi$, $i \in J_r$ and $q \in J^0_2$.

For each $i \in J_r$, we fix a closed point $x_i \in U_i$ and
  apply \lemref{lem:Type-2-9} to find admissible hypersurface sections
  $\sX_i, \sW_i$ of $\sX$ of degrees at least $n$ (where
  $n$ is given by \corref{cor:SC-change-4}) containing $x_i$
  and elements $t_i, \pi_i \in \sO_{\sX,x_i}$
such that $Y_i$ is locally defined at $x_i$ by $\pi_i$ and
$\sX_i$ (resp. $\sW_i$) is defined by $t_i$ (resp. $t_i + \pi_i$).
Let $X_i$ (resp. $W_i$) denote the generic fiber of $\sX_i$
  (resp. $\sW_i$).
  Let $\sX' \subset \sX$ be any hypersurface section
  for which items (1) and (2) of \corref{cor:SC-change-4} hold.
  We now have the following key lemma.

\begin{lem}\label{lem:SC-change-6}
The restriction map
    \begin{equation}\label{eqn:SC-change-6-0}
    \frac{\Br(X)}{\Fil'_{0} \Br(X)} \to \
    \frac{\Br(X')}{\Fil'_{0} \Br(X')} \bigoplus
    \left(\stackrel{r}{\underset{i =1}\bigoplus}
    \frac{\Br(X_i)}{\Fil'_{0} \Br(X_i)}\right) \bigoplus
    \left(\stackrel{r}{\underset{i =1}\bigoplus}
    \frac{\Br(W_i)}{\Fil'_{0} \Br(W_i)}\right)
  \end{equation}
  is injective.
  \end{lem}
 \begin{proof}
 Since $\Br(X)\{p'\} \subset \Fil_{0} \Br(X) = \Fil'_{0} \Br(X)$
  (cf. \cite[Cor.~2.5]{Kato-89}), the lemma is actually about the
  $p$-primary torsion Brauer classes. So let $\chi \in \Br(X)\{p\}$ such that
  $\chi \notin \Fil_0 \Br(X)$.  Let $m \ge 1$ be the smallest
  positive integer such that $\chi \in \Br(X)[p^m]$. We shall prove by induction
  on $m$ that $\chi|_Z$ does not lie in $\Fil_0 \Br(Z)$ for some
  $Z \in \sS$, where we let $\sS = \{X', X_i, W_i| i \in J_r\}$.
  This will prove the lemma. 

By \lemref{lem:Br-fil-0}, we choose $\omega \in H^2_{p^m}(X)$ such that
  $\kappa_X(\omega) = \chi$. We consider $\omega$ as an element of
  $H^2(X)$ via the canonical map $H^2_{p^m}(X) \to H^2(X)$
(note that the image of $\omega$ in $H^2(X)$ is not zero because its image
$\chi$ in $\Br(X)$ is not zero)
  and let $n_i = \Sw_{Y_i}(\omega)$. We let
  $D_\chi = \stackrel{r}{\underset{i =1}\sum} n_i Y_i$ be the Swan divisor of
  $\omega$ so that $\omega \in \Fil_{D_\chi} H^2_{p^m}(X)$.
   We let $\ov{\omega}$ (resp. $\ov{\chi}$) denote the image of $\omega$
   (resp. $\chi$) in  $\frac{\Fil_{D_\chi} H^2_{p^m}(X)}{\Fil_{0} H^2_{p^m}(X)}$
   (resp. $\frac{\Fil_{D_\chi} \Br(X)[p^m]}{\Fil_{0} \Br(X)[p^m]}$).
   Our assumption implies that $\ov{\omega} \neq 0$.
    We now look at the map
  $R_* \colon \frac{\Fil_{D_\chi} H^2_{p^m}(X)}{\Fil_{0} H^2_{p^m}(X)} \to
    \frac{\Fil_{{D_\chi}/p} H^2_{p^{m-1}}(X)}{\Fil_{0} H^2_{p^{m-1}}(X)}$
(cf. ~\eqref{eqn:Rest**}).

    \vskip.2cm

{\bf{Step~1:}} As the first step in the proof, we shall show that
    $\chi|_Z \notin \Fil_0 \Br(Z)$
    for some $Z \in \sS$ if we are given that $R_*(\ov{\omega}) = 0$.
    
To prove this, we write $n_i = q_i m_i + s_i$ with
$s_i \in J^0_{m_i-1}$ for $i \in J_r$, where recall that
$\sX_s = \stackrel{r}{\underset{i =1}\sum} m_i Y_i$.
    We let $q = {\rm max}\{q_1, \ldots , q_r\}$.
We can assume (after a permutation) that $q = q_1$. We let
  $D^o_{\omega} = \stackrel{r}{\underset{i =1}\sum} s_i Y_i$ and
 $D^\dagger_\omega = \stackrel{r}{\underset{i =1}\sum} (q_1m_i + s_i) Y_i =
q_1 \sX_s + D^o_{\omega}$. We then have that
$D^\dagger_\omega$ is $(\omega,1)$-admissible and $D^o_{\omega} \in \Phi$.
In particular, $\omega \in \Fil_{D^\dagger_\omega} H^2_{p^m}(X)$.
We let $D^{\heartsuit}_\omega = D^\dagger_\omega - Y_1$.
We consider $\omega$ as an element of $\Fil_{D^\dagger_\omega} H^2(X)$
and denote its image
in  $\frac{\Fil_{D^\dagger_\omega} H^2(X)}{\Fil_{D^{\heartsuit}_\omega} H^2(X)}$ by
$\ov{\omega}$.
 We now divide the proof of Step~1 into several cases and subcases.

\vskip.2cm

     {\bf{Case~1:}} Suppose that $q_1 \ge 1$.

 In this case, we have $n_1 \ge 1$.
  In particular, $\omega \notin \Fil_{D^{\heartsuit}_\omega} H^2(X)$
  (as $\Sw_{Y_1}(\omega) > n_1 -1 \ge 0$). That is, 
  $\ov{\omega} \in \frac{\Fil_{D^\dagger_\omega} H^2(X)}{\Fil_{D^{\heartsuit}_\omega} H^2(X)}$
  is not zero. 

\vskip.2cm

{\bf{Subcase~1:}} Assume in Case~1 that $\omega$ has type I at $Y_1$.

  In this case,  $\ov{\omega}$ does not die under the composite map
  (cf. \lemref{lem:RSW-spl})
  \[
  \frac{\Fil_{D^\dagger_\omega} H^2(X)}{\Fil_{D^{\heartsuit}_\omega} H^2(X)} \xrightarrow{\Rsw}
  H^0_\zar(Y_1, \Omega^2_{\sX}(\log Y)(D^\dagger_\omega) \otimes_{\sO_{\sX}} \sO_{Y_1})
  \xrightarrow{\nu^0_1} H^0_\zar(Y_1, \Omega^1_{Y_1}(\log F_1)(D^\dagger_\omega)).
  \]

We now look at the diagram
\begin{equation}\label{eqn:SC-change-6-1}
  \xymatrix@C1.2pc{
    H^0_\zar(Y_1, \Omega^1_{Y_1}(\log F_1)(D^\dagger_\omega))
    \ar[r]^-{\theta_{Y_1}}_-{\cong} \ar[d] &
    H^0_\zar(Y_1, \Omega^1_{Y_1}(\log F_1)(D^o_{\omega})) \ar[d] \\
    H^0_\zar(Y'_1, \Omega^1_{Y'_1}(\log F'_{1})(D^\dagger_\omega))
    \ar[r]^-{\theta_{Y'_1}}_-{\cong}  &
    H^0_\zar(Y'_1, \Omega^1_{Y'_1}(\log F'_1)(D^o_{\omega})),}
\end{equation}
where the horizontal arrows are induced by the canonical isomorphism
$\theta_{\sX} \colon \sO_{\sX} \xrightarrow{\cong} \sO_{\sX}(q_1\sX_s)$ which sends
the unit element to $\pi^{-q_1}$. It is then clear that this diagram is commutative
and its horizontal arrows are bijective.
As the right vertical arrow is injective by \corref{cor:SC-change-4},
it follows that the left vertical arrow is also injective.

By \corref{cor:RSW-gen-1}, we have a commutative diagram 
\begin{equation}\label{eqn:SC-change-6-2}
  \xymatrix@C.8pc{
    \frac{\Fil_{D^\dagger_\omega} H^2(X)}{\Fil_{D^{\heartsuit}_\omega} H^2(X)} \ar[r]^-{\Rsw}
    \ar[d] &
H^0_\zar(Y_1, \Omega^2_{\sX}(\log Y)(D^\dagger_\omega) \otimes_{\sO_{\sX}} \sO_{Y_1})
\ar[r]^-{\nu^0_1} \ar[d] & H^0_\zar(Y_1, \Omega^1_{Y_1}(\log F_1)(D^\dagger_\omega))
\ar[d] \\
\frac{\Fil_{D'_\omega} H^2(X')}{\Fil_{D'^{\heartsuit}_\omega} H^2(X')} \ar[r]^-{\Rsw} &
H^0_\zar(Y'_1, \Omega^2_{\sX'}(\log (Y'))(D'_\omega) \otimes_{\sO_{\sX'}} \sO_{Y'_1})
\ar[r]^-{\nu^0_1}  & H^0_\zar(Y'_1, \Omega^1_{Y'_1}(\log F'_1)(D'_\omega)),}
\end{equation}
where $D'_\omega = D^\dagger_\omega|_{\sX'}$ and $D'^{\heartsuit}_\omega = D'_\omega - Y'_1$.
We conclude from this diagram that
$\ov{\omega}$ does not die in
$\frac{\Fil_{D'_\omega} H^2(X')}{\Fil_{D'^{\heartsuit}_\omega} H^2(X')}$.
Since $\Fil_0 H^2(X') \subset \Fil_{D'^{\heartsuit}_\omega} H^2(X')$
(note that $D'^{\heartsuit}_\omega$ is effective on $\sX'$), it follows that
  $\omega$ does not die in $\frac{\Fil_{D'_\omega} H^2(X')}{\Fil_0H^2(X')}$.
  In particular, it does not die in $\frac{H^2(X')}{\Fil_0H^2(X')}$.
Equivalently, $\chi|_{X'} \notin \Fil_0 \Br(X')$.

\vskip.2cm

{\bf{Subcase~2:}} Assume in Case~1 that $\omega$ has type II at $Y_1$ and
$d_{\sX} \ge 4$.

In this case, $\Rsw(\ov{\omega})$ lies in
$H^0_\zar(Y_1, \Omega^2_{Y_1}(\log F_1)(D^\dagger_\omega))
\subset
H^0_\zar(Y_1, \Omega^2_{\sX}(\log Y)(D^\dagger_\omega) \otimes_{\sO_{\sX}} \sO_{Y_1})$
  (cf. \lemref{lem:RSW-spl}). 
 The same argument as in Subcase~1 shows that the map
  $H^0_\zar(Y_1, \Omega^2_{Y_1}(\log F_1)(D^\dagger_\omega)) \to
 H^0_\zar(Y'_1, \Omega^2_{Y'_1}(\log F'_{1})(D'_\omega))$ is injective.
 Together with the injectivity of the horizontal arrows in the commutative diagram
\begin{equation}\label{eqn:SC-change-6-3}
  \xymatrix@C.8pc{
    H^0_\zar(Y_1, \Omega^2_{Y_1}(\log F_1)(D^\dagger_\omega)) \ar[r]^-{\nu'_1} \ar[d] &
    H^0_\zar(Y_1, \Omega^2_{\sX}(\log Y)(D^\dagger_\omega) \otimes_{\sO_{\sX}} \sO_{Y_1})
    \ar[d] \\
     H^0_\zar(Y'_1, \Omega^2_{Y'_1}(\log F'_1)(D'_\omega)) \ar[r]^-{\nu'_1} &
     H^0_\zar(Y'_1, \Omega^2_{\sX'}(\log Y')(D'_\omega) \otimes_{\sO_{\sX'}} \sO_{Y'_1})}
  \end{equation}
  and \corref{cor:RSW-gen-1}, this implies that
  $\Rsw(\ov{\omega}|_{X'}) \neq 0$ in
  $H^0_\zar(Y'_1, \Omega^2_{\sX'}(\log Y')(D'_\omega) \otimes_{\sO_{\sX'}} \sO_{Y'_1})$.
In particular, $\ov{\omega}|_{X'} \neq 0$ in
$\frac{\Fil_{D'_\omega} H^2(X')}{\Fil_{D'^{\heartsuit}_\omega} H^2(X')}$. As in Subcase~1,
this implies that $\omega$ does not die in $\frac{H^2(X')}{\Fil_0H^2(X')}$.
 Equivalently, $\chi|_{X'} \notin \Fil_0 \Br(X')$. 

 \vskip.2cm

{\bf{Subcase~3:}} Assume in Case~1 that $\omega$ has type II at $Y_1$ and
 $d_\sX = 3$.

In this case, we have $n_1 \ge 2$ by \cite[Thm.~3.2(3)]{Kato-89} and
  $\Rsw(\ov{\omega})$ lies in the subgroup
  $H^0_\zar(Y_1, \Omega^2_{Y_1}(\log F_1)(D^\dagger_\omega))$  of
  $H^0_\zar(Y_1, \Omega^2_{\sX}(\log Y)(D^\dagger_\omega) \otimes_{\sO_{\sX}} \sO_{Y_1})$,
  as in Subcase~2. In particular, $\Rsw(\ov{\omega})$ defines a non-zero element of
  $H^0_\zar(Y_1, \Omega^2_{Y_1}(\log F_1)(D^o_\omega))$ via the isomorphism
  $\theta_{Y_1} \colon H^0_\zar(Y_1, \Omega^2_{Y_1}(\log F_1)(D^\dagger_\omega))
  \xrightarrow{\cong} H^0_\zar(Y_1, \Omega^2_{Y_1}(\log F_1)(D^o_\omega))$, induced by
  $\theta_{\sX}$.

Using the commutative diagram
  \begin{equation}\label{eqn:SC-change-6-4}
  \xymatrix@C.8pc{
    H^0_\zar(Y_1, \Omega^2_{Y_1}(\log F_1)(D^\dagger_\omega)) \ar[r]
    \ar[d]_-{\theta_{Y_1}} &
    \Omega^2_{Y_1}(\log F_1)(D^\dagger_\omega) \otimes_{\sO_{Y_1}} k(x_1)
    \ar[d]^-{\theta_{Y_1}} \\
    H^0_\zar(Y_1, \Omega^2_{Y_1}(\log F_1)(D^o_\omega)) \ar[r] &
    \Omega^2_{Y_1}(\log F_1)(D^o_\omega) \otimes_{\sO_{Y_1}} k(x_1)}
  \end{equation}
  (where the horizontal arrows are the canonical restrictions)
  and ~\eqref{eqn:SC-change-5-1}, it follows that the image of $\Rsw(\ov{\omega})$
in $\Omega^2_{Y_1}(\log F_1)(D^\dagger_\omega) \otimes_{\sO_{Y_1}} k(x_1)$ is not zero.

By ~\corref{cor:Type-2-9-key} (for $p\neq 2$ or $n_1 \neq 2$)
or by \corref{cor:Type-2-12-6} (for $p=n_1=2$, recalling that
$R_*(\ov{\omega}) = 0$), we get that
either we have $\Sw_{\sX_1 \cap Y_1}(\omega|_{X_1}) \ge n_1-1 \ge 1$ or we have
  $\Sw_{\sW_1 \cap Y_1}(\omega|_{W_1}) \ge n_1-1 \ge 1$.   In particular,
$\omega|_{Z}$ does not die in $\frac{H^2(Z)}{\Fil_{0} H^2(Z)}$ for
some $Z \in \{X_1, W_1\}$. Equivalently,  $\chi|_{Z}$ does not lie in
$\Fil_0 \Br(Z)$ for some  $Z \in \{X_1, W_1\}$.

\vskip.2cm

{\bf{Case~2:}} Suppose that $q_1 = 0$.
In this case, we have $n_i=s_i$  for all $i \in J_r$. Since
$\omega \notin \Fil_{0} H^2(X)$, we must have $n_i \ge 1$ for some $i \in J_r$.
We may assume that $n_1 \ge 1$ after a permutation. Now, the proof is identical
to the proof of Case~1, the only difference being that we take
  $\theta_\sX$ to be the identity map of $\sO_\sX$.

\vskip.2cm

{\bf{Step~2:}} We shall now show in general by induction on $m \ge 1$ that
$\chi|_Z \notin \Fil_0 \Br(Z)$ for some $Z \in \sS$.

If $m =1$, then $R_*(\ov{\omega})$ is necessarily zero and we are done by
Step~1. If $m \ge 2$ and $R_*(\ov{\omega}) =0$, we are again done by
Step~1. We now assume that $m \ge 2$ and $R_*(\ov{\omega}) \neq 0$.

In this case, we must have
    $R_*(\omega) \neq 0$ in $\Fil_{{D_\chi}/p} H^2_{p^{m-1}}(X)$. Letting $\omega' =
    R_*(\omega)$, we get that $\omega' \in \Fil_{{D_\chi}/p} H^2_{p^{m-1}}(X) \setminus
    \Fil_{0} H^2_{p^{m-1}}(X)$. Letting $\chi' = \kappa_X (\omega')$, we get
    an element $\chi' \in \Fil_{{D_\chi}/p} \Br(X)[p^{m-1}]$ such that
    $\chi'$ does not lie in $\Fil_0 \Br(X)[p^{m-1}]$ by  \lemref{lem:Br-fil-0}.
    We conclude by induction on $m$ that $\chi'|_Z$ does not lie in $\Fil_0 \Br(Z)$
    for some $Z \in \sS$. Let $\iota \colon Z \inj X$ be the inclusion and let
    $D'_\chi = {D_\chi}|_Z$.

We now look at the commutative diagram (cf. \corref{cor:Fil-functorial})
    \begin{equation}\label{eqn:SC-change-6-5}
 \xymatrix@C1.5pc{
   \frac{\Fil_{D_\chi} H^2_{p^m}(X)}{\Fil_{0} H^2_{p^m}(X)} \ar[r]^-{\iota^*}
   \ar[d]_-{R_*} &
   \frac{\Fil_{D'_\chi} H^2_{p^m}(Z)}{\Fil_{0} H^2_{p^m}(Z)} \ar[d]^-{R_*} \\
   \frac{\Fil_{{D_\chi}/p} H^2_{p^{m-1}}(X)}{\Fil_{0} H^2_{p^{m-1}}(X)} \ar[r]^-{\iota^*} &
   \frac{\Fil_{{D'_\chi}/p}H^2_{p^{m-1}}(Z)}{\Fil_{0} H^2_{p^{m-1}}(Z)}.}
    \end{equation}
    As $ \chi'|_Z \notin \Fil_0 \Br(Z)$, we get that
    $\iota^*(R_*(\ov{\omega})) \neq 0$. A diagram chase shows that
    $\iota^*(\ov{\omega}) \neq 0$. Equivalently, $\chi|_Z$ does not lie in
    $\Fil_0 \Br(Z)$, as desired. This completes the induction step, and hence,
   the proof of the lemma.
 \end{proof}

We let $\wt{\Br}(X) := \frac{\Br(X)}{\Br(\sX) + \Br_0(X)} \cong
\frac{\Br(X)}{\Fil'_{-1} \Br(X)}$ (cf. ~\eqref{eqn:Fil-Br-spl}).
We shall use this notation for the remainder of this paper.

\begin{cor}\label{cor:SC-change-7}
  There are admissible hypersurface sections $\sX', \{\sX_i, \sW_i\}_{i \in J_r}$
  of $\sX$ such that the restriction map
  \begin{equation}\label{eqn:SC-change-7-0}
   \wt{\Br}(X) \to \
    \wt{\Br}(X') \bigoplus
    \left(\stackrel{r}{\underset{i =1}\bigoplus}
    \wt{\Br}(X_i)\right) \bigoplus
    \left(\stackrel{r}{\underset{i =1}\bigoplus} \wt{\Br}(W_i)\right)
  \end{equation}
  is injective.
  \end{cor}
\begin{proof}
We let
  \[
    A = \frac{\Fil'_{0} \Br(X')}{\Fil'_{-1} \Br(X')} \bigoplus
    \left(\stackrel{r}{\underset{i =1}\bigoplus}
    \frac{\Fil'_{0} \Br(X_i)}{\Fil'_{-1} \Br(X_i)}\right) \bigoplus
    \left(\stackrel{r}{\underset{i =1}\bigoplus}
    \frac{\Fil'_{0} \Br(W_i)}{\Fil'_{-1} \Br(W_i)}\right); 
    \]
    \[
B = \wt{\Br}(X') \bigoplus
    \left(\stackrel{r}{\underset{i =1}\bigoplus}
    \wt{\Br}(X_i)\right) \bigoplus
    \left(\stackrel{r}{\underset{i =1}\bigoplus} \wt{\Br}(W_i)\right);
    \ \ \mbox{and}
    \]
    \[
  C = \frac{\Br(X')}{\Fil'_{0} \Br(X')} \bigoplus
    \left(\stackrel{r}{\underset{i =1}\bigoplus}
    \frac{\Br(X_i)}{\Fil'_{0} \Br(X_i)}\right) \bigoplus
    \left(\stackrel{r}{\underset{i =1}\bigoplus}
    \frac{\Br(W_i)}{\Fil'_{0} \Br(W_i)}\right).
    \]

We now look at the commutative diagram
    \begin{equation}\label{eqn:SC-change-7-1}
      \xymatrix@C1pc{
        0 \ar[r] & \frac{\Fil'_0 \Br(X)}{\Fil'_{-1} \Br(X)} \ar[d] \ar[r] &
        \wt{\Br}(X)  \ar[r] \ar[d] & 
        \frac{\Br(X)}{\Fil'_{0} \Br(X)} \ar[r] \ar[d] & 0 \\
        0 \ar[r] & A \ar[r] & B \ar[r] & C \ar[r] & 0,}
    \end{equation}
    where the vertical arrows are the restriction maps.
    The left vertical arrow is injective by \lemref{lem:SC-change-5} while
    the right vertical arrow is injective by \lemref{lem:SC-change-6}.
    It follows that the middle vertical arrow is injective.
\end{proof}

We record the following result, which is clear from \corref{cor:SC-change-7}
using induction on $\dim(X)$.

\begin{cor}\label{cor:SC-change-8}
  There are admissible relative curves $\sC_1, \ldots , \sC_n$
  inside $\sX$ such that the restriction map
 $\wt{\Br}(X) \to \ \stackrel{n}{\underset{i =1}\bigoplus}\wt{\Br}(C_i)$
is injective.
  \end{cor}

\section{Brauer--Manin pairing II}\label{sec:BMP-2}
In this section, we continue in the setting of \S~\ref{sec:BMP} and complete the
proof of \thmref{thm:Main-3}. We then prove Theorems~\ref{thm:Main-2}
and~\ref{thm:Main-6}.
For an abelian group $A$, we shall write 
$A^\star = \Hom(A,  {\Q}/{\Z})$. For a topological abelian group $A$, we let
$A^\vee = \Hom_\cont(A, {\Q}/{\Z})$, where ${\Q}/{\Z}$ has the discrete topology.

\subsection{Extension of a theorem of Saito--Sato}\label{sec:M3-prf}
In this subsection, we state and prove our main result concerning the
Brauer--Manin pairing. The analogous result over
characteristic zero local fields was conjectured by Colliot-Th{\'e}l{\`e}ne
\cite[Conj.~1.4]{CT-Bordeaux}, and proven by Colliot-Th{\'e}l{\`e}ne--Saito
\cite{CTSaito} (prime-to-$p$ part) and Saito--Sato \cite{Saito-Sato-ENS}
(in general).

\begin{thm}\label{thm:SS-char-p}
  Let $R$ be a complete discrete valuation ring of characteristic $p$ with finite
  residue field and let
  $\sX$ be a projective semistable $R$-scheme with generic fiber $X$.  
  Then the Brauer--Manin pairing for $X$ induces an injective homomorphism
  $\frac{\Br(X)}{\Br(\sX)} \to \Hom(\CH_0(X), {\Q}/{\Z})$ and a surjective
  homomorphism
  $A_0(X) \to \wt{\Br}(X)^\star$.
  \end{thm}
\begin{proof}
We shall use the notation of \S~\ref{sec:KEF} throughout the proof
  (in particular, $k = Q(R)$).
We first reduce the proof to the case when $X$ is geometrically integral
over $k$. To do this, we use the Stein factorization theorem which tells us
that the structure map $f \colon \sX \to S$ has a canonical factorization
$\sX \xrightarrow{f'} \Spec(A) \xrightarrow{g} S$,
where $A = H^0_\zar(\sX, \sO_\sX)$ and $f'$ has geometrically connected fibers.
Since $S$ is excellent, $\sX$ is integral and normal,
and $f$ is projective and flat, it follows that $\sO_k \to A$ is a finite local
ring homomorphism of henselian discrete valuation rings.
We let $k' = Q(A)$ so that $A = \sO_{k'}$.
We next note that if $\pi'$ is a uniformizer of
$k'$, then $\pi = \pi'^e$, where $e = e({k'}/k)$. This implies that the generic
fiber (resp. reduced special fiber) of $f$ coincides with the
generic fiber (resp. reduced special fiber) of $f'$.

Using the flat base change property
for the Zariski cohomology of the structure sheaves, we get
that $H^0_\zar(X, \sO_X) = k'$. In other words, $X \to \Spec(k')$ is the Stein
factorization of $X \to \Spec(k)$.
We finally note that as $X$ is a smooth $k$-scheme,
it is geometrically reduced over $k$. This forces ${k'}/k$ to be a finite separable
field extension (cf. \cite[Tag~0366]{SP}) and $X/{k'}$ to be smooth (since
this is an irreducible component of $X_{k'}$ which is smooth over $k'$).
Since $X$ is geometrically connected over $k'$, it is actually geometrically
integral over $k'$. Putting everything together and letting
$S' = \Spec(\sO_{k'})$, we deduce that $\sO_{k'}$ is complete, $\sX$ is a semistable
$S'$-scheme with geometrically integral generic fiber $X$ and geometrically
connected reduced special fiber $Y$.

To prove the first claim of the theorem, we note that 
$\frac{\Br(X)}{\Br(\sX)} \to \Hom(\CH_0(X), {\Q}/{\Z})$ does not see whether
$\sX$ is an $S$-scheme or an $S'$-scheme because the push-forward map
$g_* \colon \Br(k') \to \Br(k)$ is bijective and the evaluation map for the
Brauer classes of the $k'$-scheme $X$ factors through $g_*$.
We can therefore replace $k$ by $k'$ to prove the first claim.

To prove the second claim, we let
$A'_0(X) = \Ker(f'_* \colon \CH_0(X) \to \CH_0(k'))$ and
look at the diagram
\begin{equation}\label{eqn:SS-char-p-0}
  \xymatrix@C1.2pc{
    A'_0(X) \ar[r] \ar[d]_-{g_*}^-{\cong} &
    \left(\frac{\Br(X)}{\Br(\sX) + \Br_0(X')}\right)^\star
    \ar[d]^-{(g^*)^\vee}
    \\
    A_0(X) \ar[r] & \wt{\Br}(X)^\star,}
\end{equation}
where $\Br_0(X')$ is the image of the pull-back map $\Br(k') \to \Br(X')$.
One checks easily that the left vertical arrow is actually the identity map
of $A_0(X)$.

This diagram is commutative by the functoriality of the Brauer--Manin pairing with
respect to proper maps. Assuming the theorem holds for $f'$, the top horizontal
arrow is surjective. To show that the same holds for $f$, it
remains to show that the right vertical arrow is surjective. We shall show
that this is in fact bijective.

Indeed, in the commutative diagram 
\begin{equation}\label{eqn:SS-char-p-1}
  \xymatrix@C1.2pc{
    \Br(\sX) + \Br(k) \ar[r] \ar[d]_-{g^*} & \Br(X) \ar[r] \ar[d]^-{\id} &
     \wt{\Br}(X) \ar[d]^-{g^*} \ar[r] & 0 \\
    \Br(\sX) + \Br(k') \ar[r] & \Br(X) \ar[r] &
    \frac{\Br(X)}{\Br(\sX) + \Br_0(X')} \ar[r] & 0}
\end{equation}
with exact rows, $g^* \colon\Br(k) \to \Br(k')$ is identified with the
map ${\Q}/{\Z} \xrightarrow{[k':k]} {\Q}/{\Z}$. In
 particular, it is surjective. It follows that
 left vertical arrow in ~\eqref{eqn:SS-char-p-1} is surjective.
 A diagram chase shows that the right vertical
arrow is bijective. This completes the proof of the reduction step.
We shall assume in the rest of the proof that $X$ is geometrically integral over
$k$.

To show that the map
$\frac{\Br(X)}{\Br(\sX)} \to \Hom(\CH_0(X), {\Q}/{\Z})$ is injective, we let
$\sA \in \Br(X)$ be such that $\<\sA, \alpha\> = 0$ for every $\alpha \in \CH_0(X)$.
This implies in particular that $\ev_{\sA}(P) = 0$ for every $P \in X^o_\fin$.
Equivalently, $\sA \in \Ev_{-2} \Br(X)$. We conclude from \thmref{thm:Main-1}
(cf. \lemref{lem:K-Ev-6}) that $\sA \in \Br(\sX)$. This proves the desired
injectivity, and hence the first claim of the theorem.

We shall prove the surjectivity of the map $A_0(X) \to \wt{\Br}(X)^\star$
by induction on $d_{\sX} = \dim(\sX)$.
The case $d_{\sX} = 2$ follows from \cite[Thm.~9.2]{Saito-Invent}. To prove the
general case, we choose admissible hypersurface sections
$\{\sX', \sX_i, \sW_i\}_{i \in J_r}$ of $\sX$ as in \corref{cor:SC-change-7} and look
at the diagram

\begin{equation}\label{eqn:SS-char-p-2}
  \xymatrix@C1pc{
    A_0(X') \bigoplus \left(\stackrel{r}{\underset{i =1}\bigoplus} A_0(X_i)\right)
    \bigoplus \left(\stackrel{r}{\underset{i =1}\bigoplus}
    A_0(W_i)\right) \ar[r] \ar[d] & A_0(X) \ar[d] \\
   \wt{\Br}(X')^\star \bigoplus
   \left(\stackrel{r}{\underset{i =1}\bigoplus}
    \wt{\Br}(X_i)\right)^\star \bigoplus
    \left(\stackrel{r}{\underset{i =1}\bigoplus}
    \wt{\Br}(W_i)\right)^\star \ar[r] & \wt{\Br}(X')^\star.}
\end{equation}
The top horizontal arrow is the sum of the push-forward
maps and the bottom horizontal arrow is the sum of the duals of the pull-back maps.

The above diagram is commutative by the functoriality of Brauer--Manin pairing with
respect to proper maps. The left vertical arrow of this diagram is surjective
by induction on $d_\sX$ since taking the discrete dual commutes with finite
direct sums. The bottom horizontal arrow is surjective by \corref{cor:SC-change-7}.
A diagram chase shows that the right vertical arrow is surjective. This completes
the proof of the theorem.
\end{proof}

The following is clear from \thmref{thm:SS-char-p}.

\begin{cor}\label{cor:SS-char-p-3}
  If $X$ is an in \thmref{thm:SS-char-p}, then the Brauer--Manin pairing for $X$
  induces an injective homomorphism $\wt{\Br}(X) \to \Hom(A_0(X), {\Q}/{\Z})$
  and an exact sequence
  \[
  \CH_0(X) \to \Hom\left(\frac{\Br(X)}{\Br(\sX)}, {\Q}/{\Z}\right) \to
     {\wh{\Z}}/{\Z} \to 0.
     \]
    \end{cor}

\vskip .2cm

\subsection{Application to the structure of $\wt{\Br}(X)$}
\label{sec:Br-Structure}
In the course of proving \thmref{thm:SS-char-p}, we established 
the following result. The analogous result over characteristic zero local fields
was
proven by Colliot-Th{\'e}l{\`e}ne--Saito \cite{CTSaito} (prime-to-$p$ case) and
Saito--Sato \cite{Saito-Sato-ENS} (general case).

\begin{thm}\label{thm}\label{thm:SS-char-p-fin}
  Let $X$ be as in \thmref{thm:SS-char-p}. There is an inclusion of abelian groups
  \begin{equation}\label{eqn:SS-char-p-fin-0}
    \wt{\Br}(X) \inj \left(\stackrel{\infty}{\underset{i =1}\bigoplus}
       {\Q_p}/{\Z_p}\right) \bigoplus T,
  \end{equation}
  where $T$ is finite. In particular, one has the following.
  \begin{enumerate}
 \item
    $\wt{\Br}(X)\{p'\}$ is finite.
  \item
    $\wt{\Br}(X) \otimes_\Z {\Z}/n$ is finite for every
    integer $n$ prime to $p$.
    \end{enumerate}
\end{thm}
\begin{proof}
To prove ~\eqref{eqn:SS-char-p-fin-0}, we let $\sX \xrightarrow{f'} S' \to S$ be
  the Stein factorization of
  $\sX \xrightarrow{f} S$ as in the proof of \thmref{thm:SS-char-p}. We showed in
  the proof of the theorem that
  the canonical map $\wt{\Br}(X) \to
  \frac{\Br(X)}{\Br(\sX) + \Br'_0(X)}$ is bijective if we let
  $\Br'_0(X)$ be the image of the map $\Br(k') \to \Br(X)$. We can therefore assume
  that $X$ is geometrically integral over $k$.
  We shall now prove ~\eqref{eqn:SS-char-p-fin-0} by first considering
  a special case.

  Suppose that $\dim(X) = 1$. In this case,
  we have an inclusion $A_0(X) \inj J_X(k)$, and hence a
  surjection $J_X(k)^\vee  \surj A_0(X)^\vee$
  by the Pontryagin duality. 
  By \cite[Thm.~21]{Clark-Lacy} and \cite[Prop.~4.9]{Folland}, one has
  either $J_X(k)^\vee = 0$ or 
  $J_X(k)^\vee \cong \left(\stackrel{\infty}{\underset{i =1}\bigoplus}
  {\Q_p}/{\Z_p}\right) \bigoplus T'$, where $T'$ is finite.
  Combining this with \cite[Thm.~23.1]{Fuchs}, we get  $A_0(X)^\vee \inj
  \left(\stackrel{\infty}{\underset{i =1}\bigoplus} {\Q_p}/{\Z_p}\right)
  \bigoplus T$, where $T$ is finite. 
  We now apply \thmref{thm:SS-char-p} to conclude  the proof of
  ~\eqref{eqn:SS-char-p-fin-0} when $\dim(X) = 1$.  The general case
follows by \corref{cor:SC-change-8}.

Item (1) of the second part of the theorem follows directly from
~\eqref{eqn:SS-char-p-fin-0}.
To prove item (2), we let $A(X) = \left(\stackrel{\infty}{\underset{i =1}\bigoplus}
{\Q_p}/{\Z_p}\right) \bigoplus T$ and consider the exact sequence
\begin{equation}\label{eqn:SS-char-p-fin-4}
  0 \to \wt{\Br}(X) \to A(X) \xrightarrow{\alpha} B(X) \to 0.
\end{equation}
It follows from \cite[Thm.~23.1]{Fuchs} that $B(X) \cong T_1 \bigoplus T_2$,
where
$T_1 = \alpha \left(\stackrel{\infty}{\underset{i =1}\bigoplus} {\Q_p}/{\Z_p}\right)$
has the form $\stackrel{N}{\underset{i =1}\bigoplus} {\Q_p}/{\Z_p}$ with
$0 \le N \le \infty$ and $T_2$ is a finite group. For any integer $n \ge 1$
prime to $p$, we now consider the exact sequence
\begin{equation}\label{eqn:SS-char-p-fin-5}
  A(X)[n] \xrightarrow{\alpha} B(X)[n] \to {\wt{\Br}(X)}/n \to {A(X)}/n \to
  {B(X)}/n \to 0.
  \end{equation}
Since $\coker(\alpha)$ and ${A(X)}/n$ are finite, we are done.
\end{proof}

\begin{remk}\label{remk:SS-char-p-fin-6}
  The proof of \thmref{thm:SS-char-p-fin} shows that ~\eqref{eqn:SS-char-p-fin-0}
  is a bijection when $X$ is a curve of positive genus. It follows that item (1) of
  the theorem fails for the $p$-primary torsion subgroup. We however expect
  that item (2) holds for every $n \in \N$.
\end{remk}

\vskip.2cm

\begin{exm}\label{exm:Trivial-Br-p}
  Using \cite[Exm.~4.3]{Lazda-Skorobogatov}, one can show that the inclusion of
  abelian groups given in ~\eqref{eqn:SS-char-p-fin-0} may not be a bijection if
  $\dim(X) \ge 2$. To see this, we let $\F$ be a finite field of
  characteristic $p \neq 2$ and let $Y$ be a Kummer surface over $\F$
  described in \cite[Exm.~4.3 (iii)]{Lazda-Skorobogatov}. We let
  $k = \F((\pi))$ and let $\sX = Y \times_{\Spec(\F)} S$, where $S = \Spec(\sO_k)$.
  Then $\sX$ is a smooth projective $S$-scheme of relative dimension two whose
  fibers are geometrically integral. We let $X$ denote the generic fiber of $\sX$.

  We let $\ov{k}$ be an
  algebraic closure of $k$ and let $k_s \subset \ov{k}$ be the separable
  closure of $k$ inside $\ov{k}$. It follows from
  \cite[Prop.~4.2]{Lazda-Skorobogatov} and \cite[Thm.~5.2.5 (ii)]{CTS}
  that $\Br(X^s)\{p\} = 0$.
  In particular, $\Br_1(X) := \Ker\left(\Br(X) \to \Br(X^s)\right)$ has the
  property that $\Br_1(X)\{p\} = \Br(X)\{p\}$. In particular,
  \propref{prop:CTS-Trans*} implies that
  $\wt{\Br}(X)\{p\} = \frac{\Br(X)\{p\}}{(\Br_0(X) + \Br(\sX))\{p\}} = 0$.
  By \thmref{thm:SS-char-p-fin}, we conclude that $\wt{\Br}(X)$ is a finite
  $p$-torsion-free group.
  \end{exm}

\vskip.2cm

In light of \thmref{thm:SS-char-p-fin} and Example~\ref{exm:Trivial-Br-p},
we pose the following question.

\begin{ques}\label{ques:Br-inf}
  Let $X$ be as in \thmref{thm:SS-char-p}. Suppose that $H^1_\zar(X, \sO_X) \neq 0$.
  Is it true that $\wt{\Br}(X) \cong
  \left(\stackrel{\infty}{\underset{i =1}\bigoplus} {\Q_p}/{\Z_p}\right)
  \bigoplus T$ for some finite group $T$?
  \end{ques}

\vskip.2cm

\begin{prop}\label{prop:CTS-Trans*}
  Let $X$ be as in \thmref{thm:SS-char-p}. Assume that $H^1_\zar(X, \sO_X) = 0$
  and $\NS(\ov{X})$ is torsion-free. Then we have
  $\Br_1(X) \subset \Fil'_{-1}\Br(X)$.
\end{prop}
\begin{proof}
Since $\NS(X^s) \subset \NS(\ov{X})$, we see that $\NS(X^s)$ is torsion-free.
In this case, the proof of \cite[Lem.~2,2, Prop.~2.3]{CTS-Trans} carries over verbatim
in positive characteristic to show that the natural inclusion
$\Br^{\nr}_1(X) := \Ker\left(\Br(X) \to \Br(X^{\nr})\right) \inj \Br_1(X)$
is an isomorphism,
where $X^{\nr}$ is the base change of $X$ by the maximal unramified extension
$k_{\nr} \subset k_s$ of $k$. It remains to show that
$\Br^{\nr}_1(X) \subset \Fil'_{-1}\Br(X)$. This can be proven by
mimicking the proof of \cite[Lem.~2.2]{CTS-Trans}. Alternatively, one
can argue as follows.

If $\chi \in \Br^{\nr}_1(X)$, then \thmref{thm:RSW-gen}
implies that $\chi \in \Fil_0 \Br(X)$. Meanwhile, it follows from
the definition of $\Fil'_{-1}\Br(X)$, \propref{prop:Kato-Br-0} and
~\eqref{eqn::K-Ev-8-0} that $\Br^{\nr}_1(X) \bigcap \Fil_0 \Br(X) \subset
\Fil'_{-1}\Br(X)$. This completes the proof.
\end{proof}

\vskip.2cm

\subsection{Application to the index of $X$}\label{sec:Index-Structure}
We record another application of \thmref{thm:SS-char-p}.
  Let $\sX$ be as in \thmref{thm:SS-char-p}. Let
  $\sX_s = \stackrel{r}{\underset{i=1}\sum} e_iY_i
\in \Div(\sX)$. For each $i \in J_r$, we let $f_i = [\F_i:\F]$, where
$\F_i$ is the algebraic closure of $\F$ in $k(Y_i)$.
We consider the following numbers. \\
(1) $I_{\Br}$ is the order of $\Ker\left(\Br(k) \to \frac{\Br(X)}{\Br(\sX)}\right)$.
\\
(2) $I_X$ is the index of $X$ as a $k$-scheme, i.e., it is the g.c.d. of the
degrees of the closed points of $X$. \\
(3) $I_0$ is the g.c.d. of the $e_if_i$. \\
The following is a special case of the
more general result \cite[Thm.~5.4.1]{Saito-Sato-ENS}.\footnote[1]{Saito--Sato
  prove the result in characteristic zero but we believe that their
proof remains valid in positive characteristic.}

\begin{cor}\label{cor:Index-degree}
One has $I_{\Br} = I_X = I_0$.
\end{cor}
\begin{proof}
  Having \propref{prop:Kato-Br-0} and \thmref{thm:SS-char-p} in hand, the proof
  of the prime-to-$p$ part of the corollary given by
  Colliot-Th{\'e}l{\`e}ne--Saito in characteristic zero
  (cf. \cite[Thm.~3.1]{CTSaito}) carries over verbatim.
\end{proof}

\subsection{Extension of two theorems of Ieronymou}\label{M2-prf}
Let $k$ be an hdvr of characteristic $p$ with the finite residue field $\F$
and the excellent ring of integers $\sO_k$. We let $S = \Spec(\sO_k)$ and
let $\sX$ be a smooth and projective $S$-scheme of absolute dimension $d_\sX \ge 1$.
We let $X$ (resp. $Y$) denote the generic (resp. special) fiber of $\sX$.
In this subsection, we shall extend two theorems of Ieronymou \cite{Ieronymou}
from characteristic zero to positive characteristic.
We refer the reader to \cite[Defn.~1]{Kollar-Szabo} for the definition of a
separably rationally connected (SRC) variety over a field.
The following is an application of Theorems~\ref{thm:Main-1} and
~\ref{thm:Main-3}. In characteristic zero, this was proven by
 Ieronymou (cf. \cite[Thm.~A]{Ieronymou}).

\begin{thm}\label{thm:Ier-0}
  Assume that the special fiber $Y$ is an SRC variety over $\F$. Then
we have $\Br(X) = {\Ev}_{-1} \Br(X)$. In particular, the evaluation
  map ${\ev}_{\sA} \colon X(k) \to {\Q}/{\Z}$ is constant for every $\sA \in \Br(X)$.
\end{thm}
\begin{proof}
  We can assume that $X$ is connected. Let $\wh{X}$ be the base change of $X$ by the
  completion map $\Spec(\wh{k}) \to \Spec(k)$, and let $w \colon \wh{X} \to X$
  denote the projection.
If $P \in X_{(0)}$, then $w^{-1}(P)$ is the unique point
$\wh{P} \in \wh{X}_{(0)}$ such that $k(\wh{P}) = \wh{k(P)} \cong
k(P) \otimes_k \wh{k}$. Furthermore, 
$\Spec(k(P)) \times_X \wh{X} = \Spec(k(\wh{P}))$ and
\begin{equation}\label{eqn:SS-char-p-fin-2}
 \ev_{\sA}(P) = \inv_{k(P)}(\sA|_{k(P)}) = \inv_{k(\wh{P})}(\sA|_{k(\wh{P})}) =
 \ev_{w^*(\sA)}(\wh{P}) \ \ \forall \ \sA \in \Br(X).
\end{equation}
Using this, we can assume that $k$ is complete.

In the latter case, it follows from \cite[Thm.~5]{Kollar-Szabo} and
  \thmref{thm:Main-3} that $\Br(X) = \Fil'_{-1} \Br(X)$. We now apply
  \thmref{thm:Main-1} to conclude that $\Br(X) = {\Ev}_{-1} \Br(X)$, which
implies the theorem.
  \end{proof}

\begin{remk}\label{remk:Ier-0-1}
  In \cite[Thm.~A]{Ieronymou}, it is assumed that $k$ is complete,
  $p$ is an odd prime and that $Y$ is
  geometrically integral. However, these assumptions can be dispensed with by
  following the proof of \thmref{thm:Ier-0} above, which also works when
  $\Char(k)=0$, rather than the proof given in op. cit..
\end{remk}

Recall that a smooth projective geometrically integral surface over a field
$K$ is called an Enriques surface if its base change to $\ov{K}$ is an
Enriques surface as in \cite[Defn.~1.1.6]{CDL}.
The following result was proven by Ieronymou (cf. \cite[Prop.~11]{Ieronymou})
in characteristic zero.

\begin{thm}\label{thm:Ier-1}
  Suppose that $\sX$ is a smooth and projective $S$-scheme whose generic fiber $X$
  is an Enriques surface over $k$. Let $\sA \in \Br(X)$ have odd order. Then the
  evaluation map ${\ev}_{\sA} \colon X(k) \to {\Q}/{\Z}$ is constant.
\end{thm}
\begin{proof}
 It suffices to show that $\sA \in {\Ev}_{-1} \Br(X)$.
  Let $\ov{Y}$ denote the base change of $Y$ to an algebraic closure of ${\F}$.
  Then $\ov{Y}$ is a Enriques surface over $\ov{\F}$ by \cite[Cor.~1.2.11]{CDL}.
  In particular, $Y$ is geometrically integral.
  Suppose first that $p =2$. Then our assumption implies that all refined
  Swan conductors of $\sA$ are zero because their targets are $\F$-vector spaces
  (cf. \corref{cor:RSW-gen-1}). It follows by \thmref{thm:RSW-gen}
  that $\sA \in \Fil'_0 \Br(X)$.
  Since $\partial_X \colon \frac{\Fil'_0 \Br(X)}{\Fil'_{-1} \Br(X)}
  \to H^1_\et(\ov{Y})$ is injective (cf. ~\eqref{eqn::K-Ev-8-0}), and
  $\Fil'_{-1} \Br(X) \subset {\Ev}_{-1} \Br(X)$
  by \lemref{lem:K-Ev-8}, it suffices to show that $\partial_X(\sA) = 0$.
  But this follows from \cite[Cor.~1.4.5, ~1.3.7]{CDL}, which together
  imply that $H^1_\et(\ov{Y})$ is either trivial or a
  group of order two.

If $p \neq 2$, we have
  $H^0_\zar(Y, \Omega^i_Y(nY)) = 0$ for $i \in J_2$ and $n \in \N_0$ by
  \cite[Thm.~1.1.8, Cor.~1.4.5]{CDL}. Hence, all refined
  Swan conductors of $\sA$ are zero (cf. \lemref{lem:RSW-spl}).  
  It follows from \thmref{thm:RSW-gen} that $\sA \in \Fil'_0 \Br(X)$. Since
  $H^1_\et(\ov{Y}) \cong {\Z}/2$ by
  Cor.~1.4.5 of op. cit., we conclude, as above, that $\sA \in {\Ev}_{-1} \Br(X)$.
  \end{proof}

\section{Albanese cokernel in positive characteristic}\label{sec:Kai}
As we mentioned in \S~\ref{sec:Intro}, the Albanese map for smooth projective varieties
over a field may not be surjective, and describing its
cokernel is, in general, a difficult problem.
In this section and the next, we study this problem for higher-dimensional smooth
projective varieties over local fields of positive characteristic. 
In characteristic zero, this problem was solved by Lichtenbaum \cite{Lichtenbaum}
for curves, and by Kai \cite{Kai} in higher dimensions under the
smoothness condition on the Picard schemes of a smooth model of the given
variety (note that this condition is automatic for curves).
For curves in positive characteristic, this was solved by Saito
\cite{Saito-Invent}.

Our setting for Sections~\ref{sec:Kai} and ~\ref{sec:Kai-2} is the following.
We fix a complete discrete valuation field $k$ of characteristic $p$ with
the ring of integers $\sO_k$, uniformizer $\pi$ and finite residue field $\F$.
A choice of the uniformizer $\pi$ determines isomorphisms
$\F[[t]] \xrightarrow{\cong} \sO_k$ and $\F((t)) \xrightarrow{\cong} k$.
We let $S =\Spec(\sO_k)$ and $\eta = \Spec(k)$.
We fix a separable closure $k_s$ of $k$ and algebraic closure $\ov{\F}$ of $\F$.
We let $\sO_{k_s}$ denote the integral closure of $\sO_k$ in $k_s$. One knows
that $\sO_{k_s}$ is a one-dimensional
valuation ring with value group $\Q$ and residue field
$\ov{\F}$. We let $S_s = \Spec(\sO_{k_s})$ and let $\eta_s$ denote the generic
point of $S_s$. We let $\Sigma^s_k$ denote the set of isomorphism classes of finite
separable field extensions of $k$. We let $\Gamma$ denote the absolute Galois group
of $k$.

We let $\sX$ be a smooth projective connected
$S$-scheme of absolute dimension $d_\sX$ with the generic fiber $X$ and
special fiber $Y$. Let $j \colon X \inj \sX$ and
$\iota \colon Y \inj \sX$ denote the inclusions.
Let $\sX^s = \sX \times_S S_s, \ \ X^s = X \times_{\Spec(k)} \Spec(k_s)$ and
$\ov{Y} = Y \times_{\Spec(\F)} \Spec(\ov{\F})$.
It is then clear that $\sX^s$ is a smooth projective $S_s$-scheme with generic
fiber $X^s$ and special fiber $\ov{Y}$.  We let $j_s \colon X^s \inj \sX^s$ and
$\iota_s \colon \ov{Y} \inj \sX^s$ denote the inclusions.
We assume that $X$ is geometrically integral over $k$.

An easy application of the
Stein factorization implies that $\sX$ is integral and $Y$ is 
geometrically integral over $\F$. The same argument also shows that $\sX^s$ and
$\ov{Y}$ are integral. By \cite[Thm.~9.4.8]{Kleiman}, the Picard scheme
$\Picc(\sX)$ exists as an $S$-scheme which is separated and locally of finite type
over $S$. In particular, all local properties (e.g., smoothness)
of Noetherian $S$-schemes are defined for $\Picc(\sX)$.
 
The main result that we shall prove is the following.

\begin{thm}\label{thm:Kai-Main}
  Assume that $\Picc(\sX)$ is a smooth $S$-scheme. Then the cokernel of the
  Albanese map
  \[
  \alb_X \colon A_0(X) \to \Alb_X(k)
  \]
  is canonically Pontryagin dual to the cokernel of the canonical map
  \[
  \Pic({X}^s)^\Gamma \to \NS({X}^s)^\Gamma.
  \]
\end{thm}

An application of \thmref{thm:Kai-Main} (and its proof)
is the following finiteness result.
Over $p$-adic fields, finiteness of the Albanese cokernel was proven
unconditionally by Saito--Sujatha \cite{Saito-Sujatha}.
But their proof does not extend to positive characteristic.

\begin{cor}\label{cor:Kai-Main-0}
  Let $X$ be as in \thmref{thm:Kai-Main}. Then
  $\coker(\Pic({X}^s)^\Gamma \to \NS({X}^s)^\Gamma)$ and $\coker(\alb_X)$
  are finite $p$-groups.
  \end{cor}

We shall prove \thmref{thm:Kai-Main} in \S~\ref{sec:Kai-2}. In this section, we
establish the existence of the Albanese map $\alb_X$ in \thmref{thm:Kai-Main} and
prove several preliminary results. We also recall Milne's extension of Tate's
duality theorem for abelian varieties to local fields of positive characteristic
and establish some results concerning this duality theorem.

\subsection{Albanese map for 0-cycles over imperfect fields}
\label{sec:Prelim-Kai-1}
We recall from \cite[p.~45-46]{Lang} that for any smooth projective variety
$Y$ over a field $F$, there exists an abelian variety $\Alb_Y$,
together with a morphism $\phi \colon Y \times_F Y \to
\Alb_Y$ having a universal property with respect to maps from $Y \times_F Y$ to
abelian varieties over $F$ (cf. \cite[\S~2.1]{Kai}). Using this universal property,
one obtains a group homomorphism $\phi_* \colon Z_0(Y)_0 \to \Alb_Y(F)$, where
$Z_0(Y)_0$ is the group of 0-cycles on $Y$ of degree zero.
$\Alb_Y$ is called the Albanese variety of $Y$.

If $F$ is a perfect field, it is known that $\phi_*$ induces a
map $\alb_Y \colon A_0(Y) \to \Alb_Y(F)$ (e.g., see \cite[\S~7, 8]{Kahn-Sujatha}
or \cite{Ramachandran}). But it is not clear (at least to the authors)
how to define such a map over imperfect fields (such as $k = \F((t))$).
In this subsection, our goal is to
show that the map $\alb_X$ in \thmref{thm:Kai-Main} indeed exists and has the
expected functorial properties.

For the remainder of \S~\ref{sec:Prelim-Kai-1}, we fix a field $K$
and a geometrically integral smooth projective variety $Y$ over $K$.
For a field extension ${K'}/K$, we let $S_{K'} = \Spec(K')$.
Let $A_Y = \amalg_{i \in \Z} \Alb^i_Y$  denote the Albanese scheme of $Y$ and let
 $\alb^1_Y \colon Y \to \Alb^1_Y \subset A_Y$ denote the universal Albanese map.
Recall from \cite[\S~2.1]{Kai} (see also \cite[\S~1]{vanHamel} and
\cite[\S~7, 8]{Kahn-Sujatha}) that $A_Y$ is a group scheme over $S_K$ which is
separated and locally of finite type. Furthermore, $\Alb^0_Y = \Alb_Y$ and $\Alb^1_Y$
is an $\Alb_Y$-torsor over $S_K$ such that $x \mapsto x - \alb^1_{Y_F}(x_0)$ defines
an isomorphism of $F$-schemes $\Alb^1_{Y_{F}} \xrightarrow{\cong} \Alb_{Y_F}$ if
$x_0 \in Y$ and $F = k(x_0)$.
 From the universal property of $\alb^1_Y$, it follows that $A_Y$ commutes with the
 base change of $Y$ by all algebraic field extensions of $K$
 (cf. \cite[Thm.~3.3]{Grothendieck} and \cite[Remark after Cor.~A.5]{Wittenberg}).

There exists a commutative diagram of $K$-schemes
\begin{equation}\label{eqn:Alb-0}
  \xymatrix@C1.8pc{
    Y^n \ar[r]^-{(\alb^1_Y)^n} \ar[d] & (\Alb_Y^1)^n \ar[d] \\
    S^n_K(Y) \ar[r]^-{\alb^n_Y} & \Alb^n_Y}
\end{equation}
for every $n \in \N$, where $S^n_K(-)$ denotes the $n$-th symmetric power functor on
$\Sch_K$, the left vertical arrow is the quotient map, and the right vertical arrow
is the summation map of $A_Y$ with respect to its addition operation.
There is a functorial morphism
of presheaves with transfer $\deg \colon A_Y \to \Z$ on $\Sch_K$ which sends
$\Alb^n_Y$ onto $n$ for all $n \in \Z$. This map is surjective whose kernel is
$\Alb_Y$.

As shown in \cite[\S~6, p.~81]{Suslin-Voevodsky} (see also \cite[\S~6.3]{Deligne}),
there exists a canonical
section $\sigma^{T'}_{T} \colon T \to S^d_T(T')$ of the projection map
$S^d_T(T') \to T$ of the $T$-scheme $S^d_T(T')$ associated to a morphism of
$K$-schemes $T' \to T$, which is finite and flat of constant degree $d$.
For $d \ge 1$, let $Y^d_{(0)}$ denote the set of closed points on $Y$ of degree
$d$ over $K$. Then the maps $\sigma^{x}_{S}$ for $x \in Y_{(0)}$ induce
a canonical map $\alpha^Y_d \colon Y^d_{(0)} \to S^d_K(Y)(K)$
(cf. \cite[\S~6, p.~81]{Suslin-Voevodsky}).
For $x \in Y^d_{(0)}$, we let
\begin{equation}\label{eqn:Alb-0-0}
  \alb_Y([x]) = \alb^d_Y(\alpha^Y_d(x)) \in \Alb^d_Y(K) \subset A_Y(K).
\end{equation}
Extending linearly, we get a homomorphism $\alb_Y \colon Z_0(Y) \to A_Y(K)$.

Suppose now that ${K'}/K$ is an ffe and there is a commutative
 diagram of $K$-schemes
 \begin{equation}\label{eqn:Alb-1}
   \xymatrix@C1pc{
     S_{K'} \ar[r]^-{\iota'} \ar[d]_-{\phi} & Y' \ar[d]^-{\phi_Y} \\
     S_K \ar[r]^-{\iota} & Y,}
 \end{equation}
 where $\iota'$ is a map of $K'$-schemes, $Y' = Y \times_{S_K, \phi} S_{K'}$ and
 $\phi_Y$ is the projection.

This gives rise to the diagram of corresponding Albanese group schemes
 and Albanese maps
 \begin{equation}\label{eqn:Alb-2-0}
   \xymatrix@C4pc{
     S_{K'} \ar[rr]^-{\iota'} \ar[dr]^->>>>{\alb^1_{{K'}}} \ar[dd]_-{\phi} & & Y'
     \ar[dd]^->>>>>{\phi_Y}
\ar[dr]^-{\alb^1_{Y'}} & \\
& \Alb^1_{{K'}} \ar[rr]^->>>>>>>{\alb_{\iota'}} \ar[dd]^->>>>>{\alb_{\phi}} & &
\Alb^1_{Y'} \ar[dd]^-{\alb_{\phi_Y}} \\
S_K \ar[dr]_-{\alb^1_{K}} \ar[rr]^->>>>>>>>>>{\iota} & & Y \ar[dr]_-{\alb^1_Y} & \\
& \Alb^1_{K} \ar[rr]_-{\alb_{\iota}} & & \Alb^1_Y,}
\end{equation}
 where $\alb_{\phi}$ (resp. $\alb_{\iota}, \ \alb_{\iota'}$) is the unique map
 between the Albanese torsors induced by $\phi$ (resp. $\iota, \ \iota'$).
In particular, the left, top and bottom faces of ~\eqref{eqn:Alb-2-0} commute.

The right face of ~\eqref{eqn:Alb-2-0} commutes because $\alb^1_{Y'}$
is the base change of $\alb^1_Y$
by $\phi$ and its vertical arrows are the canonical projections induced by $\phi$. 
The back face is the commutative diagram ~\eqref{eqn:Alb-1}.
Since $\alb^1_K$ and $\alb^1_{K'}$ are isomorphisms, a diagram chase shows that the
front face of ~\eqref{eqn:Alb-2-0} also commutes.
We have thus shown the following.
 \begin{lem}\label{lem:Alb-2}
   All faces of ~\eqref{eqn:Alb-2-0} commute.
 \end{lem}

 \vskip.2cm

We next consider the diagram
\begin{equation}\label{eqn:Alb-3}
   \xymatrix@C4pc{ 
     Z_0(S_{K'}) \ar[rr]^-{\alb_{K'}} \ar[dr]_-{\iota'_*} \ar[dd]_-{\phi_*} & &
     A_{K}(K') \ar[dr]^-{\iota'_*}
     \ar[dd]^->>>>>{\phi_{A_K *}} & \\
     & Z_0(Y') \ar[rr] \ar[dd]^->>>>>>{\phi_{Y *}} \ar[rr]^->>>>>>>>>>>>>{\alb_{Y'}}
     & & A_Y(K') \ar[dd]^-{\phi_{A_Y *}} \\
     Z_0(S_K) \ar[rr]^->>>>>>>>>>>{\alb_K} \ar[dr]_-{\iota_*} & &
     A_K(K) \ar[dr]^-{\iota_*} & \\
     & Z_0(Y) \ar[rr]^-{\alb_Y} & & A_Y(K).}
   \end{equation}
Here, $\alb_{K}$ and $\alb_{K'}$ are isomorphisms because $A_K \cong \Z_K$ and
$A_{K'} \cong \Z_{K'}$ as group schemes. The arrows $\phi_{A_K *}$ and $\phi_{A_Y *}$ are
the transfer maps (cf. \cite[\S~6, p.~81]{Suslin-Voevodsky}).

\begin{lem}\label{lem:Alb-4}
 Except for the front face, all other faces of ~\eqref{eqn:Alb-3} commute.
  The front face commutes when restricted to the subgroup $\iota'_*(Z_0(S_{K'}))$.
\end{lem}
\begin{proof}
The left face of ~\eqref{eqn:Alb-3} commutes by the functoriality of proper
push-forward on the cycle groups. The top (resp. bottom) face commutes by
the universal property of Albanese schemes over $K'$ (resp. $K$).  
To show the commutativity of the back face, we let $[S_{K'}] \in Z_0(S_{K'})$
be one of the two generators of $Z_0(S_{K'})$.
By definition of $\alb_{K'}$, and the transfer maps for the Albanese group schemes,
we see that  $\phi_{A_K *} \circ \alb_{K'}([S_{K'}]) \in A_K(K)$ is the composite map
  \begin{equation}\label{eqn:Alb-4-0}
   S_K \xrightarrow{\sigma^{K'}_K}
    S^d_K(S_{K'}) \xrightarrow{S^d(\alb^1_{K'})} S^d_K(\Alb^1_{K'})
    \xrightarrow{S^d(\alb_{\phi})}
      S^d_K(\Alb^1_K) \xrightarrow{\Sigma_K} \Alb^d_{K} \inj A_K,
  \end{equation}
  where $\Sigma_K$ is the summation map of $A_K$ and
  $d = [K':K]$ (cf. \cite[Proof of Lem.~3.2]{Spiess-Szamuely}).

Since $\alb^1_{K'}$ and  $\alb^1_{K}$ are isomorphisms, the above
  composite map is the same as the composite map
  \[
  S_K \xrightarrow{\sigma^{K'}_K} S^d_K(S_{K'}) \xrightarrow{S^d(\phi)} S^d_K(S_K)
  \xrightarrow{\alb^d_{K}} \Alb^d_K \inj A_K.
  \]
  But the latter map is the same as the map which represents the element
  $\alb_K \circ \phi_*([S_{K'}])$ of $A_K(K)$. This shows that
  the back face of ~\eqref{eqn:Alb-3} commutes.
Since $\alb_{K'}$ is an isomorphism,
  the commutativity of the right face of ~\eqref{eqn:Alb-3} would now follow if we
  show that its front face commutes when restricted to the subgroup
  $\iota'_*(Z_0(S_{K'}))$ of $Z_0(Y')$.

To show the  commutativity of the front face, we look at the diagram
  \begin{equation}\label{eqn:Alb-4-1}
    \xymatrix@C3pc{
      S_K \ar[r]^-{\sigma^{K'}_K} & S^d_K(S_{K'}) \ar[r]^-{S^d(\alb^1_{K'})} &
      S^d_K(\Alb^1_{K'}) \ar[r]^-{S^d(\alb_{\phi})} \ar[d]_-{S^d(\alb_{\iota'})} &
      S^d_K(\Alb^1_K) \ar[d]^-{S^d(\alb_{\iota})} \ar[r]^-{\Sigma_K} &
      \Alb^d_K \ar[d]^-{\alb^d_{\iota}} \\
      & & S^d_K(\Alb^1_{Y'}) \ar[r]^-{S^d(\alb_{\phi_Y})} &
      S^d_K(\Alb^1_Y) \ar[r]^-{\Sigma_Y} & \Alb^d_Y,}
\end{equation}
  where the horizontal arrows in the right-end square are the summation maps
  of the Albanese group schemes.

We now note that as an element of $A_Y(K)$,
  $(\phi_{A_Y *} \circ \alb_{Y'} \circ \iota'_*)([S_{K'}])$ is the composite map
  $\Sigma_Y \circ S^d(\alb_{\phi_Y}) \circ S^d(\alb_{\iota'}) \circ S^d(\alb^1_{K'})
  \circ \sigma^{K'}_K$.
On the other hand, using the commutativity of the top, left and the bottom faces of
  ~\eqref{eqn:Alb-3}, we get that 
$(\alb_Y \circ \phi_{Y *} \circ \iota'_*)([S_{K'}])$ is the
  composite map
  $\alb^d_{\iota} \circ \Sigma_K \circ S^d(\alb_{\phi}) \circ S^d(\alb^1_{K'})
  \circ \sigma^{K'}_K$. 
It suffices therefore to show that the two squares of ~\eqref{eqn:Alb-4-1} commute.
But the right-end square of  ~\eqref{eqn:Alb-4-1}
  commutes by the universality of Albanese group schemes over $K$, and the
  left square commutes because it is obtained by applying the functor
  $S^d_K(-)$ to the front face of ~\eqref{eqn:Alb-2-0}, which commutes by
  \lemref{lem:Alb-2}.  This completes the proof.
\end{proof}

We let $\ov{K}$ denote an algebraic closure of $K$
  and let $K' \subset \ov{K}$ denote the perfect closure of $K$ in $\ov{K}$ so
  that $K'$ is a perfect field (cf. \cite[Lem.~3.16]{Karpilovsky}). 
Let $Y'$ (resp. $\ov{Y}$) denote the base change of $Y$ by $K'$ (resp. $\ov{K}$).
Let $\pi \colon S_{K'} \to S_K$ and $\phi \colon S_{\ov{K}} \to S_{K'}$ denote the
  projections. We shall denote the projections
  $\pi_Y \colon Y' \to Y$ and $\alb_{\pi_Y} \colon \Alb^n_{Y'} \to \Alb^n_Y$
  also by $\pi$. We shall similarly denote the projections
  $\phi_{Y'} \colon \ov{Y} \to Y'$ and
  $\alb_{\phi_{Y'}} \colon \Alb^n_{\ov{Y}} \to \Alb^n_{Y'}$ by $\phi$.

\begin{lem}\label{lem:Alb-5}
  The diagrams
\begin{equation}\label{eqn:Alb-5-5} 
   \xymatrix@C2pc{
     Z_0(Y) \ar[r]^-{\alb_Y} \ar[d]_-{\pi^*} & A_Y(K) \ar[d]^-{\pi^*} & &
     Z_0(Y') \ar[r]^-{\alb_{Y'}} \ar[d]_-{\phi^*} & A_{Y}(K') \ar[d]^-{\phi^*} \\
     Z_0(Y') \ar[r]^-{\alb_{Y'}} & A_Y(K') & & Z_0(\ov{Y}) \ar[r]^-{\alb_{\ov{Y}}} &
     A_Y(\ov{K})}
    \end{equation}
are commutative.
\end{lem}
\begin{proof}
To show the commutativity of the left square,
we let $d \ge 1, \ x \in Y^d_{(0)}$ and  $F = k(x)$. We let
  $T' = S_F \times_{S_K} S_{K'}$.
  Let $\iota \colon S_F \inj Y$ and $\iota' \colon {T'} \inj {Y'}$ be the
  inclusions.
By definition, $(\pi^* \circ \alb_{Y})([x])$ is the composite map of $S_{K'}$-schemes
 \begin{equation}\label{eqn:Alb-5-6} 
   {S_{K'}} \xrightarrow{\sigma^F_K} S^d_K(S_F) \times_{S_K} {S_{K'}}
   \xrightarrow{S^d(\iota)}
   S^d_K(Y) \times_{S_K} {S_{K'}}
   \xrightarrow{(\alb^d_{Y})_{K'}} \Alb^d_{Y'} \inj A_{Y'}.
 \end{equation}
 By \cite[Exp.~V, Prop.~1.9]{SGA1} (see also \cite[Cor.~5.4]{Rydh}), this map 
 agrees with the composite of the following maps of $S_{K'}$-schemes:
 \begin{equation}\label{eqn:Alb-5-7} 
   {S_{K'}} \xrightarrow{\sigma^{T'}_{{K'}}} S^d_{K'}(T')
   \xrightarrow{S^d(\iota')}
   S^d_{K'}(Y') \xrightarrow{\alb^d_{Y'}} \Alb^d_{Y'} \inj A_{Y'}.
 \end{equation}

We now note that $Z := T'_\red$ is the spectrum of a finite field extension $F'$ of
 $K'$ because $K'$ is purely inseparable over $K$. We let $d_1 = [F':K']$ and
let $d_2$ be the multiplicity of $Z$ in the cycle $[T'] \in Z_0(Y')$ so that
$[T'] = d_2[Z]$ and $d = d_1 d_2$. We let $u \colon Z \inj T'$ and
$v \colon Z \inj Y'$ be the inclusions. 
We then get a diagram
\begin{equation}\label{eqn:Alb-5-8} 
  \xymatrix@C1.5pc{
    & \stackrel{d_2}{\underset{i =1}\prod} S^{d_1}_{K'}(Z) \ar[r] \ar@{->>}[d] &
    \stackrel{d_2}{\underset{i =1}\prod} S^{d_1}_{K'}(Y') \ar@{->>}[d] \ar[r] & 
    \stackrel{d_2}{\underset{i =1}\prod} S^{d_1}_{K'}(\Alb^1_{Y'}) \ar@{->>}[d]
    \ar[r]^-{\Sigma_{A_{Y'}}} &  \stackrel{d_2}{\underset{i =1}\prod} \Alb^{d_1}_{Y'}
    \ar[d]^-{\Sigma_{A_{Y'}}} \\
    S_{K'} \ar[ur]^-{\Delta_{\sigma^{F'}_{K'}}} \ar[r]^-{\theta}
    \ar[dr]_-{\sigma^{T'}_{K'}} &
    S^d_{K'}(Z) \ar[r]^-{S^d(v)} \ar@{^{(}->}[d]^-{S^d(u)} & S^d_{K'}(Y')
    \ar[r]^-{S^{d}(\alb^1_{Y'})} &
    S^{d}_{K'}(\Alb^1_{Y'}) \ar[r]^-{\Sigma_{A_{Y'}}}  & \Alb^d_{Y'} \\
    &  S^d_{K'}(T'). \ar[ur]_-{S^d(\iota')} & & &}
\end{equation}

In the above diagram, $\Delta_{\sigma^{F'}_{K'}}$ is the unique map whose composition
with each of the projections $\stackrel{d_2}{\underset{i =1}\prod} S^{d_1}_{K'}(Z)
\to S^{d_1}_{K'}(Z)$ is $\sigma^{F'}_{K'}$. The horizontal arrows on the top row are
obtained by taking the $d_2$-fold product of each of the maps
\[
S^{d_1}_{K'}(Z) \xrightarrow{S^{d_1}(v)} S^{d_1}_{K'}(Y')
\xrightarrow{S^{d_1}(\alb^1_{Y'})}
S^{d_1}_{K'}(\Alb^1_{Y'}) \xrightarrow{\Sigma_{A_{Y'}}} \Alb^{d_1}_{Y'}.
\]
If we let $\psi$ denote the $d_2$-fold product of this composite map, then
one checks that $\Sigma_{A_{Y'}} \circ \psi \circ \Delta_{\sigma^{F'}_{K'}} \colon
S_{K'} \to \Alb^d_{Y'}$ is $(\alb_{Y'} \circ \pi^*)([x]) \in \Alb_Y(K')$.
Since the composite map
$\Sigma_{A_{Y'}} \circ S^{d}(\alb^1_{Y'}) \circ S^d(\iota') \circ \sigma^{T'}_{K'}
\colon S_{K'} \to \Alb^d_{Y'}$ is $(\pi^* \circ \alb_Y)([x])$
(cf. \eqref{eqn:Alb-5-7}), it suffices to
show that all squares and triangles in ~\eqref{eqn:Alb-5-8} commute.

Now, the surjective vertical arrows in \eqref{eqn:Alb-5-8} are the canonical maps
induced by the factorization of the quotient maps from products to
symmetric products. Indeed, for any $W \in \Sch_{K'}$ and
$d = d_1 + \cdots + d_r$, the quotient map $W^d \to S^d_{K'}(W)$ factors through
$S^{d_1}_{K'}(W) \times \cdots \times S^{d_r}_{K'}(W) \to S^d_{K'}(W)$.
The vertical arrows are these latter maps.
In particular, the left and the middle squares
on the top commute. The commutativity of the right square is clear as
$\Sigma_{A_{Y'}}$ is the summation map of $A_{Y'}$. 
It is also clear that $S^d(\iota') \circ S^d(u) = S^d(v)$.
Since $S^d(u)$ is a (split)  nilpotent extension, there is a unique inclusion
$\theta \colon S_{K'} \inj S^d_{K'}(Z)$ such that
$S^d(u) \circ \theta = \sigma^{T'}_{K'}$.
It remains to show that the top left triangle in  ~\eqref{eqn:Alb-5-8} commutes.
As the latter is a diagram in the category $\Sch_{K'}$, it is enough to show that
the structure map $q \colon  S^d_{K'}(Z) \to S_{K'}$ is an isomorphism.

To that end, we note that ${F'}/{K'}$ is a finite separable extension (since $K'$
is perfect) which implies that $\stackrel{d}{\underset{i =1}\prod} Z$
(where the product is taken in $\Sch_{K'}$)
is reduced. As $S^d_{K'}(Z)$ is a quotient of the latter scheme, it is reduced.
Thus, if $S^d_{K'}(Z)$ consists of a single point, then the fact that $\theta$ is a
section of $q$ implies that $q$ is an isomorphism.
It therefore suffices to show that the support of $S^d_{K'}(Z)$ consists of a single
point. But the latter claim follows from \lemref{lem:Alb-6} which proves the
stronger statement that $\Delta_{\sigma^{F'}_{K'}} \colon S_{K'} \to S^{d_1}_{K'}(Z)$ is
bijective. 

We now show that the right square in ~\eqref{eqn:Alb-5-5} commutes.
As in the previous case, we let $d \ge 1, \ x \in Y'^d_{(0)}$ and $F = k(x)$. We let
$\ov{T} = S_F \times_{S_{K'}} S_{\ov{K}}$. Let
$\iota \colon S_F \inj Y'$ and $\ov{\iota} \colon \ov{T} \inj \ov{Y}$
be the inclusions. As before,  $(\phi^* \circ \alb_{Y'})([x])$ is the composite
map of $S_{\ov{K}}$-schemes
\begin{equation}\label{eqn:Alb-5-9} 
   S_{\ov{K}} \xrightarrow{\sigma^{\ov{T}}_{\ov{K}}} S^d_{\ov{K}}(\ov{T})
   \xrightarrow{S^d(\ov{\iota})}
   S^d_{\ov{K}}(\ov{Y}) \xrightarrow{S^d(\alb^1_{\ov{Y}})} S^d_{\ov{K}}(\Alb^1_{\ov{Y}})
   \xrightarrow{\Sigma_{A_{\ov{Y}}}} \Alb^d_{\ov{Y}}.
 \end{equation}

Since $F/{K'}$ is separable of degree $d$, we have that $\ov{T} =
\stackrel{d}{\underset{i =1}\amalg} T_i$, where $T_i \to \Spec(\ov{K})$ is an
  isomorphism for each $i$.
  We now look at the diagram
  \begin{equation}\label{eqn:Alb-5-10} 
    \xymatrix@C3pc{
 & \stackrel{d}{\underset{i =1}\prod} T_i \ar[r] \ar@{^{(}->}[d]^-{u} &
    \ov{Y}^d \ar@{->>}[d] \ar[r] & 
    (\Alb^1_{\ov{Y}})^d \ar@{->>}[d]
    \ar[r]^-{\Sigma_{A_{\ov{Y}}}} &  \stackrel{d}{\underset{i =1}\prod} \Alb^{d}_{\ov{Y}}
    \ar[d]^-{\Sigma_{A_{\ov{Y}}}} \\
    {S}_{\ov{K}} \ar[ur]^-{\Delta} \ar[r]_-{\sigma^{\ov{T}}_{\ov{K}}} &
    S^d_{\ov{K}}(\ov{T}) \ar[r]^-{S^d(\iota')} & S^d_{\ov{K}}(\ov{Y})
    \ar[r]^-{S^{d}(\alb^1_{\ov{Y}})} &
    S^{d}_{\ov{K}}(\Alb^1_{\ov{Y}}) \ar[r]^-{\Sigma_{A_{\ov{Y}}}}  & \Alb^d_{\ov{Y}},} 
\end{equation}     
  where $\Delta$ is the diagonal (using the isomorphism $T_i \to \Spec(\ov{K})$)
  and $u$ is the canonical inclusion.
As in the previous case, the three squares in this diagram commute and
  we only need to show that the left triangle commutes. But this follows from
  \lemref{lem:Alb-7}. This completes the proof of the lemma.
\end{proof}

\begin{lem}\label{lem:Alb-6}
  Let ${E'}/E$ be a finite extension of perfect fields of degree $n$.
  Let $W = \Spec(E)$ and $W' = \Spec(E')$. Then the canonical map
  of $E$-schemes $\sigma^{E'}_E \colon W \to S^n_{E}(W')$ is an isomorphism.
\end{lem}
\begin{proof}
As in the proof of \lemref{lem:Alb-5}, it suffices to show that
$S^n_{E}(W')$ is a singleton set. We may assume that $n \ge 2$, since
there is nothing to prove otherwise. Let $\ov{E}$ be an algebraic closure of
$E$, and let $\Gamma$ denote the absolute Galois group of $E$. We can write
 $\ov{W} := \Spec(E' \otimes_E \ov{E}) =
  \stackrel{n}{\underset{i =1}\amalg} W_i$, where $W_i \to \Spec(\ov{E})$ is an
  isomorphism for each $i$. Let $\Phi = \{W_1, \ldots, W_n\}$. We then compute
  \[
  S^n_{E}(W')(E) =
  \left[\left((S^n_{E}(W')_{\ov{E}}\right)(\ov{E})\right]^{\Gamma} =
  \left[\left(S^n_{\ov{E}}(\ov{W})\right)(\ov{E})\right]^\Gamma =
  \left(S^n(\Phi)\right)^\Gamma = \left(\Phi^n/\Sigma_n\right)^\Gamma,
  \]
  where $\Sigma_n$ is the symmetric group on $n$ letters acting on $\Phi^n$ by
  permuting the coordinates, and the second equality holds by the flat base change
  property of the quotient by finite groups
  (cf. \cite[Lem.~2.2.2]{Abramovich-Vistoli}).

We claim that $\left(\Phi^n/\Sigma_n\right)^\Gamma = [(W_1, \ldots , W_n)]$,
  where we let $[(a_1, \ldots , a_n)]$ denote the $\Sigma_n$-orbit of the point
  $(a_1, \ldots , a_n) \in \Phi^n$. To prove the claim, we let
  $x = (a_1, \ldots , a_n) \in \Phi^n$ be such that $a_i = a_j$ for some $i \neq j$.
  We then have $[x] = [(W_l, W_l, a_3, \ldots , a_n)]$ for some $W_l \in \Phi$.
In this case, we can find $W_{l'} \in \Phi$ which does not lie in the set
  $\{W_l, a_3, \ldots , a_n\}$. Since ${E'}/E$ is finite separable,
  there exists $\sigma \in \Gamma$ such that $\sigma(W_l) = W_{l'}$.
We then get $\sigma([x]) = [(W_{l'}, W_{l'}, \sigma(a_3), \ldots , \sigma(a_n))]
  \neq [(W_{l}, W_{l}, a_3, \ldots , a_n)] = [x]$. That is,
  $[x] \notin \left(\Phi^n/\Sigma_n\right)^\Gamma$.
This proves the claim, and hence the lemma.
\end{proof}

\begin{lem}\label{lem:Alb-7}
  Let $E$ be any field and let $A$ denote the $n$-fold self
  product of $E$. We let $T^n(A) = (A\otimes_E \cdots \otimes_E A)^{\Sigma_n}$,
  where $\Sigma_n$ acts on the $n$-fold
  tensor product of $A$ by permuting the coordinates. Then the canonical map
  $\sigma^A_E \colon T^n(A) \to E$ of $E$-algebras
  has a unique factorization $T^n(A) \xrightarrow{\alpha}
  E \otimes_E \cdots \otimes_E E \xrightarrow{\beta} E$,
  where $\beta$ is the product map and
  $\alpha$ is obtained by precomposing the inclusion
  $T^n(A) \inj A\otimes_E \cdots \otimes_E A$ with the tensor product of the
  projections $p_i \colon A \surj E$.
\end{lem}
\begin{proof}
To prove the lemma, we write $A = Ee_1 \bigoplus \cdots \bigoplus Ee_n$ and let
  $\omega = e_1 \wedge \cdots \wedge e_n \in  \Wedge^n_E(A)$. We now recall from
  \cite[\S~6.3]{Deligne} that
  $\sigma^A_E \colon T^n(A) \to E$ is induced by letting
  $\sigma^A_E(a_1 \otimes \cdots \otimes a_n)$ be the endomorphism of
  the $E$-module $\Wedge^n_E(A) \cong E$, given by
  $(a_1 \otimes \cdots \otimes a_n)(\omega) = (a_1'e_1) \wedge
  \cdots \wedge (a_n'e_n)=(a_1'\cdots a_n')\omega$, where $a_i'=p_i(a_i)$.
 Since $\alpha$ and $\beta$ are given by 
  $\alpha(a_1 \otimes \cdots \otimes a_n) = a_1' \otimes \cdots \otimes a_n'$
  and $\beta(a_1' \otimes \cdots \otimes a_n') 
  = a_1' \cdots a_n'$, we get the desired factorization.
\end{proof}

We can now prove the existence of the Albanese morphism for the Chow group
of 0-cycles over $K$.

\begin{prop}\label{prop:Alb-8}
  The map $\alb_Y \colon Z_0(Y) \to A_Y(K)$ descends to a canonical homomorphism
  $\alb_Y \colon \CH_0(Y) \to A_Y(K)$ which induces a homomorphism
  $\alb_Y \colon A_0(Y) \to \Alb_Y(K)$.
\end{prop}
\begin{proof}
We let $R_0(Y) \subset Z_0(Y)$ be the subgroup of 0-cycles which are rationally
equivalent to zero. We let $\psi \colon \Spec(\ov{K}) \to \Spec(K)$ be the
projection map. Since $\psi^*_Y(R_0(Y)) \subset R_0(\ov{Y})$, and since the
map $\psi^* \colon A_Y(K) \to A_Y(\ov{K}) = A_{\ov{Y}}(\ov{K})$ is injective,
\lemref{lem:Alb-5} reduces the proof of the proposition
to the case when $K$ is algebraically closed.
But this case is well-known (e.g., see \cite[Lem.~3.1]{Spiess-Szamuely}).

The second part of the lemma follows by the commutativity of the diagram
  \begin{equation}\label{eqn:Alb-5-0}
    \xymatrix@C2pc{
      Z_0(Y) \ar[r]^-{\deg} \ar[d]_-{\alb_Y} & \Z \ar@{=}[d] \\
      A_Y(K) \ar[r]^-{\deg} & \Z,}
  \end{equation}
  which is immediate from ~\eqref{eqn:Alb-0-0}
  once we recall that the degree map on the bottom row is defined by letting
  all of $\Alb^n_Y$ map to $n$ for every $n \in \Z$. 
\end{proof}

\subsection{A result on Milne--Tate duality}\label{sec:Prelim-Kai-0}
We return back to our setting of schemes over the local field $k$.
We recall the duality theorem of Milne for abelian varieties over $k$
whose characteristic zero analogue is a classical result of Tate.
Let $A$ be an abelian variety over $k$ and let $A^t = \Picc^0(A)$ be the
corresponding dual abelian variety (which exists). A theorem
of Milne (cf. \cite[Thm.~7.8]{Milne-Duality}) says that there is a perfect
continuous pairing of topological abelian groups
\begin{equation}\label{eqn:MD-0}
  {\<,\>} \colon A(k) \times H^1_\et(k, A^t) \to H^2_\et(k, \G_m) \cong {\Q}/{\Z},
\end{equation}
where $A(k)$ is endowed with the adic topology (which is profinite) and
$H^1_\et(k, A^t)$ (and $H^2_\et(k, \G_m)$) with the discrete topology.
We shall need the following two naturality properties of this pairing.

\begin{prop}\label{prop:Milne-1}
  Let $f \colon A \to B$ be a morphism of abelian varieties over $k$. Then the
  diagram
  \begin{equation}\label{eqn:MD-1}
\xymatrix@C1pc{
B(k) \times H^1_\et(k, B^t)  \ar@<6ex>[d]^-{f^t_*}
\ar[r]^-{\<,\>} & \Br(k)  \ar@{=}[d] \\
 A(k) \times H^1_\et(k, A^t) \ar@<7ex>[u]^-{f_*}
 \ar[r]^-{\<,\>} & \Br(k)}
\end{equation}
  is commutative.
\end{prop}
\begin{proof}
This follows from a routine verification using the construction of the pairing
~\eqref{eqn:MD-0} and the universal property of the Poincar{\'e} sheaf for abelian
varieties. We recall the argument here. Let $\sP_A$ denote the Poincar{\'e}
sheaf on $A \times A^t$ (cf. \cite[Exc.~9.4.3]{Kleiman}).
 By definition, it defines a canonical map of fppf sheaves of abelian groups
  $\phi_{A} \colon A \to {\sE}xt^1(A^t, \G_m)$ on $\Spec(k)$. In particular,
 we get a map $\phi_A \colon A(k) \to \Hom_{\sD_\et(k)}(A^t, \G_m[1])$
 ($\sD_\et(k)$ is defined in \S~\ref{sec:p-adic-coh}).
We thus get canonical maps
  \begin{equation}\label{eqn:MD-2}
  \xymatrix@C1pc{  
   & A(k) \times H^1_\et(k, A^t) \ar[r]^-{\cong} \ar[dl]_-{\<,\>} & A(k) \times
   \Hom_{\sD_\et(k)}(\Z, A^t[1]) \ar[d]^-{\phi_A \times \id} \\
 H^2_\et(k, \G_m) & \Hom_{\sD_\et(k)}(\Z, \G_m[2]) \ar[l]_-{\cong} & 
   \Hom_{\sD_\et(k)}(A^t, \G_m[1]) \times  \Hom_{\sD_\et(k)}(\Z, A^t[1]) \ar[l],}
  \end{equation}
  where the unnamed bottom horizontal arrow takes a pair $(h,g)$ to
  $h[1] \circ g$.

We now let $u \colon \Spec(k) \to A$ be an element of $A(k)$ and let
  $v \colon \Z \to B^t[1]$ be an element of $H^1_\et(k, B^t)$. 
  From the universality of the  Poincar{\'e} sheaf, one knows that on $B^t$,
  we have $(u \times \id_{B^t})^* \circ (f \times \id_{B^t})^* (\sP_B) \cong
  (\id_{\Spec(k)} \times f^t)^* \circ (u \times \id_{A^t})^*(\sP_A) \cong
  (u \times f^t)^*(\sP_A)$.
  If we let $\phi_A(u) \colon A^t \to \G_m[1]$ be the map induced by $u$,
  then one checks using ~\eqref{eqn:MD-2} that $\<f_*(u), v\>$ is the composition
  $\Z \xrightarrow{v} B^t[1] \xrightarrow{(\phi_A(u) \circ f^t)[1]} \G_m[2]$.
  On the other hand, $\<u, f^t_*(v)\>$ is the composition
  $\Z \xrightarrow{v} B^t[1] \xrightarrow{f^t[1]} A^t[1] \xrightarrow{\phi_A(u)[1]}
  \G_m[2]$. It follows that $\<f_*(u), v\> = \<u, f^t_*(v)\>$.
\end{proof}

Let $K \subset K'$ be two finite field extensions of $k$.
As before, we let $S_K = \Spec(K), \ S_{K'} = \Spec(K')$ and 
  $\phi \colon S_{K'} \to S_K$ the projection map. Let $A$ be an abelian
  variety over $K$ and let $A^t$ denote the dual abelian variety.
  Let $A'$ denote the base change of $A$ by $\phi$. Using the universal property
  of Picard variety, one knows that the base change of $A^t$ by $\phi$  is
  canonically identified with $A'^t$.
  In particular, the finite push-forward map on Picard groups
  induces a push-forward map $\phi_* \colon A^t(K') \to A^t(K)$. This can also be
  seen using the fact that a finite push-forward map preserves algebraic equivalence
  of codimension one cycles (cf. \cite[Prop.~10.3]{Fulton}).
  
\begin{lem}\label{lem:Milne-2}
 The diagram
  \begin{equation}\label{eqn:Milne-2-0}
    \xymatrix@C1pc{
      A^t(K) \times H^1_\et(K, A)  \ar@<6ex>[d]^-{\phi^*}
\ar[r]^-{\<,\>} & \Br(K)  \\
 A^t(K') \times H^1_\et(K', A) \ar@<7ex>[u]^-{\phi_*}
 \ar[r]^-{\<,\>} & \Br(K') \ar[u]_-{\phi_*}} 
\end{equation}
  is commutative.
\end{lem}
\begin{proof}
Let $u \in  H^1_\et(K, A)$ and $v \in A^t(K')$. Then $\phi_* (\<v, \phi^*(u)\>)$
  is described as follows: $u$ represents a morphism 
  $u \colon \Z_K \to A[1]$ in $\sD_\et(k)$ and $v$ represents a $\G_m$-torsor on
  $A' = A_{K'}$. The latter induces a unique morphism
  $v \colon \phi^*(A) \to {\G_m}_{K'}[1]$. Precomposing $v[1]$ with
  $\phi^*(u) \colon \Z_{K'} \cong \phi^*(\Z_K) \to \phi^*(A)[1]$,
  we get a morphism $\Z_{K'} \xrightarrow{v[1] \circ \phi^*(u)}  {\G_m}_{K'}[2]$.
  This morphism is represented by the element $\<v, \phi^*(u)\> \in \Br(K')$.

  Applying $\phi_*$ to $v[1] \circ \phi^*(u)$, it follows that
  $\phi_* (\<v, \phi^*(u)\>)$ is the composite morphism
  \begin{equation}\label{eqn:Milne-2-1}
    \Z_K \xrightarrow{1_{\Z}} \phi_* \phi^*(\Z_{K}) \xrightarrow{\phi_* \phi^*(u)}
  \phi_* \phi^*(A)[1]
  \xrightarrow{\phi_*(v)[1]} \phi_*({\G_m}_{K'})[2] \xrightarrow{N_{{K'}/K}[2]}
  {\G_m}_K[2],
  \end{equation}
  where $N_{{K'}/K} \colon \phi_*({\G_m}_{K'}) \to {\G_m}_K$ is the norm map and
  $1_\Z$ is the unit of adjunction.

A similar computation shows that $\phi_*(v)$ represents the composite 
  $A \xrightarrow{1_A} \phi_* \phi^*(A) \xrightarrow{\phi_*(v)} \phi_*({\G_m}_{K'})[1]
  \xrightarrow{N_{{K'}/K}[1]} {\G_m}_K[1]$ in which $1_A$ is the unit of
  adjunction for $A$. It follows that $\<\phi_*(v), u\>$ represents the composite
  morphism
  \begin{equation}\label{eqn:Milne-2-2}
    \Z_K \xrightarrow{u} A[1] \xrightarrow{1_A[1]}
    \phi_* \phi^*(A)[1] \xrightarrow{\phi_*(v)[1]}
    \phi_*({\G_m}_{K'})[2] \xrightarrow{N_{{K'}/K}} {\G_m}_K[2].
  \end{equation}
By comparing ~\eqref{eqn:Milne-2-1} and ~\eqref{eqn:Milne-2-2}, and using the
  commutative diagram
 \begin{equation}\label{eqn:Milne-2-3}
   \xymatrix@C1pc{
     \Z \ar[r]^-{1_\Z} \ar[d]_-{u} & \phi_*\phi^*(\Z) \ar[d]^-{\phi_*\phi^*(u)} \\
     A[1] \ar[r]^-{1_A[1]} & \phi_*\phi^*(A)[1],}
   \end{equation}
 we get $\<\phi_*(v), u\> = \phi_* (\<v, \phi^*(u)\>)$, as desired.
\end{proof}

\section{Albanese cokernel II}\label{sec:Kai-2}
 In this section, we continue to work in the setting of \S~\ref{sec:Kai}
and prove \thmref{thm:Kai-Main}. Since we follow the strategy
employed by Kai \cite{Kai} to prove the analogous result in characteristic zero,
we briefly describe the main steps of his proof and explain
what we do about them in positive characteristic.

Kai's proof consisted of the following three main steps:
\begin{enumerate}
\item
Theorem~3.5 of op. cit., a result of Saito--Sato \cite{Saito-Sato-ENS}.
\item
  Proposition~3.4 of op. cit, which establishes the compatibility of the
  Brauer--Manin
pairing with the Tate--Milne pairing in characteristic zero.
\item
  Proposition~3.9 of op. cit., which gives an injectivity of the specialization
  map for the relative Brauer group.
\end{enumerate}

The proof of Step~3 in op. cit. remains valid in positive characteristic, but this
is not the case for Steps~1 and ~2. Step~1 in positive characteristic is
established by \thmref{thm:Main-3} of this paper. Kai's proof of Step~2 in
characteristic zero is quite intricate, and we were unable to adapt it to
positive characteristic. We therefore give an entirely different, shorter, and
characteristic-free proof of this step. This constitutes the main technical
result of this section (cf. \thmref{thm:Milne-4}), after which we complete the
proof of \thmref{thm:Kai-Main}.

We begin with some results concerning the Picard scheme of $\sX$.
These results are necessary in order to adapt Step~3 above to positive
characteristic and to prove \thmref{thm:Milne-4}.

\subsection{Picard and Brauer groups}\label{sec:Prelim-Kai}
If $F$ is a field (resp. dvr) and $W$ is a smooth projective $F$-scheme which is
geometrically integral (resp., has  geometrically integral fibers),
then one knows that the Picard scheme $\Picc(W)$ exists as a group scheme over
$F$ which is separated and locally of finite type over $F$
(cf. \cite[Thm.~9.4.8]{Kleiman}). If $F$ is a field, we let $\Picc^0(W)$ denote
the identity component of $\Picc(W)$, called the Picard variety of $W$.
We begin with the following folklore result about the
Picard scheme.

\begin{lem}\label{lem:Pic-sch-0}
  Let $F$ be any field and let $W$ be a smooth projective geometrically integral
  $k$-scheme. Then each connected component of $\Picc(W)$ is an irreducible
  projective $F$-scheme.
  \end{lem}
  \begin{proof}
    We let $G$ be any connected component of $\Picc(W)$ and let $x \in G_{(0)}$.
    By \cite[Cor.~9.4.18.3]{Kleiman}, $G$ is a quasi-projective $F$-scheme.
    We let $F' = k(x)$ and let $G_{F'} = G_1 \amalg \cdots \amalg G_r \subset
    \Picc(W)_{F'} \cong \Picc(W_{F'})$. We let $x' \in G_{F'}(F')$  map to $x$
    under the projection $\Picc(W_{F'}) \to \Picc(W)$.
    We can assume without loss of generality that $x' \in G_1$. In this case,
    the translation map $\tau_{x'} \colon \Picc(W_{F'}) \to \Picc(W_{F'})$ defines
    an isomorphism $\Picc^0(W_{F'}) \xrightarrow{\cong} G_1$.

    Since $\Picc^0(W)$ is a geometrically irreducible projective $F$-scheme by
    \cite[Thm.~9.5.4]{Kleiman},
    it follows that the same is true for $G_1$ as an $F'$-scheme.
    In particular, $G_1$ is a geometrically irreducible projective $F'$-scheme.
    Since $G_{F'} \to G$ is a finite flat map, we get that $G_1 \to G$ is a finite
    flat map of $F$-schemes. In particular, it is surjective.
    Since $G$ is a separated $F$-scheme, it follows that
    $G$ is irreducible and a projective $F$-scheme. 
\end{proof}

  Recall that $\Picc(\sX)$ is a locally of finite type smooth $S$-scheme
  under the assumption of \thmref{thm:Kai-Main}.
  Since the identity section of $\Picc(\sX) \to S$ is connected and $\Picc(\sX)$
  is locally Noetherian (and hence locally connected, see \cite[Tag 04MF]{SP}),
there is a unique connected component
  $\Picc^0(\sX)$ of $\Picc(\sX)$ which is closed and open in $\Picc(\sX)$
  (cf. \cite[Cor.~6.1.9]{EGA1}) and contains latter's identity section.

\begin{prop}\label{prop:Pic-var}
  $\Picc^0(\sX)$ is a smooth projective $S$-group scheme whose generic and special
  fiber coincide with the Picard varieties $\Picc^0(X)$ and $\Picc^0(Y)$,
  respectively.
\end{prop}
\begin{proof}
  Since $\Picc^0(\sX)$ is a connected component $\Picc(\sX)$, it is locally of
  finite type smooth $S$-scheme. In particular, its local rings
are integral domains. It follows from \cite[Prop.~6.1.10]{EGA1} that
$\Picc^0(\sX)$ is irreducible. This implies that the generic fiber of
$\Picc^0(\sX)$ is
  irreducible. Since this generic fiber contains the identity section of
  $\Picc(X)$, it follows from \cite[Prop.~9.5.3]{Kleiman} that $\Picc^0(\sX)_\eta$
  coincides with $\Picc^0(X)$. If we let
  $\{G_i\}_{i \in I}$ be the collection of connected components of $\Picc^0(\sX)_s$,
  then each $G_i$ is also a connected component of $\Picc(Y)$ and one of these
  (say, $G_1$) contains latter's identity section. In particular,
  $G_1 = \Picc^0(Y)$. If we let $Z$ be the union of all
  other components of $\Picc^0(\sX)_s$, then $Z$ is closed in $\Picc^0(\sX)$.
  Letting $\Picc^0(\sX)' = \Picc^0(\sX) \setminus Z$, we get that $\Picc^0(\sX)'$
  is an open subscheme of $\Picc^0(\sX)$ and hence of $\Picc(\sX)$.

Since $\Picc^0(\sX)'_\eta = \Picc^0(X)$ and $\Picc^0(\sX)'_s = \Picc^0(Y)$ are
  quasi-compact, it follows that $\Picc^0(\sX)'$ is quasi-compact.
  Since it is locally of finite type
  over $S$, one deduces that $\Picc^0(\sX)'$ is an open subscheme of $\Picc(\sX)$
  which is of finite type over $S$. We conclude from \cite[Exc.~9.4.11]{Kleiman}
  that $\Picc^0(\sX)'$ is an irreducible and quasi-projective smooth $S$-scheme.
  In particular,
  \[
  \dim_k(\Picc^0(X)) = \dim_k(\Picc^0(\sX)'_\eta) =
  \dim_{\F}(\Picc^0(\sX)'_s) = \dim_{\F}(\Picc^0(Y)).
  \]

Since these Picard varieties
  are smooth over their respective base fields, we can apply
  \cite[Prop.~9.5.20]{Kleiman} to conclude that $\Picc^0(\sX)'$ is an open and
  closed group subscheme of $\Picc(\sX)$ which is proper (and hence projective)
  over $S$. Since this is also irreducible, it is  in fact a connected component
  of $\Picc(\sX)$.
  As $\Picc^0(\sX)$ is a connected component of $\Picc(\sX)$ which contains
  $\Picc^0(\sX)'$, we deduce that $\Picc^0(\sX)' = \Picc^0(\sX)$. This implies,
  in particular, that $\Picc^0(\sX)$ is a connected, smooth and projective
  group subscheme of $\Picc(\sX)$ whose generic and special fibers are the
  Picard varieties of the corresponding fibers of $\sX$.
  \end{proof}

Recall that if $W$ is a geometrically integral smooth projective variety over a
field $F$, then $\Pic(W) \subset \Picc(W)(F)$, and
one lets $\Pic^0(W) := \Pic(W) \bigcap \Picc^0(W)(F)$ and
the N{\'e}ron--Severi group $\NS(W) := {\Pic(W)}/{\Pic^0(W)}$.

\begin{defn}\label{defn:Pic-sch-1}
We let $\Pic^0(\sX) = \Pic(\sX) \bigcap \Picc^0(\sX)(S) = \Picc^0(\sX)(S)$,
  where the second equality holds because $\Pic(\sX) = \Picc(\sX)(S)$
  as $\Br(S) = 0 = \Pic(S)$ (cf. \cite[\S~2.5. p.~66]{CTS}).
\end{defn}

\begin{lem}\label{lem:Pic-sch-3}
  The restriction map $\iota^* \colon \Pic(\sX^s) \to \Pic(\ov{Y})$ induces the
  specialization maps
  \begin{equation}\label{eqn:Pic-sch-3-0}
  \spc \colon \Pic(X^s) \to \Pic(\ov{Y}) \ \  \mbox{and} \ \ \spc \colon
  \NS(X^s) \to \NS(\ov{Y})
  \end{equation}
  which are morphisms of discrete $\Gamma$-modules.
\end{lem}
\begin{proof}
The first map is canonically induced by the isomorphism
  $\Pic(\sX^s) \xrightarrow{\cong} \Pic(X^s)$. To prove the existence of the
  second map in ~\eqref{eqn:Pic-sch-3-0},
  it suffices to show that the first map sends $\Pic^0(X^s)$ into
  $\Pic^0(\ov{Y})$. For this, we first note that the specialization map
  $\spc \colon \Pic(X^s) \to \Pic(\ov{Y})$ is the direct limit of the maps
  $\Pic(X_{k'}) \to \Pic(Y_{\F'})$, where $k' \in \Sigma^s_k$ and $\F'$ is the
  residue field of $k'$. We now note that any element of $\Pic^0(X^s)$ is an
  element of $\Pic(X^s)$ which is also an element of $\Picc^0(X^s)(k_s)$.
  In particular, it is an element of $\Pic(X_{k'}) \bigcap \Picc^0(X_{k'})(k')$
  for some $k' \in \Sigma^s_k$ by \lemref{lem:Pic-sch-2}.
Note here that $\Picc^0(X^s) = \Picc^0(X)^s$.
It suffices therefore to show that the specialization map
  $\spc \colon \Pic(X_{k'}) \to \Pic(Y_{\F'})$ sends $\Pic^0(X_{k'})$ to
  $\Pic^0(Y_{\F'})$ for every $k' \in \Sigma^s_k$ whose residue field is $\F'$.
To prove this, we can assume without loss of generality that $k' = k$.

We now let $\alpha \in \Pic(X) \bigcap \Picc^0(X)(k)$. By \propref{prop:Pic-var}
and the isomorphism $\Pic(\sX) \xrightarrow{\cong} \Pic(X)$, we see that
$\Pic(\sX) \bigcap \Picc^0(\sX)(S) \to \Pic(X) \bigcap \Picc^0(X)(k)$ is a
bijection. Meanwhile, \propref{prop:Pic-var} also implies that
$\iota^* \colon \Pic(\sX) = \Picc(\sX)(S) \to \Picc(Y)(\F) =
\Pic(Y)$ sends $\Picc^0(\sX)(S)$ into $\Picc^0(Y)(\F)$. This implies that
$\spc(\alpha) \in \Pic^0(Y)$, as desired. It is clear from their construction
that the two maps in ~\eqref{eqn:Pic-sch-3-0} are $\Gamma$-equivariant.
\end{proof}

From the proof of \lemref{lem:Pic-sch-3}, we get a commutative diagram
of exact sequences
\begin{equation}\label{eqn:Pic-sch-4}
  \xymatrix@C1.7pc{
    0 \ar[r] & \Pic^0(X^s) \ar[r] \ar[d] & \Pic(X^s) \ar[r] \ar[d] &
    \NS(X^s) \ar[r] \ar[d] & 0 \\
    0 \ar[r] &  \Pic^0(\ov{Y}) \ar[r] & \Pic(\ov{Y}) \ar[r] &
    \NS(\ov{Y}) \ar[r] & 0.}
\end{equation}
Furthermore, the smoothness of $\Picc^0(\sX)$ and $\Picc(\sX)$ over $S$
(cf. \propref{prop:Pic-var}) implies that the left and the
middle vertical arrows (and hence the right vertical arrow) in this diagram are
surjective.

\begin{lem}\label{lem:Pic-sch-2}
  The cohomological and Azumaya Brauer group of $\sX^s$ coincide and each of them
  is canonically isomorphic to ${\varinjlim}_{k' \in \Sigma^s_k} \Br(\sX_{\sO_{k'}})$.
  This isomorphism also holds for $\Pic(\sX^s)$.
\end{lem}
\begin{proof}
  The first statement is well-known and the other statements follow
  if we show that the canonical map ${\varinjlim}_{k' \in \Sigma^s_k}
  H^i_\et(\sX_{\sO_{k'}}, \G_m) \to H^i_\et(\sX^s, \G_m)$ is an isomorphism for every
  $i \in \N_0$. But this is a result of EGA (cf. \cite[\S~2.2.2]{CTS}).
\end{proof}

\subsection{Construction of the key map}\label{sec:Key-map}
In this subsection, we construct the key exact exact sequence which will be
used in the proof of \thmref{thm:Kai-Main}.

For any map of schemes $f \colon W_1 \to W_2$, we let
$\Br({W_1}/{W_2}) = \Ker(f^* \colon \Br(W_2) \to \Br(W_1))$. If $W_2 = \Spec(A)$,
we shall write $\Br({W_1}/{W_2})$ also as $\Br({W_1}/{A})$. Using the
Hochschild-Serre spectral sequence
\begin{equation}\label{eqn:HS-SS}
  E^{i,j}_2 = H^i_\et(\Gamma, H^j_\et(X^s, \G_m)) \Rightarrow H^{i+j}_\et(X, \G_m),
  \end{equation}
and the fact that $H^3_\et(k, \G_m) = 0$ (cf. \cite[Cor.~7.2.2]{NSW}), one gets
natural exact sequences
\begin{equation}\label{eqn:Pic-sch-HS-SS-0}
  0 \to \Br({X^s}/X) \to \Br(X) \to \Br(X^s)^\Gamma \to H^2_\et(k, \Pic(X^s));
\end{equation}
\begin{equation}\label{eqn:Pic-sch-HS-SS-1}
0 \to \Br_0(X)  \to \Br({X^s}/X) \to H^1_\et(k, \Pic(X^s)) \to 0.
\end{equation}

From the second exact sequence, one obtains a canonical isomorphism
\begin{equation}\label{eqn:Pic-sch-HS-SS-2}
  \phi_{{X^s}/X} \colon \frac{\Br({X^s}/X)}{\Br_0(X)} \xrightarrow{\cong}
  H^1_\et(k, \Pic(X^s)).
\end{equation}
We let $\psi_{{X^s}/X} = \phi^{-1}_{{X^s}/X}$ and let $\delta_X$ denote the
composite 
\begin{equation}\label{eqn:Pic-sch-HS-SS-3}
  H^1_\et(k, \Pic^0(X^s)) \xrightarrow{\kappa_X}
  H^1_\et(k, \Pic(X^s)) \xrightarrow{\psi_{{X^s}/X}}
  \frac{\Br({X^s}/X)}{\Br_0(X)} \to  \wt{\Br}(X).
\end{equation}

In the commutative diagram (cf. \propref{prop:Kato-Br-0})
\begin{equation}\label{eqn:Pic-sch-HS-SS-4}
  \xymatrix@C1.2pc{
    & & \Br(k) \ar[r]^-{\inv_k}_-{\cong} \ar[d] & H^1_\et(\F)
    \ar@{^{(}->}[d]^-{\wt{g}^*} \\
    0 \ar[r] & \Br(\sX) \ar[r] & \Fil_0 \Br(X) \ar[r]^-{\partial_X} & H^1_\et(Y),}
  \end{equation}
the right vertical arrow is injective because $\sX$ is smooth and
$Y$ is geometrically integral.
It follows that $\Br(k) \to \Br(X)$ is injective and its image $\Br_0(X)$
has trivial intersection with $\Br(\sX)$. In particular,
the kernel of the last map in ~\eqref{eqn:Pic-sch-HS-SS-3} is
$\Br({X^s}/{\sX}) = \Br({\sX^s}/{\sX})$, where the latter equality follows from
\lemref{lem:Pic-sch-2} since the map $\Br(\sX_{\sO_{k'}}) \to \Br(X_{k'})$ is
injective for every $k' \in \Sigma^s_k$.

Combining the above discussion with the isomorphism
\[
\coker\left(\Pic(X^s)^\Gamma \to \NS(X^s)^\Gamma\right) \xrightarrow{\cong}
\Ker\left(H^1_\et(k, \Pic^0(X^s)) \to  H^1_\et(k, \Pic(X^s))\right),
\]
obtained by applying the
Galois cohomology to the exact sequence of $\Gamma$-modules
\begin{equation}\label{eqn:Pic-sch-HS-SS-5}
  0 \to \Pic^0(X^s) \to \Pic(X^s) \xrightarrow{\tau_X} \NS(X^s) \to 0,
  \end{equation}
~\eqref{eqn:Pic-sch-HS-SS-3} yields
an exact sequence
\begin{equation}\label{eqn:Pic-sch-HS-SS-6}
0 \to \coker(\Pic(X^s)^\Gamma \to \NS(X^s)^\Gamma) \to
\Ker(\delta_X) \xrightarrow{\theta} \Br({\sX^s}/{\sX}).
\end{equation}

We now look at the diagram
\begin{equation}\label{eqn:Pic-sch-HS-SS-7}
  \xymatrix@C1pc{
    H^1_\et(k, \Pic^0(X^s)) \ar[rr]^-{\kappa_X} \ar@{}[drr] |{{\rm{A}}} &  &
    H^1_\et(k, \Pic(X^s)) \ar[r] & H^1_\et(k, \NS(X^s)) \ar[d]^-{\spc} \\
    \Ker(\delta_X) \ar@{^{(}->}[u] \ar[r]^-{\theta} & \Br({\sX^s}/{\sX})
    \ar[d]_-{\iota^*} 
    \ar@{^{(}->}[r] & \frac{\Br({X^s}/X)}{\Br_0(X)} \ar[u]^-{\phi_{{X^s}/X}}_-{\cong}
    \ar@{}[dr] |{{\rm{B}}} &  H^1_\et(k, \NS(\ov{Y})) \\
    & \Br({\ov{Y}}/Y) \ar[r]^-{\phi_{{\ov{Y}}/Y}} & H^1_\et(\F, \Pic(\ov{Y}))
    \ar[r] & H^1_\et(\F, \NS(\ov{Y})) \ar[u]_-{\rho^*},}
  \end{equation}
where the right vertical arrow on the top level is induced by
\lemref{lem:Pic-sch-3}, and on the bottom level is the restriction map
$\rho^* \colon  H^1_\et(\wh{\Z}, \NS(\ov{Y})) \to  H^1_\et(\Gamma, \NS(\ov{Y}))$,
induced by the quotient map $\rho \colon \Gamma \surj \wh{\Z}$.
The top row is induced by ~\eqref{eqn:Pic-sch-HS-SS-5}.

\begin{lem}\label{lem:Pic-sch-4}
  The above diagram is commutative.
\end{lem}
\begin{proof}
Having the results of \S~\ref{sec:Prelim-Kai} at hand, the proof of the lemma
  is identical to that of \cite[Prop.~3.7]{Kai}, and we therefore omit it.
  \end{proof}

The following is Step~3 that we described in the beginning of \S~\ref{sec:Kai-2}.

\begin{lem}\label{lem:Br-res-spl}
  The map  ${\iota^*} \colon \Br({\sX^s}/{\sX}) \to \Br({\ov{Y}}/Y)$
  is injective.
\end{lem}
\begin{proof}
  The proof is identical to that of \cite[Prop.~3.9]{Kai}, and we therefore
  omit it.
  \end{proof}

\subsection{Compatibility of Brauer--Manin and Tate--Milne pairings}
\label{sec:Key-map-0}
We shall now show the compatibility between the Brauer--Manin pairing
(cf. ~\eqref{eqn:BMP-1}) and the Tate--Milne pairing (cf. ~\eqref{eqn:MD-0}).
This is the last technical step for proving \thmref{thm:Kai-Main}.
We need the following lemma.

Let $K \subset K'$ be two finite field extensions of $k$ and let
$\phi \colon S_{K'} = \Spec(K') \to \Spec(K) = S_K$ be the projection. Let
$Y = X_K$ and $Y' = X_{K'}$. We let $\Br'(Y) = {\Br(Y)}/{\Br_0(Y)}$.
    We let $G(Y) = \Hom(H^1_\et(K, \Picc^0(Y)), \Br(K))$ and
    $G(Y') = \Hom(H^1_\et(K', \Picc^0(Y')), \Br(K'))$.
    Let $(\phi^*)^\vee \colon \Hom(\Br'(Y'), \Br(K')) \to \Hom(\Br'(Y), \Br(K))$
    be the map which sends $f$ to the composite 
    $\Br'(Y) \xrightarrow{\phi^*} \Br'(Y') \xrightarrow{f} \Br(K')
      \xrightarrow{\phi_*} \Br(K)$.
      Let $(\phi^*)^\vee \colon G(Y') \to G(Y)$ be the map which sends $g$
      to the composite $H^1_\et(K, \Picc^0(Y)) \xrightarrow{\phi^*}
      H^1_\et(K', \Picc^0(Y')) \xrightarrow{g} \Br(K') \xrightarrow{\phi_*}
      \Br(K)$.

We consider the diagram
 \begin{equation}\label{eqn:Milne-2-4}
   \xymatrix@C.6pc{   
     A_0(Y') \ar[rr] \ar[dd] \ar[dr]_-{\phi_*}
     & & \Hom(\Br'(Y'), \Br(K')) \ar[dr]^-{(\phi^*)^{\vee}} \ar[dd] & \\
     & A_0(Y) \ar[dd] \ar[rr] & & \Hom(\Br'(Y), \Br(K)) \ar[dd] \\
     \Alb_Y(K') \ar[rr] \ar[dr]_-{\phi_*} & &
     G(Y') \ar[dr]^-{(\phi^*)^{\vee}} & \\
     & \Alb_Y(K) \ar[rr] & & G(Y)}
 \end{equation}
 in which the horizontal arrows are induced by the Brauer--Manin pairing
 ~\eqref{eqn:BMP-1}  on the top face and the Tate--Milne pairing ~\eqref{eqn:MD-0}
 on the bottom face. The vertical
 arrows on the left face are induced by the Albanese maps
 while on the right face, they are given by
 $f \mapsto f \circ \psi_{{Z^s}/Z} \circ \kappa_Z$ for $Z \in \{Y, Y'\}$
 (cf. ~\eqref{eqn:Pic-sch-HS-SS-3}). The slanted arrows are
 induced by $\phi$. The map $\phi_* \colon \Alb_Y(K') \to \Alb_Y(K)$
 is the transfer map for the abelian variety $\Alb_Y$ defined above.

\begin{lem}\label{lem:Milne-3}
 All faces of the diagram~\eqref{eqn:Milne-2-4} except possibly the front and the
 back faces commute.
\end{lem}
\begin{proof}
The top face of ~\eqref{eqn:Milne-2-4}
  commutes by the functoriality of the Brauer--Manin pairing with
  respect to proper maps (cf. \S~\ref{sec:BMP}).
  The right face commutes because ~\eqref{eqn:Pic-sch-HS-SS-2} is clearly
  functorial with respect to the pull-back maps induced by $\phi$. The bottom face
  commutes by applying \lemref{lem:Milne-2} to $A = \Picc^0(Y)$.
  To show that the left face commutes, it suffices to show that the diagram
  \begin{equation}\label{eqn:Milne-3-0}
   \xymatrix@C1pc{ 
     Z_0(Y') \ar[r]^-{\alb_{Y'}} \ar[d]_-{\phi_*} & A_Y(K') \ar[d]^-{\phi_*} \\
     Z_0(Y) \ar[r]^-{\alb_Y} & A_Y(K)}
  \end{equation}
  commutes, where $A_Y$ is the (total) Albanese group scheme of $Y$ and
  $\phi_*$ is the transfer map of $A_Y$ (cf. \cite[\S~2.1]{Kai}).

We let $x' \in Y'_{(0)}$ and let $x$ be its image in $Y$. We need to
  show that $\phi_* \circ \alb_{Y'}([x']) = \alb_Y \circ \phi_*([x'])$ in
  ~\eqref{eqn:Milne-3-0}. To that end, we let $F'$ (resp. $F$)
  denote the residue field of $x'$ (resp. $x$), and let $T = \Spec(F)$ and
  $T' = \Spec(F')$. We let $\iota \colon T \inj Y_F$ and
  $\iota' \colon T' \inj Y_{F'}$ denote the inclusions.
We consider the diagram
  \begin{equation}\label{eqn:Milne-3-1}
   \xymatrix@C1pc{ 
     Z_0(T') \ar[rr]^-{\alb_{T'}}_-{\cong} \ar[dr] \ar[dd] & & A_{T'}(F')
     \ar[dr]^-{\psi_{x'}}
     \ar[dd] & \\
     & Z_0(Y') \ar[rr] \ar[dd]^->>>>>>{\phi_*} & & A_Y(K') \ar[dd]^-{\phi_*} \\
     Z_0(T) \ar[rr]^->>>>>>{\alb_T} \ar[dr] & & A_T(F) \ar[dr]^-{\psi_x} & \\
     & Z_0(Y) \ar[rr]^-{\alb_Y} & & A_Y(K),}
   \end{equation}
  in which the vertical and the slanted arrows are the push-forward or the transfer
  maps and the horizontal arrows are the Albanese homomorphisms.

The left face of ~\eqref{eqn:Milne-3-1}
  commutes by the functoriality of proper push-forward
  on the cycle groups. The right face is the composition of the two squares in the
  diagram
  \begin{equation}\label{eqn:Milne-3-2}
   \xymatrix@C1pc{
     A_{T'}(F')  \ar[r]^-{\iota'_*} \ar[d] & A_Y(F') \ar[r]
     \ar[d] & A_Y(K') \ar[d] \\
     A_T(F) \ar[r]^-{\iota_*} & A_Y(F) \ar[r] & A_Y(K),}
  \end{equation}
  whose all arrows are transfer and push-forward maps. 
  The left square of this diagram commutes by \lemref{lem:Alb-4}
  since it is a special case of the right face of the diagram
  ~\eqref{eqn:Alb-3}. The right square of ~\eqref{eqn:Milne-3-2} commutes by the
  composition law of the transfer map for the total Albanese group scheme.
It follows that the right face of ~\eqref{eqn:Milne-3-1} commutes.
  
The top face  of ~\eqref{eqn:Milne-3-1} commutes by the definition of $\alb_{Y'}$
  and $\psi_{x'}$. By the same reason, the bottom face commutes.
  To show that $\phi_* \circ \alb_{Y'}([x']) = \alb_Y \circ \phi_*([x'])$ in
  ~\eqref{eqn:Milne-3-0}, it remains therefore to show
  that the back face of ~\eqref{eqn:Milne-3-1} commutes. But this follows from
  \lemref{lem:Alb-4} since it is a special case of the back face of the diagram
  ~\eqref{eqn:Alb-3}. This completes the proof.
\end{proof}

The heart of the proof of \thmref{thm:Kai-Main} lies in \thmref{thm:Main-3},
and the following result, which is of independent interest.

\begin{thm}\label{thm:Milne-4}
  With the notation of \thmref{thm:Kai-Main}, the diagram
  \begin{equation}\label{eqn:Milne-4-0}
    \xymatrix@C1pc{
      \Alb_X(k) \times  H^1_\et(k, \Picc^0(X)) \ar@<6ex>[d]^-{\delta_X}
\ar[r]^-{{\rm TM}} & \Br(k)  \ar@{=}[d] \\
 A_0(X) \ \ \ \times \ \ \ \wt{\Br}(X) \ar@<7ex>[u]^-{\alb_X}
 \ar[r]^-{{\rm BM}} & \Br(k)}
\end{equation}
of Brauer--Manin (BM) and Tate--Milne (TM) pairings is commutative.
\end{thm}
\begin{proof}
It suffices to prove the commutativity of ~\eqref{eqn:Milne-4-0} after replacing
  $\wt{\Br}(X)$ by $\Br'(X)$. Letting $H(X) = \Hom(\Br'(X), \Br(k))$, 
  it is enough to show that the diagram
\begin{equation}\label{eqn:Milne-4-1}
    \xymatrix@C1.5pc{
      A_0(X) \ar[r]^-{\alpha_X} \ar[d]_-{\alb_X} & H(X) \ar[d]^-{\delta^\vee_X} \\
      \Alb_X(k) \ar[r]^-{\beta_X} & G(X)}
\end{equation}
is commutative, where $\alpha_X$ (resp. $\beta_X$) is induced by the Brauer--Manin
(resp. Tate--Milne) pairing for $X$.
We shall prove the commutativity of ~\eqref{eqn:Milne-4-1} in several steps.

\vskip .2cm

{\bf{Step~1:}} Assume $\dim(X) \le 1$.
In this case, the lemma is already shown in \cite[\S~4]{Lichtenbaum}.

We now assume that $\dim(X) \ge 2$. We fix $\alpha \in A_0(X)$.

{\bf{Step~2:}} We assume $\alpha = [Q] - [P]$, where
$P, Q \in X(k)$. Using \cite[Thm.~7]{Altman-Kleiman}, we can find a smooth
projective curve $C \subset X$ containing $\{P,Q\}$ such that $C$ is a
scheme-theoretic intersection of very ample smooth divisors on $X$. In particular,
$C$ is geometrically integral. We let $\iota \colon C \inj X$ be the
inclusion. We then have a cycle $\alpha' = [Q] - [P] \in A_0(C)$ such that
$\alpha = \iota_*(\alpha')$. We now consider the diagram
\begin{equation}\label{eqn:Milne-4-2}
    \xymatrix@C1.5pc{
 A_0(C) \ar[rr] \ar[dr]_-{\iota_*} \ar[dd] & & H(C) \ar[dd] \ar[dr] & \\
      & A_0(X) \ar[rr] \ar[dd] & & H(X) \ar[dd]^-{\delta^\vee_X} \\
      \Alb_C(k) \ar[rr] \ar[dr] & & G(C) \ar[dr] & \\
      & \Alb_X(k) \ar[rr] & & G(X),}
    \end{equation}
where the slanted arrows are the push-forward map induced by $\iota$.

The top and the left faces of ~\eqref{eqn:Milne-4-2} commute by the
functoriality of
the Brauer--Manin pairing and the Albanese morphisms with respect to proper maps.
The right face commutes because ~\eqref{eqn:Pic-sch-HS-SS-2} is clearly
functorial with respect to the pull-back maps
$\iota^* \colon H^1_\et(k, \Pic(X^s)) \to  H^1_\et(k, \Pic(C^s))$ and
$\iota^* \colon \Br'(X) \to \Br'(C)$. The bottom face commutes by
\propref{prop:Milne-1}. The back face commutes because $\dim(C) =1$.
As $\alpha = \iota_*(\alpha')$, a diagram chase shows that
$(\delta^\vee_{X} \circ \alpha_{X})(\alpha) =
(\beta_{X} \circ \alb_{X})(\alpha)$.

{\bf{Step~3:}} We assume $X(k) \neq \emptyset$ and choose a point $P \in X(k)$.
We write $\alpha = \stackrel{n}{\underset{i =1}\sum} n_i[P_i]$. We let 
$e_i = [k(P_i):k]$ and $\alpha_i = [P_i] - e_i[P]$ so that $\alpha_i
\in A_0(X)$ for each $i$ and
$\alpha =  \stackrel{n}{\underset{i =1}\sum} n_i \alpha_i$ (note that
$\deg(\alpha) = 0$). It suffices to show  for each $i$ that
$(\delta^\vee_{X} \circ \alpha_{X})(\alpha_i) =
(\beta_{X} \circ \alb_{X})(\alpha_i)$. We can therefore assume that
$\alpha = [Q] - e[P] \in A_0(X)$, where $Q \in X_{(0)}$ and $e = [k(Q):k]$.

We let $k' = k(Q)$ and $X' = X_{k'}$. We let $\phi \colon S_{k'} \to S_k$ denote
the projection. We denote the base change map
$X' \to X$ also by $\phi$. We then have $P, Q \in X(k')$ and $\alpha =
\phi_*(\alpha')$, where $\alpha' := [Q] - [P] \in A_0(X')$.
We now look at the diagram
\begin{equation}\label{eqn:Milne-4-3}
    \xymatrix@C1.5pc{
      A_0(X') \ar[rr] \ar[dr]_-{\phi_*} \ar[dd] & & H(X') \ar[dd] \ar[dr] & \\
      & A_0(X) \ar[rr] \ar[dd] & & H(X) \ar[dd] \\
      \Alb_X(k') \ar[rr] \ar[dr] & & G(X') \ar[dr] & \\
      & \Alb_X(k) \ar[rr] & & G(X),}
    \end{equation}
where the slanted arrows are the push-forward or transfer maps induced by $\phi$.

By \lemref{lem:Milne-3}, all faces of ~\eqref{eqn:Milne-4-3} except possibly the
front and the back faces commute. Since $\alpha = \phi_*(\alpha')$, it suffices to
show that the back face of ~\eqref{eqn:Milne-4-3} commutes. But this follows from
Step~2.

{\bf{Step~4:}} We now prove the general case.
We let $\psi_X = \delta^\vee_{X} \circ \alpha_{X} - \beta_{X} \circ \alb_{X}
\in \Hom(A_0(X), G(X))$.
Since $X$ is geometrically integral, we can apply the Lang--Weil estimate to the
special fiber $Y$ of $\sX$ to find finite unramified extensions $k_1, k_2$ of $k$
such that $X(k_i) \neq \emptyset$ for $i = 1,2$ and $([k_1:k], [k_2:k]) =1$.
We let $d_i = [k_i:k]$ and $X_i = X_{k_i}$. Let $\phi_i \colon S_{k_i} \to S_k$
denote the projection. We denote the base change map
$X_i \to X$ also by $\phi_i$. 
If we let $\alpha_i = \phi^*_i(\alpha)$, then it follows from Step~2 that
$(\delta^\vee_{X_i} \circ \alpha_{X_i})(\alpha_i) =
(\beta_{X_i} \circ \alb_{X_i})(\alpha_i)$ for $i = 1,2$.
Using the diagram ~\eqref{eqn:Milne-4-3} (with $X'$ replaced by $X_i$),
we get 
$d_i (\psi_X(\alpha)) = \psi_X(d_i\alpha) = \psi_X(\phi_{i *}(\alpha_i)) = 0$ for
$i =1,2$. Since $(d_1, d_2) = 1$, this implies that $\psi_X(\alpha) = 0$.
This completes the proof.
\end{proof}

\subsection{Proofs of the main results}\label{sec:Prf-Kai-fin}
We shall now finish the proofs of \thmref{thm:Kai-Main} and
\corref{cor:Kai-Main-0}. We begin with the following result.

\begin{lem}\label{lem:Torsion-property}
  $\coker(\alb_X)$ is a torsion group and 
$\coker\left(\Pic(X^s)^\Gamma \to \ns(X^s)^\Gamma\right)$ is a finite $p$-group.
\end{lem}
\begin{proof}
Since the Albanese map is known to be surjective over algebraically closed
fields, and it is functorial with respect to base change by
field extension (cf. proof of \propref{prop:Alb-8}), any element 
$\alpha \in \coker(\alb_X)$ dies over a finite field extension of $k$. Since the
Albanese map is also functorial with respect to proper push-forward
(cf. \lemref{lem:Milne-3}), it follows that $\alpha$ is a torsion element.

To prove the second part, we follow the strategy of \cite[Rem.~3.3]{Kai}.
We let $L\Pic(X^s)$ denote the kernel of the
specialization map $\spc \colon \Pic(X^s) \to \Pic(\ov{Y})$ and let
$L\NS(X^s)$ denote the kernel of the map $\spc \colon \NS(X^s) \to \NS(\ov{Y})$
(cf. \lemref{lem:Pic-sch-3}). We fix an integer $n \in I_p$.

\vskip.2cm

{\bf Claim~1:} ${L\NS(X^s)}/n = 0$.

To prove the claim, we first note that the map
$\spc \colon L\Pic(X^s) \to L\NS(X^s)$ is surjective, as one checks using
~\eqref{eqn:Pic-sch-4}. It suffices therefore to show that
${L\Pic(X^s)}/n = 0$. To prove the latter, we use the commutative diagrams
\begin{equation}\label{eqn:Torsion-property-0}
  \xymatrix@C1.2pc{
    H^1_\et(X^s, {\Z}/n(1)) \ar[d]_-{\spc}^-{\cong} \ar@{->>}[r] &
    \Pic(X^s)[n] \ar[d]^-{\spc} & &
    \Pic(X^s)/n \ar@{^{(}->}[r] \ar[d]_-{\spc} & H^2_\et(X^s, {\Z}/n(1))
    \ar[d]^-{\spc}_-{\cong} \\
    H^1_\et(\ov{Y}, {\Z}/n(1)) \ar@{->>}[r] & \Pic(\ov{Y})[n] & &
    \Pic(\ov{Y})/n \ar@{^{(}->}[r] & H^2_\et(\ov{Y}, {\Z}/n(1)),}
  \end{equation}
obtained from the Kummer sequence.

The horizontal arrows are surjective in the first square and
injective in the second square by the Kummer sequence. The left vertical arrow
in the first square and the right vertical arrow in the second square are
bijective by the smooth proper base change theorem
(cf. \cite[Chap.~VI, Cor.~4.2]{Milne-EC}). It follows that the right vertical
arrow in the first square is surjective and the left vertical arrow in the
second square is injective.

We now use the long exact sequence
resulting from applying the functor $(-) \otimes {\Z}/n$ to the exact sequence
\[
0 \to L\Pic(X^s) \to \Pic(X^s) \to \Pic(\ov{Y}) \to 0,
\]
we conclude the proof of Claim~1.

Since $\NS(X^s) \subset \NS(\ov{X})$ by definition, and the latter
is a finitely generated abelian group (cf. \cite[\S~5.1, p.~122]{CTS}),
it follows that $\NS(X^s)$ is a finitely generated abelian group.
Combining this with the above claim, we deduce that $L\NS(X^s)$ is a finite
$p$-group. In particular, $\Ker\left(\NS(X^s)^\Gamma \to \NS(\ov{Y})^\Gamma\right)$
is a  finite $p$-group. Our next claim is the following.

\vskip.2cm

{\bf Claim~2:} The canonical map $\Pic({Y}) \to \NS(\ov{Y})^\Gamma$ is surjective.

To prove this claim, let $\Gamma_0$ denote the Galois group of $\F$.
We observe using the Hochschild-Serre spectral
sequence for $\ov{Y} \to Y$ (cf. \cite[(5.20)]{CTS} that the map $\Pic({Y}) \to
\Pic(\ov{Y})^{\Gamma_0}$ is an isomorphism. On the other hand, applying the
Galois cohomology sequence to ~\eqref{eqn:Pic-sch-HS-SS-5} (with
$X^s$ replaced by $\ov{Y}$), and using Lang's theorem which implies that
$H^1_\et(\F, \Picc^0(Y)) = 0$, we get that the map
$\Pic(\ov{Y})^{\Gamma_0} \to  \NS(\ov{Y})^{\Gamma_0} \cong \NS(\ov{Y})^{\Gamma}$
is surjective. Putting everything together, Claim~2 follows.

Using the two claims and the commutative diagram
\begin{equation}\label{eqn:Torsion-property-1}
  \xymatrix@C1.7pc{
    \Pic(X) \ar[r] \ar@{->>}[d] & \NS(X^s)^\Gamma \ar[d] \\
    \Pic(Y) \ar@{->>}[r] & \NS(\ov{Y})^\Gamma,}
\end{equation}
it follows that $\coker\left(\Pic(X) \to \NS(X^s)^\Gamma\right)$ is a finite
$p$-group. But this implies that
$\coker\left(\Pic(X^s)^\Gamma \to \NS(X^s)^\Gamma\right)$ is a finite $p$-group.
This completes the proof.
\end{proof}


\vskip .2cm

{{\bf{Proof of \thmref{thm:Kai-Main}:}}
  Lemma~\ref{lem:Torsion-property} implies that the two groups in question
  are discrete torsion groups.
A combination of \thmref{thm:Main-3}, \cite[Thm.~7.8]{Milne-Duality}
(cf. ~\eqref{eqn:MD-0}) and \thmref{thm:Milne-4} implies that that there is a
canonical isomorphism
\begin{equation}\label{eqn:Kai-Main-0-0}
  \beta_X \colon  \coker(\alb_X) \xrightarrow{\cong} (\Ker(\delta_X))^\vee.
\end{equation}
Combining this with ~\eqref{eqn:Pic-sch-HS-SS-6},
it suffices to show that $\theta$ is the zero map.

To that end, we look at the diagram ~\eqref{eqn:Pic-sch-HS-SS-7}.
Since this diagram is commutative by \lemref{lem:Pic-sch-4}, and
since the composite top horizontal arrow in this diagram is clearly zero, 
it remains to show that the composite map
\begin{equation}\label{eqn:Kai-Main-0-1}
\Br({\sX^s}/{\sX}) \xrightarrow{\iota^*} \Br({\ov{Y}}/Y)
\xrightarrow{\phi_{{\ov{Y}}/Y}} H^1_\et(\F, \Pic(\ov{Y})) \xrightarrow{\tau^*_Y}
H^1_\et(\F, \ns(\ov{Y})) \xrightarrow{\rho^*} \H^1_\et(k, \ns(\ov{Y}))
\end{equation}
in ~\eqref{eqn:Pic-sch-HS-SS-7} is injective.

Now, $\rho^*$ is injective because it the same as the inflation
map ${\rm inf} \colon H^1_\et(\F, \ns(\ov{Y})) \to \H^1_\et(k, \ns(\ov{Y}))$,
which is known to be injective. The map $\tau^*_Y$ is injective by Lang's theorem
used earlier. The map $\phi_{{\ov{Y}}/Y}$
is an isomorphism by using the analogue of ~\eqref{eqn:Pic-sch-HS-SS-2} for $Y$
and using that $\Br_0(Y) = 0$. Finally, $\iota^*$ is injective by
\lemref{lem:Br-res-spl}. This completes the proof.
\qed

\vskip .2cm

{{\bf{Proof of \corref{cor:Kai-Main-0}:}}
This is a direct consequence of  Lemma~\ref{lem:Torsion-property} and
  \thmref{thm:Kai-Main}.
  \qed

\section{A Bertini theorem for semistable schemes over a dvr}\label{sec:Bertini}
In this section, we prove a Bertini theorem for semistable schemes over a discrete
valuation ring. Special cases of this theorem were used in
Sections~\ref{sec:Esp} and ~\ref{sec:Esp-0}, and other places in this
paper. Our proofs are based on \cite{GK-JLMS},
and we adopt the notation and terminology introduced there throughout this section.
Unless otherwise stated, all schemes considered in this section are assumed to be
Noetherian, separated, excellent, and of finite Krull dimension.

\subsection{Bertini theorem for normal crossing schemes}\label{sec:NCS}
Let $Z$ be a scheme and let $X \subset Z$ be a closed subscheme. Recall that
$X$ is said to be an snc divisor on $Z$ if the latter is a regular scheme, and at
every closed point $x \in X$, the kernel of the canonical surjection $\sO_{Z,x} \surj
\sO_{X,x}$ is an ideal of the form $(x_1\cdots x_r)$, where $(x_1, \cdots , x_{m+1})$
is the maximal ideal of $\sO_{Z,x}$ and $r \in J_{m+1}$.

A scheme $X$  is called a simple normal crossing (snc) scheme
if it is reduced with irreducible components $X_1, \ldots , X_r$, each of dimension
$m$ such that for all $\emptyset \neq J \subseteq J_r$, the
scheme-theoretic intersection $X_J := {\underset{i \in J}\bigcap} X_i$ is either
empty or regular of pure dimension $m+1-|J|$. By a scheme of negative dimension, we
shall mean the empty scheme. We let $X_\emptyset = X$. A closed subscheme $W$ of
$X$ is called a stratum of $X$ if there exists $J \subset J_r$ such that $W$ is
an irreducible component of $X_J$. For $J \subset J_r$, we let
$Z_J = {\bigcup}_{J' \supsetneq J} X_{J'}$ and $U_J = X_J \setminus Z_J$.
We shall need the following property of snc schemes.

\begin{lem}\label{lem:SNC-0}
  Let $X$ be an equidimensional scheme of dimension $m$ with irreducible
  components
  $X_1, \ldots , X_r$. Let $Y \subset X$ be a reduced locally principal closed
  subscheme such that $X_i \bigcap Y$ is connected for every $i \in J_r$ and
  $U_J \bigcap Y$ is either empty or a regular equidimensional scheme of dimension
  $m - |J|$ for every $J \subset J_r$. Then we have the following.
  \begin{enumerate}
    \item
    If there is a closed immersion $X \inj X'$ of equidimensional schemes such that
    $X'$ is regular of dimension $m+1$, then $X$ is an snc scheme if and only if it
    is an snc divisor on $X'$.
    \item
    If $X$ is an snc scheme and $J \subset J' \subset J_r$ are such that
    $|J'| = |J| +1$, then $X_{J'}$ is a regular Cartier divisor on $X_J$.
  \item
    If $X$ is an snc scheme and $J \subset J_r$, then $Z_J$ is an snc divisor on
    $X_J$ whose irreducible components are the $(m - |J|)$-dimensional strata of $X$
    contained in $X_J$.
  \item
    If $X$ is an snc scheme, then $Y$ is an snc scheme of dimension $m - 1$ with
    irreducible components $X_1 \bigcap Y, \ldots , X_r \bigcap Y$.
  \item
    If $X$ is projective over a field in (4) and $Y$ is an ample divisor on $X$,
    then for every $d \ge 1$, the $d$-dimensional strata of $Y$ are of the form
    $W \bigcap Y$, where $W$ runs through $(d+1)$-dimensional strata of $X$.
    \end{enumerate}
\end{lem}
\begin{proof}
Item (1) follows from \cite[Tag~0BIA]{SP}, (2) is obvious and (3) follows 
  from (1) because one checks directly from the definition that $Z_J$ is an
  snc scheme of dimension $m- |J|$. To prove (4), we let $Y_i = X_i \bigcap Y$
  and first claim that each $Y_i$ is irreducible and regular of dimension $m-1$.
  To prove this claim, we can assume without loss of
  generality that $i =1$. We let $X^j_1 = X_1 \bigcap \cdots \bigcap X_j$ for
  $j \in J_r$ and $Y^j_1 = X^j_1 \bigcap Y$ so that we have a chain of closed
  embeddings $Z^r_1 \subsetneq \cdots \subsetneq Z^2_1 \subsetneq Z^1_1 = Z_1$ for
  $Z \in \{X_1, Y_1\}$. We shall show more generally by a descending induction on $j$
  that each $Y^j_1$ is regular of pure dimension $m - j$. We are already given that
  $Y^r_1$ is regular of pure dimension $m-r$ (note that $X_{J_r} = U_{J_r}$).

We now let $j < r$ and write $X^j_1 = W_1 \amalg \cdots \amalg W_s$ as a disjoint
  sum of connected regular schemes of dimension $m+1 -j$ such that
  $V_i := W_i \setminus X_{j+1}$ is dense in $W_i$ for
  every $i$. It follows that $Y^j_1 = (W_1 \bigcap Y) \amalg \cdots \amalg
  (W_s \bigcap Y)$ and
  $Y^j_1 \setminus Y^{j+1}_1 = (V_1 \bigcap Y) \amalg \cdots \amalg (V_s \bigcap Y)$.
  Let $x \in Y^j_1$ be a closed point. If $x \in V_i$ for some $i$, then
  $\sO_{Y^j_1, x}$ is a regular local ring of dimension $m-j$ by our assumption.
  Suppose now that $x \in W_i \bigcap Y^{j+1}_1$ for some $i$. Then we can write
  $\sO_{Y^j_1, x} = {\sO_{X^j_1,x}}/{(f)}$ and
  $\sO_{Y^{j+1}_1, x} =  {\sO_{X^{j+1}_1,x}}/{(f)} = {\sO_{X^j_1,x}}/{(f,g)}$.

We know that
  $\sO_{X^j_1,x}$ is a regular local ring of dimension $m+1-j$. The induction
  hypothesis meanwhile implies that $\sO_{Y^{j+1}_1, x}$ is a regular local ring of
  dimension $m-1-j$. It follows that if $x_1, \ldots , x_{m-1-j} \in \sO_{X^j_1,x}$
  are such that their images in $\sO_{Y^{j+1}_1, x}$ generate latter's maximal ideal,
  then $(f,g, x_1, \ldots , x_{m-1-j})$ is the maximal ideal of $\sO_{X^j_1,x}$.
  This implies that ${\sO_{X^j_1,x}}/{(f)} = \sO_{Y^j_1, x}$ is a regular local ring of
    dimension $m-j$. We have thus shown that $Y^j_1$ is a regular scheme of pure
    dimension $m-j$. Since $\dim(Y^{j+1}_1) = m-1-j$, no component of $Y^j_1$ can
    be contained in $Y^{j+1}_1$. In particular, $Y^j_1 \setminus Y^{j+1}_1$ is
    furthermore dense in $Y^j_1$. Finally, to prove our claim, it remains to
    show that each $Y_i$ is irreducible. But this is clear since $Y_i$ is  regular
    and connected.

To finish the proof of (4), we let $J \subset J_r$ and define $Y_J$ similarly to
    $X_J$. We need to show that $Y_J$ is either empty or a regular scheme of
    pure dimension $m-|J|$ if $J \neq \emptyset$. We shall prove this by a descending
    induction on $J$. If $|J| =r$, this is already clear from our assumption.
    If $0 < |J| < r$, we choose $i \in J_r \setminus J$ and let
    $J' = J \bigcup \{i\}$ so that
    $Y_{J'} = Y_J \bigcap X_i = (X_J \bigcap X_i) \bigcap Y$. We let $W$ be an
    irreducible component of $X_J$ and let $x \in (W \bigcap X_i) \bigcap Y$ be a
    closed point. We can then write $\sO_{W \bigcap Y,x} = {\sO_{W,x}}/{(f)}$ and
    $\sO_{Y_{J'},x} = {\sO_{W \bigcap X_i},x}/{(f)} = {\sO_{W,x}}/{(f,g)}$, where the last
    identity holds because $W$ is regular of dimension $m+1-|J|$ and
    $W \bigcap X_i$ is regular of dimension $m-|J|$. In particular, $W \bigcap X_i$
    is locally principal in $W$.

We now know that ${\sO_{W,x}}$ is a regular local ring of dimension $m+1 -|J|$.
    The induction hypothesis implies that $\sO_{Y_{J'},x}$ is a regular local ring of
    dimension $m-1-|J|$. It follows that if $x_1, \ldots , x_{m-|J'|} \in \sO_{W,x}$
  are such that their images in $\sO_{Y_{J'}, x}$ generate latter's maximal ideal,
  then $(f,g, x_1, \ldots , x_{m-|J'|})$ is the maximal ideal of $\sO_{W,x}$.
  This implies that ${\sO_{W,x}}/{(f)} = \sO_{Y_{J}, x}$ is a regular local ring of
  dimension $m-|J|$. We have thus shown that $Y_J$ is a regular scheme of pure
  dimension $m-|J|$. This completes the proof of (4).

 If $X$ is as in (5) and $W$ is any stratum of $X$ of dimension at least 2, then
 $W \bigcap Y$ is an ample divisor on $W$ and hence connected because
 $\dim(W) \ge 2$ (this is an easy consequence of \cite[Lem.~5.1]{GK-Jussieu}).
  Since $W \bigcap Y$ is regular as shown in (4), it must be irreducible.
  This completes the proof of the lemma.
\end{proof}

\begin{lem}\label{lem:Bertini-field}
  Let $k$ be a field of arbitrary characteristic
  and let $X \inj \P^N_k$ be a closed snc subscheme of
  dimension $m \ge 2$ with irreducible components $X_1, \ldots , X_r$.
  Let $T, Z \subset \P^N_k$ be 0-dimensional
  closed subschemes such that the following hold.
  \begin{enumerate}
  \item
    $T \bigcap X_i \neq \emptyset$ for every irreducible component $X_i$ of $X$.
  \item
    $T \bigcap Z = \emptyset$.
  \item
    $Z \bigcap X \subset X_\reg$.
  \item
    $Z$ is reduced.
 \item
    $X$ is Cohen-Macaulay.
 \end{enumerate}

  If $k$ is infinite, then for all $d \gg 0$, there exists a
  dense open $U \subset |\sI_Z(d)|$ such that every $H \in U(k)$
  has the property that $T \bigcap H = \emptyset$ and $X \bigcap H$
  is an snc scheme of
  dimension $m-1$ with irreducible components
  $X_1 \bigcap H, \ldots , X_r \bigcap H$. If $k$ is finite and we let
  \[
  \sP = \left\{f \in I^Z_\homg|X \cap H_f \ \mbox{is  an  snc  scheme  and} \
  T \cap H_f = \emptyset\right\},
  \]
  then there exists $\sP' \subseteq \sP$ such that
  $\mu_Z(\sP') > 0$. For $f \in \sP$, the irreducible components of $X \bigcap H_f$
  are $X_1 \bigcap H_f, \ldots , X_r \bigcap H_f$.
\end{lem}
\begin{proof}
We let $X_1, \ldots , X_r$ be the irreducible components of $X$.
  Assume first that $k$ is infinite and let
  \[
  \sP(d) = \{H \in |\sI_Z(d)|| T \cap H = \emptyset \ \mbox{and} \
  X_J \cap H \ \mbox{is  empty  or  regular  of  dimension} 
\]
\hspace*{1.4cm} $\ m - |J| \ \forall \ J \subset J_r\}$,
\\
where $|\sI_Z(d)| \subseteq |\sO_{\P^N_k}(d)|$ is the linear system defined by the
global sections of $\sI_Z(d)$.

Since each $X_J$ is regular, it follows from \cite[Thm.~3.4]{GK-JLMS} that
for all $d \gg 0$, there exists a dense open $U \subset |\sI_Z(d)|$
such that every $H \in U(k)$ has the property that $T \bigcap H = \emptyset$ and
$H \in \sP(d)$. We pick $H \in U(k)$ and let $Y = X \bigcap H$.
To show that $Y$ is an snc scheme of dimension $m-1$, it suffices to show
using \lemref{lem:SNC-0}(4) that $Y$ is reduced and $X_i \bigcap H$ is connected
for all $i \in J_r$.

Since $X$ is reduced, its associated points coincide with the generic points
(cf. \cite[Lem.~3.3]{GK-JAG}). It follows from condition (1) of the
lemma that $H$ does not meet any associated point of $X$. In particular,
$Y$ is a Cartier divisor on $X$. The condition (5) now implies that
$Y$ is also Cohen-Macaulay. Furthermore, $Y$ is a scheme which satisfies
Serre's $R_0$-condition as $H \in \sP(d)$. These two properties imply
that $Y$ is reduced. The claim that $X_i \bigcap H$ is connected
for any $i \in J_r$ is an easy consequence of \cite[Lem.~5.1]{GK-Jussieu}
because $X_i$ is an integral normal projective $k$-scheme of dimension at least two. 
This completes the proof of the lemma when $k$ is infinite.

We now assume $k$ to be finite and let
\[
\sP' = \left\{f \in I^Z_\homg|T \cap H_f = \emptyset  \ \mbox{and} \ X_J
\cap H_f \ 
\mbox{is empty or regular of dimension} \ m - |J| \ \forall \ 
J\right\}.
\]
Then \cite[Prop.~4.8]{GK-JLMS} implies that $\mu_Z(\sP') > 0$. It therefore remains
to show that $\sP' \subset \sP$. This follows by applying \lemref{lem:SNC-0}(4) and
repeating the argument used in the infinite field case above. The final assertion
is also proved in the same way as in the infinite field case.
\end{proof}

\subsection{Bertini theorem for semistable schemes}\label{sec:BSS}
Let $R$ be a discrete valuation ring with quotient field $k$, residue field $\F$
and let $\fm = (\pi)$ denote the maximal ideal of $R$. We let $S = \Spec(R)$.
Unlike the other sections of this paper, we allow $R$ to have mixed or
equal characteristic. Let $\sX \subset \P^N_R$ be a connected and regular closed
subscheme which is flat over $S$ and has relative dimension $m \ge 2$ over $S$.
Let $\sX_\eta$ (resp. $\sX_s$) denote the generic (resp. special) fiber of $\sX$ and
let $Y = (\sX_s)_\red$. We assume that $\sX_\eta$ is smooth over $k$ and $Y$ is an snc
divisor on $\sX$. Let $Y_1, \ldots , Y_r$ denote the irreducible components of $Y$.

The following is the main result of this section.

\begin{thm}\label{thm:Bertini-dvr}
 Let $T, Z_1, Z_2 \subset \P^N_R$ be  closed subschemes such that the following hold.
  \begin{enumerate}
  \item
    $T$ is finite over $S$.
  \item
    Either $Z_1 = \emptyset$ or
    $Z_1 \subset Y_\reg$ is reduced and finite over $\Spec(\F)$. 
  \item
    Either $Z_2 = \emptyset$ or $Z_2$ is finite {\'e}tale over $S$.
  \item
    $T \bigcap Y_i \neq \emptyset$ for every $i$.
  \item
    $T \bigcap (Z_1 \bigcup Z_2) = \emptyset$.
 \end{enumerate}
  Then for all $d \gg 0$, there are infinitely many hypersurfaces $H \subset \P^N_R$
  of degree $d$ defined over $R$ such that letting $\sX' = \sX \bigcap H$, we have
  the following.
  \begin{enumerate}
  \item
    $\sX'$ is a connected and regular scheme which is flat of relative
    dimension $m-1$ over $S$ .
  \item
    $\sX'_\eta$ is smooth over $k$.
  \item
    $\sX'_\eta$ is geometrically integral over $k$ if $\sX_\eta$ is so.
 \item
   $\sX' \bigcap Y = (\sX'_s)_\red$ is an snc divisor on $\sX'$ whose irreducible
   components are $Y_1 \bigcap H, \ldots , Y_r \bigcap H$.
  \item
    $\sX' \bigcap T = \emptyset$.
  \item
    $Z_1, Z_2 \subset H$.
  \end{enumerate}
\end{thm}
\begin{proof}
By \cite[Lem.~5.1]{GK-Jussieu}, there exists $d_1 \gg 0$ such that
  $H^i_\zar(\sX_\eta, \sO_{\sX_\eta}(-j)) = 0$ for all
  $j \ge d_1$ and $0 \le i \le 1$. It follows that
  $H^0_\zar(\sX_\eta, \sO_{\sX_\eta}) \to H^0_\zar(\sX_\eta \bigcap G, \sO_{\sX_\eta \cap G})$
  is an isomorphism for every $G \in |\sO_{\P^N_K}(j)|(k)$ whenever $j \ge d_1$.
 In particular, $\sX_\eta \bigcap G$ is geometrically connected over $k$ if
 $\sX_\eta$ is geometrically connected over $k$ (cf. \cite[Tag 03GX]{SP}).
  Combining this with \cite[Thm.~3.10]{GK-JLMS},  we get that there exists
 $d \gg 0$ such that for all $j \ge d$, we can find a dense open
  subset $U \subset |\sI_{(Z_2)_\eta}(j)|$ such that every $G \in U(k)$
  has the property that $T_\eta \bigcap G = \emptyset$, $(Z_2)_\eta \subset G$ and
  $\sX_\eta \bigcap G$ is a smooth and connected $k$-scheme of dimension $m-1$.
Furthermore, $\sX_\eta \bigcap G$ is geometrically connected over $k$ if $\sX_\eta$ is
  so.

Since $Y$ is an snc divisor on a regular scheme, it is Cohen-Macaulay.
  Under our assumptions, \lemref{lem:Bertini-field} therefore implies that for all
  $d \gg 0$, there are infinitely many hypersurfaces $F \subset \P^N_\F$ having
  degrees at least $d$ defined over $\F$ such that $Z_1 \subset F$,
  $T \bigcap F = \emptyset$ and $Y \bigcap F$ is an snc $\F$-scheme of
  dimension $m-1$. We choose such an $F$ of degree $d'$.

We now apply \cite[Lem.~6.4]{GK-JLMS} which says that if ${\rm sp} \colon
\P^{q}_k(k) \to \P^{q}_\F(\F)$ denotes the specialization map
(where $q = \dim_k(|\sI_{Z_\eta}(d')|)$) , then
  ${\rm sp}^{-1}(F) \bigcap U(k)$ is an infinite set. We pick any
  $H' \in {\rm sp}^{-1}(F) \bigcap U(k)$ and let $H$ denote the closure of $H'$ in
  $\P^{N}_R$. It is then clear that $H \subset \P^N_R$ is defined over $R$ such
  that $H_\eta = H'$ and $H_s = F$. We let $\sX' = \sX \bigcap H$ so that
  $\sX'_\eta = \sX_\eta \bigcap H'$ and $\sX'_s = \sX_s \bigcap F$.
  Since $Y \bigcap F$ is reduced, letting $Y' := (\sX'_s)_\red$, we get
  $Y' =  Y \bigcap F$. This immediately implies properties (2), (5)
and (6). This also implies that $Y'$ is an snc $\F$-scheme of dimension $m-1$
  with irreducible components as claimed in (4). 

To establish property (1), we let $x \in Y'$ be a
  closed point so that $x \in Y'_i := Y_i \bigcap F$ for some $i \in J_r$.
  Since $Y_i$ is a regular Cartier divisor on
  $\sX$ and $\sX' = \sX \bigcap H$ is a locally principal closed subscheme of $\sX$,
  there are
  non-zero elements $f, g \in \sO_{\sX,x}$ such that
  $\sO_{\sX',x} = {\sO_{\sX,x}}/{(f)}$,
  $\sO_{Y_i,x} = {\sO_{\sX,x}}/{(g)}$ and $\sO_{Y'_i,x} = {\sO_{Y_i,x}}/{(f)} =
  {\sO_{\sX,x}}/{(f,g)}$. Since $\sO_{Y'_i,x}$ is a regular local ring of
  dimension $m-1$, we can find elements $x_1, \ldots , x_{m-1} \in \sO_{\sX,x}$ whose
  images in $\sO_{Y'_i,x}$ generate latter's maximal ideal. This implies that
  $(f,g, x_1, \ldots , x_{m-1})$ is the maximal ideal of $\sO_{\sX,x}$. We conclude
  from this that $\sO_{\sX',x} = {\sO_{\sX,x}}/{(f)}$ is a regular local ring of
  dimension $m$. We have thus shown that $\sX'$ is regular such that
  $\dim_k(\sX'_\eta) = \dim_\F(\sX'_s) = m-1$.

\enlargethispage{20pt}
  
We next claim that $\sX'$ is connected. If not, then connectedness of $\sX_\eta$
  implies that $\sX'$ has a connected component which is faithfully flat over $S$
  and any other connected component must be contained in $\sX'_s$. In particular,
  it has dimension at most $m-1$. But this is a contradiction since the dimension
  of $\sX'$ at each of its closed point is $m$ by \cite[Thm.~14.2]{Matsumura}.
  It follows that $\sX'$ is connected (hence integral) and flat over $S$. This
  proves (1). If $\sX_\eta$ is geometrically integral (in particular, geometrically
  connected) over $k$, then it follows from the choice of $U$ that
  $\sX'_\eta = \sX_\eta \bigcap H_\eta = \sX_\eta \bigcap H'$ is geometrically
  connected over $k$. In combination with (2), this implies that $\sX'_\eta$ is
  geometrically integral over $k$. This proves (3).
Finally, we apply \lemref{lem:SNC-0}(1) to conclude that
  $Y'$ is an snc divisor on $\sX'$. This completes the proof.
\end{proof}

\vskip .4cm

\noindent\emph{Acknowledgements.}
Part of this work was carried out while the first author was visiting the Institute
of Mathematical Sciences (IMSc), Chennai, during the summer of 2025. He would like to
thank the institute for the invitation, support, and excellent working conditions.
The second author was supported by a SERB-CRG research grant from the Government of
India while working on part of this project at ISI Bangalore.

\end{document}